# Properties of Eigenvalues on Riemann Surfaces with Large Symmetry Groups

by

Joseph Cook

A Doctoral Thesis

Submitted in partial fulfillment of the requirements

for the award of

Doctor of Philosophy of Loughborough University

December 2, 2018

# Abstract

On compact Riemann surfaces, the Laplacian $\Delta$ has a discrete, non-negative spectrum of eigenvalues $\{\lambda_i\}$ of finite multiplicity. The spectrum is intrinsically linked to the geometry of the surface. In this work, we consider surfaces of constant negative curvature with a large symmetry group. It is not possible to explicitly calculate the eigenvalues for surfaces in this class, so we combine group theoretic and analytical methods to derive results about the spectrum.

In particular, we focus on the Bolza surface and the Klein quartic. These have the highest order symmetry groups among compact Riemann surfaces of genera 2 and 3 respectively. The full automorphism group of the Bolza surface is isomorphic to $\mathrm{GL}_2(\mathbb{Z}_3) \rtimes \mathbb{Z}_2$; in Chapter 3 we analyze the irreducible representations of this group and prove that the multiplicity of $\lambda_1$ is 3, building on the work in [41], and identify the irreducible representation that corresponds to this eigenspace. This proof relies on a certain conjecture, for which we give substantial numerical evidence and a hopeful method for proving. We go on to show that $\lambda_2$ has multiplicity 4.

The full automorphism group of the Klein quartic is isomorphic to $\mathrm{PGL}_2(\mathbb{Z}_7) \rtimes \mathbb{Z}_2$. In Chapter 4 we give some preliminary results about the surface, including a calculation of its systole length and Fenchel-Nielsen coordinates, and prove that the multiplicity of $\lambda_1$ for the Klein quartic is 6, 7, or 8. In the appendices, we include numerical computations of eigenvalues for the Bolza surface, the Klein quartic, and other interesting hyperbolic surfaces, such as the Fermat quartic.

# Acknowledgements

I am indebted to my supervisor Alexander Strohmaier for providing me with an interesting project that I have enjoyed from the start, and for introducing me to several fascinating areas of mathematics. His encouragement, expertise, and patience in explaining and re-explaining key areas of subject matter have been invaluable in the production of this work. I am grateful to Eugenie Hunsicker for the beneficial advice and feedback she has given me at yearly progress review meetings, which has been imperative in shaping this thesis. I would like to thank Claudia Garetto for agreeing to co-supervise this project, and for checking the final draft.

I am grateful for the support and friendship of my fellow PhD students, both at Loughborough University and the University of Leeds. The discussions and opportunities to present my work to my peers have kept me motivated throughout my time as a postgraduate, and have vastly improved my confidence in speaking about mathematics.

I would like to thank the EPSRC for providing financial support throughout my period of study, as well as the Loughborough University School of Science, and University of Leeds School of Mathematics, for facilitating and financing my trips to various conferences over the past three and a half years. Furthermore, I would like to thank the School of Mathematics at the University of Leeds for hosting me when I moved here with my supervisor. They have given me all the benefits of a full University of Leeds PhD student, for which I am beholden.

I would like to thank the providers of free software: notably, FreeFEM++ and GAP. These programs have vastly enriched the theoretical work I have done. In particular, images produced in FreeFEM++ have been an asset when presenting my work at conferences and seminars.

Finally, I would like to acknowledge that the completion of this research project would not have been possible without the continued love and support my family and close friends.

*Dedicated to Alice*

# Contents

# Chapter 1

# Introduction

**Paris, February 1809**

The court of Napoleon Bonaparte fills with the sound of a Haydn melody, as our performer teases a curious rotating glass cylinder with felt-covered rods. The performer is Ernst Chladni, the "father of acoustics". Upon finishing his concert, he proceeds to demonstrate his experiments with acoustics to the emperor. These involve bowing a fixed plate at different frequencies, and observing how the patterns made by sand on the plate change as the frequency increases. The patterns that he recorded on a square plate were documented in his seminal work, *Akusik*.

Whilst the music of the clavicylinder is of little interest to the contemporary mathematician, the investigation of the patterns of sand on a vibrating plate, that is, the nodal lines of Dirichlet eigenfunctions, remains an interesting problem. Napoleon himself remarked to Chladni that if one could apply some of the calculations involved in the mathematical formulation of acoustics to areas curved in more than one direction, it "could be useful for applications to other subjects as well" [76].

What do we understand by the term "curved in more than one direction"? Broadly speaking, the curvature of a surface $S$ at a point $p$ measures how much the surface differs from its tangent plane $T_pS$, that is, how far away it is from being "flat". If the curvature is constant, then this measure takes the same value for every $p \in S$. For a simple example of a surface of constant curvature, one can consider the sphere. There are several different notions of curvature; in this work we predominantly consider Gaussian curvature. The normal vector $\mathbf{n}$ forms a plane with each vector $v$ in $T_p\mathbb{S}^2$, with each direction $v$ corresponding the direction of a curve on the surface. The principal curvatures $\kappa_i$ are the (orthogonal) directions that give the maximum and minimum "bending" of $S$ at $p$. The Gaussian curvature is then defined at a point $p$ as the product of the principal curvatures, that is,
$$K = \kappa_1(p)\kappa_2(p).$$

Figure 1.1: Principal directions at a point $p$ on $\mathbb{S}^2$

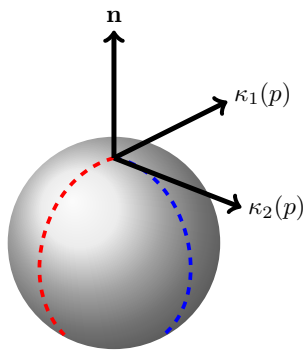

When dealing with surfaces of constant Gaussian curvature, there are three main cases: when $K$ is positive, when it is negative, and when it is zero. We call these surfaces spherical, flat or hyperbolic. Having positive Gaussian curvature means that the principal curvatures at $p$ curve "in the same way", so the area around $p$ resembles a cap. Conversely, negative curvature implies that the principal curvatures curve in different directions, so $p$ is both a maximum and a minimum. Every point on a hyperbolic surface resembles a saddle, or indeed, a Pringle. Hitherto, we have alluded to the notion of surfaces, but have not said precisely what we mean. We take this opportunity to formally introduce the setting for the analysis in this work.

**Definition 1.1 (Riemann surface)** *A Riemann surface is a two (real) dimensional, orientable Riemannian manifold. The prefix "Riemannian" indicates that a smooth manifold has a Riemannian metric, that is, a smooth family of inner products $g_p(\,\cdot\,,\,\cdot\,)$ such that*
$$g_p : T_pM \times T_pM \to \mathbb{R}$$
*depends smoothly on $p \in M$.*

If $X$ and $Y$ are $C^\infty$ vector fields on $M$, then their inner product with respect to the Riemannian metric is a $C^\infty$ real valued function on $M$. The Riemannian metric gives us a consistent way of measuring angles on the surface, and thus endows the surface with a conformal structure.

The property of being orientable means that one can make a consistent choice of normal vector on the surface, that is, we can move a normal vector along a closed loop continuously, without it changing sign. Figure 1.2 shows the Möbius strip, which is an example of a manifold that is not orientable. To see this, move the outward pointing normal vector around the closed loop. When it returns to its original position, it will be inward pointing, so the Möbius strip is not orientable. If we fix a an orientation (normal vector) on a surface, we say that it is oriented. Where a surface $S$ has boundary, we define an orientation on its interior, and this in turn induces an orientation on the boundary $\partial S$.

These additional conditions on a smooth manifold of two real dimensions allow us to give the manifold a complex structure. This leads us to an alternative definition of a Riemann surface as a connected, one (complex) dimensional complex manifold. The complex structure allows us to define holomorphic functions between Riemann surfaces; in particular, the transition maps are holomorphic.

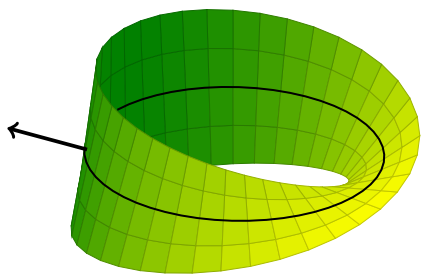

Figure 1.2: The Möbius strip is not orientable

The Bolza surface and Klein quartic are interesting because they have maximal symmetry among hyperbolic surfaces in their respective topological classes. The topological invariant that we use to classify such surfaces is the genus.

**Definition 1.2 (Genus)** *The genus of a connected, orientable surface is the maximum number of simple, closed, non-intersecting curves that one can cut along such that the resulting surface remains connected.*

For example, the genus of the torus is 1, because one can cut along a simple, closed curve to obtain a cylinder, but another such cut would result in two cylinders. The genus may also be thought of as the number of handles glued onto a sphere; one sees straightaway that the sphere has genus 0. The Bolza surface has genus 2, and the Klein quartic has genus 3.

The following theorem provides a link between the geometry of a surface and its topology (see, for example, [16]).

**Theorem 1.3 (Gauss-Bonnet)** *Let $S$ be a compact Riemann surface without boundary, with Gaussian curvature $K$. Then*

$$\int_S K \, dA = 2\pi \chi(S),$$

*where dA is the area measure of S, and $\chi(S)$ is its Euler characteristic. $\chi(S)$ is another topological invariant, and can be calculated from the genus g using the formula*

$$\chi(S) = 2 - 2g.$$

**Remark 1.4** *For compact Riemann surfaces of constant curvature, the uniformization theorem (introduced in Section 2.2.1) enables us to normalize the Gaussian curvature to $K = \pm 1$ or $K = 0$, depending respectively on whether the surface is positively or negatively curved, or flat. On surfaces of constant negative curvature, the Gauss-Bonnet theorem then implies*

$$\text{Area}(S) = 4\pi(g - 1),$$

*where g is the genus of the surface.*

Napoleon showed great foresight in his comment regarding the study of curved surfaces. He was (inadvertently, one supposes) predicting an arena for the study of: the Helmholtz wave equation

$$\nabla^2 A - k^2 A = 0,$$

where $\nabla^2$ is the Laplacian, $A$ is the amplitude and $k$ is the wave number; the stationary wave equation

$$\left(\nabla^2 - \frac{1}{c^2}\frac{\partial^2}{\partial t^2}\right)u(\mathbf{x}, t) = 0;$$

and the Schrödinger equation

$$H\Psi = E\Psi, \tag{1.1}$$

where $H$ is the Hamiltonian and $E$ is the energy of the state $\Psi$. We will examine this third equation shortly in the context of the Bolza surface.

Although spectral theory in a Euclidean setting was commonplace in physics prior to the 1950s, during this decade the geometry and spectra of hyperbolic surfaces became an avenue of interest for quantum physicists. These surfaces provide a model for studying quantum chaos; Selberg's trace formula helped to bridge the gap between classical and quantum mechanics.

We will study the Laplace-Beltrami operator acting on functions on compact Riemann surfaces of constant negative curvature, in particular, those $\phi_i$ that solve the equation

$$\Delta \phi_i - \lambda_i \phi_i = 0, \tag{1.2}$$

where the eigenvalues $\lambda_i$ form a discrete spectrum

$$0 = \lambda_0 < \lambda_1 \leq \lambda_2 \leq \ldots.$$

This is a heuristic overview of the problem. In the following section, we give details of the different boundary value problems on the surfaces that we consider. We go into greater detail about the formulation of Equation (1.2) in Section 2.1, where we define the Laplacian and its domain, and prove that the eigenfunctions $\phi_i$ exist and form an orthonormal basis of the $L^2$ space of functions on the surface.

The interplay between hyperbolic surfaces and Fuchsian groups gives us additional tools to analyze their spectral theory. In Chapter 2, we give an outline of hyperbolic geometry, and the analytical methods used to investigate the spectra of these surfaces. We introduce Selberg's trace formula (Theorem 2.47), and give a generalization to cases where an isometry group is acting.

## Eigenvalue Problems

A topic of central interest in this work is the eigenvalue spectrum of the Laplacian corresponding to a particular closed or boundary value problem, which will take the form of one of the following. We follow the definitions of [15]. The Riemann surfaces that we consider have no boundary, so their eigenvalues and eigenfunctions satisfy the closed eigenvalue problem:

**Definition 1.5 (Closed eigenvalue problem)** *For a connected, compact surface $M$, find all $\lambda \in \mathbb{R}$ such that*
$$\Delta\phi - \lambda\phi = 0$$
*has nontrivial solutions $\phi \in C^2(M)$.*

In practice, we will not tackle the closed eigenvalue problem head on. It is not possible to calculate such eigenvalues explicitly, and the surface as a whole is often too complicated to analyze. Rather, we will use the symmetry of the surface to break the closed problem down into one of the following, where we consider a boundary value problem on a subspace of the surface. We define the following problems on the more general setting of a connected, smooth, oriented Riemannian manifold $M$. Additional conditions will be specified. Giving the definitions on general manifolds will allow us to consider these problems on sub-manifolds later.

**Definition 1.6 (Dirichlet eigenvalue problem)** *For $\partial M \neq \emptyset$, $\overline{M}$ compact and connected, find all $\lambda \in \mathbb{R}$ such that*
$$\Delta\phi - \lambda\phi = 0$$
*has nontrivial solutions $\phi \in C^2(M) \cap C^0(\overline{M})$ that satisfy the boundary condition*
$$\phi = 0$$
*on $\partial M$.*

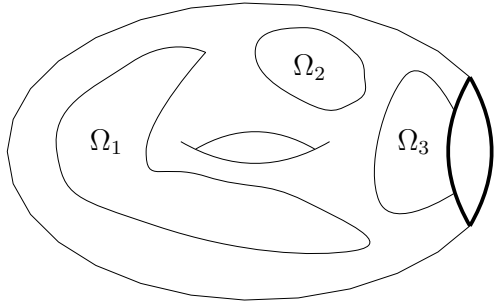

Figure 1.3: Domain monotonicity for a torus with a disk removed

**Remark 1.7 (Domain monotonicity, Dirichlet data)** *Unless otherwise stated, "domain" refers to an open, connected subset in $M$ with compact closure. Ideally, when $\Omega$ is a subset of a compact Riemann surface, its boundary $\partial\Omega$ will be a smooth curve. In practice, this is not always the case; some of the domains we wish to consider have piecewise smooth boundary, by which we mean that the boundary is formed as the finite union of smooth curves. Let*

$$\Omega_1, \ldots, \Omega_m$$

*be pairwise disjoint domains with piecewise smooth boundary and compact closure in $M$. If $M$ has boundary, the boundary of $\Omega_i$, $i = 1, \ldots, m$, meets $\partial M$ tranversally. For example, in Figure 1.3, $M$ is a torus with a disk removed, $\partial M$ is a circle, and we have three pairwise disjoint domains. For $i = 1, \ldots, m$, consider the Dirichlet eigenvalue problem on $\partial\Omega_i \cap M$ (where $\partial M \neq \emptyset$, we do not impose a constraint on $\partial\Omega_i \cap M$). We arrange all the eigenvalues of these domains into the increasing sequence*

$$0 \leq \nu_1 \leq \nu_2 \leq \ldots$$

*repeated according to multiplicity. Then for $\lambda_k$ in the Laplace spectrum of $M$, we have*

$$\lambda_k \leq \nu_k$$

*for all $k \in \mathbb{N}$.*

*A useful corollary of this is that for a Dirichlet eigenvalue problem on $\Omega \subset M$ and any eigenvalue problem on $M$, we have*

$$\lambda_k(\Omega) \geq \lambda_k(M)$$

*for all $k \in \mathbb{N}$. For example, in Figure 1.3 the first (and indeed $k^{th}$) Dirichlet eigenvalues of $\Omega_1$, $\Omega_2$, and $\Omega_3$ will all be larger than the first (respectively $k^{th}$) eigenvalue of the manifold in which they are contained, that is, the torus with a disk removed, where we impose Dirichlet conditions on the boundary circle.*

Our survey of spectral theory in Section 2.1 will focus on the "Dirichlet Laplacian", that is, the Laplace operator associated with the above problem. Another important type of boundary condition is the Neumann problem.

**Definition 1.8 (Neumann eigenvalue problem)** *For $\partial M \neq \emptyset$, $\overline{M}$ compact and connected, find all $\lambda \in \mathbb{R}$ such that*

$$\Delta \phi - \lambda \phi = 0$$

*has nontrivial solutions $\phi \in C^2(M) \cap C^1(\overline{M})$ that satisfy the boundary condition*

$$\frac{\partial \phi}{\partial \boldsymbol{n}} = 0$$

*on $\partial M$, where $\frac{\partial}{\partial \boldsymbol{n}}$ is the outward unit normal vector field on $\partial M$.*

**Remark 1.9 (Domain monotonicity, Neumann data)** *Let*

$$\Omega_1, \ldots, \Omega_m$$

*be as above. Considering now the Neumann eigenvalue problem, where the normal derivative of the eigenfunction is required to be 0 on $\partial \Omega \cap M$, we arrange all the eigenvalues of these domains into the increasing sequence*

$$0 \leq \mu_1 \leq \mu_2 \leq \ldots$$

*repeated according to multiplicity. Then for $\lambda_k$ in the Laplace spectrum of $M$, we have*

$$\mu_k \leq \lambda_k$$

*for all $k \in \mathbb{N}$.*

Finally, it is possible to have a boundary value problem that requires both Dirichlet and Neumann conditions.

**Definition 1.10 (Mixed eigenvalue problem)** *For $\partial M \neq \emptyset$, $\overline{M}$ compact and connected compact and connected, $N$ an open submanifold of $\partial M$, find all $\lambda \in \mathbb{R}$ such that*

$$\Delta \phi - \lambda \phi = 0$$

*has nontrivial solutions $\phi \in C^2(M) \cap C^1(M \cup N) \cap C^0(\overline{M})$ that satisfy the boundary conditions*

$$\phi = 0 \text{ on } \partial M - N, \quad \frac{\partial \phi}{\partial \boldsymbol{n}} = 0 \text{ on } N.$$

Using the symmetry group of a hyperbolic surface, we can reduce the closed problem on the surface to a simpler Dirichlet, Neumann or mixed problem on a fundamental domain of the group action on the surface. Useful to this approach is the concept of nodal lines and nodal domains, that is, the points on which the eigenfunction is zero and how they partition the domain. This has an obvious link to Dirichlet boundary problems, and we hope to exploit the connection in Section 3.4.

**Definition 1.11 (Nodal domain)** *Let $S$ be a closed Riemann surface and*

$$C^0(S) \ni f : S \to \mathbb{R}.$$

*The nodal set of $f$ is $Z_f = \{f^{-1}(0)\}$, and each connected component of*

$$\overline{S} \setminus \{f^{-1}(0)\}$$

*is a nodal domain of $f$.*

**Theorem 1.12 (Courant's nodal domain theorem, [19])** *Consider a Neumann eigenvalue problem on $M$, $\partial M \neq 0$, $\overline{M}$ compact and connected with eigenvalues*

$$0 = \lambda_0 \leq \lambda_1 \leq \lambda_2 \leq \ldots$$

*corresponding to the complete orthonormal basis of eigenfunctions $\{\phi_0, \phi_1, \ldots\}$ of $L^2(M)$. The number of nodal domains of $\phi_k$ is $\leq k+1$.*

*Now consider a Dirichlet eigenvalue problem on the same region, with eigenvalues*

$$0 < \lambda_1 \leq \lambda_2 \leq \ldots$$

*corresponding to the complete orthonormal basis of eigenfunctions $\{\phi_1, \phi_2, \ldots\}$ of $L^2(M)$. Then the number of nodal domains of $\phi_k$ is $\leq k$.*

A classic paper on the nodal sets of eigenfunctions is [18], where it is proved that for the Laplace operator acting on functions on a Riemannian manifold, the nodal sets of eigenfunctions are smooth sub-manifolds, apart from a closed set of lower dimension. Cheng also generalizes Theorem 1.12 to manifolds of higher dimension, and proves that for a Riemann surface $S$ and $V \in C^\infty(S)$, a solution $\varphi$ of the Schödinger equation

$$(\Delta + V)\varphi = 0$$

with nodal set $Z_\varphi$ has isolated critical points on $Z_\varphi$, and the meeting points of nodal lines form an equiangular system.

A survey of historical results relating to the nodal sets of the Laplacian on surfaces of variable curvature is given in [85]. In more recent years, there have been some interesting results linking the topology of nodal sets to the spectrum of $\Delta$; Otal and Rosas proved an initial result in 2009 [66] that for a compact, orientable, hyperbolic surface $S$ of genus $g$,

$$\lambda_{2g-2}(S) > \frac{1}{4}. \tag{1.3}$$

The value $2g - 2$ is of course equal to $-\chi(S)$, the Euler characteristic of $S$. This is the finite area version of the classical result by Buser [12], that for every hyperbolic surface of genus $g \geq 2$

$$\lambda_{4g-2}(S) > \frac{1}{4}.$$

The proof of Equation (1.3) uses a novel technique in which the topology of the nodal sets of eigenfunctions is considered, and the analysis is done separately for nodal sets that are topologically equivalent to a cross cap, an annulus or a disk. The constant $\frac{1}{4}$ is significant, because in many applications of spectral theory it is important to distinguish the two cases when an eigenvalue $\lambda \geq \frac{1}{4}$ and when $\lambda < \frac{1}{4}$. This comes from the fact the the function

$$x \mapsto \sin\left(x\left(\lambda_n - \frac{1}{4}\right)^{\frac{1}{2}}\right)$$

is bounded in the first case but grows exponentially in the second, as $x \to \infty$ [13]. $\frac{1}{4}$ is well known to be the bottom of the spectrum of the hyperbolic plane. For this reason, eigenvalues in the interval $\left[0, \frac{1}{4}\right]$ have the special title of *small eigenvalues*. For a compact Riemann surface $S$, the small eigenvalues of $\Delta(S)$ influence the asymptotic length distribution of its closed geodesics (see, for example, [21, 73]). In the same paper, Otal and Rosas prove that any finite area hyperbolic surface of genus $g$ with $n$ holes has at most $2g - 2 + n$ small eigenvalues.

Following the techniques of [66], Mondal proved several results that use the topology of nodal sets of Laplace eigenfunctions to give bounds on small eigenvalues of the spectrum that depend on the systole of the surface, that is, the length of the shortest closed geodesic (not necessarily unique). The systole of a Riemann surface $S$ will be denoted $s(S)$. It has a multiplicity corresponding to the number of non-homotopic geodesics of minimal length that cannot be contracted. Starting with the main result (in this area) of his thesis [57], published in [58], he proved that for $S$ closed and hyperbolic,

$$\lambda_{2g-2}(S) > \frac{1}{4} + \epsilon_0(S),$$

where

$$\epsilon_0(S) < \min\left\{\frac{1}{4(g-1)}, \frac{1}{4}\left(\left(\frac{\cosh(\rho_0)}{\sinh(\rho_0)}\right)^2 - 1\right)\right\}$$

and

$$2s(S)\sinh(\rho_0) = \text{Area}(S).$$

More results of a similar nature can be found in papers such as [58–60]. In joint work with Ballmann and Matthiesen [4], he proved that a closed Riemann surface $S$ with negative Euler characteristic has at most $-\chi(S)$ small eigenvalues. This paper generalizes earlier results, including Equation (1.3), by removing the restriction that the surface must be of negatively curved, that is, it is true for arbitrary Riemannian metrics.

## Comparing eigenvalues

In general, we cannot compute eigenvalues of surfaces, or domains within surfaces, explicitly. There are several analytical tools that use the geometry of the surface to create both upper and lower bounds on particular eigenvalues. To bound eigenvalues from below, we can compare them to the eigenvalues of disks having equivalent area; these are known to minimize the eigenvalue provided a certain condition is satisfied, namely the isoperimetric inequality:

**Definition 1.13 (Isoperimetric inequality)** *Let $D$ be an open disk in the hyperbolic plane $\mathbb{H}$ (see Section 2.2 for a discussion on hyperbolic geometry). An open domain $\Omega$ in a hyperbolic surface $S$ is said to satisfy the isoperimetric inequality if*

$$A(\Omega) = A(D) \implies l(\partial\Omega) \geq l(\partial D),$$

*with equality if and only if $\Omega$ is isometric to $D$, where $A(\Omega)$ denotes the area of $\Omega$, and $l(\partial\Omega)$ denotes the length of its boundary.*

Caution must be used with this definition, since there exists a whole class of results called isoperimetric inequalities that relate the area of a domain to its perimeter, and do not take into account physical properties such as torsional rigidity. Isoperimetric inequalities are said to have their origins in ancient Carthage, where Queen Dido was tasked with delimiting the boundary of her new city, such that it would occupy a region that could be encompassed by a single ox-hide. Legend has it that she cut the hide into long thin strips and chose a circular boundary, having conjectured that this would maximize the area of land. More succinctly we have the famous formula

$$L^2 \geq 4\pi A,$$

where a curve $C$ of length $L$ encloses a region of area $A$, and equality holds if and only if $C$ is a circle. We take this form of the inequality from [65], where Osserman gives an interesting historical exposition of significant results following on from this classical theorem.

The term isoperimetric inequality is used interchangeably in the literature for problems of the above type, and those of more relevance to this work: isoperimetric problems, where (often) variational calculus is used to show that in a certain class of domains, a physical property is shown to be extremal for radially symmetric regions. In [67], Payne gives a similar treatise to that of Osserman, covering the extension of isoperimetric inequalities to physical properties of surfaces. Whilst being a thorough survey of key research up to 1967, the fact that the spectral theory of surfaces is still a hot topic means that this reference is good introduction to early work in the field but is far from complete. A slightly more up to date survey article is given by Ashbaugh [1]. This mainly pertains to the Laplacian on Euclidean domains, but as one sees in [15], many of these results can be generalized to different geometries, including hyperbolic.

Of course, here the physical property we are interested in is eigenvalues; in particular, $\lambda_1$. The most important result is the following, which gives a lower bound on the first positive eigenvalue of a domain $\Omega$ satisfying Definition 1.13.

**Theorem 1.14 (Faber-Krahn inequality)** *Let $\Omega$ be an open domain in a hyperbolic surface $S$, and $D$ be a disk in hyperbolic space (more information on such disks is given in Section 2.3.1). If $\Omega$ satisfies the isoperimetric inequality, then*

$$A(\Omega) = A(D) \implies \lambda_1(\Omega) \geq \lambda_1(D),$$

*where $\lambda_1$ is the first Dirichlet eigenvalue.*

To prove the Faber-Krahn inequality, the levels sets of eigenfunctions on an open domain are radially symmetrized so that the first Dirichlet eigenvalue can be compared to that of a geodesic disk with equal area. This theorem can be stated more generally for $n$-dimensional manifolds satisfying an equivalent isoperimetric property; details can be found in [15]. While Theorem 1.14 is a powerful tool for analyzing domains that satisfy Definition 1.13, not every domain that we will study has this property. In particular, domains in surfaces that resemble annuli can have a boundary that is shorter than that of the disk with equal area.

Another tool for creating lower bounds on $\lambda_1$ is Cheeger's inequality. For the cases investigated in this work, it is not as strong as the Faber-Krahn inequality, so we will not actually use it. However, we state it as a key result in the theory of isoperimetric problems.

**Definition 1.15 (Cheeger constant)** *Let $M$ be a non-compact Riemann manifold of dimension $n \geq 2$, possibly with boundary, and possibly having compact closure. The Cheeger constant $\mathfrak{h}(M)$ is defined as the following infimum:*

$$\mathfrak{h}(M) = \inf_{\Omega} \frac{l(\partial \Omega)}{\text{Area}(\Omega)},$$

*where $\Omega$ ranges over all open domains of $S$ with compact closure and smooth boundary and $l(\partial \Omega)$ denotes the length of the boundary of $\Omega$.*

**Theorem 1.16 (Cheeger's inequality, [17])** *For $\Omega \subset S$ as above, we have*

$$\lambda_1(\Omega) \geq \frac{\mathfrak{h}^2(\Omega)}{4}.$$

We will also be interested in creating upper bounds on eigenvalues. In his classic treatise on acoustics [69], Rayleigh gave a variational characterization of the first eigenvalue for the Dirichlet problem. This involves taking the gradient of a function. We stated earlier that the eigenfunctions of the Laplacian on a compact Riemann surface $S$ form an orthonormal basis of $L^2(S)$. Recall that $L^2(S)$ is the space of square integrable functions on $S$, that is, a function $f$ on $S$ is in $L^2(S)$ if

$$\int_S |f|^2 dS \leq \infty.$$

where $dS$ is the Riemannian measure on $S$. $L^2(S)$ is special among the $L^p(S)$ spaces because its inner product

$$\langle f, g \rangle_{L^2} = \int_S f \overline{g} \, dS$$

makes it a Hilbert space. The inner product induces the norm

$$\|f\|^2 = \langle f, f \rangle_{L^2}.$$

Note that the $L^2$ norm gives no indication about how smooth a function is. When we apply the Laplacian, we may therefore uncover huge (or indeed, undefined) second derivatives, even if the $L^2$ norm of a function is small. The Laplacian is unbounded on $L^2(S)$ and operations such as taking the gradient of $f \in L^2(S)$ simply do not make sense. Therefore, we need to formulate a different space of functions before we can introduce Rayleigh's principle. Rather than considering the whole space of square integrable functions, we will work with those functions in $L^2(S)$ for which we can define a reasonable notion of differentiability. First consider the multi-index

$$\alpha = (\alpha_1, \ldots, \alpha_d), \qquad \alpha_j \in \mathbb{N}$$

of order

$$|\alpha| := \sum_{i=1}^d \alpha_i,$$

and for $\psi(x) : \mathbb{R}^d \to \mathbb{C}$ sufficiently smooth, define

$$D^\alpha \psi(x) := \frac{\partial^{|\alpha|} \psi(x)}{\partial x_1^{\alpha_1} \ldots \partial x_d^{\alpha_d}}.$$

For convenience, we define $D^1\psi(x)$ to be the gradient $\nabla \psi(x)$ and $D^2\psi(x)$ to be the Hessian matrix for $\psi(x)$, the trace of which is equal to $\Delta \psi(x)$. We are now ready to define our notion of differentiability.

**Definition 1.17 (Weak derivative)** *Let $\psi$ and $g$ be locally integrable functions on a domain $\Omega \subset S$. We say that $g$ is the $\alpha^{\text{th}}$-weak partial derivative of $\psi$ if for all $\varphi \in C_0^\infty(\Omega)$,*

$$\int_\Omega \varphi(x) g(x) dx = (-1)^{|\alpha|} \int_\Omega D^\alpha \varphi(x) \psi(x) dx.$$

With the weak derivative, we can now define the spaces of functions on which we consider our eigenvalue problems, namely, Sobolev spaces. In Section 2.1, we will see how to extend the Laplacian on these spaces to a self-adjoint operator on $L^2(S)$, but for now, we simply give enough detail to define the subsequent "variational" theorems.

**Definition 1.18 (Sobolev space)** *Let $\Omega$ be as above. The Sobolev spaces, indexed by $p \in [1, \infty)$ and $k \in \mathbb{N}$, are given by*

$$W^{k,p} := \{\psi \in L^p(\Omega) \,|\, D^\alpha \psi \in L^p(\Omega),\, 0 \leq |\alpha| \leq k\},$$

*and are complete with respect to the Sobolev norm*

$$\|\psi\|_{W^{k,p}} := \left( \sum_{0 \leq |\alpha| \leq k} \|D^\alpha \psi\|_{L^p(\Omega)}^p \right)^{\frac{1}{p}}$$

*where $D^\alpha \psi$ is the weak partial derivative. We may also define $W_0^{k,p}(\Omega)$ as the closure of $C_0^\infty(\Omega)$ with respect to the Sobolev norm.*

**Remark 1.19** *Again, we are most interested in the case where $p = 2$, since $W^{k,2}(\Omega)$ is a separable Hilbert space with inner product*

$$\langle \phi, \psi \rangle_{W^{k,2}} = \sum_{0 \leq |\alpha| \leq k} \langle D^\alpha \phi, D^\alpha \psi \rangle_{L^2}.$$

The variational characterization of the Laplacian (and its associated eigenvalue problems) is concerned with the validity of the formula

$$\langle \Delta \phi, f \rangle_{L^2} = D[\phi, f],$$

where the $\phi$ is typically an eigenfunction, and the quantity on the right is the Dirichlet integral

$$D[f, h] := \langle D^1 f, D^1 h \rangle_{L^2}$$

for $f, h \in W^{1,2}(\Omega)$ in the case when we consider a closed or Neumann problem and $f, h \in W_0^{1,2}(\Omega)$ when we consider a Dirichlet problem. $D^1$ denotes the first weak derivative. We can now state the following theorems as in [15].

**Theorem 1.20 (Rayleigh's theorem)** *Let $S$ be a closed and connected Riemann surface, $\Omega \subset S$ be a domain with compact closure and non-empty piecewise smooth boundary. Consider a fixed eigenvalue problem on $\Omega$ with eigenvalues*

$$\lambda_1 \leq \lambda_2 \leq \ldots$$

*repeated according to multiplicity. For any $0 \neq f \in W^{1,2}(\Omega)$ in the case when we consider a closed or Neumann problem, and $0 \neq f \in W_0^{1,2}(\Omega)$ when we consider a Dirichlet problem,*

$$\lambda_1(\Omega) \leq \frac{D[f,f]}{\|f\|^2},$$

*with equality if and only if $f$ eigenfunction corresponding to $\lambda_1$. The object on the right is called the Rayleigh quotient, and will be denoted $\mathcal{R}(f)$.*

We can use Rayleigh's quotient to characterize $\lambda_k$ for any $k$. In particular, we have the following pair of theorems, taken from [6]:

**Theorem 1.21 (Max-min)** *Let $M$ be a smooth, connected, compact manifold, possibly with boundary $\partial M$. For the Neumann eigenvalue problem on $M$, or closed eigenvalue problem where $\partial M = \emptyset$, we have*

$$\lambda_k = \sup_{V_k} \inf\{\mathcal{R}(f) \,|\, f \neq 0,\, f \perp V_k\},$$

*where $V_k$ runs through the $k$-dimensional subspaces of $W^{1,2}(M)$. We consider the $(k-1)$-dimensional subspaces of $W_0^{1,2}(M)$ for the Dirichlet eigenvalue problem on $M$.*

**Theorem 1.22 (Min-max)** *Let $M$ be as above. For the Neumann eigenvalue problem on $M$, or closed eigenvalue problem where $\partial M = \emptyset$, we have*

$$\lambda_k = \inf_{V_k} \sup\{\mathcal{R}(f) \,|\, f \neq 0,\, f \in V_k\},$$

*where $V_k$ runs through the $(k+1)$-dimensional subspaces of $W^{1,2}(M)$. We consider the $k$-dimensional subspaces of $W_0^{1,2}(M)$ for the Dirichlet eigenvalue problem on $M$.*

**Remark 1.23** *In this work, we use Theorem 1.22 in Section 3.3 to produce an upper bound on the second positive eigenvalue of the Bolza surface. We do not actually use its sister theorem at any point, but include it for completion.*

The eigenvalues in the spectrum of the Laplacian on a compact Riemann surface are known to have finite multiplicity; it is a reasonable question to ask what this multiplicity is, and and whether it is possible to bound the multiplicity for a given eigenvalue either from below or above. We will see in Section 2.4 that the multiplicities of a compact

Riemann surface are connected to the irreducible representations of its group of automorphisms. Well established results are that $0 = \lambda_0$ has simple multiplicity for the Laplacian on a closed, connected surface, and for a Schrödinger operator with a smooth potential, the first eigenvalue, whilst not necessarily 0, also has multiplicity one [5, 15].

For a Schrödinger operator $H$ on a surface $S$, Cheng [18] bounded the multiplicity $m_1$ of $\lambda_1$ by a quantity that is quadratic in $\chi(S)$, with the corollary that for the Laplacian on the sphere $\mathbb{S}^2$ with the standard round metric, $m_1(\mathbb{S}^2) = 3$. Besson improved this to a bound that is linear in $\chi(S)$ [7], highlighting that for the projective plane $\mathbb{P}^2$ with round metric, and the flat torus $\mathbb{T}^2$, we have $m_1(\mathbb{P}^2) = 5$ and $m_1(\mathbb{T}^2) = 6$ respectively. Nadirashvili slightly improved the linear bound to

$$m_1(S) = 5 - 2\chi(S)$$

for surfaces $S$ with negative Euler characteristic [62]. Sévannec [75] improved this bound for surfaces for surfaces of negative Euler characteristic, asymptotically by a factor of 2. In particular, he proved the following

**Theorem 1.24** *Let $S$ be a closed surface. If $\chi(S) < 0$, then $m_1(S) < 5 - \chi(S)$.*

**Example 1.25** *Theorem 1.24 already provides upper bounds for the multiplicities of the first eigenvalues of the Bolza surface and Klein quartic. Recall that the Euler characteristic of a surface $S$ may be given in terms of the genus $g$ as*

$$\chi(S) = 2 - 2g.$$

*Hence, for surfaces of genus 2, $m_1(S) < 7$ and for surfaces of genus 3, $m_1(S) < 9$.*

## The Bolza surface

In Chapter 3, we will look at the specific case of the Bolza surface, denoted by $\mathcal{B}$, which can be obtained by identifying opposite edges of a regular octagon in the Poincare disk $\mathbb{D}$. The vertices are at the points

$$2^{-\frac{1}{4}} e^{\frac{\pi i k}{4}}, \qquad k = 0, 1, \ldots, 7,$$

and are joined by geodesic segments of equal length. The group of (orientation preserving) symmetries of $\mathcal{B}$ is $\mathrm{GL}_2(\mathbb{Z}_3)$, which has order 48 (see, for example, [46, 51]). The full group of isometries, including reflections, is the semi-direct product $\mathrm{GL}_2(\mathbb{Z}_3) \rtimes \mathbb{Z}_2$. This is the highest order of a symmetry group for any genus 2 surface [51]. The Bolza surface takes its name from Oskar Bolza, a student of Klein, who investigated algebraic curves; in particular binary sextic equations with linear transformation into themselves [8]. The algebraic curve for the Bolza surface is

$$w^2 = z(z^4 - 1), \tag{1.4}$$

where $w$ and $z$ are two functions.

Figure 1.4: The Bolza surface in $\mathbb{D}$

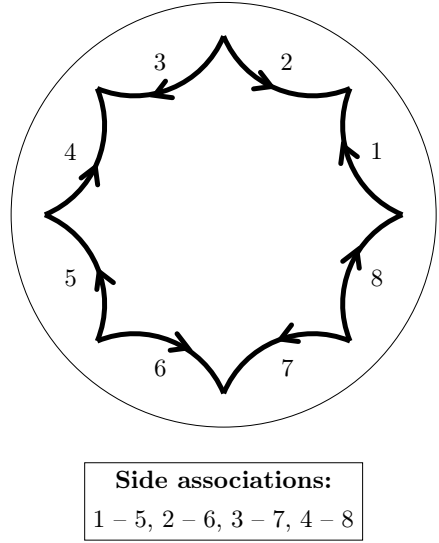

Side associations:
$1 - 5, \ 2 - 6, \ 3 - 7, \ 4 - 8$

The Bolza surface has attracted the attention of spectral theorists, as it is conjectured to be the surface that maximizes the first positive eigenvalue of the Laplace operator among all compact Riemann surfaces of genus two with constant negative curvature. The conjecture is given in [77], and is based on numerical calculations carried out by the authors. This is thought to be related to the fact that the Bolza surface has the largest automorphism group in this class of surfaces. As well as having maximal symmetry in this class, the Bolza surface also maximizes the systole [72], which Jenni [41] showed to be

$$s(\mathcal{B}) = 2 \operatorname{arccosh}\left(1 + \sqrt{2}\right)$$

and have multiplicity 12. These are the four associated edges in Figure 1.4, the four diameters of $\mathbb{D}$ connecting their midpoints, and the curves that connect the midpoint of one edge with the midpoint connecting the adjoining edge. The formula for the systole follows immediately from the side lengths of the (2, 3, 8) triangles that tessellate the surface (see, for example, Figure 3.3). The systole of the Bolza surface was considered in [47], where it was conjectured to be optimal for the systolic ratio

$$SR(S) = \sup_G \frac{s(G)^2}{\operatorname{Area}_G(S)}$$

of CAT(0) metrics $G$ on surfaces $S$ of genus 2 (that is, metrics of negative curvature including those that are singular).

It is the same work by Jenni [41] that inspired the method for the investigation of the Bolza surface in this work. He investigated the class of Riemann surfaces of genus $g$, labelled $F_g$, obtained by associating the opposite sides of a regular $4n$-gon with angles $\frac{\pi}{2g}$. In this class, all the surfaces have one rotational and two reflective symmetries; the Bolza surface stands out because it has an extra rotational symmetry, due to the fact that the regular octagon can be tessellated by equilateral triangles (see Figure 3.1); the curves that make up this tessellation are the 12 geodesics with length equal to the systole of the surface. As well as his results on the Bolza surface, he also proved that for $g > 2$, $F_g$ has $2g$ such curves of length

$$s(F_g) = 4\operatorname{arccosh}\left(\sqrt{2}\cos\left(\frac{\pi}{4g}\right)\right).$$

For $g \geq 20$, he identified the group representation that appears in the first eigenspace of $F_g$, proved that the first eigenvalue has multiplicity 2, and also proved that

$$\lim_{g \to \infty} g \cdot \lambda_1(F_g) = 2.$$

Using the same methodology, he proved similar results on the first eigenvalues of a surface of genus 2 in [42], where the surface varies by changing the length $l(\mu)$ of one of the sides $\mu$ that make up its fundamental domain $F_\mu$. In particular,

$$\lim_{l(\mu) \to 0} \frac{l(\mu)}{\lambda_1(F_\mu)} = \pi^2.$$

This surface does not have any obvious maximal properties.

The Bolza surface has been investigated in the context of eigenvalue maximization in [39]. Here the authors look at the surface as a ramified double cover of the sphere where

$$\Pi : \mathcal{P} \to \mathbb{S}^2$$

is a degree 2 branched cover of the sphere, and $\mathcal{P}$ is a surface of genus 2 with the conformal structure of the Bolza surface, that is,

$$\mathcal{P} := \left\{(z, w) \in \mathbb{C} \,|\, w^2 = F(z) := z\frac{(z-1)(z-i)}{(z+1)(z+i)}\right\}.$$

It has been shown by Yang and Yau [84] that for non-constant holomorphic maps of degree $d$ from a closed, orientable Riemann surface $S$ to $\mathbb{S}^2$, the first positive eigenvalue $\lambda_1$ of $S$ satisfies the following inequality

$$\lambda_1(S)\operatorname{Area}(S) \leq 8\pi d. \tag{1.5}$$

This built on Hersch's result [34] that for any metric $g$ on $\mathbb{S}^2$,

$$\lambda_1(g)\text{Area}_g(\mathbb{S}^2) \leq 8\pi,$$

where $\lambda_1(g)$ and $\text{Area}_g$ refer to the first eigenvalue and area with respect to the metric $g$. It was extended to non-orientable surfaces by Li and Yau in [54], where they showed that for and metric $g$ on the real projective plane $\mathbb{R}P^2$

$$\lambda_1(g)\text{Area}_g(\mathbb{R}P^2) \leq 12\pi.$$

Similar inequalities hold for an orientable surface $S$ of genus $\gamma$ immersed into $\mathbb{R}^{2+p}$, that is, $M = X(\bar{M})$ where $\bar{M}$ is a compact orientable surface and $X: \bar{M} \to \mathbb{R}^{2+p}$ is an immersion. In this case, Ilias and El Soufi [22] improved the bound of Yang and Yau to

$$\lambda_1(S)\text{Area}(S) \leq 8\pi \left\lfloor \frac{\gamma+3}{2} \right\rfloor, \tag{1.6}$$

where $\lfloor \cdot \rfloor$ is the floor function. For example, the Bolza surface has area $4\pi$ by the Gauss-Bonnet theorem, so we may use Equation (1.6) to state

$$\lambda_1(\mathcal{B}) \leq 4,$$

when we view the Bolza surface as an immersed surface in $\mathbb{R}^{2+p}$.

For surfaces of genus 2, we can take $d$ in Equation (1.5) to be 2; further information can be found in [29], which also gives background information on many of the concepts in this section relating to algebraic geometry, that will not reappear in the main body of this work. The Bolza surface has 6 ramification points, by the Riemann-Hurwitz formula. These are called the Weierstrass points of the surface and they occur at

$$\{0, \infty, \pm 1, \pm i\}.$$

Note that 0 and $\infty$ are fixed points of Equation (1.4), and 0, $\pm 1$, and $\pm i$ are its roots when $w = 0$. These points are singularities of the metric $g_0$ on $\mathcal{P}$ that is inherited from $\mathbb{S}^2$ by pulling back the round metric. In [39], the following conjecture is made

**Conjecture 1.26** *The metric $g_0$ on $\mathcal{P}$ gives equality in Equation* (1.5), *that is*

$$\lambda_1(\mathcal{P}, g_0)\text{Area}((\mathcal{P}, g_0)) = 16\pi.$$

As part of their investigation into this conjecture, the authors reduce a spectral problem on $\mathbb{S}^2$ to a mixed boundary problem on the semicircle. Of particular interest is the fact that the two problems in Figure 1.5 are isospectral under an exchange of Dirichlet and Neumann conditions. This phenomenon, known as "Dirichlet-Neumann swap isospectrality", is investigated further in [40], where the authors explain how to construct an isospectral domain. It was encountered in this work during the investigation of the 4-dimensional irreducible representations of the Bolza surface (see Section 3.4), where numerical results indicate that the hyperbolic pentagon with mixed boundary conditions, shown in Figure 1.6, also has the property of swap isospectrality. The relevant finite element code is included in Appendix C.

Figure 1.5: Boundary value problems on a semicircle, [39] (dashed lines indicate a Neumann boundary condition, and solid, Dirichlet)

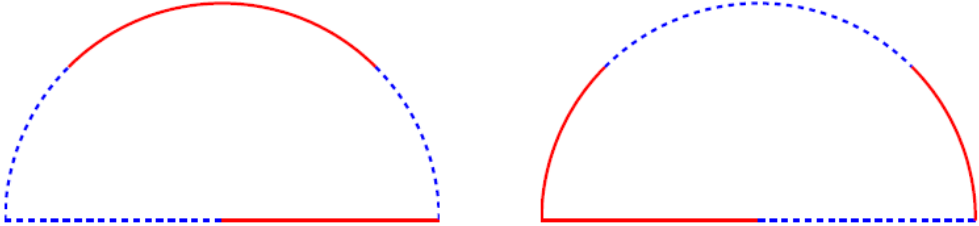

Figure 1.6: Boundary value problems on a pentagon

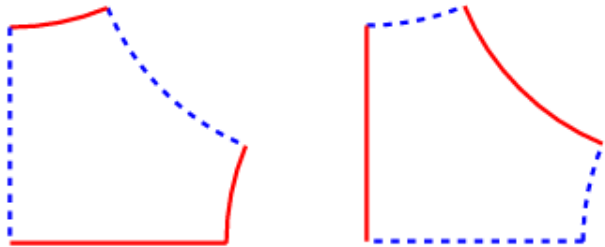

In the above case, the conformal structure of the Bolza surface was investigated as a double cover of the sphere $\mathbb{S}^2$, however, genus 2 surfaces are also ramified double covers of the Riemann sphere [24]. The double cover consists of two sheets, and there is a function, $J$, called the *hyperelliptic involution* that switches between the sheets. In the case of the Bolza surface, the involution is the centre of the full isometry group (see Proposition 3.5). The hyperelliptic involution leaves geodesics invariant, and its fixed points are the Weierstrass points of the surface [31]. For the Bolza surface, we have seen that the Weierstrass points are $\{0, \pm 1, \pm i, \infty\}$.

Suppose we start in the Riemann sphere and draw the geodesics that connect our Weierstrass points to get a cross (see Figure 1.7). The line from 1 to $\infty$ is dashed; the complement of the cross joining the other five points is simply connected, but this fails to be the case when we add $\infty$. Thus, $\infty$ is a branch point when we lift to the cover. When we lift, angles are preserved, except at Weierstrass points where they are halved. On the picture on the left of Figure 1.7, imagine taking the following journey: start at $\infty$, walk to 1, carry on to walk to 0, turn by $\frac{\pi}{2}$, walk to $i$, turn right around and walk back to 0. Carry on moving around the cross until arriving back at 1, and then carry on to $\infty$.

When we lift to the surface, the angles are halved at Weierstrass points, so the same journey becomes: start at $\infty$, walk to 1, turn by $\frac{\pi}{2}$, walk to 0, turn by $\frac{\pi}{4}$, walk to $i$, carry on and walk to 0. Continuing the journey takes us around the upper half on the octagon on the right hand side of Figure 1.7. This is one sheet of the cover; the other is its conjugate. They are glued together along the dotted line.

Figure 1.7: Creating a double cover of the Riemann sphere

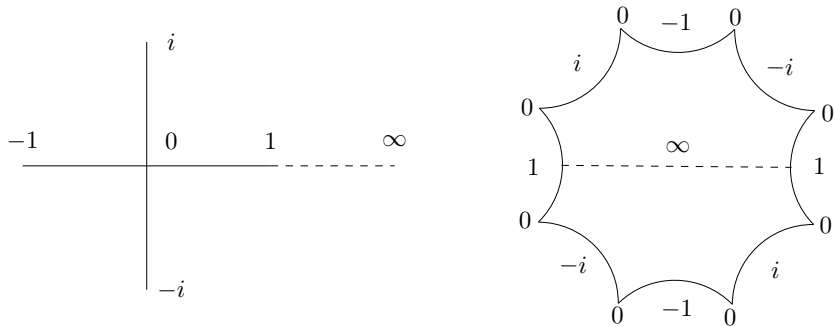

We will see in Section 3.4 that an eigenfunction in the first (non-trivial) eigenspace has as its nodal lines one line of a $(4, 4, 4)$ triangle, and its rotations around the centre of the fundamental polygon (the notation for such a triangle is explained in Section 2.2.3). When we take the quotient of the Bolza surface by the hyperelliptic involution $J$, the centre points of the 16 $(4, 4, 4)$ triangles will map to 8 points on the Riemann sphere. These 8 points give a cubical subdivision of the quotient surface that is dual to the subdivision obtained from the 6 ramification (or Weierstrass) points that induce an octahedral symmetry. Using the numerical results in Appendix C, we can plot the values of the three eigenfunctions corresponding to $\lambda_1$ in $\mathbb{R}^3$; this is shown in Figure 1.8. One can clearly see the octahedral symmetry, as well as 5 of the 6 singular points that come from the ramifications (the other is at the rear). It appears that the eigenfunctions corresponding to $\lambda_1$ provide an immersion into $\mathbb{R}^3$; for a symmetric image such as the one in Figure 1.8, one should take the eigenfunctions to be orthonormal.

The spectral theory of the Bolza surface has been studied in the physics literature as the quantum mechanical Hadamard-Gutzwiller model for quantum chaos [3, 30, 63].

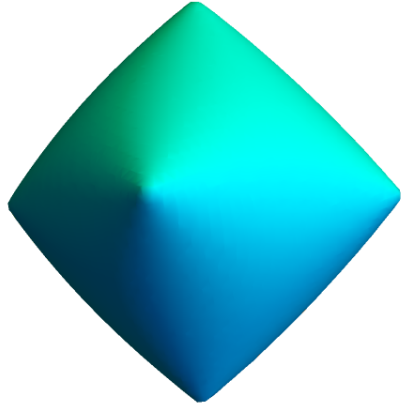

Figure 1.8: Plotting the values of the functions in the first eigenspace in $\mathbb{R}^3$

This is essentially billiards on a hyperbolic surface. The advantage of using a hyperbolic surface over a Euclidean model is that Selberg's trace formula is exact, as opposed to the semiclassical approximations given by the Euclidean case. The Bolza surface is one of the simplest surfaces on which to study this dynamical system, due to its high symmetry, low genus, and the fact that much information is already known about it from the context of algebraic geometry. It can be described as the quotient of the Poincaré disk by the so called "octagon group" (see Section 2.2 for more details on quotient surfaces). This is the Fuchsian group with four generators [63]

$$g_k = \begin{pmatrix} 1+\sqrt{2} & (2+\sqrt{2})\alpha \exp\left(\frac{ik\pi}{4}\right) \\ (2+\sqrt{2})\alpha \exp\left(\frac{-ik\pi}{4}\right) & 1+\sqrt{2} \end{pmatrix},$$

where

$$\alpha := \sqrt{\sqrt{2}-1}, \qquad k = 0, \ldots, 3,$$

and the inverses are given by conjugation with

$$\begin{pmatrix} \exp\left(\frac{i\pi}{2}\right) & 0 \\ 0 & \exp\left(\frac{-i\pi}{2}\right) \end{pmatrix}.$$

They satisfy the relation

$$g_0 g_1^{-1} g_2 g_3^{-1} g_0^{-1} g_1 g_2^{-1} g_3 = I_2,$$

where $I_2$ is the $2 \times 2$ identity matrix. The quotient of the Poincaré disk by this group tessellates the unit disk with regular octagons.

The two theorems we prove in Chapter 3 are to do with the multiplicities of low-lying eigenvalues in the spectrum of $\mathcal{B}$. In particular

**Theorem 1.27** *The dimension of the first eigenspace of $\Delta(\mathcal{B})$ is 3.*

This result was stated by Jenni in [42] after being proved in his thesis [41]. Upon careful inspection, a gap was found in this proof, relating to the analysis of 3 and 4 dimensional irreducible representations. In this work, we give a proof, utilizing the entire automorphism group as opposed to a smaller subgroup. We give the associated irreducible representation explicitly. The proof relies on a conjecture, for which numerical evidence is given in Appendix C. We also suggest how this conjecture may be proved in future work. We also prove

**Theorem 1.28** *The dimension of the second eigenspace of $\Delta(\mathcal{B})$ is 4.*

**Remark 1.29** *Both of these theorems are consistent with the numerical calculations of eigenvalues listed in Table C.1. We observed in the numerics that there are 3 eigenvalues roughly equal to 3.84 and 4 eigenvalues roughly equal to 5.35.*

## The Klein quartic

The Klein quartic, also called Klein's curve, has the largest automorphism group of all compact Riemann surfaces of genus 3. Klein studied the quotient of the hyperbolic plane by Fuchsian groups, that is, discrete subgroups of $\mathrm{PSL}_2(\mathbb{R})$. He found the group

$$\Gamma(7) = \mathrm{PSL}_2(\mathbb{Z}_7)$$

to be of particular interest, and highlighted some of its special properties [48]. Since then, it has been studied extensively both as a Riemann surface and as an algebraic curve; the best introduction is the book "The Eightfold Way" [53], written to accompany the installation of a statue of the surface at the MSRI, Berkeley. The book covers the geometry of the surface, a number theoretical approach to the surface, and many topics in algebra that highlight the extensive interest in the curve. So far, however, the spectral theory of the surface remains relatively untouched; the large symmetry group suggests that surface may maximize $\lambda_1$ among genus 3 compact Riemann surfaces.

We obtain the surface by associating sides $2n+1$ and $2n+6 \mod (14)$, $n = 0, \ldots, 6$ in Figure 1.9. It is possible to tessellate this fundamental polygon with 336 triangles (see Figure 4.2). This corresponds to the full isometry group of the surface, $\mathrm{PGL}_2(\mathbb{Z}_7)$, which has order 336. The group of orientation preserving isometries is half the size of this, that is, 168. This means that the Klein quartic not only has the largest automorphism group for genus 3, it is also optimal, by the following classical theorem.

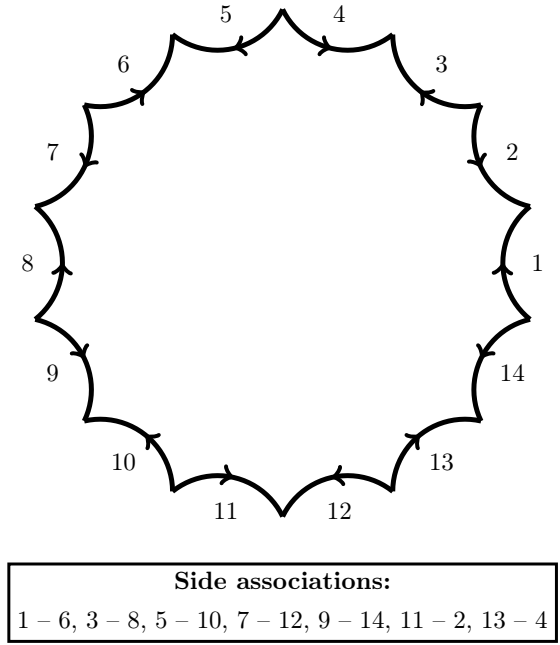

Figure 1.9: The Klein quartic

Side associations:
1 – 6, 3 – 8, 5 – 10, 7 – 12, 9 – 14, 11 – 2, 13 – 4

**Theorem 1.30 (Hurwitz's automorphism theorem, [36])** *The number of automorphisms (not including reflections) of a compact Riemann surface of genus $g > 1$ cannot exceed $84(g - 1)$.*

We see that for genus 3, this bound is 168. The symmetry group of the Klein quartic is therefore maximal. Any surface that achieves this bound is referred to as a Hurwitz surface. A Hurwitz surface does not occur for every single genus $g > 1$; the next Hurwitz surface is the Macbeath surface of genus 7, which has 504 orientation preserving isometries [55]. The group is isomorphic to the simple group $\mathrm{PSL}_2(\mathbb{Z}_8)$ [55] and is further investigated in [83]. This surface is also known as Fricke's surface, as Fricke investigated the corresponding isometry group in [25]. The fundamental polygon is shown in [25], along with the relevant side associations.

The Klein quartic, as its name suggests, can be obtained as the locus of a projective quartic curve, namely

$$x^3y + y^3z + z^3x = 0. \tag{1.7}$$

The study of this equation accounts for much of the interest in Klein's curve in number theory. Stating the key results would involve introducing many terms that are irrelevant to the rest of this work - the interested reader can find a thorough exposition in [23] instead.

The Fermat quartic is a similar surface that also comes from a projective equation; in fact there is a whole class of Fermat surfaces that are obtained as the locus of

$$x^k + y^k + z^k = 0, \qquad k \in \mathbb{N}, \quad k > 2.$$

For the genus 3 version, we set $k = 4$. The fundamental polygon for the Fermat quartic is a 16-gon, shown in Figure D.2 along with the correct side identifications. Similar to the Bolza surface, the Fermat quartic can be tessellated by $(2, 3, 8)$ triangles; as it is a genus 3 surface, there are 192 of them. The Fermat quartic therefore has an automorphism group of order 192. A presentation of the group of orientation preserving elements is given in [52]. We can extend this by adding a (non-commutative) $\mathbb{Z}_2$ action to get the following presentation of the whole group:

$$\langle a, b, c \, | \, a^8 = b^3 = c^2 = (ab)^2 = (a^2b^2)^3 = (a^4b^2)^3 = acac = bcbc = e \rangle,$$

where $e$ is the identity element. In Remark D.1 of Appendix D we will comment further on this group, and how it relates to numerical computations of the eigenvalues of the Fermat quartic.

In Chapter 4, we introduce the notion of a pants decomposition and give the Fenchel-Nielsen parameters of the Klein quartic. We show how to calculate a closed formula for its systole, and that the multiplicity of the systole is 21 (this is not new information, see [72]). We finally show that its first eigenspace is either 6, 7, or 8 dimensional. Of course, this is not a complete result as in the case of the Bolza surface, however it uses the same mechanisms, and a proof along the same lines as that of Theorem 3.9 is theoretically possible. Based on our findings in Table D.1, we can at least conjecture that the multiplicity of the first positive eigenvalue of the Klein quartic is 8.

# Chapter 2

# Background

## 2.1 Spectral theory

The Laplace-Beltrami operator $\Delta$, on a Riemannian manifold $M$, is given in local coordinates by the formula

$$\Delta = -\sum_{i,\,k=1}^{n} \frac{1}{\sqrt{g}} \frac{\partial}{\partial x_i} \sqrt{g} g^{ik} \frac{\partial}{\partial x_k},$$

where $g_{ik}$ is the metric tensor, $g$ its determinant, and $g^{ik}$ are the components of the dual metric on the cotangent bundle of $M$.

There exists an orthonormal basis $\{\phi_i \,|\, i \in \mathbb{N}_0\}$ in $L^2(M)$ consisting of eigenfunctions

$$\Delta \phi_i = \lambda_i \phi_i \tag{2.1}$$

such that

$$0 = \lambda_0 < \lambda_1 \leq \lambda_2 \leq \ldots \tag{2.2}$$

form an ascending sequence with $\infty$ as the only accumulation point. We understand Equation (2.1) in terms of the weak derivative.

The sequence in Equation (2.2) is called the spectrum of $\Delta$, and the individual $\lambda_i$s are referred to as eigenvalues. In this work, $M$ will be either a compact Riemann surface, or a subspace of a compact Riemann surface on which we consider a Neumann, Dirichlet of mixed boundary problem. Recall the definition of these eigenvalue problems from Chapter 1. In particular, note that we take $\Omega$ to have compact closure and piecewise smooth boundary. The property of compactness, along with the smoothness condition where we consider a bounded $\Omega$, ensure that the spectrum is discrete.

This result is entirely non-trivial; to see where it comes from, we will consider the example of the Dirichlet Laplacian $\Delta_D$ on a domain $\Omega$ in a Riemann surface (recall Definition 1.6). The approach we take is via sesquilinear forms, which will be defined shortly. The spectral theory of $\Delta$ can also be explained through heat kernel asymptotics, details of which are given in [13] and [28]. We start by introducing some classical theory of unbounded operators, following [50] and [64]. Another good reference for the theory of unbounded operators is [78].

Let $\mathcal{H}$ be a separable, complex Hilbert space. Let $A$ be a linear operator in $\mathcal{H}$. A bounded linear operator $A$ can be defined for all $u \in \mathcal{H}$; it satisfies the inequality

$$\|Au\| \leq C\|u\|,$$

where the least possible value of $C$ defines the norm $\|A\|$. Conversely, differential operators, in particular $\Delta$, are unbounded and as a result do not have this property. They are however defined for some $u \in \mathcal{H}$, so to properly define an unbounded operator $A$, we must also give its domain $D(A) \subset \mathcal{H}$. Often we will write $(A, D(A))$. The graph of the operator is used to define $A$ as

$$\mathrm{gr}(A) := \{(u, Au) \,|\, u \in D(A)\} \subset \mathcal{H} \times \mathcal{H},$$

and comes with a norm (called the graph norm):

$$\|u\|_A = \|u\| + \|Au\|.$$

**Definition 2.1 (Densely defined operator)** *We say that $A$ is densely defined if $D(A)$ is dense in $\mathcal{H}$.*

**Example 2.2** *$\Delta^\Omega$, that is, the Laplacian acting on $\Omega$, is defined on $C_0^\infty(\Omega)$, and satisfies Dirichlet boundary conditions trivially due to the property of compact support. Note that it is also a non-negative operator since for $\phi \in C_0^\infty(\Omega)$, an integration by parts shows that*

$$(\phi, \Delta\phi) = \|\nabla\phi\|^2 \geq 0.$$

*$C_0^\infty(\Omega)$ is dense in $L^2(\Omega)$, so $(\Delta^\Omega, C_0^\infty(\Omega))$ is densely defined on $L^2(\Omega)$. Using Theorem 2.13, we will create an extension of $P_D^\Omega = (\Delta^\Omega, C_0^\infty(\Omega))$ to define $\Delta_D^\Omega$ on $L^2(\Omega)$.*

Being densely defined is an essential property for defining the adjoint of an unbounded operator. For our unbounded, densely defined operator $A$, we want to construct an extension of $A$ to a continuous operator on all of $\mathcal{H}$, that is, a linear operator $\tilde{A}$ such that $D(A) \subset D(\tilde{A})$ and for all $u \in D(A)$, $Au = \tilde{A}u$. In the case of bounded operators, we have the bounded linear transformation (or so-called 'BLT') theorem that allows us to do this (see, for example, [70] for details of this and other named but unacknowledged theorems in this section). On the other hand, for unbounded operators

$$\sup_{0 \neq u \in D(A)} \frac{\|Au\|}{\|u\|} = \infty,$$

so a continuous extension is not possible. The problem then, is to create a natural extension of $D(A)$ in a way that allows us to replicate, or at least approximate, the spectral theory of bounded operators. First, we need to introduce the concept of closability.

**Definition 2.3 (Closed operator)** *A is called closed if for all $u, v \in \mathcal{H}$ and a sequence $\{u_n\}_{n \in \mathbb{N}} \subset D(A)$, the conditions*

$$u_n \to u \quad \text{and} \quad Au_n \to v$$

*as $n \to \infty$ imply*

$$u \in D(A) \quad \text{and} \quad Au = v.$$

*In other words A is closed if $D(A)$ is a complete vector space with respect to the graph norm.*

Any bounded operator is closed by the closed graph theorem. Unbounded operators are not necessarily closed, but we have a similar notion: $A$ is called closable if it has a closed extension, that is, if its graph remains the graph of an operator $\bar{A}$ when its closure $\overline{\text{gr}(A)}$ is taken. Once again, question of domains makes defining the adjoint harder for unbounded operators; for bounded operators, it follows from the Riesz representation theorem.

**Definition 2.4 (Adjoint of an unbounded operator)** *Let A be a densely defined operator in $\mathcal{H}$. Its adjoint $A^*$ is uniquely defined as follows*

$$D(A^*) := \{u \in \mathcal{H} \,|\, \exists u^* \in \mathcal{H} \text{ s.t. } \forall v \in D(A), (u, Av) = (u^*, v)\},$$
$$A^* u := u^*$$

.

**Proposition 2.5** *$A^*$ is always closed, regardless of whether or not A is.*

To see this, we can look at an alternative definition of the adjoint, taking advantage of the fact that we can describe an operator by its graph. For $A^*$,

$$\text{gr}(A^*) = R(\text{gr}(A)^\perp),$$

where $R$ acts like a rotation on $\mathcal{H} \times \mathcal{H}$, namely $R(u_1, u_2) = (u_2, -u_1)$.

**Definition 2.6 (Symmetric operator)** *A is symmetric if it is densely defined and its adjoint $A^*$ is an extension of A.*

An easy consequence of Proposition 2.5 is that if $A$ is symmetric, then it is closable. The preceding definitions have laid the groundwork for us to define self-adjointness for unbounded operators. These will be crucial for developing a functional calculus of operators $f(A)$, $f : \mathbb{R} \to \mathbb{C}$, whereby we will be able to use function analytic methods in our spectral analysis of operators.

**Definition 2.7 (Self-adjoint operator)** *$A$ is self-adjoint if it is densely defined, $A^* = A$ and $D(A) = D(A^*)$. A symmetric operator $A$ is essentially self-adjoint if its closure $\bar{A}$ is self-adjoint.*

We now consider sesquilinear forms. The aim is to establish a one to one correspondence between these forms and self-adjoint operators; this will allow us to create self-adjoint extensions of densely defined, symmetric, non-negative operators.

**Definition 2.8 (Sesquilinear form)** *A sesquilinear form in $\mathcal{H}$ is a map*

$$q : D(q) \times D(q) \to \mathbb{C},$$

*where $D(q) \subset \mathcal{H}$, such that*

$$(u, v) \to q(u, v)$$

*is linear in the second entry and semi-linear in the first. The mapping*

$$D(q) \to \mathbb{C}$$
$$u \mapsto q[u] := q(u, u)$$

*is called the quadratic form associated to $q$. $q$ is semi-bounded below if*

$$q[u] \geq \alpha \|u\|^2$$

*for some $\alpha \in \mathbb{R}$. $q$ is called elliptic if*

$$|q[u]| \geq c\|u\|^2$$

*for some $c > 0$. Finally, $q$ is closed if it is semi-bounded below and $D(q)$ is a complete Hilbert space with respect to the norm*

$$\|u\|_q^2 = q(u, u) + (|\alpha| + 1)\|u\|^2.$$

**Example 2.9** *For our study of the Dirichlet Laplacian, we will introduce a sesquilinear form $Q_D^\Omega$ such that*

$$D(Q_D^\Omega) := D(P_D^\Omega)$$

*and for $\phi \in W_0^{1,2}$,*

$$Q_D^\Omega[\phi] := (\phi, \Delta\phi) = \|\nabla\phi\|^2.$$

**Definition 2.10 (Adjoint of a sesquilinear form)** *The adjoint $q^*$ of $q$ is defined as*

$$D(q^*) := D(q),$$
$$q^*(\phi, \psi) := \overline{q(\psi, \phi)}.$$

*$q$ is symmetric if $q^* = q$.*

Now let $q$ be a symmetric, densely defined sesquilinear form that is bounded from below. As before, $q$ is closable if it has a closed extension, that is, if for all $\{\psi_n\}_{n \in \mathbb{N}} \subset D(q)$,

$$D(q) \ni \psi_n \xrightarrow{n \to \infty} 0$$
$$q[\psi_n - \psi_k] \xrightarrow{n, k \to \infty} 0$$

imply $q[\psi_n] \xrightarrow{n \to \infty} 0$, then $q$ has closure $\bar{q}$ defined as

$$D(\bar{q}) := \left\{ \psi \in \mathcal{H} \mid \exists \{\psi_n\} \subset D(q) \text{ s.t. if } \psi_n \xrightarrow{n \to \infty} \psi, \ q[\psi_n - \psi_k] \xrightarrow{n, k \to \infty} 0 \right\},$$
$$\bar{q}[\psi] := \lim_{n \to \infty} q[\psi_n].$$

**Example 2.11** *The form $Q_D^\Omega$ is not closed on $C_0^\infty$, but note that the norm*

$$\|\phi\|_{Q_D^\Omega}^2 = Q_D^\Omega[\phi] + (|\alpha| + 1)\|u\|^2$$

*is the norm on the Sobolev space $W_0^{1,2}(\Omega)$ (compare Definition 1.18). The domain of the closure $\tilde{Q}_D^\Omega = \overline{Q_D^\Omega}$ is therefore $W_0^{1,2}(\Omega)$, and for $\phi \in W_0^{1,2}(\Omega)$,*

$$\tilde{Q}_D^\Omega[\phi] = \|\nabla \phi\|^2,$$

*where we interpret $\nabla \phi$ in a distributional sense.*

From the sesquilinear form $q$, we wish to induce an operator $A_q$ that is self-adjoint and bounded from below. The correspondence is as follows. Suppose we have $A_q$ as above. Then

$$\dot{q}[\psi] := (\psi, A_q \psi), \quad D(\dot{q}) := D(A_q)$$

is densely defined, symmetric and bounded from below. It is closable, so we take the closure $q := \bar{\dot{q}}$ as our desired form. On the other hand, suppose we have $q$ as above. Then we will use the following to define the corresponding operator:

**Theorem 2.12 (Representation theorem, [45])** *Let $q$ be densely defined, symmetric, closed, and bounded from below in $\mathcal{H}$. Then*

$$D(A_q) := \{\psi \in D(q) \mid \exists \eta \in \mathcal{H} \text{ s.t. } \forall \phi \in D(q), \ q(\phi, \psi) = (\phi, \eta)\},$$
$$A_q \psi := \eta$$

*is self-adjoint and bounded from below. $D(q)$ is called the form domain of the operator $A_q$.*

We now introduce a special construction that produces the self-adjoint operator with the smallest form domain, that is, the domain of $q$ is contained in the form domain of any other self-adjoint extension $\tilde{A}_q$ of $A_q$. We will use this to define the Dirichlet Laplacian as the self-adjoint extension of a minimal operator that has the smooth, compactly supported functions as its domain.

**Theorem 2.13 (Friedrichs extension)** *Let $A$ be a symmetric operator, bounded from below. Let $q_A$ be the sesquilinear form such that*

$$q_A[\psi] := (\psi, A\psi),$$

*and*

$$D(q_A) := D(A)$$

*is densely defined, symmetric and bounded from below. If $q_A$ is closable, we take the closure*

$$\tilde{q}_A := \overline{q}_A.$$

*Finally, we can associate to $\tilde{q}_A$ the self-adjoint linear operator $A_{\tilde{q}}$ via Theorem 2.12.*

**Example 2.14** *Using Theorem 2.12 and Theorem 2.13, we are now ready to define $\Delta_D^\Omega$ on $L^2(\Omega)$ as a self-adjoint extension of $P_D^\Omega$ via the quadratic form $\tilde{Q}_D^\Omega$. The definition in the Representation Theorem gives us*

$$D(\Delta_D^\Omega) = \left\{ \psi \in W_0^{1,2}(\Omega) \,|\, \exists \eta \in L^2(\Omega) \text{ s.t. } \forall \phi \in W_0^{1,2}(\Omega), (\nabla \phi, \nabla \psi) = (\phi, \eta) \right\}, \quad (2.3)$$
$$\Delta_D^\Omega \psi = \eta.$$

We may abbreviate this definition slightly, using the knowledge that $C_0^\infty(\Omega)$ is dense in $W_0^{1,2}(\Omega)$, and that the condition in 2.3 defines the distributional Laplacian, to state

$$D(\Delta_D^\Omega) = \left\{ \psi \in W_0^{1,2}(\Omega) \,|\, \Delta \psi \in L^2(\Omega) \right\},$$
$$\Delta_D^\Omega \psi = \Delta \psi.$$

We define the spectrum of an unbounded operator via the inverse, $(A - z)^{-1}$ of a self-adjoint operator in $\mathcal{H}$. It exists for all values of $z$ that are not eigenvalues of $A$.

**Definition 2.15 (Resolvent)** *The resolvent set of $A$,*

$$\rho(A) := \{ z \in \mathbb{C} \,|\, (z - A) : D(A) \to \mathcal{H} \text{ is bijective and } (z - A)^{-1} \text{ is bounded} \}$$

*is an open subset of $\mathbb{C}$. The resolvent operator*

$$\rho(z) := (z - A)^{-1}$$

*depends analytically on $z$ inside $\rho(A)$.*

**Example 2.16** *The resolvent operator for the Dirichlet Laplacian, $\rho(-1) = (\Delta_D^\Omega + 1)^{-1}$ is a bounded map*

$$\rho(-1) : L^2(\Omega) \to W_0^{1,2}(\Omega).$$

**Definition 2.17 (Spectrum)** *The spectrum of $A$ is the complement*

$$\sigma(A) := \mathbb{C} \setminus \rho(A).$$

**Remark 2.18** *We are interested in the discrete spectrum $\sigma_{dis}$, which is the set of isolated eigenvalues of finite multiplicity, that is*

$$\sigma_{dis}(A) = \{\lambda \in \mathbb{C} \,|\, \exists \psi \in D(A),\, \|\psi\| = 1,\, A\psi = \lambda\psi\}.$$

*If $\lambda \in \sigma_{dis}(A)$, then its multiplicity is given by $\dim(\ker(A - \lambda)) < \infty$.*

Now that we have defined the resolvent operator, we just need to show that it is compact, and then we can apply the spectral theorem for compact, self-adjoint operators. First, we define what it means to be a compact operator.

**Definition 2.19 (Compact operator)** *A bounded operator $A$ is compact if for a bounded sequence $(u_n)_{n \in \mathbb{N}}$, $(Au_n)_n$ contains a converging sub-sequence.*

There are alternative criteria for a bounded operator to be compact, for example, we have

**Proposition 2.20** *Let $A$ be bounded operator on a Hilbert space $\mathcal{H}$. Then the following are equivalent [64]:*

*(i) $A$ is compact;*

*(ii) $A$ is the norm limit of a sequence of finite rank operators;*

*(iii) for any weakly convergent sequence $u_n \rightharpoonup u$, $Au_n \to u$ strongly;*

*(iv) for any finite orthonormal system $(e_n)$, $Ae_n \to 0$.*

**Theorem 2.21 (Rellich-Kondrachov embedding theorem)** *Let $\Omega$ be a bounded domain of a compact Riemann surface $S$. The the embedding*

$$W^{1,2}(\Omega) \subset L^2(\Omega)$$

*is compact.*

**Example 2.22** *We have seen that $\rho(-1)$ is bounded from $L^2(\Omega)$ to $W_0^{1,2}(\Omega)$. By the Rellich-Kondrachov embedding theorem, $W_0^{1,2}(\Omega)$ has a compact embedding into $L^2(\Omega)$, since $\Omega$ is compact, as is the surface $S$ for which it is a subset.*

The following theorem guarantees that any complex number in the spectrum of a compact operator is an eigenvalue.

**Theorem 2.23 (Fredholm alternative)** *Let $A$ be a compact operator. For $z_0 \neq 0$,*

*(i) either $(z_0 - A)^{-1}$ exists and is a bijection $\mathcal{H} \to \mathcal{H}$,*

*(ii) or $z_0$ is an eigenvalue of finite multiplicity.*

Therefore, the spectrum of the compact resolvent $\rho(-1)$ consists of $\sigma_{dis} \cup \{0\}$, and its eigenvalues have finite multiplicity. Moreover the only accumulation point of its spectrum is zero.

The upshot of Theorem 2.23 is that we can "invert" the property of a compact resolvent $\rho(z)$ of $A$ to get information about what we are really interested in - the spectrum of $A$. In particular, $A$ has a purely discrete spectrum that has $\infty$ as its only accumulation point. To see this, note that for the eigensystem $(\mu_n, u_n)$ of $\rho_z$,

$$(z_0 - A)^{-1} u_n = \mu_n u_n \iff (1 - z_0 \mu_n) u_n = -\mu_n A u_n.$$

Rearranging the right hand sides tells us that the eigenvalues of $A$ are

$$\lambda_n = \frac{z_0 \mu_n - 1}{\mu_n}.$$

One can see that an accumulation of the resolvent $\mu_n \to 0$ will correspond to an accumulation at infinity for the spectrum of $A$. Returning to $\Delta_D^\Omega$, we have seen that it has a compact resolvent, and therefore a discrete spectrum. We finish the section by stating the general spectral theorem for self adjoint operators.

**Theorem 2.24** *If $A$ is compact and self-adjoint, then it admits an orthonormal eigenbasis $(u_n, \lambda_n)$ with $\lambda_n \to 0$.*

## 2.2 Hyperbolic geometry

We now want to look at the Laplacian in two specific models of hyperbolic geometry. First, we will consider the Poincaré disk model. This is the setting through which we view the Riemann surfaces. Next, we will consider the Laplacian in geodesic polar coordinates on the Poincaré upper half plane, $\mathbb{H}$. Harmonic analysis on $\mathbb{H}$ will give us a method to calculate the eigenvalues of hyperbolic disks. As a result, we can then use the Faber-Krahn inequality as a comparison tool for Dirichlet boundary problems on our surfaces.

### 2.2.1 The upper half plane

As a model for hyperbolic geometry, we consider the upper half plane

$$\mathbb{H} = \{z \in \mathbb{C} \mid \Im z > 0\},$$

with metric

$$ds_\mathbb{H}^2 = y^{-2}\left(dx^2 + dy^2\right).$$

The Laplacian with this metric takes the form

$$\Delta_\mathbb{H} = -y^2\left(\frac{\partial^2}{\partial x^2} + \frac{\partial^2}{\partial y^2}\right).$$

The distance function is

$$\rho(z,w) = \log\left(\frac{|z-\overline{w}| + |z-w|}{|z-\overline{w}| - |z-w|}\right).$$

Geodesics in $\mathbb{H}$ are semicircles and half lines perpendicular to the real line. The isometries of hyperbolic space are given by Möbius transformations

$$t: z \mapsto \frac{az+b}{cz+d},$$

such that $a, b, c, d \in \mathbb{R}$ and $ad - bc = 1$. With these conditions, a Möbius transformation $t$ determines a matrix of the form

$$\begin{pmatrix} a & b \\ c & d \end{pmatrix},$$

up to a sign. Since the identity transformation can be given by $a = d = \pm 1$ (and $b = c = 0$), the group of orientation preserving isometries $\mathrm{Isom}^+(\mathbb{H})$ is

$$\mathrm{PSL}_2(\mathbb{R}) = \mathrm{SL}_2(\mathbb{R})/(\pm 1).$$

In this thesis, we will be working with the full group of isometries $\mathrm{Isom}(\mathbb{H})$; this is given by the Möbius transformations and the reflection $z \mapsto -\overline{z}$. The function

$$\cosh(\rho(z,w)) = 1 + 2u(z,w),$$

where

$$u(z,w) = \frac{|z-w|^2}{4\Im z \Im w},$$

defines a point pair invariant, that is, for $t \in \mathrm{Isom}(\mathbb{H})$,

$$\cosh\left(\rho(tz, tw)\right) = \cosh(\rho(z,w)).$$

Any element in $\mathrm{Isom}^+(\mathbb{H})$, not equal to the identity, is conjugate to a unique matrix of one of the following forms:

$$\begin{pmatrix} \gamma & 0 \\ 0 & 1 \end{pmatrix}, \quad \begin{pmatrix} \cos(\theta) & -\sin(\theta) \\ \sin(\theta) & \cos(\theta) \end{pmatrix}, \quad \text{or} \quad \begin{pmatrix} 1 & 1 \\ 1 & 0 \end{pmatrix},$$

where $\gamma > 1$ and $0 < \theta \leq \frac{\pi}{2}$. These three cases respectively describe hyperbolic, elliptic, and parabolic transformations.

Where we have an orientation reversing isometry $t \in \text{Isom}^-(\mathbb{H}) = \text{Isom}(\mathbb{H}) \setminus \text{Isom}^+(\mathbb{H})$,

$$t : z \mapsto \frac{a\bar{z} + b}{c\bar{z} + d},$$

such that $ad - bc = -1$, $t$ will be conjugate to a matrix of the form

$$\begin{pmatrix} \gamma & 0 \\ 0 & 1 \end{pmatrix} \quad \text{or} \quad \begin{pmatrix} -1 & 0 \\ 0 & 1 \end{pmatrix},$$

where $\gamma < -1$. The first case describes a glide reflection, and the second, a reflection.

**Definition 2.25 (Fuchsian group)** *A subgroup $\Gamma \subset \text{PSL}_2(\mathbb{R})$ is called Fuchsian if it acts on $\mathbb{H}$ discontinuously. A Fuchsian group $\Gamma$ is said to be of the first kind if its limit set, that is, the limit points of orbits of $\Gamma$ is $\partial \mathbb{H}$*

**Definition 2.26 (Fundamental domain)** *A set $F \subset \mathbb{H}$ is a fundamental domain for $\Gamma$ if $F$ is a closed domain in $\mathbb{H}$, and the translations of $F$ under $\Gamma$ tessellate $\mathbb{H}$.*

Every finite area hyperbolic surface $S$ can be represented as the quotient of the $\mathbb{H}$ by the action of a Fuchsian group $\Gamma$ of the first kind, that is

$$S \cong \Gamma \setminus \mathbb{H}.$$

If $\Gamma$ contains only hyperbolic elements, then the surface will be smooth and compact, with genus $g \geq 2$. This neat isomorphism comes from the uniformization theorem, which we state in two parts (see, for example, [24]).

**Theorem 2.27 (Koebe, Poincaré)** *For any smooth Riemannian metric on a surface, there exists a conformally equivalent metric of constant curvature.*

We can use a global rescaling to restrict ourselves to curvatures of $+1$, $0$, and $-1$. Then we have the following

**Theorem 2.28 (Hopf)** *Up to isometry, the only smooth, complete, simply connected surfaces of constant curvature are $\mathbb{S}^2$, $\mathbb{R}^2$, and $\mathbb{H}^2$.*

The complex analogue of this is that every simply connected Riemann surface is holomorphically equivalent to $\mathbb{C} \cup \infty$, $\mathbb{C}$ or $\mathbb{H}$.

Since Fuchsian groups act discontinuously, a fundamental domain $F$ will share a border with a finite number of translates of itself in the associated tessellation. If $\Gamma$ is a Fuchsian group of the first kind, its fundamental domain can be chosen to be a convex polygon given by

$$D(w) = \{z \in \mathbb{H} \mid \rho(z, w) < \rho(z, \gamma w) \, \forall \gamma \in \Gamma, \, \gamma \neq I\},$$

for $w \in \mathbb{H}$ not fixed by such a $\gamma$, where $I$ is the identity action. This object is referred to as the Dirichlet polygon. Side associations can be given on $D(w)$ such that the group actions that "glue" a side to its pair, also generate $\Gamma$.

### 2.2.2 The Poincaré disk

The Poincaré disk is a conformal model of hyperbolic geometry; the property of preserving Euclidean angles, and radial symmetry, means that it is a convenient model in which to view symmetric hyperbolic surfaces. It is the open unit disk in the complex plane

$$\mathbb{D} = \{z \in \mathbb{C} \,|\, |z| < 1\},$$

with Riemannian measure

$$ds^2 = \frac{4dzd\bar{z}}{(1-|z|^2)^2}.$$

Since the metric is radially symmetric about the origin, it means that the regular polygons we will consider actually 'look' regular, as opposed to the upper half plane model, where they are skewed towards the real line. This makes it much easier to see, for example, the rotational symmetry of the surfaces.

Geodesics in this model are diameters of the unit circle, and circles meeting the boundary of the disk perpendicularly. All the lines in the tessellations shown in Chapter 3 and Chapter 4 are made up of geodesics.

Figure 2.1: Geodesics in the Poincare disk

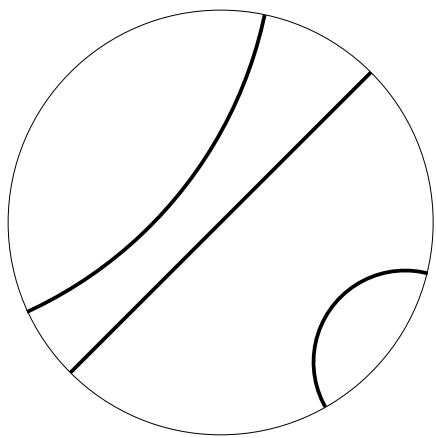

### 2.2.3 Triangle groups

The fundamental polygons for the surfaces we consider in this paper can be tessellated by hyperbolic triangles: for the Bolza surface, the fundamental octagon is tessellated with 96 triangles having angles $(\frac{\pi}{2}, \frac{\pi}{3}, \frac{\pi}{8})$; the 14-gon of the Klein quartic contains 336 triangles with angles $(\frac{\pi}{2}, \frac{\pi}{3}, \frac{\pi}{7})$; and the surface M3 [72] has 96 $(\frac{\pi}{2}, \frac{\pi}{3}, \frac{\pi}{12})$ triangles.

**Definition 2.29 (Triangle group)** *A triangle group is a subgroup $\Gamma \subset \mathrm{PSL}_2(\mathbb{R})$ generated by reflections in the sides of an $(l, m, n)$ triangle (a triangle with angles $\frac{\pi}{l}$, $\frac{\pi}{m}$ and $\frac{\pi}{n}$, where $\frac{1}{l} + \frac{1}{m} + \frac{1}{n} < 1$ and $l, m, n \in \mathbb{N}$). For an $(l, m, n)$ triangle with side reflections given by $x, y, z$, the corresponding triangle group $T(l, m, n)$ has the presentation*

$$T(l, m, n) = \langle x, y, z \mid x^2 = y^2 = z^2 = (xy)^l = (yz)^m = (zx)^n = e \rangle.$$

Figure 2.2: Generators of an $(l, m, n)$ triangle group

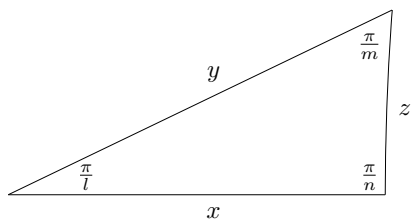

The area of a right angled hyperbolic triangle is

$$\frac{\pi}{2} - \alpha - \beta,$$

where $\alpha$ and $\beta$ are the sizes of the other two angles. This formula is known as the Gauss defect, and comes from the more general Gauss-Bonnet theorem. Hence a $(2, 3, 7)$ triangle has area $\frac{\pi}{42}$, a $(2, 3, 8)$ triangle has area $\frac{\pi}{24}$ and a $(2, 3, 12)$ triangle has area $\frac{\pi}{12}$. We can use the following formulae from [13] to calculate the hyperbolic lengths of the sides, based on the fact that we know all the angles.

**Theorem 2.30** *For any right angled hyperbolic triangle, the following hold:*

$$\cosh(c) = \cot(\alpha) \cot(\beta),$$

$$\cos(\alpha) = \cosh(a) \sin(\beta),$$

*where $c$ is the hypotenuse, $a$ is the side with opposite angle $\alpha$, and $b$ is the side with opposite angle $\beta$.*

The hyperbolic side lengths for the $(2, 3, 7)$ triangle are:

$$a = \operatorname{arccosh}\left(\frac{\csc\left(\frac{\pi}{7}\right)}{2}\right) \approx 0.545,$$

$$b = \operatorname{arccosh}\left(\frac{2\cos\left(\frac{\pi}{7}\right)}{\sqrt{3}}\right) \approx 0.283,$$

$$c = \text{arccosh}\left(\frac{\cot\left(\frac{\pi}{7}\right)}{\sqrt{3}}\right) \approx 0.621.$$

For the (2, 3, 8) triangle, they are:

$$L = \text{arccosh}\left(\frac{\csc\left(\frac{\pi}{8}\right)}{2}\right) \approx 0.7642,$$

$$M = \text{arccosh}\left(\frac{2\cos\left(\frac{\pi}{8}\right)}{\sqrt{3}}\right) \approx 0.3635,$$

$$N = \text{arccosh}\left(\frac{\cot\left(\frac{\pi}{8}\right)}{\sqrt{3}}\right) \approx 0.8607.$$

We will not need the side lengths for a (2, 3, 12) triangle in this report, but clearly they will take a similar form.

The isometries of the surfaces come from the movements of a triangle around the tessellation by reflections in one or more of its sides. The relations between these reflections form the triangle group, and where the polygons have side associations, additional relations will arise in the presentation of the group.

## 2.3 Harmonic analysis on $\mathbb{H}$

### 2.3.1 Disks

Here we take "disk" to mean a general disk in hyperbolic space rather than the Poincaré disk model. Due to the radial symmetry of the disk, it is possible to give an explicit formula for the eigenfunctions. From these, we can calculate the eigenvalues of the disk and use them to bound the eigenvalues of a domain in our surface from below. First, we introduce a more convenient polar co-ordinate system, known as the geodesic polar co-ordinates. Recall the different forms that a group action in $\text{Isom}^+(\mathbb{H})$ can take (for example, reflection or rotation matrices). An element $g \in \text{Isom}^+(\mathbb{H})$ can be written as $g = k(\varphi)a(e^{-r})k(\theta)$, where

$$k(\varphi) = \begin{pmatrix} \cos(\varphi) & \sin(\varphi) \\ -\sin(\varphi) & \cos(\varphi) \end{pmatrix}, \quad 0 \leq \varphi < \pi,$$

$$a(e^{-r}) = \begin{pmatrix} e^{-r/2} & 0 \\ 0 & e^{r/2} \end{pmatrix}, \quad r \geq 0,$$

and $k(\theta)$ is a rotation by angle $2\theta$ around $i$. This way of expressing $g$ is known as Cartan's decomposition [38]. Given a point $z \in \mathbb{C}$ such that $z = x + yi$, in geodesic polar co-ordinates, $z = k(\varphi)e^{-r}i$ and we use the following to define the change of co-ordinates:

$$y = (\cosh(r) + \sinh(r)\cos(2\varphi))^{-1},$$

$$x = y \sinh(r) \sin(2\varphi).$$

The Riemannian measure and length element in terms of $(r, \varphi)$ are now

$$ds^2 = dr^2 + (2\sinh(r))^2 d\varphi^2,$$

$$dz = (2\sinh r) dr\, d\varphi.$$

These give us formulae for the area of a disk D:

$$\text{Area}(D) = 4\pi \sinh^2\left(\frac{r}{2}\right),$$

and circumference of its boundary

$$l(\partial D) = 2\pi \sinh(r),$$

where we normalize the disk to be centred at $i$ and take $r$ as above to be its radius. The Laplacian takes the form

$$\Delta = -\left(\frac{\partial^2}{\partial r^2} + \frac{1}{\tanh r}\frac{\partial}{\partial r} + \frac{1}{(2\sinh r)^2}\frac{\partial^2}{\partial \varphi^2}\right),$$

and eigenfunctions are Legendre p functions

$$P_{-s}(\cosh(r)) = {}_2\mathbf{F}_1\left(s, 1-s; 1; \frac{1-\cosh(r)}{2}\right),$$

where

$${}_2\mathbf{F}_1(a, b; c; z) = \sum_{k=0}^{\infty} \frac{(a)_k (b)_k}{(c)_k k!} z^k, \qquad (s)_k = s \cdots (s+k-1)$$

is the hypergeometric function and $\lambda = s(1-s)$ are the eigenvalues. The full derivation of these eigenfunctions is given in [38].

### 2.3.2 Cylinders

Not all domains in hyperbolic surfaces satisfy the isoperimetric inequality. It may still be desirable to investigate the eigenvalues of such regions; we will see a particular case of this when we look at the Bolza surface in Section 3.4. One region that does not, in general, satisfy the isoperimetric inequality is the hyperbolic cylinder. We take the definition from [57].

**Definition 2.31 (Hyperbolic Cylinder)** *A hyperbolic cylinder $\mathcal{C}$, with core geodesic $\gamma$ of length $l$, is the quotient of $\mathbb{H}$ by a hyperbolic isometry*

$$\tau = \begin{pmatrix} e^{l/2} & 0 \\ 0 & e^{-l/2} \end{pmatrix}, \quad l > 0.$$

*The core geodesic is the image of the axis of $\tau$ under the quotient map $\mathbb{H}/\langle \tau \rangle$.*

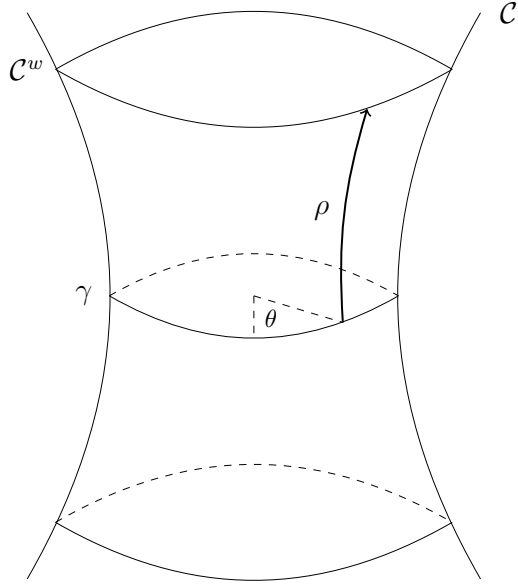

Figure 2.3: Collar $\mathcal{C}^w$ contained in hyperbolic cylinder $\mathcal{C}$

We use Fermi coordinates to describe a point $p \in \mathcal{C}$. Denote by $l$ the length of $\gamma$, parametrized with constant speed. $\rho$ is the signed distance of $p$ from $\gamma$ and $\theta$ is the projection of $p$ onto $\gamma$. For example, in Figure 2.3, the Fermi coordinates of the point on the boundary of $\mathcal{C}^w$ at the end of the arrow are $(\rho, \theta)$.

**Definition 2.32 (Collar)** *The collar $\mathcal{C}^w$, for $w \geq l$ is*

$$\mathcal{C}^w = \{(\rho, \theta) \in \mathcal{C} \,|\, l \cosh(\rho) < w,\, 0 \leq \theta \leq 2\pi\}.$$

**Remark 2.33** *Note that $\mathcal{C}^w$ is diffeomorphic to a annulus $\Omega$ of constant width $\rho$, with two boundary components of length $w$. In particular, for such an annulus, we have the following formulae for its area and the length of its boundary:*

$$A(\Omega) = 2l \sinh(\rho),$$
$$l(\partial \Omega) = 2l \cosh(\rho),$$

*where again, $l$ is the length of the core geodesic, and $\rho$ is the constant distance from the core geodesic to the boundary.*

With respect to the Fermi coordinates, $\mathcal{C}$ has the metric

$$ds^2 = d\rho^2 + \frac{l^2 \cosh^2 \rho}{4\pi^2} d\theta^2.$$

The Laplacian on $\mathcal{C}$ is given by

$$\Delta = -\left(\frac{d^2}{d\rho^2} + \tanh \rho \frac{d}{d\rho} + \frac{4\pi^2}{l^2 \cosh^2(\rho)} \frac{d^2}{d\theta^2}\right).$$

There is a fundamental system of even and odd eigenfunctions; these are given and analyzed in [11]. With respect to the core geodesic $\gamma$, the even function is

$$\phi_k^{even}(s;\rho) = (\cosh(\rho))^{\frac{2\pi i k}{l}} {}_2\mathbf{F}_1\left(\frac{s}{2} + \frac{\pi i k}{l}, \frac{1-s}{2} + \frac{\pi i k}{l}; \frac{1}{2}; -\sinh^2(\rho)\right)$$

and the odd function is

$$\phi_k^{odd}(s;\rho) = \sinh(\rho)(\cosh(\rho))^{\frac{2\pi i k}{l}} {}_2\mathbf{F}_1\left(\frac{1+s}{2} + \frac{\pi i k}{l}, \frac{2-s}{2} + \frac{\pi i k}{l}; \frac{3}{2}; -\sinh^2(\rho)\right)$$

where, once again, the corresponding eigenvalues are $\lambda = s(1-s)$.

Although annuli in hyperbolic cylinders do not satisfy the isoperimetric inequality in Definition 1.13, Mondal proved that hyperbolic cylinders satisfy a similar geometric isoperimetric inequality. First we need the notion of symmetrization:

**Definition 2.34 (Symmetrization)** *Let $\Omega$ be an essential annulus with core geodesic $\gamma$. The symmetrization of $\Omega$ is the annulus $\Omega_0$ of constant width that is symmetric with respect to the core geodesic $\gamma$ and has the same area as $\Omega$.*

Then we have the following:

**Lemma 2.35** *[57, Lemma 2.2.3] Let $\Omega$ be an essential annulus with piecewise smooth boundary and $\Omega_0$ be the symmetrisation of $\Omega$. Then $l(\partial\Omega) \geq l(\partial\Omega_0)$.*

**Proof:** The proof is a partial reproduction of that in [57]. Let $\Omega$ be an essential annulus in a cylinder $\mathcal{C}$ with core geodesic $\gamma$ of length $k$, and $\Omega_0$ its symmetrisation. Lemma 2.35 is proved in two parts; first, when $\partial\Omega$ is graph over $\gamma$, and the more general case when $\partial\Omega$ is arbitrary. The two parts of the proof are similar, so here we just reproduce the first part.

Since $\partial\Omega$ is a graph over $\gamma$, we can take two piecewise functions $r_1$ and $r_2$ to map $\gamma$ to the boundary in the positive and negative directions respectively (recall that the distance in Fermi coordinates is signed). Formally, for $i \in \{1, 2\}$,

$$r_i : [0, k] \to \mathbb{R}, \quad r_i(0) = r_i(k),$$

and $\partial\Omega$ is parametrised in Fermi coordinates as

$$\{(t, r_i(t)) \mid t \in [0, k]\}.$$

Without loss of generality, we can assume $r_1(t) > r_2(t)$ for $t \in [0, k]$. We integrate with respect to the hyperbolic cylinder metric (see Section 2.3.2) to calculate

$$A(\Omega) = \int_0^k \int_{r_1(t)}^{r_2(t)} \cosh(r) dr dt = \int_0^k (\sinh(r_2(t)) - \sinh(r_1(t))) dt$$

and
$$l(\partial\Omega) = \int_0^k (r_1'(t) + 1)^{\frac{1}{2}} \cosh(r_1(t))dt + \int_0^k (r_2'(t) + 1)^{\frac{1}{2}} \cosh(r_2(t))dt.$$

Note that
$$l(\partial\Omega) \geq \int_0^k (\cosh(r_1(t)) + \cosh(r_2(t)))dt,$$

and let $L$ denote the expression on the right hand side of this inequality. For the symmetrisation $\Omega_0$, the functions from the core geodesic to the boundary are constant, for example equal to $\pm\rho$, where $\rho = \sinh^{-1}\left(\frac{A(\Omega)}{2k}\right)$. Taking such a $\rho$ ensures that $A(\Omega) = A(\Omega_0)$. Here, recall the simpler formulae
$$l(\partial\Omega_0) = 2|k|\cosh(\rho)$$

and
$$A(\Omega_0) = 2|k|\sinh(\rho).$$

To prove the lemma, we show that
$$l(\partial\Omega)^2 - A(\Omega)^2 \geq l(\partial\Omega_0)^2 - A(\Omega_0)^2.$$

The result follows from the equality of the surface areas. First, note that
$$l(\partial\Omega)^2 - A(\Omega)^2 \geq (L + A(\Omega))(L - A(\Omega)).$$

Let us consider the two factors on the right hand side of this inequality:
$$L + A(\Omega) = \int_0^k ((\cosh(r_2(t)) + \sinh(r_2(t))) + (\cosh(r_1(t)) - \sinh(r_1(t))))dt$$
$$= \int_0^k (\exp(r_2(t)) + \exp(-r_1(t)))dt$$

and
$$L - A(\Omega) = \int_0^k (\exp(-r_2(t)) + \exp(r_1(t)))dt.$$

By Hölder's inequality,
$$(L + A(\Omega))(L - A(\Omega)) \geq$$
$$\left(\int_0^k \left((\exp(r_2(t)) + \exp(-r_1(t)))^{\frac{1}{2}}(\exp(-r_2(t)) + \exp(r_1(t)))^{\frac{1}{2}}\right)dt\right)^2$$
$$= \left(\int_0^k (2 + 2\cosh(r_1(t) + r_2(t)))^{\frac{1}{2}} dt\right)^2$$
$$\geq 4k^2,$$

since $\cosh(t) \geq 1$. On the other hand, $4k^2$ is precisely the value of
$$l(\partial\Omega_0)^2 - A(\Omega_0)^2.$$

Therefore, Lemma 2.35 is proved for $\Omega$ with 'nice' boundary. □

## 2.4 Representation theory

**Definition 2.36 (Representation)** *Let $V$ be an $n$-dimensional vector space over $\mathbb{C}$. Then $\mathrm{GL}(V)$ denotes the the group of isomorphisms $V \to V$. If $V$ has a finite basis of vectors $\{e_1, \ldots, e_n\}$, then each linear map $a : V \to V$ is defined by an $n \times n$ matrix $(a_{ij})$ such that*

$$a(e_j) = \sum_{i=j}^{n} a_{ij} e_i.$$

*A representation of a finite group $G$ is a group homomorphism*

$$\rho : G \to \mathrm{GL}(V).$$

*For each $g, h \in G$, we have a map $\rho(g) : V \to V$ such that*

(i) $\rho(gh) = \rho(g)\rho(h)$,

(ii) $\rho(1) = 1$,

(iii) $\rho(g^{-1}) = \rho(g)^{-1}$.

Denote by $\dim(V)$ the dimension of the representation. A sub-representation of $\rho$ is a sub-vector space $W \subseteq V$ such that $\rho(g)w \in W$ for all $w \in W$ and $g \in G$. A subspace with this property is called $G$-invariant (or just invariant).

**Definition 2.37 (Irreducible representation)** *Let $\rho : G \to \mathrm{GL}(V)$ be a representation. We say that $\rho$ is irreducible if $V \neq 0$ and if no subspace of $V$ is $G$-invariant, except $\{0\}$ and $V$, that is, only "trivial" subspaces are $G$-invariant.*

**Definition 2.38 (Unitary representation)** *Where $V$ is a complex Hilbert space, a representation $\rho$ of $G$ is unitary if $\rho(g) = U_g$ is a unitary operator for every $g \in G$, that is, $U_g$ is surjective and preserves the inner product on $V$. We have*

$$U_g^* U_g = U_g U_g^* = I,$$

*that is, $U_g$ is the inverse of its adjoint.*

**Definition 2.39 (Character)** *The character of a representation $\rho$ of a group $G$, as above, is a function $\chi_\rho : G \to F$, where $F$ is the field over which the vector space $V$ is defined. For $g \in G$, it is given by*

$$\chi_\rho(g) = \mathrm{Tr}(\rho(g)).$$

**Remark 2.40** *We will not always include a subscript when talking about characters; it will depend on the context. Characters take a constant value on each conjugacy class. They will be used in the next section, where we introduce a general version of Selberg's trace formula that uses a character to restrict to a particular irreducible representation.*

In Chapter 3 and Chapter 4, we look at the irreducible representations of the isometry groups of specific surfaces. We use these to investigate the orthonormal basis of eigenfunctions of $L^2(\Gamma \setminus \mathbb{H})$; our first step towards achieving this is the Peter-Weyl theorem.

**Theorem 2.41 (Peter-Weyl Theorem, taken from [49])** *Let $\mathcal{H}$ be the Hilbert space of square integrable functions over $G$, that is $L^2(G)$. Let $\rho$ be a unitary representation of $G$ that acts on $\mathcal{H}$ as*

$$\rho(h)f(g) = f(h^{-1}g).$$

*Then $L^2(G)$ decomposes as an orthogonal direct sum of unitary irreducible representations of $G$, where the multiplicity of a representation is equal to its degree. We have*

$$L^2(G) = \bigoplus_\rho E_\rho^{\oplus \dim(E_\rho)},$$

*where the sum runs over the total spaces $E_\rho$ of the representations $\rho$.*

Next, we have the following property of the Laplacian:

**Lemma 2.42** *The Laplace operator commutes with the group of isometries of the surface (see, for example, [14]).*

For a isometry $\phi : S \to S$ of a compact Riemann surface $S$, and an eigenfunction of the Laplacian $f : S \to \mathbb{R}$ with eigenvalue $\lambda$, this means

$$\Delta(f \circ \phi) = (\Delta f) \circ \phi = (\lambda f) \circ \phi = \lambda(f \circ \phi).$$

Consider the eigenspaces of the Laplacian

$$\mathcal{E}_i = \mathcal{E}_i(S) = \{\phi_i \in L^2(S) \,|\, \Delta\phi_i = \lambda_i\phi_i\}.$$

We can use the Peter-Weyl theorem along with the invariance of the Laplacian under isometries to decompose the eigenspace $\mathcal{E}_i$ as a direct sum of eigenspaces of irreducible representations (these can appear more than once). Huber showed that for every irreducible representation, there are infinitely many eigenspaces in which this representation is involved [35].

## 2.5 Selberg's trace formula

We start by following the basic presentation of the trace formula given in [10, 38, 56]. Recall that a hyperbolic surface $S$ can be represented as a quotient of the hyperbolic plane by a Fuchsian group that only has hyperbolic elements (with the exception of the identity), that is, $S \cong \Gamma \setminus \mathbb{H}$. The space of square integrable functions on $S$ comes from the functions $f \in L^2(\mathbb{H})$ that are invariant under the action of an isometry

$$T_\gamma(z)f = f(\gamma^{-1}z) = f(z), \quad \gamma \in \Gamma,$$

and are square integrable over any fundamental domain $\mathcal{F}_\Gamma$ of $\Gamma$ (recall Definition 2.26). We obtain the Hilbert space $L^2(\Gamma \setminus \mathbb{H})$ by defining the obvious inner product

$$\langle f_1, f_2 \rangle = \int_{\mathcal{F}_\Gamma} f_1 \overline{f_2} d\mu$$

on these functions, where $d\mu$ is the Riemannian measure. In a similar way, we can identify functions in $C^\infty(\Gamma \setminus \mathbb{H})$ with the functions in $C^\infty(\mathbb{H})$ that are invariant under the action of $T_\gamma$, $\gamma \in \Gamma$. We have seen that the Laplacian commutes with isometries; as a result it commutes with $T_\gamma$ and thus maps $C^\infty(\mathbb{H}) \to C^\infty(\mathbb{H})$.

**Definition 2.43 (Point pair invariant)** *A point pair invariant for a Fuchsian group $\Gamma$ is a function $k : \mathbb{H} \times \mathbb{H} \to \mathbb{C}$ that satisfies the following relations:*

1. *$k(\gamma z, \gamma w) = k(z, w)$ for all $\gamma \in \Gamma$, $z, w \in \mathbb{H}$;*

2. *$k(z, w) = k(w, z)$ for all $z, w \in \mathbb{H}$.*

Such a function $k(z, w)$ depends on the hyperbolic distance $d(z, w) = u$ between the two points. In this way, we can view $k$ as a function of one non-negative variable $k(u)$. If we take $k(z, w)$ to be the kernel of an integral operator

$$(Lf)(z) = \int_{\mathbb{H}} k(z, w) f(w) d\mu,$$

then $L$ is invariant under the action of $\text{Isom}(\mathbb{H})$. Both the Laplacian and $L$ are invariant under isometries, and in particular, we have the following theorem (see [38, Theorem 1.9]):

**Theorem 2.44** *The invariant integral operators commute with the Laplace operator.*

This commutative property gives a relation between eigenfunctions of the Laplacian and those of the invariant integral operators; in particular, if $f$ is a solution to

$$(\Delta f - \lambda f) = 0,$$

where $\lambda = \frac{1}{2} + t^2$, then it is also a solution to

$$Lf = h(t)f,$$

where $h(t)$ can be computed from $k(z, w)$ using the Selberg/Harish-Chandra transform [38]:

$$q(v) = \int_v^\infty k(u)(u-v)^{-\frac{1}{2}} du,$$

$$g(r) = 2q\left(\left(\sinh\left(\frac{r}{2}\right)\right)^2\right),$$

$$h(t) = \int_{-\infty}^{\infty} e^{irt} g(r) dr.$$

To study the spectrum of the Laplacian on $\Gamma \setminus \mathbb{H}$, note that the point pair invariant on $\mathbb{H} \times \mathbb{H}$ induces a point pair invariant on $\Gamma \setminus \mathbb{H} \times \Gamma \setminus \mathbb{H}$ with

$$k_\Gamma(z, w) = \sum_{\gamma \in \Gamma} (\gamma z, w).$$

This allows us to define an invariant integral operator

$$(Kf)(z) = \int_{\Gamma \setminus \mathbb{H}} k_\Gamma(z, w) f(w) d\mu$$

on functions $f \in C^\infty(\Gamma \setminus \mathbb{H})$. If such an $f$ is an eigenfunction of the Laplacian, then it again satisfies

$$Kf = h(t)f,$$

since

$$\int_{\mathcal{F}_\Gamma} k_\Gamma(z, w) f(w) d\mu = \int_{\mathbb{H}} k(z, w) f(w) d\mu,$$

where $\mathcal{F}_\Gamma$ is a fundamental domain for $\Gamma \setminus \mathbb{H}$.

**Remark 2.45** *In this presentation, we have given little attention to the convergence of the various series and integrals considered - the finer details can be found in [56]. The relation between eigenfunctions of the Laplacian and of the invariant integral operators depends on k being smooth, and of compact support, on $\mathbb{H}$ (or $\Gamma \setminus \mathbb{H}$). It is more common to specify conditions on the spectral function h(t), since this function, along with its Fourier transform, appear in Selberg's trace formula below. In particular, we require the following conditions:*

*(i) h is analytic on $|\Im(t)| \leq \sigma$ for some $\sigma > \frac{1}{2}$;*

*(ii) $h(t) = h(-t)$;*

*(iii) $|h(t)| \ll (1 + \Re(t))^{-2-\delta}$ for fixed $\delta > 0$, uniformly $\forall\, t$ in $|\Im(t)| \leq \sigma$.*

Let $\{\varphi_0, \varphi_1, \ldots\} \in C^\infty(\Gamma \setminus \mathbb{H})$ be eigenfunctions of $\Delta$, and recall that these form an orthonormal basis of $L^2(\Gamma \setminus \mathbb{H})$. In particular, for $f \in C^2(\Gamma \setminus \mathbb{H})$, we have the spectral decomposition

$$f(z) = \sum_j \langle f, \varphi_j \rangle \varphi_j(z).$$

Now for $h(t)$ satisfying the conditions in Remark 2.45, we have the eigenvalue equation

$$K\varphi_j = h(t_j)\varphi_j,$$

and the decomposition of our point pair invariant (see [56, Proposition 9])

$$k_\Gamma(z, w) = \sum_{j=0}^{\infty} h(t_j)\varphi_j(z)\overline{\varphi_j(w)}.$$

**Definition 2.46 (Trace class)** *Let $\mathcal{F}_\Gamma$ be a fundamental domain for $\Gamma \backslash \mathbb{H}$ as above. The integral operator with $K$ with $k_\Gamma(z, w)$ as its kernel are of trace class if $k_\Gamma(z, w)$ is absolutely integrable for $z = w$. The trace of $K$ is*

$$\mathrm{Tr} K = \int_{\mathcal{F}_\Gamma} k_\Gamma(z, z) d\mu.$$

Selberg's trace formula allows us to compute $\mathrm{Tr} K$ in two different ways; one using the spectral decomposition of $k_\Gamma$, and the other using the geometry of $\Gamma$. First, we integrate over the spectral decomposition of $k_\Gamma$, that is

$$k_\Gamma(z, z) = \sum_j h(t_j)|\varphi_j(z)|^2,$$

to get

$$\mathrm{Tr} K = \sum_j h(t_j),$$

since the $\varphi_j$ are orthonormal. This gives us the spectral side of the formula. Now recall

$$k_\Gamma(z, z) = \sum_{\gamma \in \Gamma} k(\gamma z, z).$$

Taking the trace of $K$, we get

$$\mathrm{Tr} K = \sum_{\gamma \in \Gamma} \int_{\mathcal{F}_\Gamma} k_\Gamma(\gamma z, z) d\mu.$$

Let $[\gamma]$ denote the conjugacy class of an element of $\Gamma$, that is

$$[\gamma] = \left\{\tilde{\gamma} \in \Gamma \,|\, \tilde{\gamma} = \tau\gamma\tau^{-1} \text{ for some } \tau \in \Gamma\right\}.$$

We define the partial kernel over a conjugacy class to be

$$k_{[\gamma]}(z, z) = \sum_{\gamma \in [\gamma]} k(\gamma z, z).$$

Note that

$$k_\Gamma(z, z) = \sum_{[\gamma]} k_{[\gamma]}(z, z) \quad \text{and} \quad \mathrm{Tr} K = \sum_{[\gamma]} \mathrm{Tr} K_{[\gamma]},$$

where $K_{[\gamma]}$ is the integral operator that has kernel $k_{[\gamma]}$. The centralizer of $\gamma$ is

$$Z_\gamma = \left\{\tau \in \Gamma \,|\, \tau\gamma = \gamma\tau\right\}.$$

Two elements $\tau, \tau' \in \Gamma$ give the same conjugate of $\gamma$ if and only if $\tau'\tau^{-1}$ belongs to the centralizer of $\gamma$. The upshot is that the fundamental domain of the centralizer of a hyperbolic element is in general much simpler than that of the fundamental domain $\mathcal{F}_\Gamma$ (and therefore the integrals over such a domain are much simpler). Such a fundamental domain is shown in Figure 2.4.

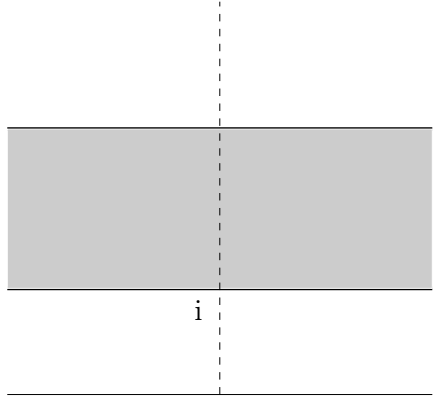

Figure 2.4: Fundamental domain for the centralizer of a hyperbolic element

For $\gamma \in [\gamma]$

$$\mathrm{Tr} K_{[\gamma]} = \sum_{\tau \in Z_\gamma \backslash \Gamma} \int_{\mathcal{F}_\Gamma} k(\tau_{-1}\gamma\tau z, z) d\mu$$

$$= \int_{Z_\gamma \backslash \mathbb{H}} k(\gamma z, z) d\mu.$$

This depends only on $[\gamma]$, that is, if $\gamma' = g^{-1}\gamma g$, $g \in \mathrm{SL}_2(\mathbb{R})$, then

$$\mathrm{Tr} K_{[\gamma]} = \int_{g^{-1}Z(\gamma)g \backslash \mathbb{H}} k(g^{-1}\gamma g z, z) d\mu = \mathrm{Tr} K_{[\gamma']}.$$

The integral depends only on the distance $d(z, \gamma z)$, and can be explicitly computed (see [38, Sections 10.4 and 10.5]). The most common version of the trace formula is for compact Riemann surfaces, where there are only hyperbolic conjugacy classes (and the identity motion). There is an exact correspondence between hyperbolic conjugacy classes and the lengths of primitive geodesics (multiple loops are not included), so this is the most interesting version for seeing the relation between the energy and length spectra.

**Theorem 2.47 (Selberg's trace formula, [74])** *Let $h(t)$ be as above, and let $g(t)$ be its Fourier transform. Let $S$ be a compact Riemann surface. Then*

$$\sum_{j=0}^\infty h(t_j) = \frac{\mathrm{Area}(S)}{4\pi} \int_{-\infty}^\infty h(t) \tanh(\pi t) t \, dt + \sum_{[\gamma]} \sum_{n=1}^\infty \frac{l_\gamma g(n l_\gamma)}{2\sinh(n l_\gamma / 2)}.$$

One interesting application of this theorem is to prove Weyl's law, the famous formula that describes the asymptotic properties of Laplace eigenvalues. We define the set of eigenvalues below a certain value as

$$N(\lambda) = \#j\{\lambda_j \leq \lambda\}.$$

**Theorem 2.48 (Weyl's law, [82])** *Let $S$ be a compact Riemann surface. Then*

$$N(\lambda) \sim \frac{\text{Area}(S)}{4\pi}\lambda, \quad \lambda \to \infty.$$

This can be proved by using the heat kernel

$$h(t) = e^{-\beta t^2}$$

with Fourier transform

$$g(s) = \frac{e^{-s^2/2\beta}}{\sqrt{4\pi\beta}}$$

in Selberg's trace formula, where $\beta > 0$. In particular, with $\lambda_j = t_j^2 + \frac{1}{4}$, we take the limit

$$\lim_{\beta \to 0} \sum_{j=0}^{\infty} e^{-\beta\lambda_j},$$

and the result follows from a Tauberian theorem of Wiener-Ikehara (see Theorem 9.4.7 in [13]).

If one has information about the irreducible representations of the automorphism group of a surface, one can consider a more general version of Selberg's trace formula, given in [9]. The advantage of this is that one can restrict the formula to a single unitary character $\chi$ on $\Gamma$, that is, consider the space $L^2(\Gamma \setminus \mathbb{H}, \chi)$ rather than simply $L^2(\Gamma \setminus \mathbb{H})$. $L^2(\Gamma \setminus \mathbb{H}, \chi)$ is the Hilbert space of functions $f : \mathbb{H} \to \mathbb{C}$ such that

$$f(Tz) \equiv \chi(T)f(z)$$

for all $T \in \Gamma$ and $f \in L^2(\Gamma \setminus \mathbb{H})$.

**Theorem 2.49 ( [9])** *Let $M$ be the set $T \in \Gamma$ such that $T$ does not have a cusp of $[\Gamma^+]$ as a fixed point. Let $I$ denote the identity action, and $M_I = M \cup [I]$. Let $\mathcal{F}$ be a fundamental domain for $\Gamma \setminus \mathbb{H}$. For each non-elliptic $T \in M$, there exists a hyperbolic element or glide reflection $T_0 \in Z(T)$ such that the infinite cyclic group*

$$[T_0] = \{T_0^n \,|\, n \in \mathbb{Z}\}$$

has index 1, 2, or 4 in $Z(T)$. Then for each admissible $h(t)$ we have

$$\sum_{T \in M_I} \int_{\Gamma \backslash \mathbb{H}} \chi(T) k(z, Tz) d\mu z = \frac{\text{Area}(\mathcal{F})}{4\pi} \int_{-\infty}^{\infty} h(t) \tanh(\pi t) t dt$$
$$+ \sum_{\substack{T \in M \\ non-ell}} \frac{\log N(T_0)}{[Z_\Gamma(T) : [T_0]]} \frac{\chi(T) g(\log N(T))}{N(T)^{\frac{1}{2}} - \text{sgn}(\det T) N(T)^{-\frac{1}{2}}}$$
$$+ \sum_{\substack{T \in M \\ elliptic}} \frac{\chi(T)}{2|Z_\Gamma(T)| \sin(\theta(T))} \int_{-\infty}^{\infty} \frac{\exp(-2\theta(T)t)}{1 + \exp(-2\pi t)} h(t) dt,$$

where the $[T]$ sums are taken over distinct conjugacy classes of $\Gamma$, and the non-elliptic sum consists of conjugacy classes where $T$ is a reflection, hyperbolic transformation or glide reflection (see Section 2.2 on hyperbolic geometry; also for $\theta$). The sums and integrals are all convergent with good majorants. The elliptic sum is finite. Note the sgn function in the non-elliptic sum, which distinguishes between orientation preserving and orientation reversing isometries.

This theorem is interesting in that it includes orientation reversing isometries. The disadvantage is that it only applies to unitary characters. For characters of higher dimensional irreducible representations, Theorem 4.1 in Chapter 8 of [33] provides a similar formula, although this does not deal with reflections. One might conjecture that it is possible to combine the two formulae to get a more general version of Selberg's trace formula for characters of higher dimensional representations of an automorphism group that includes reflections.

# Chapter 3

# The Bolza Surface

## 3.1 Symmetry group of the Bolza surface

In his thesis, Jenni [41] lists the generators of the group of isometries, and the relations between them. The whole isometry group is generated by 4 actions, which are shown in Figure 3.1:

- $R$ - rotation of order 8 around the origin (centre of the Bolza surface) that maps $Q_i \mapsto Q_{i+1}$ (modulo 8), and preserves $P_1$ and $P_2$;

- $S$ - reflection in the real line that maps $Q_i \mapsto Q_{9-i}$, and preserves $P_1$ and $P_2$;

- $T$ - reflection that switches the centre point with the outmost, that is, maps $P_1 \mapsto P_2$ and preserves the $Q_i$s;

- $U$ - rotation of order 3 around the centre of a (4, 4, 4) triangle, 16 of which tessellate the surface. This is the only action which permutes the $P$ and $Q$ points; for example, it maps $P_1 \mapsto Q_8$, $Q_8 \mapsto Q_1$, and $Q_1 \mapsto P_1$.

**Remark 3.1** *Each element of the group of isometries can be written uniquely in the canonical form*

$$R^i S^j T^k U^l, \quad i \in \{0, \ldots, 7\}, \quad j, k \in \{0, 1\}, \quad l \in \{0, 1, 2\}.$$

*An easy computation confirms that there will be $8 \cdot 2 \cdot 2 \cdot 3 = 96$ elements [41].*

**Remark 3.2** *$R$ and $U$ are orientation preserving actions, whereas $S$ and $T$ are orientation reversing. Therefore the orientation preserving elements of $\mathrm{Isom}(\mathcal{B})$, denoted $\mathrm{Isom}^+(\mathcal{B})$, will have the form $R^i(ST)^j U^l$. The orientation reversing elements will take the form either $R^i S U^l$ or $R^i T U^l$.*

Figure 3.1: Symmetries of the Bolza surface

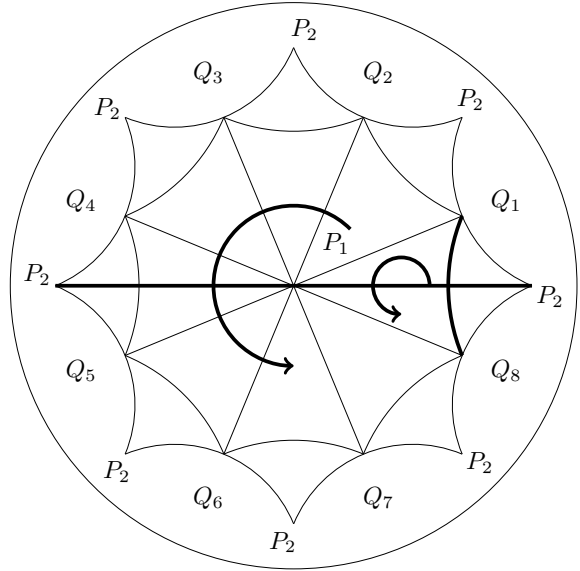

Calculations in GAP [26] show that the relations given in Satz $3_2$ of Jenni's thesis give a finitely presented group of order 12. This implies that there is a problem with at least one of his relations. After scrutinizing the actions of $R$, $S$, $T$, and $U$ on the tessellation regular octagons, we offer the following corrected theorem:

**Theorem 3.3 (Correction to Satz $3_2$ in [41])** *The generators of the Bolza surface satisfy the following relations:*

$$R^8 = S^2 = T^2 = U^3 = RSRS = STST = RTR^3T = e,$$

$$UR = R^7U^2, \ U^2R = STU, \ US = SU^2, \ UT = RSU,$$

*where $e$ is the trivial (identity) action.*

**Remark 3.4** $UR = R^7U^2$ *replaces the erroneous* $UR = RU^2$.

We check that these relations are satisfied by considering the actions of $R$, $S$, $T$, and $U$ on the 16 equilateral triangles shown in Figure 3.2, where the vertices $P_i$, $i \in \{1, 2\}$, and $Q_j$, $j \in \{1, \ldots, 4\}$, represent the Weierstrass points of the surface; the bar notation denotes that a point is in the other sheet of the cover (compare Figure 1.7).

We tabulate this information in Table 3.1, and from this point it is straightforward to verify that the relations hold by "chasing" the triangle around the table. For example, if we want to know where $RS$ sends $P_1Q_1Q_2$, we observe

$$RS(P_1Q_1Q_2) = R(S(P_1Q_1Q_2))$$

Figure 3.2: The Bolza surface tessellated by 16 equilateral triangles

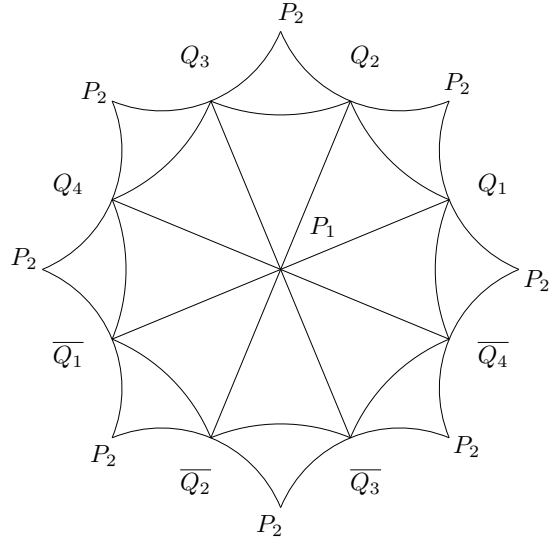

$$= R(P_1\overline{Q_3Q_4})$$
$$= P_1\overline{Q_4}Q_1.$$

Table 3.2 records the orbits of the triangles under the action of the generators.

**Proposition 3.5** *The centre of the isometry group is $R^4$.*

**Proof:** To prove this, we first have to show that $R^4$ commutes with the three other isometries. We start with $S$. Using $RS = SR^7$,

$$\begin{aligned} R^4 S &= R^3 S R^7 \\ &= S R^5 R^7 \\ &= S R^{12} \\ &= S R^4. \end{aligned}$$

For $T$, we will use $R^3 T = T R^7$; this comes from multiplying both sides of $RTR^3T = e$ by $TR^7$. Using $RTR^3 = T$ and $RT = TR^5$,

$$\begin{aligned} R^4 T &= R^3 T R^5 \\ &= T R^7 R^5 \\ &= T R^4. \end{aligned}$$

And for $U$, we use $UR = R^7 U^2$, $RUR = U^2$, and $RU = U^2 R^7$,

Table 3.1: Actions of isometries on equilateral triangles

| Triangle | $R$ | $S$ | $T$ | $U$ |
|---|---|---|---|---|
| $P_1\overline{Q_4Q_1}$ | $P_1Q_1Q_2$ | $P_1\overline{Q_4}Q_1$ | $P_2\overline{Q_4}Q_1$ | $P_1\overline{Q_4}Q_1$ |
| $P_1Q_1Q_2$ | $P_1Q_2Q_3$ | $P_1\overline{Q_3Q_4}$ | $P_2\overline{Q_1}Q_2$ | $P_1\overline{Q_3Q_4}$ |
| $P_1Q_2Q_3$ | $P_1Q_3Q_4$ | $P_1\overline{Q_2Q_3}$ | $P_2Q_2Q_3$ | $P_2\overline{Q_3Q_4}$ |
| $P_1Q_3Q_4$ | $P_1Q_4\overline{Q_1}$ | $P_1\overline{Q_1}Q_2$ | $P_2\overline{Q_3Q_4}$ | $P_2Q_4\overline{Q_1}$ |
| $P_1Q_4\overline{Q_1}$ | $P_1\overline{Q_1}Q_2$ | $P_1Q_4\overline{Q_1}$ | $P_2Q_4\overline{Q_1}$ | $P_1Q_4\overline{Q_1}$ |
| $P_1\overline{Q_1}Q_2$ | $P_1\overline{Q_2Q_3}$ | $P_1Q_3Q_4$ | $P_2\overline{Q_1}Q_2$ | $P_1Q_3Q_4$ |
| $P_1\overline{Q_2Q_3}$ | $P_1\overline{Q_3Q_4}$ | $P_1Q_2Q_3$ | $P_2\overline{Q_2Q_3}$ | $P_2Q_3Q_4$ |
| $P_1\overline{Q_3Q_4}$ | $P_1\overline{Q_4}Q_1$ | $P_1Q_1Q_2$ | $P_2Q_3Q_4$ | $P_2\overline{Q_4}Q_1$ |
| $P_2\overline{Q_4}Q_1$ | $P_2Q_1Q_2$ | $P_2\overline{Q_4}Q_1$ | $P_1\overline{Q_4}Q_1$ | $P_1Q_1Q_2$ |
| $P_2Q_1Q_2$ | $P_2Q_2Q_3$ | $P_2\overline{Q_3Q_4}$ | $P_1\overline{Q_1}Q_2$ | $P_1\overline{Q_2Q_3}$ |
| $P_2Q_2Q_3$ | $P_2Q_3Q_4$ | $P_2\overline{Q_2Q_3}$ | $P_1Q_2Q_3$ | $P_2Q_2Q_3$ |
| $P_2Q_3Q_4$ | $P_2Q_4\overline{Q_1}$ | $P_2\overline{Q_1}Q_2$ | $P_1\overline{Q_3Q_4}$ | $P_2Q_1Q_2$ |
| $P_2Q_4\overline{Q_1}$ | $P_2\overline{Q_1}Q_2$ | $P_2Q_4\overline{Q_1}$ | $P_1Q_4\overline{Q_1}$ | $P_1\overline{Q_1}Q_2$ |
| $P_2\overline{Q_1}Q_2$ | $P_2\overline{Q_2Q_3}$ | $P_2Q_3Q_4$ | $P_1Q_1Q_2$ | $P_1Q_2Q_3$ |
| $P_2\overline{Q_2Q_3}$ | $P_2\overline{Q_3Q_4}$ | $P_2Q_2Q_3$ | $P_1\overline{Q_2Q_3}$ | $P_2Q_2Q_3$ |
| $P_2\overline{Q_3Q_4}$ | $P_2\overline{Q_4}Q_1$ | $P_2Q_1Q_2$ | $P_1Q_3Q_4$ | $P_2\overline{Q_1}Q_2$ |

$$R^4 U = R^3 U^2 R^7$$
$$= U R^5 R^7$$
$$= U R^4.$$

Obviously $R^4$ commutes with $R$. Since it commutes with each of the basic actions it will commute with every isometry of the form $R^i S^j T^k U^l$. Now we just need to show that there are no more elements in the centre. We know from analysis in GAP [26] that the centre is isomorphic to $\mathbb{Z}_2$, and must contain the identity element. Hence $R^4$ is the only non-trivial element that can appear. $(R^4)^2 = R^8 = e$, so when we create the isomorphism, we map $R^4$ to the non-identity element in $\mathbb{Z}_2$.

□

The Bolza surface can also be tessellated by 96 (2, 3, 8) triangles, each with area $\pi/24$. A (2, 3, 8) triangle is fundamental domain for $\text{Isom}(\mathcal{B})$. The domains that we consider in Section 3.4 will be constructed using these triangles; for example, the equilateral (4, 4, 4) triangles of Figure 3.1 can be barycentrically subdivided into 6 (2, 3, 8) triangles. Reflections in the geodesic lines in Figure 3.3 can be described as a combination of $R$, $S$, $T$, and $U$. In the figures of Section 3.4, a bold and a dashed line will respectively denote Dirichlet and Neumann boundary conditions. The subscripts on $R$, $S$, $T$, and $U$ make it clear which irreducible representation we are using, and distinguishes the group

Table 3.2: Checking the relations between the generators of the Bolza surface

|  | $RSRS$ | $STST$ | $RTR^3T$ | $UR$ | $R^7U$ | $U^2R$ | $STU$ | $US$ | $SU^2$ | $UT$ | $RSU$ |
|---|---|---|---|---|---|---|---|---|---|---|---|
| $P_1\overline{Q_4}Q_1$ | $P_1\overline{Q_4}Q_1$ | $P_1\overline{Q_4}Q_1$ | $P_1\overline{Q_4}Q_1$ | $P_1\overline{Q_3Q_4}$ | $P_1\overline{Q_3Q_4}$ | $P_2\overline{Q_4}Q_1$ | $P_2\overline{Q_4}Q_1$ | $P_1\overline{Q_4}Q_1$ | $P_1\overline{Q_4}Q_1$ | $P_1Q_1Q_2$ | $P_1Q_1Q_2$ |
| $P_1Q_1Q_2$ | $P_1Q_1Q_2$ | $P_1Q_1Q_2$ | $P_1Q_1Q_2$ | $P_2\overline{Q_3Q_4}$ | $P_2\overline{Q_3Q_4}$ | $P_2\overline{Q_1Q_2}$ | $P_2\overline{Q_1Q_2}$ | $P_2\overline{Q_4}Q_1$ | $P_2\overline{Q_4}Q_1$ | $P_1Q_2Q_3$ | $P_1Q_2Q_3$ |
| $P_1Q_2Q_3$ | $P_1Q_2Q_3$ | $P_1Q_2Q_3$ | $P_1Q_2Q_3$ | $P_2Q_4\overline{Q_1}$ | $P_2Q_4\overline{Q_1}$ | $P_1\overline{Q_1Q_2}$ | $P_1\overline{Q_1Q_2}$ | $P_2Q_3Q_4$ | $P_2Q_3Q_4$ | $P_2Q_2Q_3$ | $P_2Q_2Q_3$ |
| $P_1Q_3Q_4$ | $P_1Q_3Q_4$ | $P_1Q_3Q_4$ | $P_1Q_3Q_4$ | $P_1Q_4\overline{Q_1}$ | $P_1Q_4\overline{Q_1}$ | $P_1Q_4\overline{Q_1}$ | $P_1Q_4\overline{Q_1}$ | $P_1Q_3Q_4$ | $P_1Q_3Q_4$ | $P_2\overline{Q_1Q_2}$ | $P_2\overline{Q_1Q_2}$ |
| $P_1Q_4\overline{Q_1}$ | $P_1Q_4\overline{Q_1}$ | $P_1Q_4\overline{Q_1}$ | $P_1Q_4\overline{Q_1}$ | $P_1Q_3Q_4$ | $P_1Q_3Q_4$ | $P_2Q_4\overline{Q_1}$ | $P_2Q_4\overline{Q_1}$ | $P_1Q_4\overline{Q_1}$ | $P_1Q_4\overline{Q_1}$ | $P_1\overline{Q_1Q_2}$ | $P_1\overline{Q_1Q_2}$ |
| $P_1\overline{Q_1Q_2}$ | $P_1\overline{Q_1Q_2}$ | $P_1\overline{Q_1Q_2}$ | $P_1\overline{Q_1Q_2}$ | $P_2Q_3Q_4$ | $P_2Q_3Q_4$ | $P_2Q_1Q_2$ | $P_2Q_1Q_2$ | $P_2Q_4\overline{Q_1}$ | $P_2Q_4\overline{Q_1}$ | $P_1\overline{Q_2Q_3}$ | $P_1\overline{Q_2Q_3}$ |
| $P_1\overline{Q_2Q_3}$ | $P_1\overline{Q_2Q_3}$ | $P_1\overline{Q_2Q_3}$ | $P_1\overline{Q_2Q_3}$ | $P_2\overline{Q_4}Q_1$ | $P_2\overline{Q_4}Q_1$ | $P_1Q_1Q_2$ | $P_1Q_1Q_2$ | $P_2\overline{Q_3Q_4}$ | $P_2\overline{Q_3Q_4}$ | $P_2Q_2Q_3$ | $P_2Q_2Q_3$ |
| $P_1\overline{Q_3Q_4}$ | $P_1\overline{Q_3Q_4}$ | $P_1\overline{Q_3Q_4}$ | $P_1\overline{Q_3Q_4}$ | $P_1\overline{Q_4}Q_1$ | $P_1\overline{Q_4}Q_1$ | $P_1\overline{Q_4}Q_1$ | $P_1\overline{Q_4}Q_1$ | $P_1\overline{Q_3Q_4}$ | $P_1\overline{Q_3Q_4}$ | $P_2Q_1Q_2$ | $P_2Q_1Q_2$ |
| $P_2\overline{Q_4}Q_1$ | $P_2\overline{Q_4}Q_1$ | $P_2\overline{Q_4}Q_1$ | $P_2\overline{Q_4}Q_1$ | $P_1\overline{Q_2Q_3}$ | $P_1\overline{Q_2Q_3}$ | $P_2Q_3Q_4$ | $P_2Q_3Q_4$ | $P_1Q_1Q_2$ | $P_1Q_1Q_2$ | $P_1\overline{Q_4}Q_1$ | $P_1\overline{Q_4}Q_1$ |
| $P_2Q_1Q_2$ | $P_2Q_1Q_2$ | $P_2Q_1Q_2$ | $P_2Q_1Q_2$ | $P_2Q_2Q_3$ | $P_2Q_2Q_3$ | $P_2Q_2Q_3$ | $P_2Q_2Q_3$ | $P_2\overline{Q_1Q_2}$ | $P_2\overline{Q_1Q_2}$ | $P_1Q_3Q_4$ | $P_1Q_3Q_4$ |
| $P_2Q_2Q_3$ | $P_2Q_2Q_3$ | $P_2Q_2Q_3$ | $P_2Q_2Q_3$ | $P_2Q_1Q_2$ | $P_2Q_1Q_2$ | $P_1\overline{Q_2Q_3}$ | $P_1\overline{Q_2Q_3}$ | $P_2\overline{Q_2Q_3}$ | $P_2\overline{Q_2Q_3}$ | $P_2\overline{Q_3Q_4}$ | $P_2\overline{Q_3Q_4}$ |
| $P_2Q_3Q_4$ | $P_2Q_3Q_4$ | $P_2Q_3Q_4$ | $P_2Q_3Q_4$ | $P_1\overline{Q_1Q_2}$ | $P_1\overline{Q_1Q_2}$ | $P_1Q_3Q_4$ | $P_1Q_3Q_4$ | $P_1Q_2Q_3$ | $P_1Q_2Q_3$ | $P_2\overline{Q_4}Q_1$ | $P_2\overline{Q_4}Q_1$ |
| $P_2Q_4\overline{Q_1}$ | $P_2Q_4\overline{Q_1}$ | $P_2Q_4\overline{Q_1}$ | $P_2Q_4\overline{Q_1}$ | $P_1Q_2Q_3$ | $P_1Q_2Q_3$ | $P_2\overline{Q_3Q_4}$ | $P_2\overline{Q_3Q_4}$ | $P_1\overline{Q_1Q_2}$ | $P_1\overline{Q_1Q_2}$ | $P_1Q_4\overline{Q_1}$ | $P_1Q_4\overline{Q_1}$ |
| $P_2\overline{Q_1Q_2}$ | $P_2\overline{Q_1Q_2}$ | $P_2\overline{Q_1Q_2}$ | $P_2\overline{Q_1Q_2}$ | $P_2\overline{Q_2Q_3}$ | $P_2\overline{Q_2Q_3}$ | $P_2\overline{Q_2Q_3}$ | $P_2\overline{Q_2Q_3}$ | $P_2Q_1Q_2$ | $P_2Q_1Q_2$ | $P_1\overline{Q_3Q_4}$ | $P_1\overline{Q_3Q_4}$ |
| $P_2\overline{Q_2Q_3}$ | $P_2\overline{Q_2Q_3}$ | $P_2\overline{Q_2Q_3}$ | $P_2\overline{Q_2Q_3}$ | $P_2\overline{Q_1Q_2}$ | $P_2\overline{Q_1Q_2}$ | $P_1Q_2Q_3$ | $P_1Q_2Q_3$ | $P_2Q_2Q_3$ | $P_2Q_2Q_3$ | $P_2Q_3Q_4$ | $P_2Q_3Q_4$ |
| $P_2\overline{Q_3Q_4}$ | $P_2\overline{Q_3Q_4}$ | $P_2\overline{Q_3Q_4}$ | $P_2\overline{Q_3Q_4}$ | $P_1Q_1Q_2$ | $P_1Q_1Q_2$ | $P_1\overline{Q_3Q_4}$ | $P_1\overline{Q_3Q_4}$ | $P_1\overline{Q_2Q_3}$ | $P_1\overline{Q_2Q_3}$ | $P_2Q_4\overline{Q_1}$ | $P_2Q_4\overline{Q_1}$ |

generators from their irreducible representations.

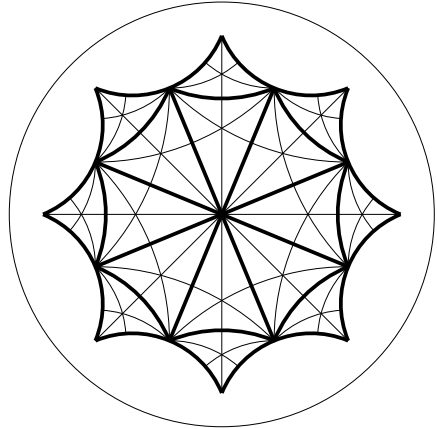

Figure 3.3: Tessellation of $\mathcal{B}$ by (2, 3, 8) and (4, 4, 4) (in bold) triangles

**Remark 3.6** *Recall the side lengths of (2, 3, 8) triangles calculated in Section 2.2.3, and the comment in the introduction regarding the systoles of the Bolza surface. From Figure 3.3, we can see that the systole of the Bolza surface consists of four "medium length" (2, 3, 8) triangle sides, that is,*

$$s(\mathcal{B}) = 4\operatorname{arccosh}\left(\frac{\csc\left(\frac{\pi}{8}\right)}{2}\right).$$

*This is equivalent to the formula*

$$2\operatorname{arccosh}\left(1+\sqrt{2}\right)$$

*from [41], given in the introduction.*

## 3.2 Length spectrum of the Bolza surface

An expression for the lengths of primitive geodesics of the Bolza surface in terms of algebraic numbers has been given in [2]. We have seen the form of the generators in the introduction; an equivalent form is given by

$$g_k = \begin{pmatrix} \cosh\left(\frac{l_0}{2}\right) & \exp\left(i\left(\frac{k\pi}{4}\right)\right)\sinh\left(\frac{l_0}{2}\right) \\ \exp\left(-i\left(\frac{k\pi}{4}\right)\right)\sinh\left(\frac{l_0}{2}\right) & \cosh\left(\frac{l_0}{2}\right) \end{pmatrix}, \quad k=0,\ldots,3,$$

where

$$\cosh\left(\frac{l_0}{2}\right) = \cot\left(\frac{\pi}{8}\right) = 1+\sqrt{2}$$

describes the length of the first geodesic. In fact, the periodic orbit matrices for all geodesics in the spectrum take a similar form:

$$g = \begin{pmatrix} \alpha & \beta \\ \beta^* & \alpha^* \end{pmatrix},$$

where

$$\alpha = a + \sqrt{2}b$$

and

$$\beta = \sqrt{1+\sqrt{2}}\left(c + \sqrt{2}d\right), \quad a,\,b,\,c,\,d \in \mathbb{Z}[i].$$

The lengths of the corresponding geodesics are determined by

$$|\Re(\alpha)| = m + \sqrt{2}n, \quad m,\,n \in \mathbb{N}.$$

We learn from [2] that $n$ runs through all the natural numbers, apart from $n=48$ and $n=72$, where there are no orbits. $m$ is the odd natural number that minimizes

$$|m/n - \sqrt{2}|.$$

Once we have obtained a sequence of $(m,n)$ pairs from the above rule, we have a handy way to compute the lengths of primitive geodesics in closed form, in particular

$$\cosh\left(\frac{l_{n-1}}{2}\right) = m + \sqrt{2}n, \quad m,\,n \in \mathbb{N}.$$

Having an exact formula for the length spectrum will be helpful later in the chapter, when we create bounds on eigenvalue multiplicity using Selberg's trace formula.

## 3.3 Upper bounds on eigenvalues of the Bolza surface

We use the Rayleigh quotient (Theorem 1.20) to give upper bounds on the first two distinct positive eigenvalues of $\Delta(\mathcal{B})$. We have calculated these numerically and listed them in Table C.1; in particular, the first and second positive eigenvalues are roughly equal to 3.84 and 5.35 respectively. Since the Rayleigh quotient is an equality for the first eigenfunction, the closer the test function is to the first eigenfunction, the better the bound will be. Rather than simply trying to guess a function, we can use the calculation for $\chi_8$ in Section 3.4 to tell us where the function should be zero (see Figure 3.13).

We know that the first eigenfunction satisfies mixed boundary conditions on a (4, 4, 4) equilateral triangle; Dirichlet on one side and Neumann on the other two. We also know that the eigenfunction should be radially symmetric with respect to the origin of the unit disk, since the (4, 4, 4) triangle is rotated about the centre of the Poincare disk by $R$. We will work in the upper half plane model when performing calculations, since the metric is slightly easier to work with. Transformations between the Poincaré disk and the upper half plane are given by the Cayley transform:

$$f(z) = \frac{z+i}{iz+1},$$

where the inverse function naturally maps $\mathbb{H}$ to $\mathbb{D}$. Under this transformation, a (4, 4, 4) triangle that is symmetric with respect to the real line is now symmetric with respect to the $\Im$ axis, including our manufactured boundary conditions. This means we can simplify our calculations further by integrating over the (2, 4, 8) triangle in Figure 3.4, denoted $\Omega_T$, where the equations of the boundary curves are given by

$$c_1(x) = \sqrt{q^2 - x^2}$$

and

$$c_2(x) = \sqrt{\csc^2\left(\frac{\pi}{8}\right) - \left(x + \cot\left(\frac{\pi}{8}\right)\right)^2},$$

where $q = \sqrt{3 + 2\sqrt{2} - 2\sqrt{4 + 3\sqrt{2}}}$.

After some experimentation with various test functions, we settled on the function

$$\varphi_1(x, y) = \left(y - \sqrt{1 - 22x^2} - 1\right)\left(y - \sqrt{q^2 - x^2}\right),$$

which is equal to zero on $c_1$ due to the $y - \sqrt{q^2 - x^2}$ factor. Note that we only enforce a (Dirichlet) boundary condition along $c_1$; the function is "free" along the other two sides of the domain in Figure 3.4 (recall from Theorem 1.20 that test functions for a Neumann boundary condition in the variational characterization used for the min-max principle only have to be in $W^{1,2}$ and not identically 0). The factor of $y - \sqrt{1 - 22x^2} - 1$ represents the circular symmetry that we observed when plotting the eigenfunction in FreeFEM++. We started with 1 as the coefficient of $x^2$, and after several numerical experiments, found 22 to be among the optimal integer coefficients for producing a good bound.

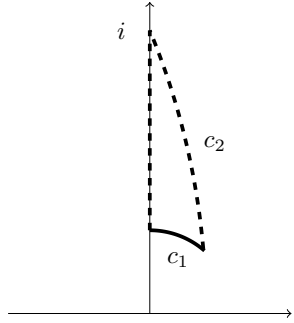

Figure 3.4: The domain of integration $\Omega_T$ in $\mathbb{H}$

We wish to evaluate the $L^2$ norms of both this function and its gradient over $\Omega_T$. Unfortunately, these cannot be explicitly calculated, however we can use Mathematica [37] to integrate explicitly with respect to $y$ and then approximate the final integral using Simpson's rule, as we will then be integrating a continuous curve over an interval in one variable. We are aware of the limitations of performing symbolic integration using computer software, however the integrands that we integrate with respect to $y$ consist of terms that are integer powers of $y$. We have chosen Mathematica to perform our integrations, but have checked them in SAGE [79], which is open source, and have also checked some of the less cumbersome ones by hand.

First, we consider $\varphi_1(x, y)$. Integrating with respect to $y$ allows us to define a function in $x$,

$$h(x) = \int_{c_1(x)}^{c_2(x)} \frac{\left(y - \sqrt{1 - 22x^2} - 1\right)^2 \left(y - \sqrt{q^2 - x^2}\right)^2}{y^2} dy.$$

We do not display the closed evaluated form of $h(x)$ here, since it is long, unpleasant, and tells us little information. It is given as a sum of polynomial terms of $x$, as well as their square roots and logs, and can be calculated in Mathematica with the following line of code:

```
q = Sqrt[3 + 2 Sqrt[2] - 2 Sqrt[4 + 3 Sqrt[2]]];
u[x_, y_] = (y - 1 - Sqrt[1 - 22 x^2]) (y - Sqrt[q^2 - x^2]);
Integrate[(u[x, y] u[x, y])/y^2, {y, Sqrt[q^2 - x^2],
  Sqrt[Csc[Pi/8]^2 - (x + Cot[Pi/8])^2]}]
```

We now need to approximate the integral

$$\int_0^{2^{\frac{1}{4}}-1} h(x)dx.$$

Simpson's rule allows us to control the accuracy of our approximation by selecting the step size, and comes with a formula for a rigorous error bound. Moreover, both the error bound and the approximation can be given symbolically. Recall that for an even value of $n$, Simpson's rule approximates an integral as

$$\int_a^b f(x)dx \approx S_n,$$

where

$$S_n = \frac{\Delta x}{3}\left(f(x_0) + 4f(x_1) + 2f(x_2) + 4f(x_3) + \ldots + 2f(x_{n-2}) + 4f(x_{n-1}) + f(x_n)\right),$$

$$\Delta x = \frac{b-a}{n}, \quad \text{and} \quad x_i = a + i\Delta x, \quad i \in \{0, \ldots, n\}.$$

The error for such an approximation depends on the maximum value of the fourth derivative on the interval. In particular, suppose that $\left|f^{(4)}(x)\right| \leq K$ for $a \leq x \leq b$ and $K \in \mathbb{R}$. Then

$$\left|\int_a^b f(x)dx - S_n\right| \leq \frac{K(b-a)^5}{180n^4}.$$

We will chose $n$ such that our error will be less than $10^{-6}$ every time we approximate with Simpson's rule, that is, we take the smallest even value of $n$ such that

$$\left(10^6 \left(\max\left|f^{(4)}(x)\right|\right) \frac{(2^{\frac{1}{4}} - 1)^5}{180}\right)^{\frac{1}{4}} < n.$$

For $h(x)$, the maximum value will occur at the upper end of the integral, that is, at $x = 2^{\frac{1}{4}} - 1$. This is due to the fact that both $h(x)$ function and its fourth derivative include the term $\log\left(1 - x^2 - 2x \cot\left(\frac{\pi}{8}\right)\right)$, which has a singularity just outside the interval.

The following line of code gives us the closed form of the fourth derivative evaluated at $2^{\frac{1}{4}} - 1$, and then substituted into our error formula.

```
h[x_]=Integrate[(u[x, y] u[x, y])/y^2, {y, Sqrt[q^2 - x^2],
   Sqrt[Csc[Pi/8]^2 - (x + Cot[Pi/8])^2]}];
adh[x_] = Abs[D[h[x], {x, 4}]];
((10^6) (adh[2^(1/4) - 1] (2^(1/4) - 1)^5)/180)^(1/4)
```

Again, the symbolic output is too long to include here, but it is approximately equal to 89.7433. Therefore, we take our $n$ to be the next even integer, 90, to guarantee the required level of accuracy. We can now approximate the integral of $h(x)$, and thus the $L^2$ norm of $\varphi_1(x, y)$, using $S_{90}$. The following input in Mathematica gives us our closed form.

```
phi = ((2^(1/4) - 1)/(270)) Total[{h[0], h[2^(1/4) - 1],
    Total[Table[
      2 h[i], {i, 2 (2^(1/4) - 1)/90,
       88 (2^(1/4) - 1)/90, (2^(1/4) - 1)/45}]],
```

```
Total[Table[
    4 h[i], {i, (2^(1/4) - 1)/90,
    89 (2^(1/4) - 1)/90, (2^(1/4) - 1)/45}]]
```

Once again, we do not display the output due to its length. It is approximately equal to 0.0280905. This gives $\|\varphi_1\|^2$ correct to within $10^{-6}$. Since we divide by $\|\varphi_1\|^2$ in the Rayleigh quotient, and want an upper bound on the first eigenvalue, we want a lower bound on $\|\varphi_1\|^2$. We take
$$\frac{28089}{1000000} < \|\varphi_1\|^2.$$

Now we need to bound the $\|\nabla \varphi_1\|^2$, this time from above. This time, since we are working with the gradient, we must analyze
$$\|\nabla \varphi_1\|^2 = \int_0^{2^{\frac{1}{4}}-1} \int_{c_1(x)}^{c_2(x)} \frac{\langle \nabla \varphi, \nabla \varphi \rangle}{y^2} dy dx.$$

Note that we have the following formula for taking the scalar product of a vector $\mathbf{X} = X^n \frac{\partial}{\partial x^n}$ of dimension $n$ with itself, with respect to the Riemannian metric:
$$\langle \mathbf{X}, \mathbf{X} \rangle = g_{ij} X^i X^j,$$
where we use the Einstein summation convention. We also have that for the gradient of a function $f(x_1, \ldots, x_n)$,
$$\nabla f = \sum_{j=1}^{n} \sum_{i=1}^{n} g^{ij} \frac{\partial f}{\partial x^i} \frac{\partial}{\partial x^j}.$$

For a function $f(x, y)$ in $\mathbb{H}$, we have
$$\langle \nabla f, \nabla f \rangle = g_{11}(g^{11} \partial_x f)^2 + g_{22}(g^{22} \partial_y f)^2$$
$$= g^{11}(\partial_x f)^2 + g^{22}(\partial_y f)^2$$
$$= y^2 \left((\partial_x f)^2 + (\partial_y f)^2\right).$$

Therefore, when we integrate with respect to the hyperbolic metric, the $y^2$ term cancels. We are essentially integrating the Euclidean gradient of a function of two variable with respect to the Euclidean metric. As before, we can integrate explicitly with respect to the $y$ variable to define a function in $x$:
$$g(x) = \int_{c_1(x)}^{c_2(x)} \left(\left(-\sqrt{1-22x^2} - \sqrt{q^2-x^2} + 2y - 1\right)^2 \right.$$
$$\left. + \left(\frac{x\left(-\sqrt{1-22x^2} + y - 1\right)}{\sqrt{q^2-x^2}} + \frac{22x\left(y - \sqrt{q^2-x^2}\right)}{\sqrt{1-22x^2}}\right)^2\right) dy.$$

This can be evaluated using the following in Mathematica; again we suppress the lengthy output.

```
q = Sqrt[3 + 2 Sqrt[2] - 2 Sqrt[4 + 3 Sqrt[2]]];
u[x_, y_] = (y - 1 - Sqrt[1 - 22 x^2]) (y - Sqrt[q^2 - x^2]);
Integrate[
 Grad[u[x, y], {x, y}].Grad[u[x, y], {x, y}], {x, 0, 2^(1/4) - 1}, {y,
   Sqrt[q^2 - x^2], Sqrt[Csc[Pi/8]^2 - (x + Cot[Pi/8])^2]}]
```

As before, we find a value for $n$, based on the maximum value of $|g^{(4)}(x)|$, to ensure $S_n$ is accurate enough. Once again, we find that we take the maximum value of $x$ in our interval, that is, $2^{\frac{1}{4}} - 1$, to achieve this maximum. Our error formula is

```
g[x_]=Integrate[
 Grad[u[x, y], {x, y}].Grad[u[x, y], {x, y}], {x, 0, 2^(1/4) - 1}, {y,
   Sqrt[q^2 - x^2], Sqrt[Csc[Pi/8]^2 - (x + Cot[Pi/8])^2]}]
adg[x_] = Abs[D[g[x], {x, 4}]];
((10^6) (adg[2^(1/4) - 1] (2^(1/4) - 1)^5)/180)^(1/4)
```

This time we find that we must take $n = 150$ in Simpson's formula. Our approximation is calculated with the following in Mathematica:

```
gradphi=((2^(1/4) - 1)/(450)) Total[{g[0], g[2^(1/4) - 1],
   Total[Table[
     2 g[i], {i, 2 (2^(1/4) - 1)/150,
      148 (2^(1/4) - 1)/150, (2^(1/4) - 1)/75}]],
   Total[Table[
     4 g[i], {i, (2^(1/4) - 1)/150,
      149 (2^(1/4) - 1)/150, (2^(1/4) - 1)/75}]]}];
```

The symbolic form of this is approximately equal to 0.1164682, correct again to within $10^{-6}$. We need an upper bound for $\|\nabla\varphi_1\|^2$, so we take

$$\|\nabla\varphi_1\|^2 < \frac{116469}{1000000}.$$

We can finally bound the first positive eigenvalue of the Bolza surface by

$$\lambda_1 \leq \mathcal{R}(\varphi_1) < \frac{116469}{28089} \approx 4.1464274.$$

It is possible to bound higher eigenvalues using the min-max principle. We will need an upper bound for $\lambda_2$ in order to prove that it has multiplicity 4. The Rayleigh quotient of a function $\varphi_2$ that is orthonormal to $\varphi_1$ will give us such a bound, if we chose the function appropriately (that is, $\{\varphi_1, \varphi_2\}$ form a 2 dimensional subspace of $W_0^{1,2}$ as required by Theorem 1.22).

**Remark 3.7** *Why is it sufficient to only consider two test functions? There are three functions corresponding to the first eigenvalue on the entire surface, but here we have restricted ourselves to a fundamental domain $\Omega$ of a subgroup of the isometry group $\mathrm{Isom}(\mathcal{B})$. By reflection, we can extend an eigenfunction on this domain to one on the entire surface (see Remark 3.11). In particular, when we reflect the domain along its "free" edges, that is, the edges with a Neumann boundary condition, we get a larger region with closed*

Dirichlet boundary. It is well known that the first Dirichlet eigenvalue of a connected domain is simple (has multiplicity one) and strictly positive (see, for example, [27]), so there will only be one corresponding function whose Rayleigh quotient, when evaluated over $\Omega$, gives this eigenvalue. Therefore, we only need to consider a two dimensional subspace of $W_0^{1,2}$ to get an upper bound on the second positive eigenvalue of $\Omega$.

The test function we would like to work with is

$$\varphi_2(x,\,y) = \left( \frac{11y}{2} - \frac{\left(\frac{11y}{2}\right)^3}{3!} + \frac{\left(\frac{11y}{2}\right)^5}{5!} - \frac{\left(\frac{11y}{2}\right)^7}{7!} + \frac{\left(\frac{11y}{2}\right)^9}{9!} - \frac{\left(\frac{11y}{2}\right)^{11}}{11!} \right) \left( 1 - \frac{x^2}{100} \right).$$

Here the $1 - \frac{x^2}{100}$ factor is to represent the behaviour of an eigenfunction of $\lambda_2$ observed in a numerical plot. We started with $1 - x^2$, and found experimentally that reducing the size of the $x^2$ coefficient gave a better bound. Note that the first term in the product is clearly an approximation of $\sin\left(\frac{11y}{2}\right)$. We initially constructed an orthonormal set of test functions using trigonometric functions as an obvious choice (again experimenting to find an "optimal" coefficient; although we couldn't produce accurate enough bounds with these alone, we found that by modifying them with a factor derived from our numerical insight, we were able to produce good enough bounds. We use the Taylor series approximation so that we can evaluate $\langle \varphi_1, \varphi_2 \rangle$ (integrating a polynomial function is much easier with respect to the hyperbolic metric over this domain). Fortunately, the approximation gives a slightly better bound than the function itself.

We use the Gram-Schmidt procedure on $\varphi_2$ and $\varphi_1$ to construct a function $\tilde{\varphi}_2$ that is orthogonal to $\varphi_1$, by subtracting a "projection" of $\varphi_1$ from $\varphi_2$:

$$\tilde{\varphi}_2(x,\,y) = \varphi_2(x,\,y) - \frac{\langle \varphi_1, \varphi_2 \rangle}{\|\varphi_1\|^2} \varphi_1(x,\,y)$$

Let us consider the inner product in the numerator of the projection. As we have alluded, we may integrate it explicitly with respect to $y$:

```
q = Sqrt[3 + 2 Sqrt[2] - 2 Sqrt[4 + 3 Sqrt[2]]];
u[x_, y_] = (y - 1 - Sqrt[1 - 22 x^2]) (y - Sqrt[q^2 - x^2]);
v[x_, y_] = (11 y/2 - ((11 y/2)^3)/3! + ((11 y/2)^5)/
    5! - ((11 y/2)^7)/7! + ((11 y/2)^9)/9! - ((11 y/2)^11)/11!) (1 -
    x^2/100);
Integrate[(u[x, y] v[x, y])/y^2, {y, Sqrt[q^2 - x^2],
  Sqrt[Csc[Pi/8]^2 - (x + Cot[Pi/8])^2]}]
```

We define a function $p(x)$, $x \in [0, 2^{\frac{1}{4}} - 1]$, as the output of this integral and once again use Simpson's rule to approximate it. Similar computations to before show us that we need $n = 146$ for our desired degree of accuracy, and we proceed in the same vein:

```
p[x_]=Integrate[(u[x, y] v[x, y])/y^2, {y, Sqrt[q^2 - x^2],
  Sqrt[Csc[Pi/8]^2 - (x + Cot[Pi/8])^2]}];
P = ((2^(1/4) - 1)/(3 146)) Total[{p[0], p[2^(1/4) - 1],
```

```
    Total[Table[
        2 p[i], {i, 2 (2^(1/4) - 1)/146,
        144 (2^(1/4) - 1)/146, (2^(1/4) - 1)/73}]],
    Total[Table[
        4 p[i], {i, (2^(1/4) - 1)/146,
        145 (2^(1/4) - 1)/146, (2^(1/4) - 1)/73}]]}];
```

This constant $P$ is correct to with $10^{-6}$. We get an upper and lower bound bound on $\langle \varphi_1, \varphi_2 \rangle$:
$$-\frac{9662}{1000000} < \langle \varphi_1, \varphi_2 \rangle < -\frac{9660}{1000000}.$$
We have already seen the lower bound on $\|\varphi_1\|^2$, and we can easily obtain an upper bound from our knowledge of the error. We get
$$\frac{28089}{1000000} < \|\varphi_1\|^2 < \frac{28091}{1000000}.$$
Therefore we get can bound the projection coefficient from above and below by
$$-\frac{9662}{28089} < \frac{\langle \varphi_1, \varphi_2 \rangle}{\|\varphi_1\|^2} < -\frac{9660}{28091}.$$
We will also need to consider terms involving the square of the projection coefficient, so we include now the relevant bounds:
$$\frac{93315600}{789104281} < \left(\frac{\langle \varphi_1, \varphi_2 \rangle}{\|\varphi_1\|^2}\right)^2 < \frac{93354244}{788991921}.$$
We are now ready to bound the $L^2$ norm of $\tilde{\varphi}_2$ from below. We have
$$\|\tilde{\varphi}_2\|^2 = \int_0^{2^{\frac{1}{4}}} \int_{c_1(x)}^{c_2(x)} \frac{1}{y^2} \left( \varphi_2^2 - 2\frac{\langle \varphi_1, \varphi_2 \rangle}{\|\varphi_1\|^2} \varphi_1 \varphi_2 + \left(\frac{\langle \varphi_1, \varphi_2 \rangle}{\|\varphi_1\|^2}\right)^2 \varphi_1^2 \right) dy\, dx$$
$$\geq \int_0^{2^{\frac{1}{4}}} \left( \int_{c_1(x)}^{c_2(x)} \left(\frac{\varphi_2}{y}\right)^2 dy + 2 \left(\frac{9660}{28091}\right) p(x) + \frac{93315600}{789104281} h(x) \right) dx.$$
As before, we can integrate
$$\int_{c_1(x)}^{c_2(x)} \left(\frac{\varphi_2}{y}\right)^2 dy$$
explicitly using the following line of code in Mathematica, where the output defines a function $c(x)$:

```
q = Sqrt[3 + 2 Sqrt[2] - 2 Sqrt[4 + 3 Sqrt[2]]];
v[x_, y_] = (11 y/2 - ((11 y/2)^3)/3! + ((11 y/2)^5)/
    5! - ((11 y/2)^7)/7! + ((11 y/2)^9)/9! - ((11 y/2)^11)/11!) (1 -
    x^2/100);
Integrate[(v[x, y] v[x, y])/y^2, {y, Sqrt[q^2 - x^2],
  Sqrt[Csc[Pi/8]^2 - (x + Cot[Pi/8])^2]}]
```

We now use Simpson's rule to approximate
$$\int_0^{2^{\frac{1}{4}}} \left( c(x) + 2\left(\frac{9660}{28091}\right) p(x) + \frac{93315600}{789104281} h(x) \right) dx. \tag{3.1}$$

We find that we need to take $n = 194$ for the required degree of accuracy, and the following code generates our approximation. As before, the symbolic output is yet longer and so we do not display it.

```
c[x_] = Integrate[(v[x, y] v[x, y])/y^2, {y, Sqrt[q^2 - x^2],
    Sqrt[Csc[Pi/8]^2 - (x + Cot[Pi/8])^2]}];
v2[x_] = c[x] + 2 (9660/28091) p[x] + (9660/28091)^2 h[x];
V = ((2^(1/4) - 1)/(3 194)) Total[{v2[0], v2[2^(1/4) - 1],
    Total[Table[
      2 v2[i], {i, 2 (2^(1/4) - 1)/194,
        192 (2^(1/4) - 1)/194, (2^(1/4) - 1)/97}]],
    Total[Table[
      4 v2[i], {i, (2^(1/4) - 1)/194,
        193 (2^(1/4) - 1)/194, (2^(1/4) - 1)/97}]]}];
```

From this, we achieve the following lower bound:

$$\|\tilde{\varphi}_2\|^2 > \frac{252552}{1000000}.$$

**Remark 3.8** *Whilst this is an accurate bound for the integral in Equation (3.1), we have lost some overall accuracy through the propagation of error through division and squaring in our bounds for the projection coefficient.*

All that remains is to bound $\|\nabla \tilde{\varphi}_2\|^2$. This time we want an upper bound, so we use the alternative bounds on our projection coefficient to get the following inequality:

$$\|\nabla \tilde{\varphi}_2\|^2 \leq \int_0^{2^{\frac{1}{4}}} \left( \int_{c_1(x)}^{c_2(x)} \langle \nabla \varphi_2, \nabla \varphi_2 \rangle dy + 2 \left(\frac{9662}{28089}\right) \int_{c_1(x)}^{c_2(x)} \langle \nabla \varphi_2, \nabla \varphi_1 \rangle dy \right.$$
$$\left. + \frac{93354244}{788991921} \int_{c_1(x)}^{c_2(x)} \langle \nabla \varphi_1, \nabla \varphi_1 \rangle dy \right) dx$$

Note that here, the angular brackets are for the scalar product of vector fields rather than the $L^2$ inner product. The integrals with respect to $y$ can all be evaluated using Mathematica; indeed we have already analyzed the third. At this stage we are doing nothing new, so can speed up the process to define our function of $x$:

```
gradv2[x_] =
  Integrate[
    Grad[v[x, y], {x, y}].Grad[v[x, y], {x, y}], {y, Sqrt[q^2 - x^2],
      Sqrt[Csc[Pi/8]^2 - (x + Cot[Pi/8])^2]}] +
  2 (9662/28089) Integrate[
    Grad[u[x, y], {x, y}].Grad[v[x, y], {x, y}], {y, Sqrt[q^2 - x^2],
      Sqrt[Csc[Pi/8]^2 - (x + Cot[Pi/8])^2]}] + (9662/
    28089)^2 Integrate[
    Grad[u[x, y], {x, y}].Grad[u[x, y], {x, y}], {y, Sqrt[q^2 - x^2],
      Sqrt[Csc[Pi/8]^2 - (x + Cot[Pi/8])^2]}]
```

We need to set $n = 192$ to reach our desired degree of accuracy, and $S_{192}$ is calculated with the following, where "gradv2zero" is the value of the function at zero..

```
gradv2zero =
  1/474393442370555512667661926400000
    (42733474853681838723362084401119381839 +
      30216456180729532565343854445094502400 Sqrt[2] -
      14884598099207334170365080185497190400 Sqrt[4 + 3 Sqrt[2]] -
      10524501996471252904537914913220198400 Sqrt[8 + 6 Sqrt[2]] -
      14701142201104864428066135243859376480655624 Sqrt[
        6 + 4 Sqrt[2] - 4 Sqrt[4 + 3 Sqrt[2]]] -
      20790554698446527074663610993142635407996047 Sqrt[
        3 + 2 Sqrt[2] - 2 Sqrt[4 + 3 Sqrt[2]]] +
      72415726188913689841313372525831677977849996 Sqrt[
        24 + 17 Sqrt[2] - 8 Sqrt[4 + 3 Sqrt[2]] -
          6 Sqrt[8 + 6 Sqrt[2]]] +
      512056515273113952388877447270207152417122 4 Sqrt[
        2 (24 + 17 Sqrt[2] - 8 Sqrt[4 + 3 Sqrt[2]] -
          6 Sqrt[8 + 6 Sqrt[2]])]);
GradV = ((2^(1/4) - 1)/(3 192)) Total[{gradv2zero,
    gradv2[2^(1/4) - 1],
    Total[Table[
      2 gradv2[i], {i, 2 (2^(1/4) - 1)/192,
        190 (2^(1/4) - 1)/192, (2^(1/4) - 1)/96}]],
    Total[Table[
      4 gradv2[i], {i, (2^(1/4) - 1)/192,
        191 (2^(1/4) - 1)/192, (2^(1/4) - 1)/96}]]}]
```

The upshot of this is that we now have the upper bound we were seeking, namely,

$$\|\nabla \tilde{\varphi}_2\|^2 < \frac{1408244}{1000000},$$

and thus, an upper bound on the second positive eigenvalue of the Bolza surface:

$$\lambda_2 \leq \mathcal{R}(\tilde{\varphi}_2) < \frac{1408244}{252552} \approx 5.57606.$$

If a Dirichlet boundary problem corresponding to an irreducible representation gives a lower bound in the Faber-Krahn inequality that is greater than $\mathcal{R}(\varphi_1)$, it cannot appear in $\mathcal{E}_1$. We aim to prove that there are only two irreducible representations that allow a lower bound that is less than $\mathcal{R}(\varphi_1)$. One will be the trivial representation, which corresponds to the trivial eigenspace $\mathcal{E}_0$, the other to $\mathcal{E}_1$. We will use the notation $\lambda_1(\chi_i)$ to denote the first positive eigenvalue of the $i^{th}$ irreducible representation, and the notation $\Omega_i$ to denote the boundary of the subspace of $\mathcal{B}$ corresponding to boundary conditions coming from the $i^{th}$ irreducible representation. We will then show that the trivial representation cannot appear.

## 3.4 First eigenspace of the Bolza surface

**Theorem 3.9** *The dimension of the first eigenspace $\mathcal{E}_1$ of $\Delta$ is 3, and the irreducible representation appearing in this eigenspace is $\chi_8$.*

In order to calculate the dimension of the first eigenspace, we must use the irreducible representations of $\text{Isom}(\mathcal{B})$ to create boundary value problems on the surface. Then, using the Faber-Krahn inequality, we can determine which irreducible representations appear in $\mathcal{E}_1$ (counting multiplicity).

Since we can describe the isometry group as a finitely presented group on 4 elements, satisfying a set of relations, we can use GAP [26] to obtain information about the representation theory of the group. This is shown in Appendix A. To come up with boundary value problems corresponding to a particular irreducible representation, we need the matrices associated to each of the four generators of $\text{Isom}(\mathcal{B})$.

Given a fundamental domain $\Omega \subset \mathcal{B}$ with respect to the action of either $\text{Isom}(\mathcal{B})$ or one of its subgroups, a reflection in one of its sides will correspond to the action of some $\Phi \in \text{Isom}(\mathcal{B})$. Such an $\Omega$ consists of one or more $(2, 3, 8)$ triangles. We consider the following three domains:

- a single $(2, 3, 8)$ triangle, which is the fundamental domain given by the action of the whole isometry group;

- half of a $(4, 4, 4)$ triangle, corresponding the the subgroup of $\text{Isom}(B)$ generated by $R$, $S$, and $T$;

- a right-angled pentagon of area $\frac{\pi}{2}$, corresponding to the subgroup of $\text{Isom}(B)$ generated by $R^4$, $S$, and $T$.

For every irreducible representation $\chi_i$ of $\text{Isom}(\mathcal{B})$, we want to see how a non-zero vector is moved by $\chi_i$ to adjacent fundamental domains. The aim is to construct a vector $v$ from the basis of eigenvectors of the representation matrices, that satisfies either $\Phi v = v$ or $\Phi v = -v$; this will imply that either the normal derivative of the vector valued function, or the function itself, will vanish on $\partial \Omega$, hence corresponding to either a Neumann or Dirichlet boundary condition respectively. Different irreducible representations will give different boundary conditions on a fundamental domain. Since our aim is to use the Faber-Krahn inequality, we must have a domain with purely Dirichlet boundary conditions. If we have a mixed boundary value problem on $\Omega$, we extend the domain by reflecting in a boundary geodesic that has a Neumann condition. In this way, we obtain closed domains with only Dirichlet boundary conditions. The isoperimetric inequality is known to be satisfied in the hyperbolic plane, so for simply connected $\Omega$ with Dirichlet boundary conditions, we may apply Faber-Krahn. Before we use this method on the complicated isometry group of the Bolza surface, we give an example using the simplest group $\mathbb{Z}_2$ to make the method clear.

**Example 3.10 ($\mathbb{Z}_2$)** $\mathbb{Z}_2$ *is the isometry group of a kite $K$ with only one line of symmetry, generated by a single reflection $S$.*

Figure 3.5: Kite $K$ with one line of symmetry

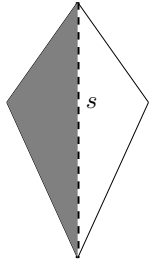

$\mathbb{Z}_2$ *has two irreducible representations, $M_1$ and $M_2$, both one dimensional. $M_1$ is the trivial representation, where $S_1$ takes the value 1. On $M_2$, $S_2$ takes the value $-1$. A fundamental domain for $\mathbb{Z}_2$ acting on $K$ is one of two scalene triangles (for example, the one shaded in Figure 3.5). We look at a vector valued function $f$ that takes value $v \neq 0$ on the shaded region, and see how it is affected by the group action. Since both irreducible representations are one dimensional, $v \in \mathbb{R}$.*

*When we map the shaded region to the white region, the value of $v$ on the white region becomes $Sv$. Restricting to $M_1$, $S_1 = 1$, so $S_1 v = v$. Hence all functions $f$ that correspond to $M_1$ are even with respect to the line of symmetry. On the other hand, restricting to $M_2$, $S_2 v = -v$, so the functions that correspond to $M_2$ are odd with respect to the line of symmetry. Functions corresponding to $M_2$ will therefore disappear on this line. If we were solving a Dirichlet boundary value problem on $K$, we could restrict the problem to the fundamental domain with Dirichlet boundary conditions in the case of $M_2$; in the case of $M_1$ we could restrict to a mixed boundary problem on the fundamental domain, with a Neumann condition on the dashed line.*

**Remark 3.11** *This technique was used in [41, 42] in precisely the same manner as we intend to use it, and was also used in [43] to investigate the relationship between geodesics on a hyperbolic surface, and the nodal lines of Laplace eigenfunctions. Given a surface $S$ with isometry group $\mathrm{Isom}(S)$, and a fundamental domain $\Omega$ for the action of a subgroup of $\mathrm{Isom}(S)$, it allows us to extend an eigenfunction $f$ of the mixed boundary problem on $\Omega$, to one on the whole surface $S$. The following sketched proof comes from [43]. For the Dirichlet case, let $f = 0$ on a geodesic boundary edge $B \in \partial \Omega$. Since the Laplacian commutes with the isometry group of $S$, for a reflection $r_B \in \mathrm{Isom}(S)$ about $B$ such that $f \circ r_B = -f$, the composition is also an eigenfunction of the Laplacian with the same eigenvalue. Let $g = f$, $g : S \to \mathbb{R}$ on one component of $S - B$ and $g = -f \circ r_B$ on the other. $g$ is now a weak solution to the eigenvalue problem, and by the elliptic* *regularity theorem, a strong solution also. $f = g$ on an open set, so $f = g$ on $S$, and $f = -f \circ r_B$. For a Neumann boundary condition, differentiating an even function gives an odd function and the hence the gradient of $v$ will vanish on this part of $\partial \Omega$, and we follow a similar argument.*

In this proof, we will systematically produce lower bounds on the first eigenspace in which the representation can appear. We will show that only $\chi_8$ can appear in the first eigenspace of $\mathcal{B}$. We will start with the 1 dimensional irreducible representations, $\chi_1, \ldots, \chi_4$. Since $\chi_1$ gives Neumann conditions on all three sides of a $(2, 3, 8)$ triangle, it will correspond to the first Neumann eigenspace, that is, the trivial eigenspace $\mathcal{E}_0$. We also have to show that it does not appear again in $\mathcal{E}_1$, but we will defer the proof of this to Lemma 3.15. In the following, we will denote by $\nu_x$ the first zero of the Legendre function

$$P_{-s}\left(\cosh\left(2\operatorname{arcsinh}\left(\sqrt{\frac{x}{4\pi}}\right)\right)\right) = P_{-s}\left(1 + \frac{x}{2\pi}\right),$$

where $x$ is fixed and we consider $P_{-s}$ as a function of $t$ with $s = \frac{1}{2} + it$. Then $\nu_x$ is the first eigenvalue of the hyperbolic disk of area $x$ with Dirichlet boundary condition (see the analysis in [38, Chapter 1]).

By the Faber-Krahn inequality, we can use $\nu_x$ as a lower bound for Dirichlet boundary problems on domains on area $x$ that satisfy the isoperimetric inequality (Definition 1.13). In turn, we need a lower bound on $\nu_x$ for each $x$. The largest domain whose first eigenvalue we bound using the Faber-Krahn inequality has area $\pi$. By the property of domain monotonicity for Dirichlet eigenvalues (Remark 1.7), a disk with area $\pi$ will have a lower first Dirichlet eigenvalue than any disk with smaller area, that is,

$$\nu_\pi < \nu_x \qquad \forall\, x < \pi.$$

Therefore, we just need to bound $\nu_\pi$ from below in order to get a lower bound on every other $\nu_x$ that we consider. Note that this bound will be far from optimal as $x$ decreases, but it suffices for it to be larger than our upper bound on the first eigenvalue of the whole surface, that is, $\frac{116469}{28089}$.

For $x = \pi$, we consider the first zero of the function $P_{-s}\left(\frac{3}{2}\right)$, where $s = \frac{1}{2} - it$. As a hypergeometric function, this is

$$f(t) = {}_2\mathbf{F}_1\left(\frac{1}{2} - it, \frac{1}{2} + it; 1; -\frac{1}{4}\right).$$

We can approximate $f(t)$ using its power series; recall the general definition from Section 2.3.1:

$${}_2\mathbf{F}_1(a, b; c; z) = \sum_{k=0}^\infty \frac{(a)_k (b)_k}{(c)_k k!} z^k,$$

and that this power series converges if $c$ is not a negative integer, and $|z| < 1$. As both of these conditions are satisfied for $f(t)$, we know that the functions given by the partial sums, that is,

$$f_n(t) = \sum_{k=0}^{n} \frac{\left(\frac{1}{2} - it\right)_k \left(\frac{1}{2} + it\right)_k}{(1)_k k!} \left(-\frac{1}{4}\right)^k,$$

approximate $f(t)$ more and more closely for each $n$. Considering the numerator

$$\left(\frac{1}{2} - it\right)_k \left(\frac{1}{2} + it\right)_k,$$

we see that for $k = 1$ we have

$$\left(\frac{1}{2} - it\right)\left(\frac{1}{2} + it\right) = \left(\frac{1}{4} + t^2\right),$$

for $k = 2$,

$$\left(\frac{1}{2} - it\right)\left(\frac{1}{2} + it\right)\left(\frac{3}{2} - it\right)\left(\frac{3}{2} + it\right) = \left(\frac{1}{4} + t^2\right)\left(\frac{9}{4} + t^2\right),$$

and so on. This means that the numerator will always be a polynomial with even powers of $t$; in particular it will be positive. Now note that the series is alternating, due to the $\left(-\frac{1}{4}\right)$ term. This means that for odd $n$, $f_n(t) < f(t)$, and for even $n$, $f_n(t) > f(t)$. Therefore, we can bound the first positive zero of $f(t)$ from below by the first positive zero of $f_n(t)$ for $n$ odd. It turns out that taking $f_3(t)$ already gives us a good enough lower bound.

We need to find the first positive real root of

$$f_3(t) = -\frac{t^6}{2304} + \frac{109t^4}{9216} - \frac{8035t^2}{36864} + \frac{15479}{16384}.$$

Note that this function is cubic in $t^2$, so we can use Cardano's formula to find a real value of $t^2$ such that $f_3(t) = 0$. Cardano's formula states that the real root of a general cubic equation

$$ax^3 + bx^2 + cx + d = 0$$

is given by

$$x = -\frac{1}{3a}\left(b + C + \frac{\Delta_0}{C}\right),$$

where

$$\Delta_0 = b^2 - 3ac,$$
$$\Delta_1 = 2b^3 - 9abc + 27a^2d, \quad \text{and}$$
$$C = \sqrt[3]{\frac{\Delta_1 \pm \sqrt{\Delta_1^2 - 4\Delta_0^3}}{2}}.$$

For $f_3(t) = 0$, we find that

$$t^2 = \frac{1}{12}\left(109 - \frac{1528}{\sqrt[3]{108\sqrt{789} - 1495}} + 8\sqrt[3]{108\sqrt{789} - 1495}\right)$$
$$= \frac{23}{4}$$

is a root, and thus $\frac{\sqrt{23}}{2}$ is the positive zero of $f_3(t)$. To check this, we have

$$f_3\left(\frac{\sqrt{23}}{2}\right) = -\frac{12167}{147456} + \frac{57661}{147456} - \frac{184805}{147456} + \frac{15479}{16384} = 0.$$

This gives us our lower bound:

$$\nu_\pi > \left(\frac{\sqrt{23}}{2}\right)^2 + \frac{1}{4} = 6.$$

For $\chi_2$, $R_2v = U_2v = v$ and $S_2v = T_2v = -v$, so any reflection will give a Dirichlet condition (since all reflections are orientation reversing). We therefore have a Dirichlet boundary problem on one of the (2, 3, 8) triangles, with area $\frac{\pi}{24}$, so we may use the Faber-Krahn inequality, which gives us

$$\lambda_1(\chi_2) \geq \nu_{\frac{\pi}{24}} > \nu_\pi > 6.$$

This rules out $\chi_2$.

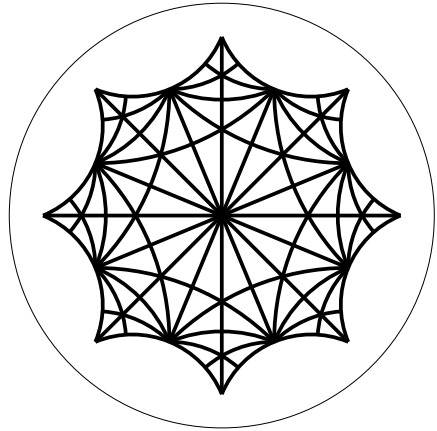

Figure 3.6: Boundary conditions for $\chi_2$

For $\chi_3$, $R_3v = T_3v = -v$ and $S_3v = U_3v = v$. The Dirichlet boundary problem is one of the equilateral (4, 4, 4) triangles, with area $\frac{\pi}{4}$. The Faber-Krahn inequality gives us

$$\lambda_1(\chi_3) \geq \nu_{\frac{\pi}{4}} > \nu_\pi > 6.$$

This rules out $\chi_3$.

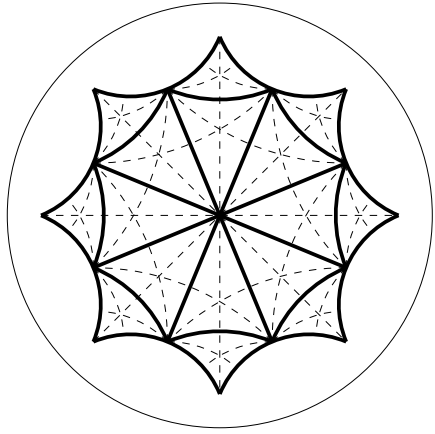

Figure 3.7: Boundary conditions for $\chi_3$

For $\chi_4$, $R_4 v = S_4 v = -v$ and $T_4 v = U_4 v = v$. In this case the Dirichlet boundary problem is the double of a (2, 3, 8) triangle (reflection along the hypotenuse). This has area $\frac{\pi}{12}$. The Faber-Krahn inequality gives us

$$\lambda_1(\chi_4) \geq \nu_{\frac{\pi}{12}} > \nu_\pi > 6.$$

This rules out $\chi_4$.

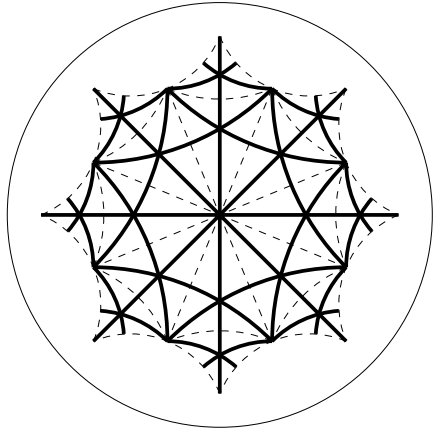

Figure 3.8: Boundary conditions for $\chi_4$

For the 2 and 3 dimensional irreducible representations, we consider the irreducible representations of the symmetry group (given in Appendix A after the first set of irreducible representations) acting on the standard basis $((1, 0)^T, (0, 1)^T$ for dimension 2, and similarly for dimension 3). These correspond to the basis of eigenfunctions in the eigenspace corresponding to a representation, for example, if $e_1$ is a function in the eigenspace corresponding to $\chi_5$, its vector representation will be $(1, 0)^T$.

For $\chi_5$, first notice that the matrix representation of $T_5$ is $-I_2$, where $I_n$ denotes the $n \times n$ identity matrix. Now see that $R_5 S_5 = -I_2$ and $R_5^2 = I_2$. We represent the function $\tilde{u} = e_1 - e_2$ as a vector $u = \begin{pmatrix} 1 \\ -1 \end{pmatrix}$, we see that $S_5 u = -u$. Therefore

$$R_5^i S_5 u = \begin{cases} R_5 S_5 u, & \text{if } i \in \{1, 3, 5, 7\}, \\ S_5 u, & \text{if } i \in \{0, 2, 4, 6\}, \end{cases}$$

and

$$R_5^i S_5 u = -u \qquad i \in \{0, \ldots, 8\}.$$

We see that the Dirichlet boundary problem is half of a $(4, 4, 4)$ triangle, with area $\frac{\pi}{8}$. The Faber-Krahn inequality gives us

$$\lambda_1(\chi_5) \geq \nu_{\frac{\pi}{8}} > \nu_\pi > 6.$$

This rules out $\chi_5$.

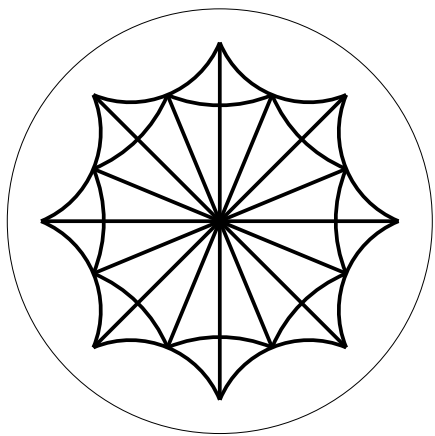

Figure 3.9: Boundary conditions for $\chi_5$

For $\chi_6$, note that $T_6 = R_6 S_6 = I_2$. We only need to consider the eigenvectors of $S_6$. Again, let $\tilde{u} = e_1 - e_2$, then $u = \begin{pmatrix} 1 \\ -1 \end{pmatrix}$, and

$$S_6 u = R_6^i S_6 u = -u, \qquad i \in \{0, 2, 4, 6\}.$$

Here we have a larger Dirichlet boundary problem; we use the side associations to see that it has area $\pi$. The Faber-Krahn inequality gives us

$$\lambda_1(\chi_6) \geq \nu_\pi > 6.$$

This rules out $\chi_6$.

Things are slightly trickier in 3 dimensions; we must consider more than one boundary value problem for each representation. We consider one of the 16 $(4, 4, 4)$ triangles that tessellates the Bolza surface, as well as its line of symmetry with respect to the real line.

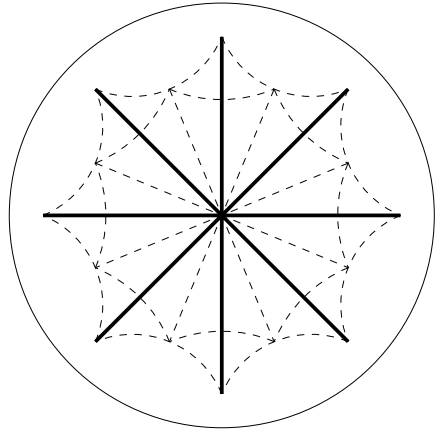

Figure 3.10: Boundary conditions for $\chi_6$

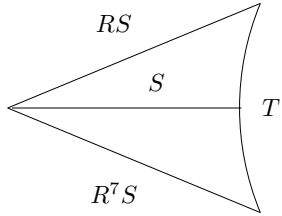

Figure 3.11: (4, 4, 4) triangle with line of symmetry

For $\chi_7$, we consider the vector $e_1 = (1,\ 0,\ 0)^T$. We have $R_7 S_7 e_1 = T_7 e_1 = e_1$, $S_7 e_1 = R_7^7 S_7 e_1 = -e_1$. For $e_2 = (0,\ 1,\ 0)^T$ we have $R_7 S_7 e_2 = S_7 e_1 = -e_2$, $R_7^7 S_7 e_2 = T_7 e_2 = e_2$, and for $e_3 = (0,\ 0,\ 1)^T$ we have $R_7 S_7 e_3 = R_7^7 S_7 e_3 = e_3$, $S_7 e_3 = T_7 e_3 = -e_3$. When closed, the Dirichlet boundary problem for each basis function, shown in Figure 3.12, contains an area of $\frac{\pi}{4}$. As above, Faber-Krahn shows us

$$\lambda_1(\chi_7) \geq \nu_{\frac{\pi}{4}} > \nu_\pi > 6.$$

This rules out $\chi_7$.

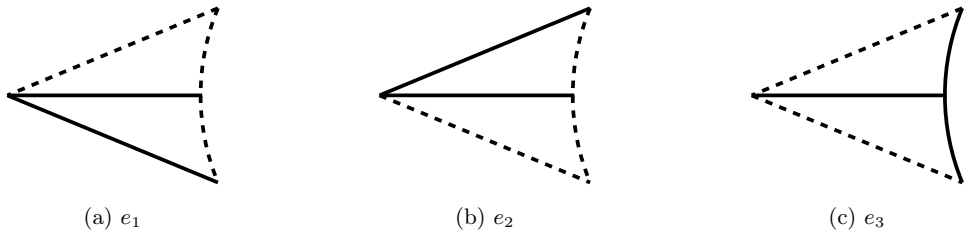

(a) $e_1$      (b) $e_2$      (c) $e_3$

Figure 3.12: Boundary conditions for $\chi_7$

For $\chi_8$, we again consider $e_1$, $e_2$, and $e_3$ as above. We have: $R_8 S_8 e_1 = S_8 e_1 = T_8 e_1 = e_1$ and $R_8^7 S_8 e_1 = -e_1$; $R_8^7 S_8 e_2 = S_8 e_2 = T_8 e_2 = e_2$ and $R_8 S_8 e_2 = -e_2$; and $R_8 S_8 e_3 = S_8 e_3 = R_8^7 S_8 e_3 = e_3$ and $T_8 e_3 = -e_3$. We observe that in each case the Dirichlet boundary contains an area of $2\pi$, that is, half of the surface. Courant's nodal domain theorem tells us that the first eigenfunction has exactly 2 nodal domains, so we now have a candidate for the first eigenspace. We would expect the Faber-Krahn inequality to produce a bound that does not rule out $\chi_8$. Indeed, it tells us that

$$\lambda_1(\chi_8) \geq \nu_{2\pi} \approx 3.66204.$$

This means that $\chi_8$ could appear in the first eigenspace.

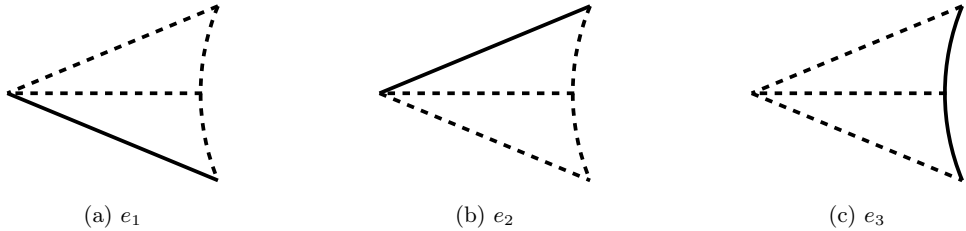

(a) $e_1$        (b) $e_2$        (c) $e_3$

Figure 3.13: Boundary conditions for $\chi_8$

For $\chi_9$, we have: $R_9 S_9 e_1 = S_9 e_1 = T_9 e_1 = -e_1$ and $R_9^7 S_9 e_1 = e_1$; $R_9^7 S_9 e_2 = S_9 e_2 = T_9 e_2 = -e_2$ and $R_9 S_9 e_2 = e_2$; and $R_9 S_9 e_3 = S_9 e_3 = R_9^7 S_9 e_3 = -e_3$ and $T_9 e_3 = e_3$. The Dirichlet boundary problems corresponding to $e_1$ and $e_2$ involve a triangle with area $\frac{\pi}{8}$, and $e_3$ corresponds to a triangle of area $\frac{\pi}{4}$. The Faber-Krahn inequality gives us

$$\lambda_1(\chi_9) \geq \nu_{\frac{\pi}{4}} > \nu_\pi > 6.$$

This rules out $\chi_9$.

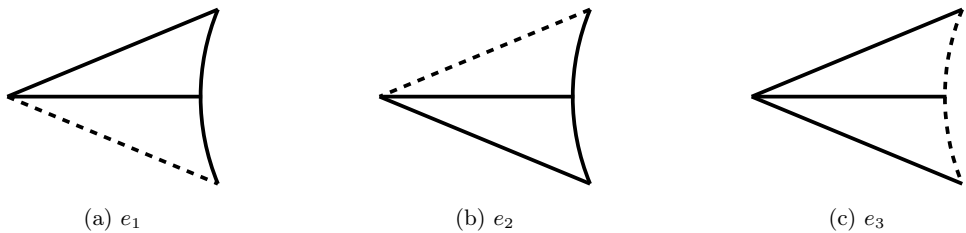

(a) $e_1$        (b) $e_2$        (c) $e_3$

Figure 3.14: Boundary conditions for $\chi_9$

Finally, for $\chi_{10}$ we have: $R_{10} S_{10} e_1 = T_{10} e_1 = -e_1$ and $R_{10}^7 S_{10} e_1 = S_{10} e_1 = e_1$; $R_{10}^7 S_{10} e_2 = T_{10} e_2 = -e_2$ and $R_{10} S_{10} e_2 = S_{10} e_2 = e_2$; and $R_{10} S_{10} e_3 = R_{10}^7 S_{10} e_3 = -e_3$ and $T_{10} e_3 = S_{10} e_3 = e_3$. In each case, the Dirichlet boundary problem corresponds to a quadrilateral with area $\frac{\pi}{2}$. The Faber-Krahn inequality yields

$$\lambda_1(\chi_{10}) \geq \nu_{\frac{\pi}{2}} > \nu_\pi > 6.$$

This rules out $\chi_{10}$.

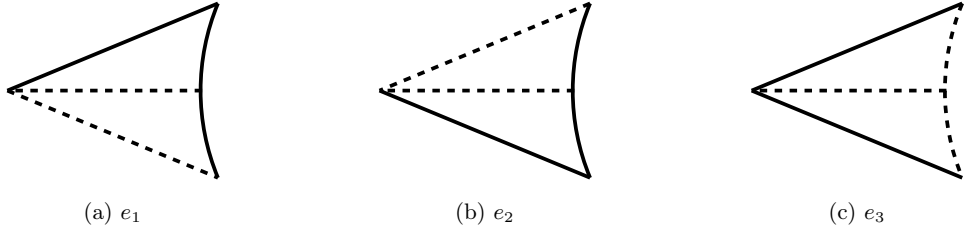

(a) $e_1$   (b) $e_2$   (c) $e_3$

Figure 3.15: Boundary conditions for $\chi_{10}$

The four dimensional irreducible representations of the isometry group of the Bolza surface correspond to the boundary conditions on a (non-regular) right-angled hyperbolic pentagon, of area $\frac{\pi}{2}$, shown in Figure 3.16. Note that due to the relations in Theorem 3.3, we have $R^2TR^6 = T$ and $RTR^7 = R^4T$. For $\chi_{11}$, the vector

$$\begin{pmatrix} \frac{1}{2}\left(-i - \sqrt{3}\right) \\ \frac{1}{2}\left((-2-i) + \sqrt{3}\right) \\ \frac{1}{4}\left((-2+i) + \sqrt{3}\right)\left(i + \sqrt{3}\right) \\ 1 \end{pmatrix}$$

satisfies the boundary conditions. For $\chi_{12}$ it is

$$\begin{pmatrix} \frac{1}{2}\left((1+2i) - i\sqrt{3}\right) \\ 1 \\ \frac{1}{2}\left(i - \sqrt{3}\right) \\ 1 \end{pmatrix},$$

and for $\chi_{13}$ it is

$$\begin{pmatrix} 1 \\ -1 \\ \left(\frac{1}{2} + \frac{i}{2}\right)\left(1 + \sqrt{3}\right) \\ 1 \end{pmatrix}.$$

Figure 3.16: Pentagon with boundary conditions

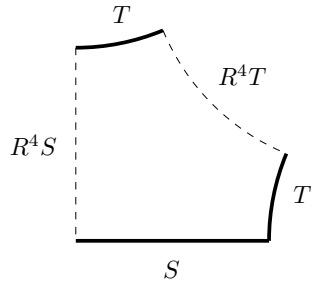

There is no way to create a region with closed Dirichlet boundary as we have done with the triangles corresponding to lower dimensional irreducible representations. We therefore cannot use the Faber-Krahn inequality to create a lower bound on the first eigenvalue of the pentagon.

On the other hand, we can reflect the pentagon along the vertical side, which has a Neumann boundary condition, using $R^4S$, to create a right-angled hyperbolic hexagon of area $\pi$. This is shown in Figure 3.17, along with its boundary conditions. Observe that the hexagon now has two Neumann boundaries of the same length, $a$ and $a'$. In particular, the side lengths of the hexagon are as follows:

$$l(a) = l(a') = l(\alpha) = 2\operatorname{arccosh}\left(\frac{\csc\left(\frac{\pi}{8}\right)}{2}\right) \approx 1.52857,$$

$$l(b) = l(b') = \operatorname{arccosh}\left(\frac{\csc\left(\frac{\pi}{8}\right)}{2}\right) \approx 0.764285,$$

$$l(c) = 2\operatorname{arccosh}\left(\frac{2\cos\left(\frac{\pi}{8}\right)}{\sqrt{3}}\right) + 2\operatorname{arccosh}\left(\frac{\cot\left(\frac{\pi}{8}\right)}{\sqrt{3}}\right) \approx 2.44845,$$

$$l(d) = \operatorname{arccosh}\left(\frac{2\cos\left(\frac{\pi}{8}\right)}{\sqrt{3}}\right) + \operatorname{arccosh}\left(\frac{\cot\left(\frac{\pi}{8}\right)}{\sqrt{3}}\right) \approx 1.22423.$$

An eigenfunction on the hexagon is even with respect to the vertical side $d$; it will take the same values on $a$ as it does on $a'$. Therefore we can "glue" $a$ to $a'$ to create an annulus $\Omega$ that lies within a hyperbolic cylinder. It has one geodesic boundary and one piecewise geodesic boundary composed of three geodesics. Eigenfunctions on this surface will be periodic and will satisfy Dirichlet conditions on the boundary.

Figure 3.17: Hexagon from reflection along the vertical

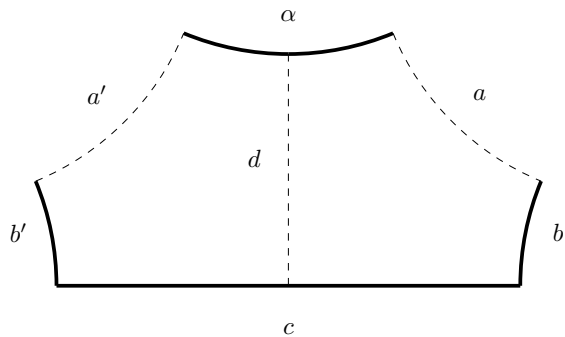

We have already seen that the circumference of a disk of area $\pi$ is

$$2\pi \sinh\left(2\operatorname{arcsinh}\left(\frac{1}{2}\right)\right) \approx 7.02481.$$

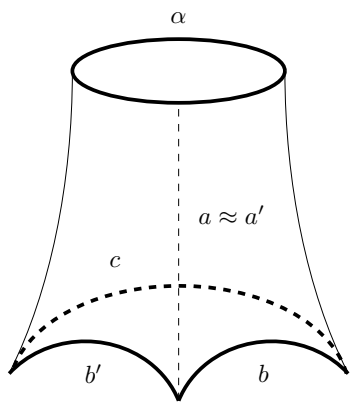

Figure 3.18: Cylinder from gluing $a$ to $a'$

On the other hand, the boundary of $\Omega$, after gluing, has length

$$2\operatorname{arccosh}\left(\frac{2\cos\left(\frac{\pi}{8}\right)}{\sqrt{3}}\right) + 2\operatorname{arccosh}\left(\frac{\cot\left(\frac{\pi}{8}\right)}{\sqrt{3}}\right) + 4\operatorname{arccosh}\left(\frac{\csc\left(\frac{\pi}{8}\right)}{2}\right) \approx 5.50559.$$

The annulus $\Omega$ does not satisfy the hypothesis of the Faber-Krahn inequality. In order to rule out the 4 dimensional irreducible representations, we need to make use of the following conjecture about the first Dirichlet eigenvalue of $\Omega$.

**Conjecture 3.12** *Let $\Omega$ be the annulus described above (alternatively, the right angled pentagon with boundary conditions as in Figure 3.16). Then $\lambda_1(\Omega) > 4.15$.*

**Remark 3.13** *There is strong numerical evidence to support this conjecture. Table C.2 shows the eigenvalues of $\Omega$ with boundary conditions as in Figure 3.16, calculated using FreeFEM++ (the corresponding code is also given in Appendix C). One can see that the first positive eigenvalue is approximately 5.35, and can compare this with the spectrum of $\mathcal{B}$ given in Table C.1, where the same eigenvalue appears with multiplicity 4. Unfortunately, we were not able to prove this conjecture rigorously, but in Section 3.4.1, we outline a method that could be used to produce the relevant bound.*

The truth of this conjecture would imply that $\lambda(\chi_{11})$, $\lambda(\chi_{12})$, and $\lambda(\chi_{13})$ do not appear in $\mathcal{E}_1$. Assuming this, we continue our proof of Theorem 3.9 by ruling out the trivial representation, $\chi_1$. We require the following corollary to Courant's nodal domain theorem (Theorem 1.12):

**Proposition 3.14** *[18, Proposition 1.1] For Neumann boundary problems, the number of nodal domains of the eigenfunction corresponding to the first positive eigenvalue is 2. For Dirichlet boundary problem problems, the number of nodal domains of the eigenfunction corresponding to the second positive eigenvalue is 2.*

**Lemma 3.15** *The trivial representation does not appear in $\mathcal{E}_1$.*

**Proof:** To prove this, we will consider the first positive Neumann eigenvalue of the (2, 3, 8) triangle. Proposition 3.14 implies that there will be a nodal line within the triangle that separates it into two domains. The line can intersect two of the boundary lines, one boundary line twice, or be contained entirely in the triangle. In the case where it is contained entirely, it already bounds a Dirichlet region in the triangle; by domain monotonicity it will certainly have a first eigenvalue larger than that of the (2, 3, 8) triangle with Dirichlet boundary conditions. In the case where it intersects one boundary twice, when we extend the eigenfunction to the entire surface, the nodal line will bound a Dirichlet region over two (2, 3, 8) triangles, and so will have a first eigenvalue larger than that of a hyperbolic disk of area $\frac{\pi}{12}$.

In the case where the nodal line intersects two edges of the triangle, we can quickly convince ourselves that the largest possible domain it can bound occurs when it intersects the two largest edges; when we extend such a nodal line to the surface, it will be contained in one of 6 small octagons, each of area $\frac{2\pi}{3}$, which tessellate the surface. By domain monotonicity, this Dirichlet region will have an eigenvalue larger than that of the $\frac{2\pi}{3}$ octagon. Using Faber-Krahn, we see that such an octagon $\Omega_O$ satisfies

$$\lambda_1(\Omega_O) \geq \nu_{\frac{2\pi}{3}} > \nu_\pi > 6.$$

Hence the trivial representation does not appear in the first eigenspace. □

The final thing we must show is that $\chi_8$ only appears once in the first eigenspace, that is, we know that the dimension of the first eigenspace is a multiple of 3, but it remains to prove that it is exactly 3. To do this, we prove that the multiplicity is strictly less than 6. We use the upper bound on the first eigenvalue (from Section 3.3) to create a bound on the multiplicity using Selberg's trace formula (see Theorem 2.47). Considering the spectral side of the formula, we have

$$h(t_0) + mh(t_1) < \sum_j h(t_j),$$

where $m$ is the multiplicity that we wish to bound. Rearranging this gives

$$m < \frac{\sum_j h(t_j) - h(t_0)}{h(t_1)},$$

and we will choose a function $h(t)$ such that

$$m < \frac{\sum_j h(t_j) - h(t_0)}{h(t_1)} < \frac{\sum_j h(t_j) - h(t_0)}{h\left(\sqrt{\frac{116469}{28089} + 1/4}\right)},$$

that is, it is decreasing on the interval containing $t_1$ and $\sqrt{\frac{116469}{28089} + 1/4}$ (recalling our bound on the first eigenvalue). We calculate the value of $\sum_j h(t_j)$ using the geometric side of Selberg's trace formula. We only have precise information about the multiplicity of the systole, so we use the test function

$$h(t) = \left(\frac{\sin(tL)}{tL}\right)^4,$$

as in [38], so that the Fourier transform

$$g(\xi) = \frac{1}{192L^4}\left((\xi - 4L)^3\text{sgn}(\xi - 4L) - 4(\xi - 2L)^3\text{sgn}(\xi - 2L) - 4(2L + \xi)^3\text{sgn}(2L + \xi)\right.$$
$$\left. + (4L + \xi)^3\text{sgn}(4L + \xi) + 6\xi^3\text{sgn}(\xi)\right)$$

is supported on $[-4L, 4L]$, where $4L$ is taken to be less than the length of the second primitive geodesic. Using the formula from Section 3.2, we know this is

$$2\operatorname{arccosh}\left(3 + 2\sqrt{2}\right) \approx 4.896904895356152.$$

In this way, we only get a contribution from the first length

$$2\operatorname{arccosh}\left(1 + \sqrt{2}\right) \approx 3.057141839.$$

We can choose any

$$L < \frac{1}{2}\operatorname{arccosh}\left(3 + 2\sqrt{2}\right),$$

however, we want our bound to be optimal. We can run a program in Mathematica (see Appendix F) to numerically plot how the bound on $m$ changes when we evaluate the geometric side of the formula for different values of $L$. In this way, we find that defining $h(t)$ as above with $L = \frac{93}{100}$ gives the strongest bound on $m$. To do our analysis, we alter Theorem 2.47 by using the equality

$$k(z, z) = -\frac{1}{2\pi}\int_0^\infty \frac{g'(t)}{\sinh\left(\frac{t}{2}\right)}dt$$

for the point pair invariant. In our case, $g(t)$ has compact support in $\left[-\frac{93}{25}, \frac{93}{25}\right]$, so the trace formula is

$$\sum_i h(\rho_i) = \int_0^{\frac{93}{25}} \frac{g'(t)}{\sinh\left(\frac{t}{2}\right)}dt + 24\left(\frac{l_1 g(l_1)}{2\sinh(l_1/2)}\right).$$

The right hand side of this equation can be evaluated explicitly using Mathematica. It is equal to the rather complicated

$$F = \frac{200}{74805201\sqrt{2+\sqrt{2}}}\left(25\sqrt{2+\sqrt{2}}\left(27900\operatorname{Li}_2\left(\frac{1}{e^{93/50}}\right) - 55800\operatorname{Li}_2\left(-\frac{1}{e^{93/100}}\right)\right.\right.$$

$$- 55800 \operatorname{Li}_2\left(\frac{1}{e^{93/100}}\right) + 5000 \operatorname{Li}_3\left(\frac{1}{e^{93/25}}\right) - 45000 \operatorname{Li}_3\left(\frac{1}{e^{93/50}}\right)$$

$$- 60000 \operatorname{Li}_3\left(-\frac{1}{e^{93/100}}\right) + 100000 \operatorname{Li}_3\left(\frac{1}{e^{93/100}}\right) - 105000\,\zeta(3)$$

$$- 25947 \log\left(e^{93/100} - 1\right) + 25947 \log\left(1 + e^{93/100}\right) + 69192 \tanh^{-1}\left(\frac{1}{e^{93/50}}\right)$$

$$- 51894 \tanh^{-1}\left(\frac{1}{e^{93/100}}\right) - 69192 \coth^{-1}\left(e^{93/100}\right) - 34596 \log\left(\tanh\left(\frac{93}{200}\right)\right)$$

$$+ 34596 \log\left(\tanh\left(\frac{93}{100}\right)\right) + 232500\sqrt{2+\sqrt{2}}\,\pi^2$$

$$+ 3217428(2^{3/4}) \cosh^{-1}\left(1+\sqrt{2}\right) - 5189400(2^{3/4}) \cosh^{-1}\left(1+\sqrt{2}\right)^2$$

$$+ 2790000\; 2^{3/4} \cosh^{-1}\left(1+\sqrt{2}\right)^3 - 500000(2^{3/4}) \cosh^{-1}\left(1+\sqrt{2}\right)^4\Bigg),$$

where

$$\operatorname{Li}_n(z) = \sum_{k=1}^{\infty} \frac{z^k}{k^n}$$

is the polylogarithm function, and

$$\zeta(s) = \sum_{k=1}^{\infty} k^{-s}$$

is the Riemann zeta function (for $\Re(s) > 1$).

The terms involving trigonometric functions, hyperbolic functions, or logs are fine, however we need to approximate the terms involving a power series in order to give an upper bound for $F$. Let us start by analyzing the terms with a dilogarithm function ($\operatorname{Li}_2(z)$). We may simplify the first two negative terms using the following functional equation for the dilogarithm,

$$\operatorname{Li}_2(z) + \operatorname{Li}_2(-z) = \frac{1}{2} \operatorname{Li}_2(z^2),$$

to get

$$-55800 \left( \operatorname{Li}_2\left(-\frac{1}{e^{93/100}}\right) + \operatorname{Li}_2\left(\frac{1}{e^{93/100}}\right) \right) = -27900 \operatorname{Li}_2\left(\frac{1}{e^{93/50}}\right),$$

which cancels with the dilogarithm term in the expression for $F$. Moving on to the trilogarithm terms, we can use the similar formula,

$$\operatorname{Li}_3(z) + \operatorname{Li}_3(-z) = \frac{1}{4} \operatorname{Li}_3(z^2),$$

to simplify

$$-45000 \operatorname{Li}_3\left(\frac{1}{e^{93/50}}\right) - 60000 \operatorname{Li}_3\left(-\frac{1}{e^{93/100}}\right) + 100000 \operatorname{Li}_3\left(\frac{1}{e^{93/100}}\right)$$

into
$$-60000 \operatorname{Li}_3\left(\frac{1}{e^{93/50}}\right) + 160000 \operatorname{Li}_3\left(\frac{1}{e^{93/100}}\right).$$

Our problem is now reduced to finding upper bounds for $\operatorname{Li}_3\left(e^{-\frac{93}{25}}\right)$ and $\operatorname{Li}_3\left(e^{-\frac{93}{100}}\right)$, and lower bounds for $\operatorname{Li}_3\left(e^{-\frac{93}{50}}\right)$ and $\zeta(3)$. For the upper bounds, note that $e^{-\frac{93}{25}}$ and $e^{-\frac{93}{100}}$ are both less than 1. The power series for the polylogarithm has a positive integer in the denominator that is greater than or equal to 1, so

$$\operatorname{Li}_3\left(e^{-\frac{93}{25}}\right) = \sum_{k=1}^{\infty} \frac{\left(e^{-\frac{93}{25}}\right)^k}{k^3} < \sum_{k=1}^{\infty} \left(e^{-\frac{93}{25}}\right)^k.$$

The term on the right is a geometric series, minus the first term, so the sum of the geometric series gives us

$$\operatorname{Li}_3\left(e^{-\frac{93}{25}}\right) < \frac{1}{1-e^{-\frac{93}{25}}} - 1 = \frac{1}{e^{\frac{93}{25}} - 1} \approx 0.0248358.$$

The same trick with $\operatorname{Li}_3\left(e^{-\frac{93}{100}}\right)$ does not give such a close bound, however we can make it as strong as we like by choosing where we start approximating with the geometric series, that is, we subtract the first $k+1$ terms of the geometric series and add the first $k$ terms of the trilogarithm series. With $k=8$, we get

$$\operatorname{Li}_3\left(e^{-\frac{93}{100}}\right) < \frac{1}{e^{93/100} - 1} - \frac{511}{512 e^{186/25}} - \frac{342}{343 e^{651/100}} - \frac{215}{216 e^{279/50}} - \frac{124}{125 e^{93/20}}$$
$$- \frac{63}{64 e^{93/25}} - \frac{26}{27 e^{279/100}} - \frac{7}{8 e^{93/50}}$$
$$\approx 0.417148.$$

For our lower bounds, we can simply truncate the power series at a point where a suitable degree of accuracy has been achieved. For $\operatorname{Li}_3\left(e^{-\frac{93}{50}}\right)$, truncating after 5 terms gives us

$$\operatorname{Li}_3\left(e^{-\frac{93}{50}}\right) > \frac{1}{125 e^{93/10}} + \frac{1}{64 e^{186/25}} + \frac{1}{27 e^{279/50}} + \frac{1}{8 e^{93/25}} + \frac{1}{e^{93/50}} \approx 0.158852,$$

and for $\zeta(3)$ summing the first 30 terms yields

$$\zeta(3) > \frac{1518061660370247564611888793148945960 3}{1263451477568240939757534871315200000 0} \approx 1.20152.$$

Our precise upper bound for $F$ is

$$\tilde{F} = \frac{200}{74805201\sqrt{2+\sqrt{2}}} \left(25\sqrt{2+\sqrt{2}} \left(5000\left(\frac{1}{e^{\frac{93}{25}}-1}\right) - 60000\left(\frac{1}{125 e^{93/10}}\right.\right.\right.$$
$$\left.\left.\left. + \frac{1}{64 e^{186/25}} + \frac{1}{27 e^{279/50}} + \frac{1}{8 e^{93/25}} + \frac{1}{e^{93/50}}\right) + 160000\left(\frac{1}{e^{93/100} - 1}\right.\right.\right.$$

$$
\begin{aligned}
&-\frac{511}{512e^{186/25}} - \frac{342}{343e^{651/100}} - \frac{215}{216e^{279/50}} - \frac{124}{125e^{93/20}} - \frac{63}{64e^{93/25}} \\
&-\frac{26}{27e^{279/100}} - \frac{7}{8e^{93/50}} \Bigg) \\
&- 105000 \left( \frac{151806166037024756461188793148945960 3}{1263451477568240939757534871315200000 0} \right) \\
&- 25947 \log\left(e^{93/100} - 1\right) + 25947 \log\left(1 + e^{93/100}\right) + 69192 \tanh^{-1}\left(\frac{1}{e^{93/50}}\right) \\
&- 51894 \tanh^{-1}\left(\frac{1}{e^{93/100}}\right) - 69192 \coth^{-1}\left(e^{93/100}\right) - 34596 \log\left(\tanh\left(\frac{93}{200}\right)\right) \\
&+ 34596 \log\left(\tanh\left(\frac{93}{100}\right)\right) + 232500 \sqrt{2 + \sqrt{2}} \pi^2 \\
&+ 3217428(2^{3/4}) \cosh^{-1}\left(1 + \sqrt{2}\right) - 5189400(2^{3/4}) \cosh^{-1}\left(1 + \sqrt{2}\right)^2 \\
&+ 2790000\ 2^{3/4} \cosh^{-1}\left(1 + \sqrt{2}\right)^3 - 500000(2^{3/4}) \cosh^{-1}\left(1 + \sqrt{2}\right)^4 \Bigg),
\end{aligned}
$$

which is approximately equal to 1.60267. We have

$$h(t_0) + mh(t_1) < \sum_j h(t_j),$$

where $m$ is the dimension of $\mathcal{E}_1$. We therefore have

$$h(t_0) + mh(t_1) < F < \tilde{F},$$

and a simple rearrangement gives

$$m < \frac{\tilde{F} - h(t_0)}{h(t_1)}.$$

Since $h(t)$ is decreasing on the interval $\left[0, \frac{93}{25}\right]$, we can apply our bound on the first eigenvalue to state

$$h(t_1) > h\left(\sqrt{\frac{116469}{28089} + 1/4}\right).$$

We conclude that

$$m < \frac{\tilde{F} - h(i/2)}{h\left(\sqrt{\frac{116469}{28089} + 1/4}\right)} \approx 5.87519,$$

where

$$h(i/2) = \frac{1600000000 \sinh^4\left(\frac{93}{200}\right)}{74805201}$$

and

$$h\left(\sqrt{\frac{116469}{28089} + 1/4}\right) = \frac{15585025600000000 \sin^4\left(\frac{93\sqrt{\frac{48643}{3121}}}{200}\right)}{176999686686876249}.$$

We have shown that $m$ is less than 6, so $\chi_8$ can only appear once in $\mathcal{E}_1$. The dimension of $\mathcal{E}_1$ is 3.

### 3.4.1 Method for proving Conjecture 3.12

Consider the annulus $\Omega$ in Figure 3.18 with Dirichlet boundary conditions. In particular, we can take $\Omega$ to be a subspace of a hyperbolic cylinder with $\alpha$ as its core geodesic. Table C.2 shows that numerically,

$$\lambda_1(\Omega) \approx 5.35.$$

Regarding the fact that 5.35 is relatively much larger than the figure of 4.15 that we want to bound it away from, the outlook for producing a rigorous bound seems good.

The method that we have in mind to prove Conjecture 3.12 is a generalization of the maximum principle. In particular, we have the following:

**Theorem 3.16** *[68, Theorem 10] Let $u(\boldsymbol{x})$ be an eigenfunction of the Laplacian with eigenvalue $\lambda$ on a compact Riemannian manifold $(M, g)$ with Lipschitz boundary $\partial M$. If there exists a function $w(\boldsymbol{x})$ that is strictly positive on $M \cup \partial M$, and*

$$(-\Delta + \lambda)u(\boldsymbol{x}) \leq 0 \qquad \text{in } M,$$

*then $u(\boldsymbol{x})/w(\boldsymbol{x})$ cannot attain a non-negative maximum in $M$ unless it is a constant. If $u(\boldsymbol{x})/w(\boldsymbol{x})$ attains its non-negative maximum at a point $p \in \partial M$ that lies on the boundary of a ball in $M$ and if $u/w$ is not constant, then*

$$\frac{\partial}{\partial \nu}\left(\frac{u}{w}\right) > 0 \qquad \text{at } p,$$

*where $\frac{\partial}{\partial \nu}$ is any outward directional derivative.*

**Corollary 3.17** *Let $(M, g)$ be as above, and assume $\partial M \neq \emptyset$. Assume there exists a strictly positive function $w \in C^\infty(M)$ with*

$$(-\Delta + \lambda)w \leq 0 \qquad \text{in } M.$$

*Then there is no Dirichlet eigenvalue in $[0, \lambda]$ in $M$.*

**Corollary 3.18** *Let $(M, g)$ be as above, and assume $\partial M \neq \emptyset$. Assume there exists a function $w \in C^\infty(M)$ that is positive on the interior of $M$ and non-negative on $M$, satisfying*

$$(-\Delta + \lambda)w \leq 0 \qquad \text{in } M.$$

*Then there is no Dirichlet eigenvalue in $[0, \lambda)$ in $M$.*

**Proof:** [Sketch of proof of Corollary 3.18] If we have a function $w$ that is strictly positive on the interior of $M$, then we can make $M$ a tiny bit smaller to obtain a region $M_\delta \subset M$.

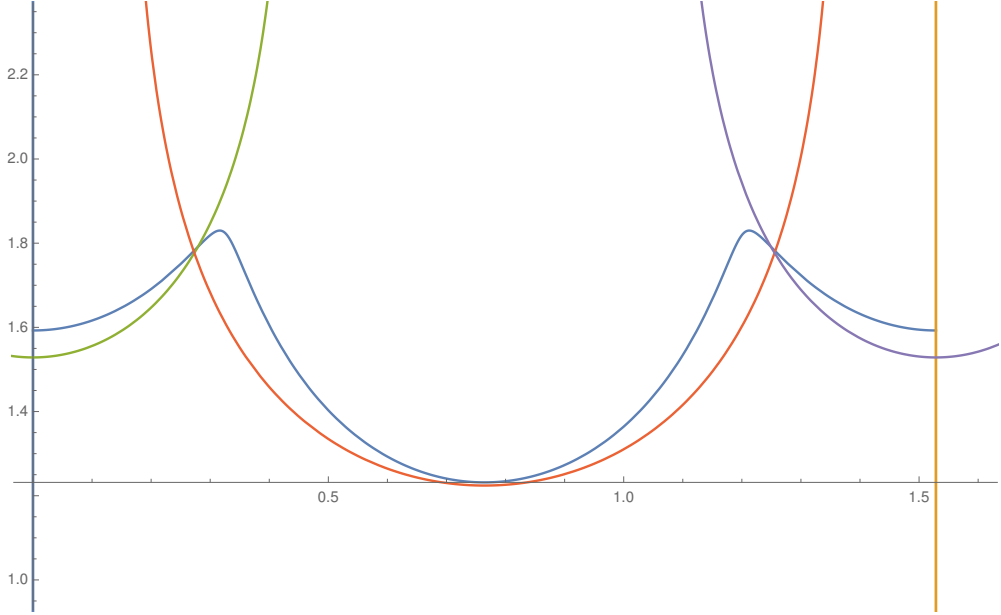

Figure 3.19: Comparing the nodal line of the test function $w$, as a function of $t$ (blue curve), to the boundary of $\Omega$ (red, lilac, green curves)

Now on $M_\delta$ we can construct a function $w_\delta$ that is equal to $w$ on $M_\delta$, and cut it off smoothly in the region $M \setminus M_\delta$ by multiplying by a cut off function $\chi \in C_0^\infty(M)$, such that $w_\delta$ is zero on $\partial M$. On $M_\delta$, $w_\delta$ is strictly positive and sub-harmonic, so there is no Dirichlet eigenvalue in $[0, \lambda]$ on $M_\delta$ by Corollary 3.17. Then we can create a sequence of domains $M_\delta$ getting progressively closer to $M$ such that the first eigenvalues of $M_\delta$ will converge to the first eigenvalue of $M$ as $\delta$ goes to zero. In particular, we use the sequence of corresponding functions $w_\delta \chi$ in the Rayleigh quotient to bound the first Dirichlet eigenvalue of $M$ from above (recall that by domain monotonicity, $\lambda_1(M) \leq \lambda_1(M_\delta)$). The gradient of $w_\delta \chi$ only differs from that of $w$ in an increasingly smaller area, and we can make its contribution to the Rayleigh quotient arbitrarily small (that is, such that there is no Dirichlet eigenvalue on $M$ in the interval $[0, \lambda + \epsilon]$ for $\epsilon > 0$). We obtain that there is no eigenvalue in $[0, \lambda)$. $\square$

Now for our annulus $\Omega$, numerical experiments show that if we take a certain combination of the first three odd eigenfunctions of the cylinder with core geodesic $\alpha$ (see Section 2.3.2) and eigenvalue $\lambda = 5$, we can find a test function that satisfies the conditions of Corollary 3.18 on $\Omega$. Note that along $\alpha$, such functions are zero by construction but the test function we create will be positive everywhere else on $\Omega$. Thus we show that there is no Dirichlet eigenvalue in the interval $[0, 5)$. In particular, the eigenfunctions are

$$\phi_0(\rho) = \sinh(\rho)\,_2\mathbf{F}_1\left(\frac{3}{4} + i\frac{\sqrt{\lambda - 1/4}}{2}, \frac{3}{4} + i\frac{\sqrt{\lambda - 1/4}}{2}; \frac{3}{2}; -\sinh^2(\rho)\right),$$

$$\phi_1(\rho) = \sinh(\rho)(\cosh(\rho))^{\frac{2\pi i}{l}} {}_2\mathbf{F}_1\left(\frac{3}{4} + i\frac{\sqrt{\lambda - 1/4}}{2} + \frac{\pi i}{l}, \frac{3}{4} + i\frac{\sqrt{\lambda - 1/4}}{2} + \frac{\pi i}{l}; \frac{3}{2}; -\sinh^2(\rho)\right),$$

$$\phi_2(\rho) = \sinh(\rho)(\cosh(\rho))^{\frac{4\pi i}{l}} {}_2\mathbf{F}_1\left(\frac{3}{4} + i\frac{\sqrt{\lambda - 1/4}}{2} + \frac{2\pi i}{l}, \frac{3}{4} + i\frac{\sqrt{\lambda - 1/4}}{2} + \frac{2\pi i}{l}; \frac{3}{2}; -\sinh^2(\rho)\right),$$

where

$$l = l(\alpha) = 2\operatorname{arccosh}\left(\frac{\csc\left(\frac{\pi}{8}\right)}{2}\right).$$

The combination that we take is

$$w(\rho, t) = \phi_0(\rho) - 0.038 \sin\left(\frac{2\pi t}{l} - \frac{\pi}{2}\right) \phi_1(\rho) - 0.0002 \cos\left(\frac{4\pi t}{l}\right) \phi_2.$$

This function is sub-harmonic over $\Omega$; in particular

$$(-\Delta + \lambda)w(\rho, t) = 0$$

over $\Omega$, since $\phi_0$, $\phi_1$, and $\phi_2$ have been constructed as eigenfunctions of the Laplacian (with respect to the cylinder metric described in Section 2.3.2). It is not immediately obvious that $w$ is strictly positive over $\Omega$, but we can plot the nodal line of $w$ as a function of $t$ to see at which values of $\rho$ the function is zero, and then compare it to the boundary of $\Omega$. This comparison is shown in Figure 3.19. The blue curve shows the values of $\rho$ where $w$ changes sign as $t$ runs over the length of the core geodesic. The green, red and lilac curves show the boundary of $\Omega$ in Fermi coordinates. Close inspection shows that the nodal line is always outside the boundary of $\Omega$. The code that produced this plot is given in Appendix F.

**Remark 3.19** *The proof of Conjecture 3.12 now amounts to showing that the test function $w(\rho, t)$ is positive in $\Omega$. One could perhaps do this by approximating the hypergeometric functions by polynomials in $\rho$, which are easier to analyze. Here, care must be taken to ensure that the first and second derivatives are approximated well enough such that the approximating function is still sub-harmonic over $\Omega$. One can see from Figure 3.19 that the test function $w$ is only just positive inside $\Omega$, so finding a sufficiently accurate polynomial approximant could prove difficult. On the other hand, we have room to manoeuver with $\lambda$; the current value of 5 is still "far away" from the value of 4.15 that would be the lowest sufficient value to prove Conjecture 3.12. Changing $\lambda$ to a slightly smaller value may give more freedom to find a polynomial approximant and show that it satisfies the hypothesis of Corollary 3.17.*

**Remark 3.20** *At one point, we thought it would be possible to prove Conjecture 3.12 using a Faber-Krahn style inequality for half cylinders. This would amount to proving* *that the first Dirichlet eigenvalue of $\Omega$ is greater than the first Dirichlet eigenvalue of the symmetric half cylinder with $\alpha$ as its core geodesic, such that its area is equal to that of $\Omega$. Proving this condition seems to be less attainable than the more direct approach described above, however, we include the alternative method in Appendix G.*

## 3.5 Second eigenspace of Bolza surface

Having calculated lower bounds for the eigenvalues of each irreducible representation, we have done the bulk of the work to make a statement about the second eigenspace, $\mathcal{E}_2$. We begin with the following

**Lemma 3.21** *The dimension of $\mathcal{E}_2$ is a multiple of 4.*

**Proof:** Recall our upper bound on $\lambda_2$ in Section 3.3. The only irreducible representations that can appear in this eigenspace are $\chi_8$ and the four dimensional ones. To rule out $\chi_8$, we will return to Proposition 3.14. The eigenfunction corresponding to the second positive eigenvalue of the domain in Figure 3.13 has 2 nodal domains. So we have two Dirichlet boundary problems each with the same eigenvalue. The entire region has area $2\pi$, so unless the nodal line splits it equally, one of the nodal domains will have area less than $\pi$. Therefore, we can take $\pi$ as an upper bound on the area of at least one of the nodal domains, and use the Faber-Krahn inequality so conclude that the eigenvalue of the nodal domain, and therefore the whole region, is greater than $\nu_\pi$ and therefore greater than 6. This means that only the four dimensional irreducible representations can appear in $\mathcal{E}_2$.
□

**Theorem 3.22** *The dimension of $\mathcal{E}_2(\mathcal{B})$ is 4.*

**Proof:** By Lemma 3.21, we already know that the multiplicity of $\lambda_2$ is a multiple of 4. To prove that it is exactly four, we will prove that the multiplicity is less than 8. We will use Selberg's trace formula in the same way as for the multiplicity bound above, using our upper bounds on the first and second positive eigenvalues. Again, we use the test function
$$h(t) = \left(\frac{\sin(tL)}{tL}\right)^4,$$
so that the Fourier transform of $h(t)$ is supported on $[-4L, 4L]$. We choose $4L$ to be less than the length of the shortest closed geodesic of the Bolza surface to avoid contribution from the length spectrum. Recall that the shortest closed geodesic has length
$$2\operatorname{arccosh}\left(1 + \sqrt{2}\right) \approx 3.057141839.$$

Choosing
$$L = 3/4,$$
Selberg's trace formula becomes
$$\sum_j h(t_j) = 2 \int_0^\infty h(t) \tanh(\pi t) t dt.$$

Figure 3.20: $h(t)$ on $[0, 4]$

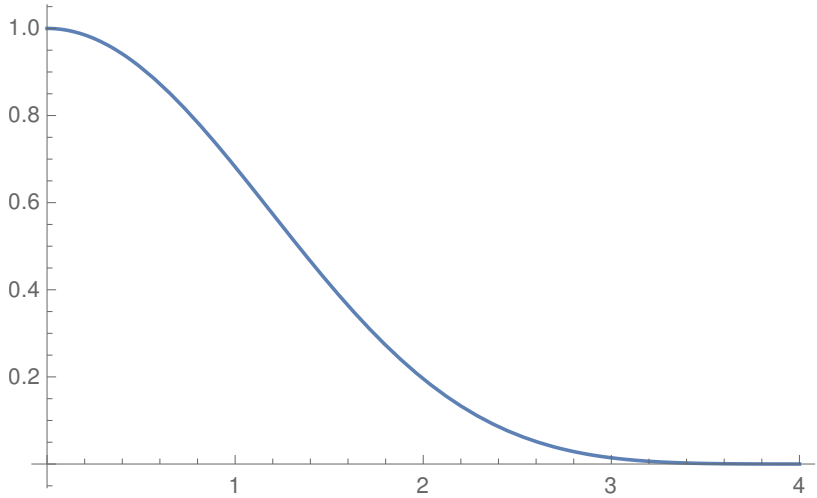

Now we use the straightforward inequality
$$\tanh(x) < 1$$
to state
$$\sum_j h(t_j) < 2 \int_0^\infty h(t) t dt.$$
The integral on the right hand side can be explicitly calculated; for $L = 3/4$, we have
$$2 \int_0^\infty h(t) t dt = \frac{32 \log(2)}{9}.$$
On the other hand, by Theorem 3.9 we have
$$h(t_0) + 3h(t_1) + mh(t_2) < \sum_j h(t_j),$$
where $m$ is the dimension of $\mathcal{E}_2$. We therefore have
$$h(t_0) + 3h(t_1) + mh(t_2) < \frac{32 \log(2)}{9},$$

and a simple rearrangement gives
$$m < \frac{\frac{32\log(2)}{9} - h(t_0) - 3h(t_1)}{h(t_2)}.$$

Since $h(t)$ is decreasing on the interval $[0, 4L]$, we can apply our bounds on the first two eigenvalues to state
$$h(t_1) > h\left(\sqrt{\frac{116469}{28089} + 1/4}\right)$$
and
$$h(t_2) > h\left(\sqrt{\frac{1408244}{252552} + 1/4}\right).$$

We conclude that
$$m < \frac{\frac{32\log(2)}{9} - h(i/2) - 3h\left(\sqrt{\frac{116469}{28089} + 1/4}\right)}{h\left(\sqrt{\frac{1408244}{252552} + 1/4}\right)}$$

$$= \frac{4070947840281 \csc^4\left(\frac{\sqrt{\frac{2017659}{10523}}}{8}\right)\left(\frac{32\log(2)}{9} - \frac{39897665536 \sin^4\left(\frac{3\sqrt{\frac{48643}{3121}}}{8}\right)}{63885819123} - \frac{4096}{81}\sinh^4\left(\frac{3}{8}\right)\right)}{453564534784}$$

$\approx 7.1091.$

$\square$

# Chapter 4

# The Klein Quartic

## 4.1 Systole of the Klein quartic

The systole of the Klein quartic is described in Klein's original paper [48] on the properties of the quartic curve, and is expounded upon in [44, 80]. We will use $\mathcal{K}$ to denote the surface. Consider a tessellation of $\mathcal{K}$ by heptagons (these can be barycentrically subdivided into 14 (2, 3, 7) triangles, see Figure 4.2). Starting with the central heptagon, join the midpoints of two adjoining edges with a geodesic. Extend this geodesic in both directions until it meets the boundary of the fundamental 14-gon. On its way, it passes through the midpoints of a further 3 and a half heptagon midpoints in either direction; a total of 8 midpoints. For this reason, Karcher and Weber coin it an "eight step geodesic". Note that the two boundary edges that it meets are actually the same edge once the side associations are made; that is, we have a simple closed geodesic. In fact, the length of this curve turns out to be the systole of the surface [44]. An order 7 rotation around the centre of the 14-gon maps this geodesic onto six other copies.

We obtain two more geodesics of length equal to the systole by creating curves that pass through adjoining edges of heptagons at varying distances from the centre; these are shown in Figure 4.2. There are 7 heptagons surrounding the central heptagon. The magenta curve passes through two edges of one of these. This time, when it reaches the boundary, it does not form a closed curve since the two sides are not associated with one another. We make the relevant side association and continue the geodesic from here, following the pattern of passing through midpoints of heptagon edges, until the geodesic is closed. Since we could have chosen any of the 7 heptagons, we get a further 7 closed geodesics. The orange (and green) curves pass through two adjoining edges of a heptagon in the second closest group from the centre. Again, there will be 7 of these. They follow the same pattern through 8 heptagon edge midpoints, thus having the same length. The multiplicity of the systole is 21.

Figure 4.1: Calculating $\frac{1}{8}$ of the systole

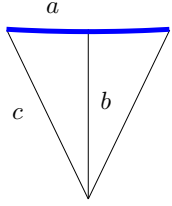

We calculate the systole using basic hyperbolic trigonometry to give a closed formula for one of the "steps", and then multiply this by 8. Similar to the formulae in Theorem 2.30 for right-angled hyperbolic triangles, we have

$$\sinh(a) = \sin(\alpha)\sinh(c) \qquad (4.1)$$

for a right-angled hyperbolic triangle where $c$ is the length of the hypotenuse and $\alpha$ is the angle opposite the side of length $a$ [13]. Also recall from Section 2.2.3 that the side of medium length in a (2, 3, 7) triangle has length

$$\operatorname{arccosh}\left(\frac{1}{2}\csc\left(\frac{\pi}{7}\right)\right).$$

Consider the blue geodesic in Figure 4.2. We can zoom in on its centre "step" to get the triangle shown in Figure 4.1. Since step is symmetric with respect to $b$, we can divide the larger triangle into two smaller right-angled triangles. Note that the medium side of a (2, 3, 7) triangle is the hypotenuse of the right-angled triangle in Figure 4.1. Let $a$ be half the length of the blue line, that is, the shortest side of the triangle on the left. We calculate the length of $a$ using Equation (4.1), where $\alpha$ is the angle opposite $a$ and is equal to $\frac{\pi}{7}$. We have

$$\sinh(a) = \sin\left(\frac{\pi}{7}\right)\sinh\left(\operatorname{arccosh}\left(\frac{1}{2}\csc\left(\frac{\pi}{7}\right)\right)\right)$$

$$= \sin\left(\frac{\pi}{7}\right)\sqrt{\left(-1+\frac{1}{2}\csc\left(\frac{\pi}{7}\right)\right)\left(1+\frac{1}{2}\csc\left(\frac{\pi}{7}\right)\right)}$$

$$= \sin\left(\frac{\pi}{7}\right)\left(\frac{1}{2}\sqrt{\csc^2\left(\frac{\pi}{7}\right)-4}\right).$$

Thus

$$l(a) = \operatorname{arcsinh}\left(\left(\frac{1}{2}\sqrt{\csc^2\left(\frac{\pi}{7}\right)-4}\right)\sin\left(\frac{\pi}{7}\right)\right).$$

This is half of the length of one step of the geodesic, and there are eight steps altogether, therefore

$$s(\mathcal{K}) = 16\operatorname{arcsinh}\left(\left(\frac{1}{2}\sqrt{\csc^2\left(\frac{\pi}{7}\right)-4}\right)\sin\left(\frac{\pi}{7}\right)\right) \approx 3.93594624883.$$

**Remark 4.1** *Note that there is more than one way to calculate a closed formula for the systole; [72] came up with*

$$8 \operatorname{arccosh}\left(\frac{3}{2} - 2\sin^2\left(\frac{\pi}{7}\right)\right) \approx 3.93594624883.$$

Figure 4.2: Geodesics of length equal to the systole

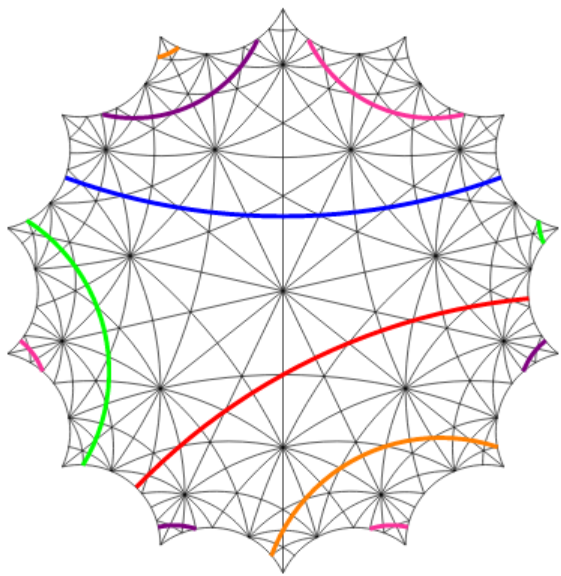

**Remark 4.2** *Numerical calculations of the other values in the geodesic length spectrum of the Klein quartic are given in the doctoral thesis of Vogeler [81], along with multiplicities. Vogeler also gives similar length spectra calculations for many other Hurwitz surfaces, including the Macbeath surface mentioned in the introduction to this work.*

## 4.2 Fenchel-Nielsen parameters of the Klein quartic

To discuss Fenchel-Nielsen parameters, we first need to introduce the notion of a pants decomposition. A pair of pants (sometimes called a Y-piece) is the simplest hyperbolic surface; topologically, it is a sphere with three disks removed. A compact Riemann surface of genus $g \geq 2$ can be decomposed into pairs of pants by cutting along $3g - 3$ simple, closed, non-intersecting geodesics. These geodesics form the boundaries of each pair of pants. A pair of pants has a specific geometry; we can form them by gluing together two right-angled (hyperbolic) hexagons, where alternate sides of each hexagon will be equal to half the length of one of the boundary geodesics of the pair of pants. The lengths of the sides in between these curves are fixed, since a geodesic joining two other geodesics perpendicularly is unique. This construction is shown in Figure 4.3; the bold curves on the pair of pants (on the left) are the boundary curves $\gamma_i$, $i \in \{1, 2, 3\}$, with length $l_i = l(\gamma_i)$. The pants are obtained by gluing two copies of the hexagon on the right along the dashed lines. The bold curves in the hexagon have length $\frac{l_i}{2}$. A complete derivation of the geometry of hyperbolic hexagons, staring from a hyperbolic triangle, is given in [13]. This is also an excellent reference for the fine details of pants decompositions.

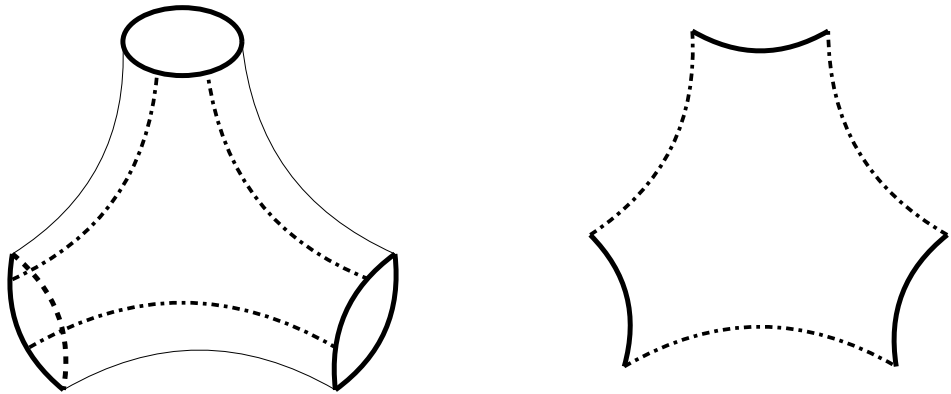

Figure 4.3: A hyperbolic pair of pants can be split into two identical hexagons

The Teichmüller space $\mathcal{T}_g$ is the space of all Riemann surfaces of genus $g$ up to marking equivalence. We will not go into depth here, but again, [13] provides a thorough introduction to this theory. The Fenchel-Nielsen parameters are a set of $6g - 6$ coordinates that completely classify a surface in Teichmüller space. Knowing the lengths of these geodesics provides half of these parameters. Additionally, for each geodesic we need to know whether it is twisted at all before being glued back together. This means that we also have $3g - 3$ twist parameters. The Fenchel-Nielsen parameters (or coordinates) are given as

$$\{l_1, t_1, \ldots, l_{3g-3}, t_{3g-3}\},$$

where $l_i$ and $t_i$, $i \in \{1, \ldots, 3g - 3\}$, respectively denote the length and twist parameters of a geodesic $\gamma_i$.

To obtain a decomposition of the Klein quartic into four pairs of pants, simply take six non-intersecting curves of length equal to the systole, as is shown in Figure 4.2. This decomposition is described in [44], along with a decomposition into hexagons that are not right-angled. The latter decomposition is interesting in that all the twist parameters are zero. An inspection of the number of regular heptagons enclosed by the curves in Figure 4.2 confirms that each pair of pants contains six heptagons. A heptagon is made up of 14 triangles of area $\frac{\pi}{42}$. Therefore each pair of pants has area

$$6(14)\left(\frac{\pi}{42}\right) = 2\pi.$$

Figure 4.4: Calculating the twist parameter $t_1$ of $\gamma_1$

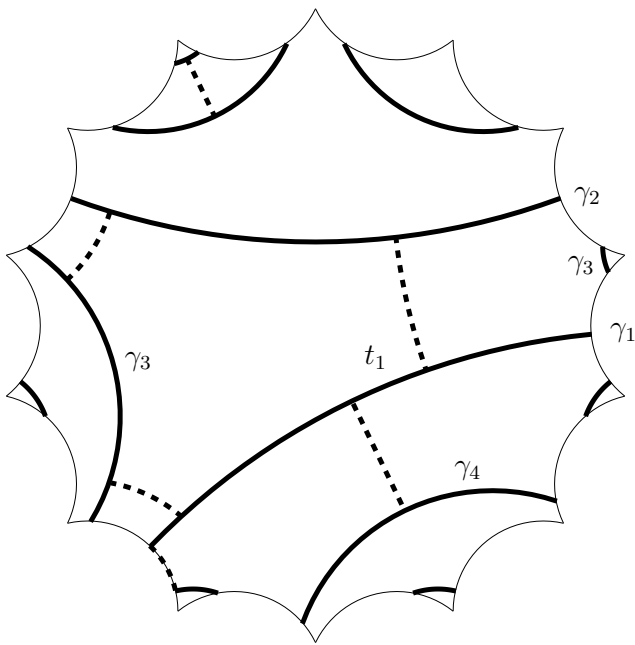

This is the correct area for a pair of pants having constant negative curvature.

We now calculate the six twist parameters. We start with the 'central' pair of pants, bounded by the geodesics $\gamma_1$, $\gamma_2$, and $\gamma_3$ in Figure 4.4. We draw the unique geodesic segment that meets $\gamma_1$ and $\gamma_2$ perpendicularly. To calculate the twist parameter $t_1$ of $\gamma_1$, consider the other pair of pants for which $\gamma_1$ is a boundary (the other two boundaries are $\gamma_4$ and $\gamma_5$). We draw the unique geodesic segment that meets $\gamma_1$ and $\gamma_4$ perpendicularly. When we twist along $\gamma_1$, we want the two perpendicular sections (dashed lines in the figure) to match up. The twist parameter for $\gamma_1$ is the distance along the $\gamma_1$ between the points where it intersects the perpendicular geodesics. Simply by examining Figure 4.2, we see that this is $\frac{1}{8}$ of the systole. Therefore

$$t_1 = \operatorname{arccosh}\left(\frac{3}{2} - 2\sin^2\left(\frac{\pi}{7}\right)\right) \approx 0.4919932811037915.$$

We can play a similar game with the other boundary curves to conclude that the twist parameters for each curve will be the same, that is, $\frac{1}{8}$ of the systole. Therefore, we have a complete set of Fenchel-Nielsen coordinates for the Klein quartic:

$$\left\{s(\mathcal{K}), \frac{s(\mathcal{K})}{8}, s(\mathcal{K}), \frac{s(\mathcal{K})}{8}, s(\mathcal{K}), \frac{s(\mathcal{K})}{8}, s(\mathcal{K}), \frac{s(\mathcal{K})}{8}, s(\mathcal{K}), \frac{s(\mathcal{K})}{8}, s(\mathcal{K}), \frac{s(\mathcal{K})}{8}\right\},$$

where $s(\mathcal{K})$ is the systole of $\mathcal{K}$.

## 4.3 Spectral theory of the Klein quartic

The symmetry group of the Klein quartic comes from the 336 (2, 3, 7) triangles that tessellate the surface. We saw the presentation of such a group in the hyperbolic geometry section; since there are side associations on our surface, we have the following presentation:

$$\langle a, b, c \,|\, a^2 = b^2 = c^2 = (ab)^2 = (ac)^7 = (bc)^3 = (abc)^8 = e \rangle,$$

where the lengths of $a$, $b$ and $c$ are given in Section 2.2. As with the Bolza surface, we can use this presentation in GAP to get information about the structure of the group, and more importantly, its irreducible representations (see Appendix B). We see that the isometry group (including reflections) of the Klein quartic has two one dimensional, three six dimensional, two seven dimensional, and two eight dimensional irreducible representations.

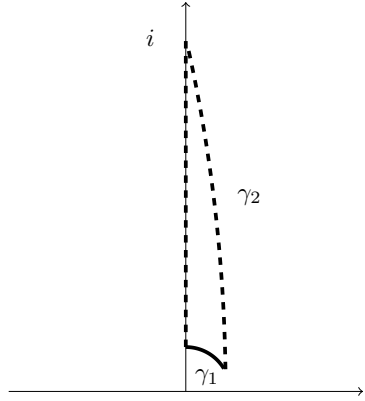

Figure 4.5: The domain of integration in $\mathbb{H}$

As with the Bolza surface, we can exploit the symmetry group and its representation theory to prove results about the spectrum of $\mathcal{K}$, which is given numerically in Table D.1. To begin with, we need an upper bound on the first eigenvalue of the surface; in Table D.1 we see that this is roughly 2.67. Following the strategy of Section 3.3, we construct a test function that is close to the first eigenfunction. This can then be used in the Rayleigh quotient. We work once again in the upper half plane model, although this time we use a slightly different triangle. Joining the vertices of the fundamental 14-gon with (Euclidean) straight lines tessellates the Klein quartic with 14 (7, 7, 7) triangles of area $\frac{4\pi}{7}$. We take one of these and map it to $\mathbb{H}$ using the Cayley transform on its vertices. We then use the fact that it has an additional line of symmetry with respect to the $\Im$ axis, and work with half the triangle, which has area $\frac{2\pi}{7}$. This triangle, denoted $\Omega_K$, is shown in Figure 4.5. The boundary curves are given by
$$\gamma_1(x) = \sqrt{r^2 - x^2}$$
and
$$\gamma_2(x) = \sqrt{\csc^2\left(\frac{\pi}{14}\right) - \left(x + \cot\left(\frac{\pi}{14}\right)\right)^2},$$
where
$$r = \sqrt{\frac{-(\cos\left(\frac{\pi}{14}\right) - 1)\exp\left(4\left(\cosh^{-1}\left(\frac{2\cos\left(\frac{\pi}{7}\right)}{\sqrt{3}}\right) + \cosh^{-1}\left(\frac{\cot\left(\frac{\pi}{7}\right)}{\sqrt{3}}\right) + \cosh^{-1}\left(\frac{1}{2}\csc\left(\frac{\pi}{7}\right)\right)\right)\right) + 1 + \cos\left(\frac{\pi}{14}\right)}{(1 + \cos\left(\frac{\pi}{14}\right))\exp\left(4\left(\cosh^{-1}\left(\frac{2\cos\left(\frac{\pi}{7}\right)}{\sqrt{3}}\right) + \cosh^{-1}\left(\frac{\cot\left(\frac{\pi}{7}\right)}{\sqrt{3}}\right) + \cosh^{-1}\left(\frac{1}{2}\csc\left(\frac{\pi}{7}\right)\right)\right)\right) + 1 - \cos\left(\frac{\pi}{14}\right)}}.$$

The test function we use is
$$\phi_1(x, y) = y - \sqrt{r^2 - x^2}.$$

Note that $\phi_1(x, y) = 0$ on $\gamma_1(x)$. Our aim is to bound its $L^2$ norm on $\Omega_K$ from below. We start by integrating it explicitly with respect to $y$, to define a function $u(x)$ as
$$u(x) = \int_{\gamma_1(x)}^{\gamma_2(x)} \left(\frac{\phi_1(x, y)}{y}\right)^2 dy$$
$$= 2\sqrt{r^2 - x^2}\log\left(\frac{\sqrt{r^2 - x^2}}{\sqrt{\csc^2\left(\frac{\pi}{14}\right) - \left(x + \cot\left(\frac{\pi}{14}\right)\right)^2}}\right)$$
$$- \frac{r^2 - x^2}{\sqrt{\csc^2\left(\frac{\pi}{14}\right) - \left(x + \cot\left(\frac{\pi}{14}\right)\right)^2}} + \sqrt{\csc^2\left(\frac{\pi}{14}\right) - \left(x + \cot\left(\frac{\pi}{14}\right)\right)^2}.$$

We now use Simpson's rule in exactly the same way as with the Bolza surface, with a maximum error of $10^{-6}$. To approximate
$$\int_0^p u(x) dx,$$
where
$$p = \frac{\sin\left(\frac{\pi}{14}\right)\left(\exp\left(4\left(\cosh^{-1}\left(\frac{2\cos\left(\frac{\pi}{7}\right)}{\sqrt{3}}\right) + \cosh^{-1}\left(\frac{\cot\left(\frac{\pi}{7}\right)}{\sqrt{3}}\right) + \cosh^{-1}\left(\frac{1}{2}\csc\left(\frac{\pi}{7}\right)\right)\right)\right) - 1\right)}{(1 + \cos\left(\frac{\pi}{14}\right))\exp\left(4\left(\cosh^{-1}\left(\frac{2\cos\left(\frac{\pi}{7}\right)}{\sqrt{3}}\right) + \cosh^{-1}\left(\frac{\cot\left(\frac{\pi}{7}\right)}{\sqrt{3}}\right) + \cosh^{-1}\left(\frac{1}{2}\csc\left(\frac{\pi}{7}\right)\right)\right)\right) + 1 - \cos\left(\frac{\pi}{14}\right)},$$
with this degree of accuracy, we need to take $n = 1648$. The following input in Mathematica [37] generates $S_{1648}$.

```
\[Phi][x_, y_] = y - Sqrt[r^2 - x^2];
u[x_] = (Sqrt[Csc[Pi/14]^2 - (x + Cot[Pi/14])^2] -
    2 Sqrt[r^2 - x^2] Log[
      Sqrt[Csc[Pi/14]^2 - (x + Cot[Pi/14])^2]] - (r^2 - x^2)/
      Sqrt[Csc[Pi/14]^2 - (x + Cot[Pi/14])^2]) - (Sqrt[r^2 - x^2] -
    2 Sqrt[r^2 - x^2] Log[Sqrt[r^2 - x^2]] - (r^2 - x^2)/
      Sqrt[r^2 - x^2]);
U=((p)/(3 1648)) Total[{u[0], u[p],
    Total[Table[2 u[i], {i, 2 (p)/1648, 1646 (p)/1648, (p)/824}]],
    Total[Table[4 u[i], {i, (p)/1648, 1647 (p)/1648, (p)/824}]]}]
```

As before, the symbolic output is too long to display here, but it gives us the bound
$$\|\phi_1\|^2 > \frac{29997}{1000000}.$$

Next, we bound $\|\nabla\phi_1\|^2$ from above. We integrate the scalar product of
$$\nabla\phi_1 = \left\{\frac{x}{\sqrt{r^2 - x^2}}, 1\right\}$$
with itself with respect to $y$ to get
$$gu(x) = \int_{\gamma_1}^{\gamma_2} \langle \nabla\phi_1, \nabla\phi_1 \rangle \, dy$$
$$= \left(\frac{x^2}{r^2 - x^2} + 1\right)\left(\sqrt{\csc^2\left(\frac{\pi}{14}\right) - \left(x + \cot\left(\frac{\pi}{14}\right)\right)^2} - \sqrt{r^2 - x^2}\right).$$

We approximate the integral of $gu(x)$ using Simpson's rule with $n = 2522$, that is,

```
gu[x_] = Integrate[
    Grad[\[Phi]][x, y], {x, y}].Grad[\[Phi]][x, y], {x, y}], {y,
    Sqrt[r^2 - x^2], Sqrt[Csc[Pi/14]^2 - (x + Cot[Pi/14])^2]}];
gU = ((p)/(3 2252)) Total[{gu[0], gu[p],
    Total[Table[
      2 gu[i], {i, 2 (p)/2252, 2250 (p)/2252, (xval)/1126}]],
    Total[Table[4 gu[i], {i, (p)/2252, 2251 (p)/2252, (p)/1126}]]}]
```

to get the following upper bound:
$$\|\nabla\phi_1\|^2 < \frac{84077}{1000000}.$$

Finally, we can take the quotient of our bounds to get an upper bound on the first positive eigenvalue of the Klein quartic, namely
$$\lambda_1 < \frac{84077}{29997} \approx 2.8027.$$

Comparedto the numerical value of 2.76, we see that this is a reasonably good bound. We get get one result quite quickly:

**Lemma 4.3** $\mathcal{E}_1$ *of the Klein quartic is not one dimensional.*

**Proof:** We start by ruling out the non-trivial one dimensional irreducible representation. We see that reflection in each side of a (2, 3, 7) triangle changes the sign of a vector as it is moved around the surface by this representation. Therefore, we have Dirichlet boundary conditions on a (2, 3, 7) triangle. Its first eigenvalue is greater than or equal to the first eigenvalue of a hyperbolic disk of area $\frac{\pi}{42}$ by the Faber-Krahn inequality. That is,
$$\lambda_1(\rho_2) \geq \nu_{\frac{\pi}{42}} > \nu_\pi > 6,$$
meaning $\rho_2$ does not appear in the first eigenspace.

Now we consider the trivial representation, and use an argument identical to that in Lemma 3.15, using Proposition 3.14 on nodal lines and domain monotonicity. We now have a heptagon $\Omega_H$ of area $\frac{\pi}{3}$, by the Faber-Krahn inequality

$$\lambda_1(\rho_1) \geq \lambda_1(\Omega_H) \geq \nu_{\frac{\pi}{3}} > \nu_\pi > 6.$$

Therefore, the trivial representation does not appear in the first eigenspace either. $\square$

We can also use Selberg's trace formula to bound the multiplicity of the first eigenvalue, like we did in Sections 3.4 and 3.5.

**Lemma 4.4** *The multiplicity of $\lambda_1$ is less than or equal to 10.*

**Proof:** As before, we use the test function

$$h(t) = \left(\frac{\sin(tL)}{tL}\right)^4,$$

so that the Fourier transform of $h(t)$ is supported on $[-4L, 4L]$. We choose $4L$ to be less than the length of $s(\mathcal{K})$ to avoid contribution from the length spectrum. Choosing

$$L = \frac{87}{100},$$

Selberg's trace formula becomes

$$\sum_j h(t_j) = 4 \int_0^\infty h(t) \tanh(\pi t) t \, dt.$$

Note that the coefficient of the integral is a 4 now rather than a 2, as it was in the case of the Bolza surface. This is because the area of the Klein quartic is twice the area of the Bolza surface. Again we use

$$\tanh(x) < 1$$

to state

$$\sum_j h(t_j) < 4 \int_0^\infty h(t) t \, dt.$$

The integral on the right hand side can be explicitly calculated; for $L = \frac{87}{100}$, we have

$$4 \int_0^\infty h(t) t \, dt = \frac{40000 \log(2)}{7569}.$$

We have

$$h(t_0) + mh(t_1) < \frac{40000 \log(2)}{7569}.$$

Rearranging this gives

$$m < \frac{\frac{40000 \log(2)}{7569} - h(t_0)}{h(t_1)}.$$

Since $h(t)$ is decreasing on the interval $[0, 4L]$, we can apply our bound on the first eigenvalue to state
$$h(t_1) > h\left(\sqrt{\frac{84077}{29997} + 1/4}\right).$$
We conclude that
$$m < \frac{\frac{40000\log(2)}{7569} - h(i/2)}{h\left(\sqrt{\frac{84077}{29997} + 1/4}\right)}$$
$$= \frac{6361650332217601 \csc^4\left(\frac{29\sqrt{\frac{306311}{3333}}}{200}\right)\left(\frac{40000\log(2)}{7569} - \frac{1600000000\sinh^4\left(\frac{87}{200}\right)}{57289761}\right)}{17774222400000000}$$
$$\approx 10.0855.$$

□

We can combine the multiplicity bound of Lemma 4.4 with Lemma 4.3 to state that any irreducible representation of $\mathrm{Isom}(\mathcal{K})$ that appears in $\mathcal{E}_1$ only does so once. This leads us to the following obvious observation:

**Corollary 4.5** *The eigenspace corresponding to the first positive eigenvalue of the Klein quartic is either 6, 7, or 8 dimensional.*

**Remark 4.6** *The inferred multiplicities from the numerical calculations of eigenvalues, shown in Table D.1, support this corollary. In particular, they suggest that the first eigenspace is 8 dimensional.*

# Chapter 5

# Concluding Remarks

The spectral theory of Riemann surfaces has been a popular and fruitful field of study for many years, and is still an active area of research. The main body of this work documented the investigation of the Laplace spectra of two specific surfaces: the Bolza surface and the Klein quartic. These are of particular interest because they have maximal symmetry among all compact Riemann surfaces with constant negative curvature, of genus 2 and 3 respectively. In particular, the Klein quartic is a Hurwitz surface.

In Chapter 3, we used the representation theory of its automorphism group to show that the first non-trivial Laplace eigenspace of the Bolza surface has dimension 3, and that this corresponds to the irreducible representation that we denoted $\chi_8$. This proof relies on a certain conjecture, which is supported by numerical evidence. We gave a potentially effective strategy for proving this conjecture, and went on to show that the second eigenspace has dimension 4. Unfortunately, we were not able to state the exact irreducible representation corresponding to this. We suspect that the refined version of the Selberg trace formula (Theorem 2.49) could be used to pinpoint the representation we are looking for.

Our work on the Bolza surface built on work done by Jenni in the 1980s, and our investigation was part of the initial step in addressing a wider open problem related to the Bolza surface, namely

**Conjecture 5.1** *Among all compact Riemann surfaces of genus 2 with constant negative curvature, the Bolza surface globally maximizes the first positive eigenvalue in the spectrum of the Laplace-Beltrami operator.*

It is thought that by analyzing the nodal lines of eigenfunctions relating to $\chi_8$ under perturbations in Teichmüller space, one would be able to show that the Bolza surface is a local maximum for $\lambda_1$. One would then be able to use this as a building block for proving that it is a global maximum.

In Chapter 4, we aimed to analyze the Klein quartic in much the same way as the Bolza surface, to prove results about the low lying eigenspaces. Analysis of the Klein quartic in this manner proved to be more complicated than the work of Chapter 3, although we proved an initial result that the dimension of the first eigenspace is not one dimensional. Combined with the bound given on the multiplicity of the first eigenvalue using Selberg's trace formula, and the dimensions of the irreducible representation of the automorphism group, we can say the the eigenspace is either 6, 7, or 8 dimensional. Also in this section, we investigated the systole and a pants decomposition of the Klein quartic, and gave Fenchel-Nielsen coordinates for the surface in Teichmüller space. As in the case of the Bolza surface, proving a result on the first eigenspace would have been an early step towards proving the following

**Conjecture 5.2** *Among all compact Riemann surfaces of genus 3 with constant negative curvature, the Klein quartic globally maximizes the first positive eigenvalue in the spectrum of the Laplace-Beltrami operator.*

In Appendix D we have carried out numerical experiments to support this conjecture. We calculated the spectrum of the surface denoted M3 in [72], which maximizes the systole length for genus 3 surfaces, and found it to have a smaller value of $\lambda_1$ than the Klein quartic. Likewise, we investigated the Fermat quartic, and found that its first positive eigenvalue is smaller than that of the Klein quartic. Of course, the natural temptation in light of the previous two conjectures is to postulate

**Conjecture 5.3** *For $n \geq 2$, among all compact Riemann surfaces of genus n with constant negative curvature, the surface with the largest automorphism group globally maximizes the first positive eigenvalue in the spectrum of the Laplace-Beltrami operator.*

However, this assumes that there is a unique surface that maximizes the automorphism group. If we restrict our attention to Hurwitz surfaces, we can conjecture that the Macbeath surface in genus 7 has this maximization property. Unfortunately, in the next genus that admits a Hurwitz surface, that is, genus 14, there is a triple of non-isomorphic Hurwitz surfaces (see, for example, [61]). It is unclear which of these would be the best candidate for eigenvalue maximization, but it seems reasonable to guess that one of them is.

Alongside the surface M3, other surfaces different to the focal two have been touched upon in this thesis. The examples in genus 2 have simpler automorphism groups than the Bolza surface, so theoretically, results about their eigenspaces could be gained without as much work. In genus 3, the Fermat quartic and M3 surface again have smaller automorphism groups than that of the Klein quartic (although they are comparable to, if not larger than, that of the Bolza surface) and could be another avenue of study to add to the collective knowledge of these surfaces.

Whilst pants decomposition are known for the Bolza surface and Klein quartic, we do not have explicit pants decompositions for the other surfaces. Again, knowledge of these, and the corresponding Fenchel-Nielsen parameters, would be helpful in carrying out accurate computations of eigenvalues, improving on the rudimentary finite element calculations given in the appendices.

# Appendix A

# GAP Code for Bolza Surface

The following output gives details of the isometry group of the Bolza surface, computed by the computer algebra package GAP [26]. We first define the group as a free group on the 4 generators introduced in Section 3.1 with the relevant relations. We confirm that the order of this group is 96, and the "StructureDescription" command tells us that it is isomorphic to $GL_2(\mathbb{Z}_3) \rtimes \mathbb{Z}_2$.

Representatives of each conjugacy class of Isom($\mathcal{B}$) are given by the command "ConjugacyClasses". Each row of the character table tell us the value of one of the characters on each of the 13 conjugacy classes of the group, where the conjugacy classes are in the same order as the output of "ConjugacyClasses". The data at the top of the character table concerns the various "power maps" of the character table, and is irrelevant to this work. The interested reader can learn how to interpret this data by consulting the GAP manual. To interpret the character values on the different conjugacy classes, one must also be aware of the following notation used by GAP:

$$E(x) := e^{\frac{2\pi i}{x}}.$$

The "IrreducibleAffordingRepresentation"'s below the first set of irreducible representations give a simpler set of representation matrices that are easier to work with, but correspond to the same character table. One may see more easily which of these correspond to rotations, and so on. Unfortunately it was not possible to calculate these for the four dimensional representations, but we make use of them for our analysis of the 2 and 3 dimensional representations in Section 3.4.

```
gap> f:=FreeGroup("r","s","t","u");
<free group on the generators [ r, s, t, u ]>
gap> g:=f/[f.1^8,f.2^2,f.3^2,f.4^3,f.1*f.2*f.1*f.2,f.2*f.3*f.2*f.3,f.1*f.3*f.1^3*f.3,
    f.1^7*f.4^2*f.1^7*f.4^2,f.4^2*f.1*f.4^2*f.3*f.2,f.4*f.2*f.4*f.2,
        f.4*f.3*f.4^2*f.2*f.1^7];
<fp group on the generators [ r, s, t, u ]>
gap> Order(g);
```

```
gap> ConjugacyClasses(g);
[ <identity ...>^G, r^2*t^G, s^G, u^G, r^3*s*t^G, r^4^G,
  r^2*s*t^G, r^2*t*u^G, r*s^G, r^5*t^G, r*s*t*u^G,
  r^3^G, r^7*t*u*r^G ]
gap> RequirePackage("repsn");
```
————————————————————————————————————

Repsn **for** Constructing Representations of Finite Groups
Version 3.0.2

Written by
Vahid Dabbaghian
————————————————————————————————————

**true**
```
gap> StructureDescription(g);
"GL(2,3):C2"
gap> Display(CharacterTable(g));
CT1

     2  5  4  3  2  4  5  3   2  4  3  2   3   2
     3  1  1  .  1  .  1  .   1  .  .  1   .   1

        1a 4a 2a 3a 4b 2b 2c 12a 2d 8a 6a 8b 12b
    2P  1a 2b 1a 3a 2b 1a 1a  6a 1a 4b 3a 4b  6a
    3P  1a 4a 2a 1a 4b 2b 2c  4a 2d 8a 2b 8b  4a
    5P  1a 4a 2a 3a 4b 2b 2c 12b 2d 8a 6a 8b 12a
    7P  1a 4a 2a 3a 4b 2b 2c 12b 2d 8a 6a 8b 12a
   11P  1a 4a 2a 3a 4b 2b 2c 12a 2d 8a 6a 8b 12b

X.1     1  1  1  1  1  1  1   1  1  1  1  1   1
X.2     1 -1 -1  1  1  1  1  -1 -1 -1  1  1  -1
X.3     1 -1  1  1  1  1 -1  -1 -1  1  1 -1  -1
X.4     1  1 -1  1  1  1 -1   1  1 -1  1 -1   1
X.5     2 -2  . -1  2  2  .   1 -2  . -1  .   1
X.6     2  2  . -1  2  2  .  -1  2  . -1  .  -1
X.7     3 -3 -1  . -1  3  1   .  1  1  . -1   .
X.8     3 -3  1  . -1  3 -1   .  1 -1  .  1   .
X.9     3  3 -1  . -1  3 -1   . -1  1  .  1   .
X.10    3  3  1  . -1  3  1   . -1 -1  . -1   .
X.11    4  .  . -2  . -4  .   .  .  .  2  .   .
X.12    4  .  .  1  . -4  .   A  .  . -1  .  -A
X.13    4  .  .  1  . -4  .  -A  .  . -1  .   A

A = -E(12)^7+E(12)^11
  = Sqrt(3) = r3
gap> IrreducibleRepresentations(g);
[ [ r, s, t, u ] -> [ [ [ 1 ] ], [ [ 1 ] ], [ [ 1 ] ], [ [ 1 ] ] ],
  [ r, s, t, u ] -> [ [ [ 1 ] ], [ [ -1 ] ], [ [ -1 ] ], [ [ 1 ] ] ],
  [ r, s, t, u ] -> [ [ [ -1 ] ], [ [ 1 ] ], [ [ -1 ] ], [ [ 1 ] ] ],
  [ r, s, t, u ] -> [ [ [ -1 ] ], [ [ -1 ] ], [ [ 1 ] ], [ [ 1 ] ] ],
  [ r, s, t, u ] -> [ [ [ 0, 1 ], [ 1, 0 ] ], [ [ 0, -1 ], [ -1, 0 ] ],
    [ [ -1, 0 ], [ 0, -1 ] ], [ [ 0, -1 ], [ 1, -1 ] ] ],
  [ r, s, t, u ] -> [ [ [ 1, 0 ], [ -1, -1 ] ], [ [ 1, 0 ], [ -1, -1 ] ],
    [ [ 1, 0 ], [ 0, 1 ] ], [ [ -1, -1 ], [ 1, 0 ] ] ],
```

```
  [ r, s, t, u ] -> [ [ [ -1, 0, -1 ], [ 1, 0, 0 ], [ 1, 1, 0 ] ],
    [ [ -1, 0, -1 ], [ 0, -1, 1 ], [ 0, 0, 1 ] ],
      [ [ 0, 1, -1 ], [ 1, 0, 1 ], [ 0, 0, 1 ] ],
      [ [ 1, 1, 0 ], [ 0, -1, 1 ], [ 0, -1, 0 ] ] ],
  [ r, s, t, u ] -> [ [ [ 1, 1, 1 ], [ -1, 0, -1 ], [ -1, -1, 0 ] ],
    [ [ 1, 0, 0 ], [ -1, 0, -1 ], [ -1, -1, 0 ] ],
      [ [ 1, 0, 0 ], [ 0, 0, 1 ], [ 0, 1, 0 ] ],
      [ [ 1, 1, 0 ], [ -1, -1, -1 ], [ -1, 0, 0 ] ] ],
  [ r, s, t, u ] -> [ [ [ 0, -1, 1 ], [ 1, 1, 0 ], [ 0, 1, 0 ] ],
    [ [ -1, 0, 0 ], [ 0, 0, -1 ], [ 0, -1, 0 ] ],
      [ [ -1, 0, 0 ], [ 1, 0, 1 ], [ 1, 1, 0 ] ],
      [ [ 1, 0, 1 ], [ 0, 0, -1 ], [ 0, 1, -1 ] ] ],
  [ r, s, t, u ] -> [ [ [ 0, -1, 0 ], [ 0, 0, -1 ], [ -1, 1, -1 ] ],
    [ [ 1, -1, 1 ], [ 0, 0, 1 ], [ 0, 1, 0 ] ],
      [ [ 0, 0, 1 ], [ 1, -1, 1 ], [ 1, 0, 0 ] ],
      [ [ 1, -1, 1 ], [ 0, -1, 1 ], [ 0, -1, 0 ] ] ],
  [ r, s, t, u ] -> [ [ [ 0, 0, -E(12)^7, 0 ], [ 0, 0, 0, -E(3) ],
    [ -E(12)^7, E(4), 0, 0 ], [ -1, -E(3), 0, 0 ] ],
      [ [ 0, 0, 0, -E(12)^7 ], [ 0, 0, -E(3), 0 ],
        [ 0, -E(3)^2, 0, 0 ], [ E(12)^11, 0, 0, 0 ] ],
      [ [ -E(12)^11, E(12)^7, 0, 0 ], [ E(12)^7, E(12)^11, 0, 0 ],
        [ 0, 0, E(12)^11, -E(3) ], [ 0, 0, E(3), -E(12)^11 ] ],
      [ [ E(3), 0, 0, 0 ], [ -1, E(3)^2, 0, 0 ], [ 0, 0, E(3), E(4) ],
        [ 0, 0, 0, E(3)^2 ] ] ],
  [ r, s, t, u ] -> [ [ [ 0, -E(3)^2, 0, 0 ], [ 1, 0, -E(3)^2, 0 ],
    [ 0, -1, 0, -E(12)^7 ], [ -E(12)^7, 0, 0, 0 ] ],
      [ [ 0, -E(4), 0, E(3) ], [ 0, 0, -E(12)^7, 0 ], [ 0, E(12)^11, 0, 0 ],
        [ E(3)^2, 0, 1, 0 ] ],
      [ [ E(12)^11, 0, -E(12)^7, 0 ], [ 0, 0, 0, -1 ],
        [ -E(12)^7, 0, -E(12)^11, 0 ], [ 0, -1, 0, 0 ] ],
      [ [ 0, 0, E(3)^2, 0 ], [ 0, -E(3)^2, 0, -E(4) ],
        [ -1, 0, -E(3), 0 ], [ 0, -E(12)^7, 0, 0 ] ] ],
  [ r, s, t, u ] -> [ [ [ 0, -E(12)^7, 0, E(3)^2 ], [ 0, 0, E(4), 0 ],
    [ 0, 0, 0, 1 ], [ -E(3)^2, 0, E(3), 0 ] ],
      [ [ 0, 1, 0, 0 ], [ 1, 0, 0, 0 ], [ 0, -E(3)^2, 0, E(12)^7 ],
        [ -E(12)^7, 0, -E(12)^11, 0 ] ],
      [ [ -E(12)^7, 0, -E(12)^11, 0 ], [ 0, 0, 0, 1 ],
        [ -E(12)^11, 0, E(12)^7, 0 ], [ 0, 1, 0, 0 ] ],
      [ [ 1, 0, E(3), 0 ], [ 0, -E(3)^2, 0, -E(4) ], [ 0, 0, E(3)^2, 0 ],
        [ 0, -E(12)^7, 0, 0 ] ] ] ]
gap> chis:=Irr(g);;
gap> IrreducibleAffordingRepresentation(chis[1]);
[ f1, f2, f3, f4 ] -> [ [ [ 1 ] ], [ [ 1 ] ], [ [ 1 ] ], [ [ 1 ] ] ]
gap> IrreducibleAffordingRepresentation(chis[2]);
[ f1, f2, f3, f4 ] -> [ [ [ 1 ] ], [ [ -1 ] ], [ [ -1 ] ], [ [ 1 ] ] ]
gap> IrreducibleAffordingRepresentation(chis[3]);
[ f1, f2, f3, f4 ] -> [ [ [ -1 ] ], [ [ 1 ] ], [ [ -1 ] ], [ [ 1 ] ] ]
gap> IrreducibleAffordingRepresentation(chis[4]);
[ f1, f2, f3, f4 ] -> [ [ [ -1 ] ], [ [ -1 ] ], [ [ 1 ] ], [ [ 1 ] ] ]
gap> IrreducibleAffordingRepresentation(chis[5]);
[ f1, f2, f3, f4 ] -> [ [ [ 0, -1 ], [ -1, 0 ] ], [ [ 0, 1 ], [ 1, 0 ] ],
  [ [ -1, 0 ], [ 0, -1 ] ], [ [ E(3)^2, 0 ], [ 0, E(3) ] ] ]
gap> IrreducibleAffordingRepresentation(chis[6]);
[ f1, f2, f3, f4 ] -> [ [ [ 0, 1 ], [ 1, 0 ] ], [ [ 0, 1 ], [ 1, 0 ] ],
```

```
          [ [ 1, 0 ], [ 0, 1 ] ] ], [ [ E(3)^2, 0 ], [ 0, E(3) ] ] ]
gap> IrreducibleAffordingRepresentation(chis[7]);
[ f1, f2, f3, f4 ] -> [ [ [ 0, -1, 0 ], [ 1, 0, 0 ], [ 0, 0, -1 ] ],
  [ [ 0, -1, 0 ], [ -1, 0, 0 ], [ 0, 0, -1 ] ],
  [ [ 1, 0, 0 ], [ 0, 1, 0 ], [ 0, 0, -1 ] ],
  [ [ 0, 1, 0 ], [ 0, 0, 1 ], [ 1, 0, 0 ] ] ]
gap> IrreducibleAffordingRepresentation(chis[8]);
[ f1, f2, f3, f4 ] -> [ [ [ 0, 1, 0 ], [ -1, 0, 0 ], [ 0, 0, 1 ] ],
  [ [ 0, 1, 0 ], [ 1, 0, 0 ], [ 0, 0, 1 ] ],
  [ [ 1, 0, 0 ], [ 0, 1, 0 ], [ 0, 0, -1 ] ],
  [ [ 0, 1, 0 ], [ 0, 0, 1 ], [ 1, 0, 0 ] ] ]
gap> IrreducibleAffordingRepresentation(chis[9]);
[ f1, f2, f3, f4 ] -> [ [ [ 0, 1, 0 ], [ -1, 0, 0 ], [ 0, 0, 1 ] ],
  [ [ 0, -1, 0 ], [ -1, 0, 0 ], [ 0, 0, -1 ] ],
  [ [ -1, 0, 0 ], [ 0, -1, 0 ], [ 0, 0, 1 ] ],
  [ [ 0, 1, 0 ], [ 0, 0, 1 ], [ 1, 0, 0 ] ] ]
gap> IrreducibleAffordingRepresentation(chis[10]);
[ f1, f2, f3, f4 ] -> [ [ [ 0, -1, 0 ], [ 1, 0, 0 ], [ 0, 0, -1 ] ],
  [ [ 0, 1, 0 ], [ 1, 0, 0 ], [ 0, 0, 1 ] ],
  [ [ -1, 0, 0 ], [ 0, -1, 0 ], [ 0, 0, 1 ] ],
  [ [ 0, 1, 0 ], [ 0, 0, 1 ], [ 1, 0, 0 ] ] ]
```

# Appendix B

# GAP Code for Klein Quartic

The GAP code in this appendix has the same format as that in Appendix A; please refer to the commentary there in order to interpret this code.

```
gap> f:=FreeGroup("a","b","c");
<free group on the generators [ a, b, c ]>
gap> g:=f/[f.1^2,f.2^2,f.3^2,(f.1*f.2)^2,(f.2*f.3)^3,(f.1*f.3)^7,(f.1*f.2*f.3)^8];
<fp group on the generators [ a, b, c ]>
gap> Order(g);
336
gap> StructureDescription(g);
"PSL(3,2):C2"
gap> ConjugacyClasses(g);
[ <identity ...>^G, b^G, b*c^G, a*b*c*a*(c*a*b)^2*(c*a)^2*c^G,
  a*c*a*(c*a*b)^2*(c*a)^2*c^G, c*a^G, b*c*a^G, (a*b*(c*a)^2*c)^2^G,
  b*(c*a*b*c*a)^2*b*c^G ]
gap> RequirePackage("repsn");
———————————————————————————————————————
Repsn for Constructing Representations of Finite Groups
                    Version 3.0.2

                    Written by
                 Vahid Dabbaghian
———————————————————————————————————————
true
gap> Display(CharacterTable(g));
CT1

     2  4  1  1  2  .  4  3  3  3
     3  1  1  1  1  .  .  .  .  .
     7  1  .  .  .  1  .  .  .  .

        1a 3a 6a 2a 7a 2b 8a 4a 8b
     2P 1a 3a 3a 1a 7a 1a 4a 2b 4a
     3P 1a 1a 2a 2a 7a 2b 8b 4a 8a
     5P 1a 3a 6a 2a 7a 2b 8b 4a 8a
     7P 1a 3a 6a 2a 1a 2b 8a 4a 8b

X.1       1  1  1  1  1  1  1  1  1
X.2       1  1 -1 -1  1  1 -1  1 -1
X.3       6  .  .  . -1 -2  .  2  .
X.4       6  .  .  . -1  2  A  . -A
X.5       6  .  .  . -1  2 -A  .  A
X.6       7  1 -1 -1  . -1  1 -1  1
X.7       7  1  1  1  . -1 -1 -1 -1
X.8       8 -1 -1  2  1  .  .  .  .
```

```
X.9       8  -1   1  -2   1   .   .   .   .

A = -E(8)+E(8)^3
  = -Sqrt(2) = -r2
gap> IrreducibleRepresentations(g);
[ [ a, b, c ] -> [ [ [ 1 ] ], [ [ 1 ] ], [ [ 1 ] ] ],
    [ a, b, c ] -> [ [ [ -1 ] ], [ [ -1 ] ], [ [ -1 ] ] ],
    [ a, b, c ] -> [ [ [ 1/2, 3/2, 0, -1/2, 0, 1 ], [ 3/2, 1/2, -1, 1/2, 0, 1 ],
            [ 0, 0, 0, -1, 0, 0 ], [ 0, 0, -1, 0, 0, 0 ], [ 0, -1, 0, 0, 1, -1 ],
            [ -3/2, -3/2, 1, -1/2, 0, -2 ] ],
        [ [ -1/2, -3/2, 0, 1/2, 0, -1 ], [ -1/2, 1/2, 0, 1/2, 0, 0 ],
            [ -3/2, -3/2, 1, -1/2, 0, -2 ], [ 0, 0, 0, 0, 0, -1 ],
            [ -1, 0, 1, -1, -1, 0 ], [ 0, 0, 0, -1, 0, 0 ] ],
        [ [ -1/2, -1/2, 1, -1/2, -1, 0 ], [ -1, 0, 1, 0, 0, -1 ], [ 0, 0, 0, 0, -1, 0 ],
            [ -1/2, 1/2, 0, 1/2, -1, 1 ],
            [ 0, 0, -1, 0, 0, 0 ], [ 1/2, -1/2, -1, 1/2, 0, 0 ] ] ],
    [ a, b, c ] -> [ [ [ -1, -E(8)+E(8)^3, -1, 0, 0, 0 ], [ 0, 1, 0, 0, 0, 0 ],
            [ 0, 0, 1, 0, 0, 0 ],
            [ -E(8)+E(8)^3, E(8)-E(8)^3, 1-E(8)+E(8)^3, -1-E(8)+E(8)^3, 2, 0 ],
            [ -1-E(8)+E(8)^3, 0, -1, -1-E(8)+E(8)^3, 1+E(8)-E(8)^3, 0 ],
            [ -2, 1-E(8)+E(8)^3, -3+E(8)-E(8)^3, -2, 2*E(8)-2*E(8)^3, -1 ] ],
        [ [ 0, 0, 0, 1, 0, 0 ], [ 1, -1, 1-E(8)+E(8)^3, 1, -E(8)+E(8)^3, 0 ],
            [ 0, 0, 1, 0, 0, 0 ], [ 1, 0, 0, 0, 0, 0 ],
            [ 0, 0, 0, 0, 1, 0 ], [ -1, 0, -1, -1, 1+E(8)-E(8)^3, -1 ] ],
        [ [ 1, 0, 0, 0, 0, 0 ], [ 1, -1, 1, 1+E(8)-E(8)^3, -1-E(8)+E(8)^3, 0 ],
            [ -1, 1, -1, -2, 1+E(8)-E(8)^3, -1 ],
            [ 1+E(8)-E(8)^3, 0, E(8)-E(8)^3, 1+E(8)-E(8)^3, -2-E(8)+E(8)^3, 0 ],
            [ 2, -1+E(8)-E(8)^3, 2-E(8)+E(8)^3, 2-E(8)+E(8)^3, -E(8)+E(8)^3, 1-E(8)+E(8)^3 ],
            [ 1, -1, 1-E(8)+E(8)^3, 1, -E(8)+E(8)^3, 0 ] ] ],
    [ a, b, c ] -> [ [ [ 1, 0, 0, 0, 0, 0 ], [ 0, 0, 1, 0, 0, 0 ], [ 0, 1, 0, 0, 0, 0 ],
            [ -3-2*E(8)+2*E(8)^3, 5+4*E(8)-4*E(8)^3, 4+3*E(8)-3*E(8)^3, 1+E(8)-E(8)^3,
                -2-E(8)+E(8)^3, 3+2*E(8)-2*E(8)^3 ],
            [ -2-E(8)+E(8)^3, 3+2*E(8)-2*E(8)^3, 2+2*E(8)-2*E(8)^3, E(8)-E(8)^3, -1-E(8)+E(8)^3,
                1+E(8)-E(8)^3 ],
            [ E(8)-E(8)^3, -1-E(8)+E(8)^3, -1-E(8)+E(8)^3, 0, 0, -1 ] ],
        [ [ 1, 0, 0, 0, 0, 0 ], [ -1-2*E(8)+2*E(8)^3, 3+3*E(8)-3*E(8)^3, 3+2*E(8)-2*E(8)^3,
                1, -2-E(8)+E(8)^3, 2+E(8)-E(8)^3 ],
            [ 4+2*E(8)-2*E(8)^3, -6-4*E(8)+4*E(8)^3, -5-3*E(8)+3*E(8)^3, -1-E(8)+E(8)^3,
                2+2*E(8)-2*E(8)^3, -3-2*E(8)+2*E(8)^3 ],
            [ 0, 0, 0, 0, 0, 1 ], [ 1+E(8)-E(8)^3, -2-2*E(8)+2*E(8)^3, -2-E(8)+E(8)^3, -1, 1,
                -1-E(8)+E(8)^3 ], [ 0, 0, 0, 1, 0, 0 ] ],
        [ [ 3+E(8)-E(8)^3, -4-2*E(8)+2*E(8)^3, -3-2*E(8)+2*E(8)^3, -E(8)+E(8)^3,
                1+E(8)-E(8)^3, -2-2*E(8)+2*E(8)^3 ],
            [ 3+2*E(8)-2*E(8)^3, -4-3*E(8)+3*E(8)^3, -3-2*E(8)+2*E(8)^3, -1-E(8)+E(8)^3,
                2+E(8)-E(8)^3, -2-2*E(8)+2*E(8)^3 ],
            [ -E(8)+E(8)^3, E(8)-E(8)^3, 0, 1, -1, 1 ], [ 3+2*E(8)-2*E(8)^3, -5-4*E(8)+4*E(8)^3,
                -5-3*E(8)+3*E(8)^3, -1, 2+E(8)-E(8)^3, -3-2*E(8)+2*E(8)^3 ],
            [ 4+2*E(8)-2*E(8)^3, -6-4*E(8)+4*E(8)^3, -5-3*E(8)+3*E(8)^3, -1-E(8)+E(8)^3,
                2+2*E(8)-2*E(8)^3, -3-2*E(8)+2*E(8)^3 ], [ -1, 1, 1, 0, 0, 0 ] ] ],
    [ a, b, c ] -> [ [ [ -3/2, -5/4, 5/4, 11/4, 1/2, -7/4, 1 ],
            [ -1/2, -3/4, 3/4, 1/4, -1/2, 3/4, 0 ], [ 0, 0, 0, 0, 0, 1, 0 ],
            [ -3/2, -5/4, 19/12, 37/12, -1/6, -29/12, 4/3 ],
            [ 3/2, -1/4, -1/12, -19/12, 1/6, 11/12, -1/3 ], [ 0, 0, 1, 0, 0, 0, 0 ],
            [ 3/2, 3/4, 1/4, -13/4, 1/2, 13/4, -2 ] ],
        [ [ 2, 1, -2/3, -14/3, 1/3, 13/3, -8/3 ], [ 0, -1, 0, 0, 0, 0, 0 ],
            [ 1/2, -3/4, 5/12, -13/12, 1/6, 5/12, -1/3 ], [ 3/2, 5/4, -19/12, -37/12, 1/6,
                29/12, -4/3 ], [ 0, 0, 0, 0, 0, 0, 1 ],
            [ 1, 1, -4/3, -4/3, 2/3, 2/3, -1/3 ], [ 0, 0, 0, 0, 1, 0, 0 ] ],
        [ [ 0, 1, 0, 0, 0, 0, 0 ], [ 1, 0, 0, 0, 0, 0, 0 ], [ 0, -1/2, 1/2, 1/2, 0, 1/2, -1 ],
            [ 3/2, 5/4, -1/4, -7/4, 1/2, 3/4, 0 ], [ 1/2, -3/4, 3/4, -7/4, 1/2, 11/4, -1 ],
            [ 3/2, 5/4, -11/12, -17/12, 5/6, 1/12, 1/3 ],
            [ 1, 1, -4/3, -4/3, 2/3, 2/3, -1/3 ] ] ],
    [ a, b, c ] -> [ [ [ 5/3, -4, -11/3, 2/3, -7/3, 3, 1 ], [ -1, 0, 1, 1, 1, -1, 0 ],
            [ 2, -3, -5, -1, -4, 5, 1 ], [ 0, 0, 0, 1, 0, 0, 0 ],
```

```
          [ 13/3, -8, -25/3, 4/3, -20/3, 8, 3 ], [ 14/3, -8, -29/3, 2/3, -25/3, 10, 3 ],
          [ -7/3, 1, 4/3, -1/3, 8/3, -2, 0 ] ],
      [ [ 0, 1, 0, 0, 0, 0, 0 ], [ 1, 0, 0, 0, 0, 0, 0 ],
          [ 5/3, -3, -14/3, -1/3, -13/3, 5, 1 ], [ 5/3, -3, -8/3, 2/3, -4/3, 2, 1 ],
          [ 1, -1, -3, 0, -2, 3, 0 ], [ 8/3, -4, -20/3, -1/3, -16/3, 7, 1 ],
          [ 7/3, -1, -4/3, 1/3, -8/3, 2, 0 ] ],
      [ [ 2, -3, -5, -1, -4, 5, 1 ], [ -11/3, 4, 11/3, -5/3, 16/3, -5, -1 ],
          [ 4/3, 0, -1/3, 1/3, -5/3, 1, -1 ], [ -7/3, 1, 4/3, -1/3, 8/3, -2, 0 ],
          [ 8/3, -4, -20/3, -1/3, -16/3, 7, 1 ], [ 2/3, 0, -5/3, -1/3, -4/3, 2, -1 ],
          [ -7/3, 3, 7/3, -4/3, 11/3, -4, -1 ] ] ],
[ a, b, c ] -> [ [ [ 1/2, 1/2, 1/2, 1, 1, -1/2, 3/2, 1 ],
          [ -3/2, -3/2, -3/2, 0, 0, -1/2, -1/2, -1 ], [ 0, 0, 1, 0, 0, 0, 0, 0 ],
          [ 0, 0, 0, 0, 0, 1, 0 ], [ 0, 0, 0, 0, 1, 0, 0, 0 ], [ 0, 0, 0, 0, 0, 1, 0, 0 ],
          [ 0, 0, 0, 1, 0, 0, 0, 0 ], [ 3/2, 1/2, 0, -2, -3/2, 1, -3/2, 0 ] ],
      [ [ 3/2, 1/2, 0, -2, -3/2, 1, -3/2, 0 ], [ -3/2, -1/2, 0, 2, 1/2, -1, 1/2, 0 ],
          [ 0, 0, 0, 0, 1, 0, 0, 0 ], [ 1, 1, 0, -1, -1, 1, -1, 0 ],
          [ 0, 0, 1, 0, 0, 0, 0, 0 ], [ 0, 0, 0, 0, 0, 1, 0, 0 ],
          [ -1, -1, -1, 0, 0, 0, 0, 0 ], [ 1/2, 1/2, 1/2, 1, 1, -1/2, 3/2, 1 ] ],
      [ [ 0, 1, 1, 2, 2, -1, 2, 2 ], [ 2, 0, 0, -3, -2, 1, -2, -1 ],
          [ -2, -1, -1, 1, 0, -1, 0, -1 ], [ -1/2, 1/2, 1, 2, 3/2, -1, 3/2, 2 ],
          [ 0, -1, -1, -1, 0, 0, -1, -1 ], [ -1, -1, -1, 0, 0, 0, 0, 0 ],
          [ 0, 0, 0, 0, 0, 0, 1, 0 ], [ 1/2, 1/2, 0, 0, -1/2, 1, -1/2, 0 ] ] ],
[ a, b, c ] -> [ [ [ -1, 1, -1, 0, 1, 0, 0, 0 ], [ 0, 1, 0, 0, 0, 0, 0, 0 ],
          [ 0, 0, 0, 0, -1, 0, 0, 0 ], [ 1, 0, 1, 0, -1, 1, 1, -1 ],
          [ 0, 0, -1, 0, 0, 0, 0, 0 ], [ 1, 0, 0, 1, 0, 0, 1, -1 ],
          [ 0, 0, 0, 0, 0, 0, -1, 0 ], [ 0, 1, 0, 0, 0, 0, 0, -1 ] ],
      [ [ -1, 0, 1, -1, 0, 1, 0, 0 ], [ 0, -1, 0, 0, 0, 0, 0, 0 ],
          [ 0, 0, 0, 0, 1, 0, 0, 0 ], [ 1, 0, 0, 1, 0, 0, 1, -1 ], [ 0, 0, 1, 0, 0, 0, 0, 0 ],
          [ 1, 0, 1, 0, -1, 1, 1, -1 ], [ 0, 0, -1, 1, 0, -1, -1, 0 ],
          [ 0, 0, 0, 0, 0, 0, 0, -1 ] ],
      [ [ 0, 0, 0, 0, 0, -1, 0, 0 ], [ 0, 0, -1, 0, 0, -1, -1, 0 ],
          [ 1, -1, 0, 0, -1, 0, 0, 0 ], [ -1, 0, 0, -1, 0, 1, 0, 0 ],
          [ 0, 0, 0, 0, 0, 0, 1, 0 ], [ -1, 0, 0, 0, 0, 0, 0, 0 ], [ 0, 0, 0, 0, 1, 0, 0, 0 ],
          [ 0, 0, 0, 0, 0, 0, 0, -1 ] ] ] ]
```

# Appendix C

# FreeFEM++ Code for Genus 2

The following box contains the FreeFEM++ code for the Bolza surface. It contains basic comments, but for more information about the various solvers and finite element spaces available in the software, please consult the manual [32].

```
int n = 30; //number of modes
int nev=100; //number of eigenvalues to be calculated

string fnm = "bolza_surface_30n_100nev.txt"; //file name for saving

real r = 0.4550898605622276; //radius of circle
real c = 1.09868411346781; //centre of circle on x-axis

border G1 (t=5*pi/4,3*pi/4) { x=c+r*cos(t); y=r*sin(t);};
border G2 (t=5*pi/4,3*pi/4) { x=cos(pi/4)*(c+r*cos(t))-sin(pi/4)*r*sin(t);
    y=sin(pi/4)*(c+r*cos(t))+cos(pi/4)*r*sin(t);};
border G3 (t=5*pi/4,3*pi/4) { x=cos(pi/2)*(c+r*cos(t))-sin(pi/2)*r*sin(t);
    y=sin(pi/2)*(c+r*cos(t))+cos(pi/2)*r*sin(t);};
border G4 (t=5*pi/4,3*pi/4) { x=cos(3*pi/4)*(c+r*cos(t))-sin(3*pi/4)*r*sin(t);
    y=sin(3*pi/4)*(c+r*cos(t))+cos(3*pi/4)*r*sin(t);};
border G5 (t=5*pi/4,3*pi/4) { x=cos(pi)*(c+r*cos(t))-sin(pi)*r*sin(t);
    y=sin(pi)*(c+r*cos(t))+cos(pi)*r*sin(t);};
border G6 (t=5*pi/4,3*pi/4) { x=cos(5*pi/4)*(c+r*cos(t))-sin(5*pi/4)*r*sin(t);
    y=sin(5*pi/4)*(c+r*cos(t))+cos(5*pi/4)*r*sin(t);};
border G7 (t=5*pi/4,3*pi/4) { x=cos(6*pi/4)*(c+r*cos(t))-sin(6*pi/4)*r*sin(t);
    y=sin(6*pi/4)*(c+r*cos(t))+cos(6*pi/4)*r*sin(t);};
border G8 (t=5*pi/4,3*pi/4) { x=cos(7*pi/4)*(c+r*cos(t))-sin(7*pi/4)*r*sin(t);
    y=sin(7*pi/4)*(c+r*cos(t))+cos(7*pi/4)*r*sin(t);};
plot(G1(n)+G2(n)+G3(n)+G4(n)+G5(n)+G6(n)+G7(n)+G8(n));

mesh Th=buildmesh(G1(n)+G2(n)+G3(n)+G4(n)+G5(n)+G6(n)+G7(n)+G8(n),fixeborder=true);
plot(Th,wait=true,fill=true);

fespace Vh(Th,P2,periodic=[[1,y],[5,y],[3,x],[7,x],[2,y-x],[6,y-x],[4,y+x],[8,y+x]]);
    //glues opposite sides such that orientation is preserved

Vh u1, u2;
real sigma = 0.00001; //value of the shift
varf op(u1,u2)= int2d(Th)(dx(u1)*dx(u2)+ dy(u1)*dy(u2)- sigma*(4*u1*u2)/(1-x^2-y^2)^2);
```

```
varf b([u1],[u2]) = int2d(Th)((4*u1*u2)/(1-x^2-y^2)^2);

matrix OP= op(Vh,Vh,solver=Crout,factorize=1);
matrix B=b(Vh,Vh,solver=CG,eps=1e-20);

real[int] ev(nev);
Vh[int] eV(nev);

int k=EigenValue(OP,B,sym=true,sigma=sigma,value=ev,vector=eV, tol=1e-10,maxit=0,ncv=0);
for (int i=0;i<k;i++) {
u1=eV[i];
real gg = int2d(Th)(dx(u1)*dx(u1) + dy(u1)*dy(u1)); real mm= int2d(Th)(u1*u1) ;

ofstream Eva(fnm, append);
Eva <<   ev[i] << "\n";

//cout<<"----"<< i<<""<<ev[i]<<"err="
//<<dx(u1)*dx(u1) + dy(u1)*dy(u1) - (ev[i])*u1*u1 << " --- "<<endl; plot(eV[i],
    cmm="Eigen_Vector_"+i+"_valeur_=" + ev[i]  ,wait=1,value=1);
//un-comment the above two lines to display level sets for each eigenfunction
}
```

Table C.1 shows the positive eigenvalues of the Bolza surface, up to approximately 100. They can be compared with the much higher accuracy eigenvalues calculated by Strohmaier and Uski using the algorithm described in [77]. The data files associated with the paper can be obtained from the personal website of Professor Strohmaier: http://www1.maths.leeds.ac.uk/~pmtast/publications/eigdata/datafile.html. Here can be found the first 1000 eigenvalues of the Bolza surface, and other genus 2 surfaces, including the other two in this appendix.

Related to the Bolza surface is the (2, 3, 8) triangle that tessellates it. We have produced a program to calculate that eigenvalues of this - please see the Bolza surface code for comments. The one comment included in this code relates to boundary conditions - it is straightforward to change between Dirichlet and Neumann boundary conditions.

```
int n = 50;
int nev = 100;

string fnm = "238tri_50n_100nev_nnn.txt";

border G1 (t=0,0.4056163087774724) { x=t; y=0;};
border G2 (t=0.33681578765748216,0) { x=t; y=t*tan(pi/8);};
border G3 (t=0.4056163087774724,0.33681578765748216) { x=t;
    y=0.5946036542842016-sqrt((1.189207445439991^2)-(t-1.435500206409523)^2);};

plot(G1(n)+G2(n)+G3(n));
mesh Th=buildmesh(G1(n)+G2(n)+G3(n));
plot(Th);

fespace Vh(Th,P2);
```

```
Vh u1, u2;
real sigma = 0.001;

varf op(u1,u2)= int2d(Th)(dx(u1)*dx(u2)+ dy(u1)*dy(u2)- sigma*(4*u1*u2)/(1-x^2-y^2)^2);
    //naturally, this gives Neumann boundary conditions on all three sides.
    //To change one or more sides to a Dirichlet condition, add the phrase
    //'+on(G1,u1=0)' after the closed parenthesis.
    //This will give a Dirichlet condition on G1.
varf b([u1],[u2]) = int2d(Th)((4*u1*u2)/(1-x^2-y^2)^2);

matrix OP= op(Vh,Vh,solver=Crout,factorize=1);
matrix B=b(Vh,Vh,solver=CG,eps=1e-20);

real[int] ev(nev);
Vh[int] eV(nev);

int k=EigenValue(OP,B,sym=true,sigma=sigma,value=ev,vector=eV, tol=1e-10,maxit=0,ncv=0);
for (int i=0;i<k;i++) {
u1=eV[i];
real gg = int2d(Th)(dx(u1)*dx(u1) + dy(u1)*dy(u1));  real mm= int2d(Th)(u1*u1) ;

ofstream Eva(fnm, append);
Eva <<   ev[i] << "\n";

}
```

Also of interest is the spectrum of the right angled pentagon of area $\frac{\pi}{2}$ with mixed boundary conditions, studied in Section 3.4. All of the multiplicity 4 eigenvalues of $\mathcal{B}$ appear in this spectrum with simple multiplicity. Recall the discussion about isospectrality in the introduction to this work - one can see for oneself that swapping the Dirichlet and Neumann conditions in the code below (see the comments in the code for the (2, 3, 8) triangle) produces the same spectrum! Table C.2 shows the eigenvalues of this problem up to 200; one can compare this to Table C.1, where all of the eigenvalues in Table C.2 appear with multiplicity 4.

```
int n=30;
int nev=100;

real r=0.643594;

string fnm = "pentagon_30n_100nev.txt"; //file name for saving

border G1 (t=0.5,1) { x=r*2*cos(pi/8)+r*cos((9-2*t)*pi/8); y=r*sin((9-2*t)*pi/8);};
border G2 (t=-1,0) { x=cos(pi/4)*r*2*cos(pi/8)+r*cos((9-2*t)*pi/8);
    y=sin(pi/4)*r*2*cos(pi/8)+r*sin((9-2*t)*pi/8);};
border G3 (t=1,0.5) { x=r*sin((9-2*t)*pi/8); y=r*2*cos(pi/8)+r*cos((9-2*t)*pi/8);};
border G4 (t=0.545613,0) { x=0; y=t;};
border G5 (t=0,0.545613) { x=t; y=0;};
plot(G1(n)+G2(n)+G3(n)+G4(n)+G5(n));

mesh Th=buildmesh(G1(n)+G2(n)+G3(n)+G4(n)+G5(n));
```

```
plot(Th);
fespace Vh(Th,P2);

Vh u1, u2;
real sigma = 0.00001;
varf op(u1,u2)= int2d(Th)(dx(u1)*dx(u2)+ dy(u1)*dy(u2)- sigma* u1*u2)
    +on(G2,u1=0)+on(G4,u1=0);
varf b([u1],[u2]) = int2d(Th)((4*u1*u2)/(1-x^2-y^2)^2);

matrix OP= op(Vh,Vh,solver=Crout,factorize=1);
matrix B=b(Vh,Vh,solver=CG,eps=1e-20);

real[int] ev(nev);
Vh[int] eV(nev);
int k=EigenValue(OP,B,sym=true,sigma=sigma,value=ev,vector=eV, tol=1e-10,maxit=0,ncv=0);
for (int i=0;i<k;i++) {
u1=eV[i];
real gg = int2d(Th)(dx(u1)*dx(u1) + dy(u1)*dy(u1));  real mm= int2d(Th)(u1*u1) ;

ofstream Eva(fnm, append);
Eva <<  ev[i] << "\n";

}
```

We have produced FreeFEM++ codes for other genus 2 surfaces during this project, and would like to include them in this appendix. First, we give the code that computes the spectrum of the surface with automorphism group $\mathbb{Z}_{10}$:

```
int n = 50;
int nev = 50;

string fnm="Z10-Eigenvalues-n50-nev50.txt"; //file name for saving

real r = 0.4133042381223985;
real c = 1.082044543098821;
real s = 0.6687403049764225;

real x01 = s;
real y01 = 0;
real x04 = s*cos(3*pi/5);
real y04 = s*sin(3*pi/5);

real x03 = s*cos(2*pi/5);
real y03 = s*sin(2*pi/5);
real x06 = s*cos(pi);
real y06 = s*sin(pi);

real x05 = s*cos(4*pi/5);
real y05 = s*sin(4*pi/5);
real x08 = s*cos(7*pi/5);
real y08 = s*sin(7*pi/5);

real x07 = s*cos(6*pi/5);
real y07 = s*sin(6*pi/5);
```

```
real x10 = s*cos(9*pi/5);
real y10 = s*sin(9*pi/5);

real x09 = s*cos(8*pi/5);
real y09 = s*sin(8*pi/5);
real x02 = s*cos(pi/5);
real y02 = s*sin(pi/5);

border G1 (t=6*pi/5,4*pi/5) { x=c+r*cos(t); y=r*sin(t);};
border G2 (t=6*pi/5,4*pi/5) { x=cos(pi/5)*(c+r*cos(t))-sin(pi/5)*r*sin(t);
    y=sin(pi/5)*(c+r*cos(t))+cos(pi/5)*r*sin(t);};
border G3 (t=6*pi/5,4*pi/5) { x=cos(2*pi/5)*(c+r*cos(t))-sin(2*pi/5)*r*sin(t);
    y=sin(2*pi/5)*(c+r*cos(t))+cos(2*pi/5)*r*sin(t);};
border G4 (t=6*pi/5,4*pi/5) { x=cos(3*pi/5)*(c+r*cos(t))-sin(3*pi/5)*r*sin(t);
    y=sin(3*pi/5)*(c+r*cos(t))+cos(3*pi/5)*r*sin(t);};
border G5 (t=6*pi/5,4*pi/5) { x=cos(4*pi/5)*(c+r*cos(t))-sin(4*pi/5)*r*sin(t);
    y=sin(4*pi/5)*(c+r*cos(t))+cos(4*pi/5)*r*sin(t);};
border G6 (t=6*pi/5,4*pi/5) { x=cos(pi)*(c+r*cos(t))-sin(pi)*r*sin(t);
    y=sin(pi)*(c+r*cos(t))+cos(pi)*r*sin(t);};
border G7 (t=6*pi/5,4*pi/5) { x=cos(6*pi/5)*(c+r*cos(t))-sin(6*pi/5)*r*sin(t);
    y=sin(6*pi/5)*(c+r*cos(t))+cos(6*pi/5)*r*sin(t);};
border G8 (t=6*pi/5,4*pi/5) { x=cos(7*pi/5)*(c+r*cos(t))-sin(7*pi/5)*r*sin(t);
    y=sin(7*pi/5)*(c+r*cos(t))+cos(7*pi/5)*r*sin(t);};
border G9 (t=6*pi/5,4*pi/5) { x=cos(8*pi/5)*(c+r*cos(t))-sin(8*pi/5)*r*sin(t);
    y=sin(8*pi/5)*(c+r*cos(t))+cos(8*pi/5)*r*sin(t);};
border G10 (t=6*pi/5,4*pi/5) { x=cos(9*pi/5)*(c+r*cos(t))-sin(9*pi/5)*r*sin(t);
    y=sin(9*pi/5)*(c+r*cos(t))+cos(9*pi/5)*r*sin(t);};
plot(G1(n)+G2(n)+G3(n)+G4(n)+G5(n)+G6(n)+G7(n)+G8(n)+G9(n)+G10(n));
mesh Th=buildmesh(G1(n)+G2(n)+G3(n)+G4(n)+G5(n)+G6(n)+G7(n)+G8(n)+G9(n)+G10(n),
    fixeborder=true);

plot(Th,wait=true,fill=true);

fespace Vh(Th,P2,periodic=[[1,y],[4,cos(2*pi/5)*(y-y04)-sin(-2*pi/5)*(x-x04)],
    [3,cos(2*pi/5)*(y-y03)-sin(2*pi/5)*(x-x03)],[6,y],
    [5,cos(2*pi/5)*(y-y05)-sin(-2*pi/5)*(x-x05)],
    [8,cos(2*pi/5)*(y-y08)-sin(-2*pi/5)*(x-x08)],
    [7,cos(2*pi/5)*(y-y07)-sin(2*pi/5)*(x-x07)],
    [10,cos(2*pi/5)*(y-y10)-sin(-2*pi/5)*(x-x10)],
    [9,cos(2*pi/5)*(y-y09)-sin(2*pi/5)*(x-x09)],
    [2,cos(2*pi/5)*(y-y02)-sin(2*pi/5)*(x-x02)]]);

Vh u1, u2;
real sigma = 0.00001;
varf op(u1,u2)= int2d(Th)(dx(u1)*dx(u2)+ dy(u1)*dy(u2)- sigma*(4*u1*u2)/(1-x^2-y^2)^2);
varf b([u1],[u2]) = int2d(Th)((4*u1*u2)/(1-x^2-y^2)^2);
matrix OP= op(Vh,Vh,solver=Crout,factorize=1);
matrix B=b(Vh,Vh,solver=CG,eps=1e-20);

real[int] ev(nev);
Vh[int] eV(nev);

int k=EigenValue(OP,B,sym=true,sigma=sigma,value=ev,vector=eV, tol=1e-10,maxit=0,ncv=0);
for (int i=0;i<k;i++) {
```

```
u1=eV[i];
real gg = int2d(Th)(dx(u1)*dx(u1) + dy(u1)*dy(u1));   real mm= int2d(Th)(u1*u1) ;

ofstream Eva(fnm, append);
Eva <<   ev[i] << "\n";

}
```

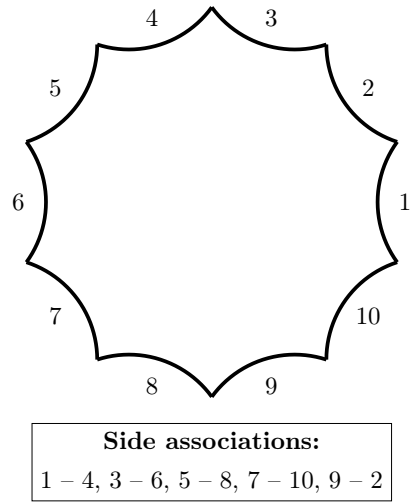

Figure C.1: Fundamental polygon for the surface with automorphism group $\mathbb{Z}_{10}$

Next, we have the surface that has a semi-regular 12-gon as its fundamental domain, and has automorphism group $(4, 6\,|\,2, 2)$ of order 24 (see [20] for a classification of this group):

```
int n = 50;
int nev = 50;

string fnm="12gon-Eigenvalues-n50-nev50.txt";  //file name for saving

real cx = 1.0206207261596574;
real cy = 0.3535533905932733;
real r  = 0.4082482904638626;

real s=0.6835992463050957;
real theta=0.22037122071230927;
real x01=s*cos(theta);
real y01=s*sin(theta);
real x02=s*cos(pi/3-theta);
real y02=s*sin(pi/3-theta);
real x03=s*cos(pi/3+theta);
real y03=s*sin(pi/3+theta);
real x04=s*cos(2*pi/3-theta);
real y04=s*sin(2*pi/3-theta);
real x05=s*cos(2*pi/3+theta);
real y05=s*sin(2*pi/3+theta);
real x06=s*cos(3*pi/3-theta);
```

```
real y06=s*sin(3*pi/3-theta);
real x07=s*cos(3*pi/3+theta);
real y07=s*sin(3*pi/3+theta);
real x08=s*cos(4*pi/3-theta);
real y08=s*sin(4*pi/3-theta);
real x09=s*cos(4*pi/3+theta);
real y09=s*sin(4*pi/3+theta);
real x10=s*cos(5*pi/3-theta);
real y10=s*sin(5*pi/3-theta);
real x11=s*cos(5*pi/3+theta);
real y11=s*sin(5*pi/3+theta);
real x12=s*cos(-theta);
real y12=s*sin(-theta);

border G1 (t=4*pi/3,pi) { x=cx+r*cos(t); y=cy+r*sin(t);};
border G2 (t=pi,2*pi/3) { x=cos(pi/3)*(cx+r*cos(t))-sin(pi/3)*(-cy+r*sin(t));
    y=sin(pi/3)*(cx+r*cos(t))+cos(pi/3)*(-cy+r*sin(t));};
border G3 (t=4*pi/3,pi) { x=cos(pi/3)*(cx+r*cos(t))-sin(pi/3)*(cy+r*sin(t));
    y=sin(pi/3)*(cx+r*cos(t))+cos(pi/3)*(cy+r*sin(t));};
border G4 (t=pi,2*pi/3) { x=cos(2*pi/3)*(cx+r*cos(t))-sin(2*pi/3)*(-cy+r*sin(t));
    y=sin(2*pi/3)*(cx+r*cos(t))+cos(2*pi/3)*(-cy+r*sin(t));};
border G5 (t=4*pi/3,pi) { x=cos(2*pi/3)*(cx+r*cos(t))-sin(2*pi/3)*(cy+r*sin(t));
    y=sin(2*pi/3)*(cx+r*cos(t))+cos(2*pi/3)*(cy+r*sin(t));};
border G6 (t=pi,2*pi/3) { x=cos(3*pi/3)*(cx+r*cos(t))-sin(3*pi/3)*(-cy+r*sin(t));
    y=sin(3*pi/3)*(cx+r*cos(t))+cos(3*pi/3)*(-cy+r*sin(t));};
border G7 (t=4*pi/3,pi) { x=cos(3*pi/3)*(cx+r*cos(t))-sin(3*pi/3)*(cy+r*sin(t));
    y=sin(3*pi/3)*(cx+r*cos(t))+cos(3*pi/3)*(cy+r*sin(t));};
border G8 (t=pi,2*pi/3) { x=cos(4*pi/3)*(cx+r*cos(t))-sin(4*pi/3)*(-cy+r*sin(t));
    y=sin(4*pi/3)*(cx+r*cos(t))+cos(4*pi/3)*(-cy+r*sin(t));};
border G9 (t=4*pi/3,pi) { x=cos(4*pi/3)*(cx+r*cos(t))-sin(4*pi/3)*(cy+r*sin(t));
    y=sin(4*pi/3)*(cx+r*cos(t))+cos(4*pi/3)*(cy+r*sin(t));};
border G10 (t=pi,2*pi/3) { x=cos(5*pi/3)*(cx+r*cos(t))-sin(5*pi/3)*(-cy+r*sin(t));
    y=sin(5*pi/3)*(cx+r*cos(t))+cos(5*pi/3)*(-cy+r*sin(t));};
border G11 (t=4*pi/3,pi) { x=cos(5*pi/3)*(cx+r*cos(t))-sin(5*pi/3)*(cy+r*sin(t));
    y=sin(5*pi/3)*(cx+r*cos(t))+cos(5*pi/3)*(cy+r*sin(t));};
border G12 (t=pi,2*pi/3) { x=cx+r*cos(t); y=-cy+r*sin(t);};

plot(G1(n)+G2(n)+G3(n)+G4(n)+G5(n)+G6(n)+G7(n)+G8(n)+G9(n)+G10(n)+G11(n)+G12(n));
mesh Th=buildmesh(G1(n)+G2(n)+G3(n)+G4(n)+G5(n)+G6(n)+G7(n)+G8(n)+G9(n)+G10(n)+G11(n)
    +G12(n),fixeborder=true);
plot(Th,wait=true,fill=true);

fespace Vh(Th,P2,periodic=[[2,(cos(pi/6))*(y-y02)-(sin(pi/6))*(x-x02)],
    [5,(cos(-pi/6))*(y-y05)-(sin(-pi/6))*(x-x05)],
    [4,(cos(pi/2))*(y-y04)-(sin(pi/2))*(x-x04)],
    [7,(cos(pi/6))*(y-y07)-(sin(pi/6))*(x-x07)],
    [6,(cos(-7*pi/6))*(y-y06)-(sin(-7*pi/6))*(x-x06)],
    [9,(cos(pi/2))*(y-y09)-(sin(pi/2))*(x-x09)],
    [8,(cos(pi/6))*(y-y08)-(sin(pi/6))*(x-x08)],
    [11,(cos(pi/6))*(y-y11)-(sin(-pi/6))*(x-x11)],
    [10,(cos(pi/2))*(y-y10)-(sin(pi/2))*(x-x10)],
    [1,(cos(pi/6))*(y-y01)-(sin(pi/6))*(x-x01)],
    [12,(cos(pi/6))*(y-y12)-(sin(-pi/6))*(x-x12)],
    [3,(cos(-pi/2))*(y-y03)-(sin(-pi/2))*(x-x03)]]);
```

```
Vh u1, u2;
real sigma = 0.00001;
varf op(u1,u2)= int2d(Th)(dx(u1)*dx(u2)+ dy(u1)*dy(u2)- sigma*(4*u1*u2)/(1-x^2-y^2)^2);
varf b([u1],[u2]) = int2d(Th)((4*u1*u2)/(1-x^2-y^2)^2);

matrix OP= op(Vh,Vh,solver=Crout,factorize=1);
matrix B=b(Vh,Vh,solver=CG,eps=1e-20);

real[int] ev(nev);
Vh[int] eV(nev);

int k=EigenValue(OP,B,sym=true,sigma=sigma,value=ev,vector=eV, tol=1e-10,maxit=0,ncv=0);
for (int i=0;i<k;i++) {
u1=eV[i];
real gg = int2d(Th)(dx(u1)*dx(u1) + dy(u1)*dy(u1));  real mm= int2d(Th)(u1*u1)  ;

ofstream Eva(fnm, append);
Eva <<  ev[i] << "\n";

}
```

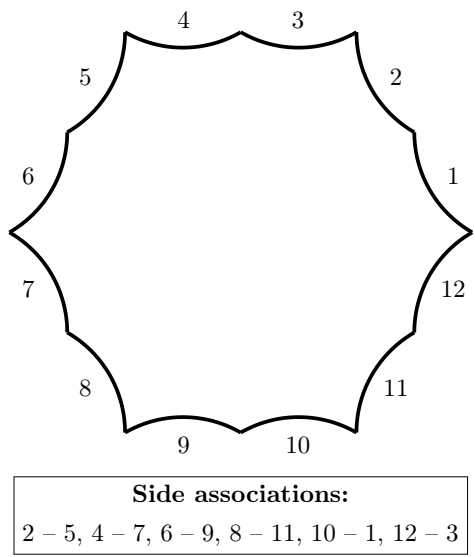

Figure C.2: Fundamental domain for the surface with order 24 automorphism group

| Eigenvalue | Numerical multiplicity |
|---:|:---:|
| 3.83618 | 3 |
| 5.35083 | 4 |
| 8.24401 | 2 |
| 14.7149 | 4 |
| 15.0386 | 3 |
| 18.645 | 3 |
| 20.5069 | 4 |
| 23.0545 | 1 |
| 28.0591 | 3 |
| 30.8081 | 4 |
| 32.6493 | 1 |
| 36.2126 | 2 |
| 38.9353 | 4 |
| 40.0913 | 3 |
| 42.8344 | 4 |
| 43.9864 | 3 |
| 50.5199 | 3 |
| 57.4478 | 4 |
| 59.4306 | 4 |
| 60.988 | 2 |
| 62.595 | 3 |
| 67.5682 | 3 |
| 71.3853 | 4 |
| 73.6255 | 2 |
| 74.8779 | 1 |
| 75.5018 | 3 |
| 75.6305 | 4 |
| 86.1472 | 4 |
| 86.6702 | 3 |
| 91.4225 | 1 |
| 93.3361 | 4 |
| 97.8419 | 3 |
| 100.675 | 3 |

Table C.1: First positive eigenvalues of the Bolza surface

| Eigenvalue | Numerical multiplicity |
|:---:|:---:|
| 5.35323 | 1 |
| 14.7246 | 1 |
| 20.5174 | 1 |
| 30.8305 | 1 |
| 38.956 | 1 |
| 42.852 | 1 |
| 57.4803 | 1 |
| 59.4601 | 1 |
| 71.4156 | 1 |
| 75.6649 | 1 |
| 86.177 | 1 |
| 93.3682 | 1 |
| 105.547 | 1 |
| 110.988 | 1 |
| 111.675 | 1 |
| 127.163 | 1 |
| 133.738 | 1 |
| 144.59 | 1 |
| 151.081 | 1 |
| 154.5 | 1 |
| 166.766 | 1 |
| 170.721 | 1 |
| 182.089 | 1 |
| 193.034 | 1 |
| 200.792 | 1 |

Table C.2: First positive eigenvalues of the mixed boundary problem on a pentagon of area $\frac{\pi}{2}$

| Eigenvalue | Numerical multiplicity |
|---|---|
| 3.48586 | 2 |
| 4.56299 | 2 |
| 5.05739 | 2 |
| 8.6528 | 2 |
| 8.96267 | 2 |
| 12.1336 | 1 |
| 14.9476 | 2 |
| 16.2255 | 2 |
| 16.9401 | 2 |
| 19.9236 | 2 |
| 20.7194 | 2 |
| 22.1597 | 1 |
| 24.0114 | 2 |
| 24.0994 | 2 |
| 28.7943 | 2 |
| 32.0614 | 2 |
| 32.3071 | 2 |
| 34.5634 | 2 |
| 36.383 | 2 |
| 36.8106 | 2 |
| 40.0014 | 1 |
| 43.8803 | 2 |
| 44.0732 | 2 |
| 45.5809 | 2 |
| 46.316 | 2 |
| 48.2176 | 2 |

Table C.3: First positive eigenvalues of $\mathbb{Z}_{10}$ surface

| Eigenvalue | Numerical multiplicity |
|---|---|
| 3.10942 | 2 |
| 4.12201 | 2 |
| 6.96949 | 1 |
| 7.78 | 2 |
| 8.98983 | 4 |
| 11.1973 | 1 |
| 15.0295 | 1 |
| 16.7819 | 2 |
| 19.5961 | 2 |
| 20.6216 | 4 |
| 20.9365 | 1 |
| 24.5947 | 2 |
| 25.0792 | 2 |
| 26.7752 | 2 |
| 28.8984 | 2 |
| 34.0612 | 4 |
| 34.7219 | 1 |
| 36.2501 | 1 |
| 39.9781 | 2 |
| 40.008 | 3 |
| 46.176 | 2 |
| 46.2604 | 4 |
| 47.5021 | 1 |

Table C.4: First positive eigenvalues of the 12-gon with order 24 symmetry group

# Appendix D

# FreeFEM++ Code for Genus 3

Similarly to Appendix C, here we give the FreeFEM++ codes for particular surfaces of genus 3. We start with the Klein quartic:

```
int n=30;
int nev=100;

string fnm="klein_quartic-Eigenvalues-n30-nev100.txt"; //file name for saving

real cx = 1.031969722;
real r = 0.254875473;

real s=0.7770942488589633;
real x06=s*cos(5*pi/7);
real y06=s*sin(5*pi/7);
real x03=s*cos(2*pi/7);
real y03=s*sin(2*pi/7);
real x05=s*cos(4*pi/7);
real y05=s*sin(4*pi/7);
real x10=s*cos(9*pi/7);
real y10=s*sin(9*pi/7);
real x07=s*cos(6*pi/7);
real y07=s*sin(6*pi/7);
real x12=s*cos(11*pi/7);
real y12=s*sin(11*pi/7);
real x09=s*cos(8*pi/7);
real y09=s*sin(8*pi/7);
real x14=s*cos(13*pi/7);
real y14=s*sin(13*pi/7);
real x11=s*cos(10*pi/7);
real y11=s*sin(10*pi/7);
real x02=s*cos(pi/7);
real y02=s*sin(pi/7);
real x13=s*cos(12*pi/7);
real y13=s*sin(12*pi/7);
real x04=s*cos(3*pi/7);
real y04=s*sin(3*pi/7);

border G1 (t=9*pi/7,5*pi/7) { x=cx+r*cos(t); y=r*sin(t);label=1;};
border G2 (t=9*pi/7,5*pi/7) { x=cos(pi/7)*(cx+r*cos(t))-sin(pi/7)*(r*sin(t));
    y=sin(pi/7)*(cx+r*cos(t))+cos(pi/7)*(r*sin(t));} label=2;;
border G3 (t=9*pi/7,5*pi/7) { x=cos(2*pi/7)*(cx+r*cos(t))-sin(2*pi/7)*r*sin(t);
    y=sin(2*pi/7)*(cx+r*cos(t))+cos(2*pi/7)*r*sin(t);} label=3;;
border G4 (t=9*pi/7,5*pi/7) { x=cos(3*pi/7)*(cx+r*cos(t))-sin(3*pi/7)*r*sin(t);
    y=sin(3*pi/7)*(cx+r*cos(t))+cos(3*pi/7)*r*sin(t);label=4;};
border G5 (t=9*pi/7,5*pi/7) { x=cos(4*pi/7)*(cx+r*cos(t))-sin(4*pi/7)*r*sin(t);
```

```
            y=sin(4*pi/7)*(cx+r*cos(t))+cos(4*pi/7)*r*sin(t);label=5;};
border G6 (t=9*pi/7,5*pi/7) { x=cos(5*pi/7)*(cx+r*cos(t))-sin(5*pi/7)*r*sin(t);
            y=sin(5*pi/7)*(cx+r*cos(t))+cos(5*pi/7)*r*sin(t);label=6;};
border G7 (t=9*pi/7,5*pi/7) { x=cos(6*pi/7)*(cx+r*cos(t))-sin(6*pi/7)*r*sin(t);
            y=sin(6*pi/7)*(cx+r*cos(t))+cos(6*pi/7)*r*sin(t);label=7;};
border G8 (t=9*pi/7,5*pi/7) { x=cos(7*pi/7)*(cx+r*cos(t))-sin(7*pi/7)*r*sin(t);
            y=sin(7*pi/7)*(cx+r*cos(t))+cos(7*pi/7)*r*sin(t);label=8;};
border G9 (t=9*pi/7,5*pi/7) { x=cos(8*pi/7)*(cx+r*cos(t))-sin(8*pi/7)*r*sin(t);
            y=sin(8*pi/7)*(cx+r*cos(t))+cos(8*pi/7)*r*sin(t);label=9;};
border G10 (t=9*pi/7,5*pi/7) { x=cos(9*pi/7)*(cx+r*cos(t))-sin(9*pi/7)*r*sin(t);
            y=sin(9*pi/7)*(cx+r*cos(t))+cos(9*pi/7)*r*sin(t);label=10;};
border G11 (t=9*pi/7,5*pi/7) { x=cos(10*pi/7)*(cx+r*cos(t))-sin(10*pi/7)*r*sin(t);
            y=sin(10*pi/7)*(cx+r*cos(t))+cos(10*pi/7)*r*sin(t);label=11;};
border G12 (t=9*pi/7,5*pi/7) { x=cos(11*pi/7)*(cx+r*cos(t))-sin(11*pi/7)*r*sin(t);
            y=sin(11*pi/7)*(cx+r*cos(t))+cos(11*pi/7)*r*sin(t);label=12;};
border G13 (t=9*pi/7,5*pi/7) { x=cos(12*pi/7)*(cx+r*cos(t))-sin(12*pi/7)*r*sin(t);
            y=sin(12*pi/7)*(cx+r*cos(t))+cos(12*pi/7)*r*sin(t);label=13;};
border G14 (t=9*pi/7,5*pi/7) { x=cos(13*pi/7)*(cx+r*cos(t))-sin(13*pi/7)*r*sin(t);
            y=sin(13*pi/7)*(cx+r*cos(t))+cos(13*pi/7)*r*sin(t);label=14;};

plot(G1(n)+G2(n)+G3(n)+G4(n)+G5(n)+G6(n)+G7(n)+G8(n)+G9(n)+G10(n)+G11(n)+G12(n)+G13(n)
    +G14(n));
mesh Th=buildmesh(G1(n)+G2(n)+G3(n)+G4(n)+G5(n)+G6(n)+G7(n)+G8(n)+G9(n)+G10(n)+G11(n)
    +G12(n)+G13(n)+G14(n),fixeborder=true);
plot(Th,wait=true,fill=true);

fespace Vh(Th,P2,periodic=[[1,y],[6,(cos(2*pi/7))*(y-y06)-(sin(-2*pi/7))*(x-x06)],
    [3,(cos(2*pi/7))*(y-y03)-(sin(2*pi/7))*(x-x03)],[8,y],
    [5,(cos(4*pi/7))*(y-y05)-(sin(4*pi/7))*(x-x05)],
    [10,(cos(2*pi/7))*(y-y10)-(sin(2*pi/7))*(x-x10)],
    [7,(cos(6*pi/7))*(y-y07)-(sin(6*pi/7))*(x-x07)],
    [12,(cos(4*pi/7))*(y-y12)-(sin(4*pi/7))*(x-x12)],
    [9,(cos(6*pi/7))*(y-y09)-(sin(-6*pi/7))*(x-x09)],
    [14,(cos(6*pi/7))*(y-y14)-(sin(6*pi/7))*(x-x14)],
    [11,(cos(4*pi/7))*(y-y11)-(sin(-4*pi/7))*(x-x11)],
    [2,(cos(6*pi/7))*(y-y02)-(sin(-6*pi/7))*(x-x02)],
    [13,(cos(2*pi/7))*(y-y13)-(sin(-2*pi/7))*(x-x13)],
    [4,(cos(4*pi/7))*(y-y04)-(sin(-4*pi/7))*(x-x04)]]);

Vh u1, u2;
real sigma = 0.00001;

varf op(u1,u2)= int2d(Th)(dx(u1)*dx(u2)+ dy(u1)*dy(u2)- sigma*(4*u1*u2)/(1-x^2-y^2)^2);
varf b([u1],[u2]) = int2d(Th)((4*u1*u2)/(1-x^2-y^2)^2);

matrix OP= op(Vh,Vh,solver=Crout,factorize=1);
matrix B=b(Vh,Vh,solver=CG,eps=1e-20);

real[int] ev(nev);
Vh[int] eV(nev);

int k=EigenValue(OP,B,sym=true,sigma=sigma,value=ev,vector=eV, tol=1e-10,maxit=0,ncv=0);
for (int i=0;i<k;i++) {
u1=eV[i];
real gg = int2d(Th)(dx(u1)*dx(u1) + dy(u1)*dy(u1)); real mm= int2d(Th)(u1*u1) ;

ofstream Eva(fnm, append);
Eva <<   ev[i] << "\n";

}
```

Of relevance to the Klein quartic (and indeed to Hurwitz surfaces in other genera) is the (2, 3, 7) triangle. We include the code for its spectrum; as with the (2, 3, 8) triangle and the pentagon in Appendix C, it is possible to add a Dirichlet condition to one or more of its sides.

```
int n = 30;
int nev=100;

string fnm="237-Eigenvalues-n30-nev100.txt"; //file name for saving

real c = 2.012192172612278;
real r = 1.74611492735219;

border G1 (t=0,0.2660772452600879) { x=t; y=0;};
border G2 (t=0.3007426187463789,0) { x=cos(pi/7)*t; y=t*sin(pi/7);};
border G3 (t=pi,3.066792828504322) { x=c+r*cos(t); y=r*sin(t);};

plot(G1(n)+G2(n)+G3(n));
mesh Th=buildmesh(G1(n)+G2(n)+G3(n));
plot(Th);

fespace Vh(Th,P2);

Vh u1, u2;
real sigma = 0.0005;

varf op(u1,u2)= int2d(Th)(dx(u1)*dx(u2)+ dy(u1)*dy(u2)- sigma*(4*u1*u2)/(1-x^2-y^2)^2);
varf b([u1],[u2]) = int2d(Th)((4*u1*u2)/(1-x^2-y^2)^2);

matrix OP= op(Vh,Vh,solver=Crout,factorize=1);
matrix B=b(Vh,Vh,solver=CG,eps=1e-20);

real[int] ev(nev);
Vh[int] eV(nev);
int k=EigenValue(OP,B,sym=true,sigma=sigma,value=ev,vector=eV, tol=1e-10,maxit=0,ncv=0);
for (int i=0;i<k;i++) {
u1=eV[i];
real gg = int2d(Th)(dx(u1)*dx(u1) + dy(u1)*dy(u1)); real mm= int2d(Th)(u1*u1) ;

ofstream Eva(fnm, append);
Eva <<  ev[i] << "\n";

}
```

Next, we include the code for the surface $M3$ in [72]:

```
int n = 30;
int nev = 100;

string fnm="m3-Eigenvalues-n30-nev100.txt"; //file name for saving

real cx = 1.0122;
real cy = 0.1566;
real r = 0.2214;

real s = 0.8045;
real theta = 0.1238;
real x01 = s*cos(theta);
real y01 = s*sin(theta);
real x02 = s*cos(pi/6-theta);
real y02 = s*sin(pi/6-theta);
real x03 = s*cos(pi/6+theta);
real y03 = s*sin(pi/6+theta);
real x04 = s*cos(2*pi/6-theta);
real y04 = s*sin(2*pi/6-theta);
real x05 = s*cos(2*pi/6+theta);
```

```
real y05 = s*sin(2*pi/6+theta);
real x06 = s*cos(3*pi/6-theta);
real y06 = s*sin(3*pi/6-theta);
real x07 = s*cos(3*pi/6+theta);
real y07 = s*sin(3*pi/6+theta);
real x08 = s*cos(4*pi/6-theta);
real y08 = s*sin(4*pi/6-theta);
real x09 = s*cos(4*pi/6+theta);
real y09 = s*sin(4*pi/6+theta);
real x10 = s*cos(5*pi/6-theta);
real y10 = s*sin(5*pi/6-theta);
real x11 = s*cos(5*pi/6+theta);
real y11 = s*sin(5*pi/6+theta);
real x12 = s*cos(6*pi/6-theta);
real y12 = s*sin(6*pi/6-theta);
real x13 = s*cos(6*pi/6+theta);
real y13 = s*sin(6*pi/6+theta);
real x14 = s*cos(7*pi/6-theta);
real y14 = s*sin(7*pi/6-theta);
real x15 = s*cos(7*pi/6+theta);
real y15 = s*sin(7*pi/6+theta);
real x16 = s*cos(8*pi/6-theta);
real y16 = s*sin(8*pi/6-theta);
real x17 = s*cos(8*pi/6+theta);
real y17 = s*sin(8*pi/6+theta);
real x18 = s*cos(9*pi/6-theta);
real y18 = s*sin(9*pi/6-theta);
real x19 = s*cos(9*pi/6+theta);
real y19 = s*sin(9*pi/6+theta);
real x20 = s*cos(10*pi/6-theta);
real y20 = s*sin(10*pi/6-theta);
real x21 = s*cos(10*pi/6+theta);
real y21 = s*sin(10*pi/6+theta);
real x22 = s*cos(11*pi/6-theta);
real y22 = s*sin(11*pi/6-theta);
real x23 = s*cos(11*pi/6+theta);
real y23 = s*sin(11*pi/6+theta);
real x24 = s*cos(theta);
real y24 = s*sin(-theta);

border G1 (t=15*pi/12,11*pi/12) { x=cx+r*cos(t); y=cy+r*sin(t);};
border G2 (t=11*pi/12,15*pi/12) { x=cos(pi/6)*(cx+r*cos(t))-sin(pi/6)*(-cy-r*sin(t));
    y=sin(pi/6)*(cx+r*cos(t))+cos(pi/6)*(-cy-r*sin(t));};
border G3 (t=15*pi/12,11*pi/12) { x=cos(pi/6)*(cx+r*cos(t))-sin(pi/6)*(cy+r*sin(t));
    y=sin(pi/6)*(cx+r*cos(t))+cos(pi/6)*(cy+r*sin(t));};
border G4 (t=11*pi/12,15*pi/12) { x=cos(2*pi/6)*(cx+r*cos(t))-sin(2*pi/6)*(-cy-r*sin(t));
    y=sin(2*pi/6)*(cx+r*cos(t))+cos(2*pi/6)*(-cy-r*sin(t));};
border G5 (t=15*pi/12,11*pi/12) { x=cos(2*pi/6)*(cx+r*cos(t))-sin(2*pi/6)*(cy+r*sin(t));
    y=sin(2*pi/6)*(cx+r*cos(t))+cos(2*pi/6)*(cy+r*sin(t));};
border G6 (t=11*pi/12,15*pi/12) { x=cos(3*pi/6)*(cx+r*cos(t))-sin(3*pi/6)*(-cy-r*sin(t));
    y=sin(3*pi/6)*(cx+r*cos(t))+cos(3*pi/6)*(-cy-r*sin(t));};
border G7 (t=15*pi/12,11*pi/12) { x=cos(3*pi/6)*(cx+r*cos(t))-sin(3*pi/6)*(cy+r*sin(t));
    y=sin(3*pi/6)*(cx+r*cos(t))+cos(3*pi/6)*(cy+r*sin(t));};
border G8 (t=11*pi/12,15*pi/12) { x=cos(4*pi/6)*(cx+r*cos(t))-sin(4*pi/6)*(-cy-r*sin(t));
    y=sin(4*pi/6)*(cx+r*cos(t))+cos(4*pi/6)*(-cy-r*sin(t));};
border G9 (t=15*pi/12,11*pi/12) { x=cos(4*pi/6)*(cx+r*cos(t))-sin(4*pi/6)*(cy+r*sin(t));
    y=sin(4*pi/6)*(cx+r*cos(t))+cos(4*pi/6)*(cy+r*sin(t));};
border G10 (t=11*pi/12,15*pi/12) { x=cos(5*pi/6)*(cx+r*cos(t))-sin(5*pi/6)*(-cy-r*sin(t));
    y=sin(5*pi/6)*(cx+r*cos(t))+cos(5*pi/6)*(-cy-r*sin(t));};
border G11 (t=15*pi/12,11*pi/12) { x=cos(5*pi/6)*(cx+r*cos(t))-sin(5*pi/6)*(cy+r*sin(t));
    y=sin(5*pi/6)*(cx+r*cos(t))+cos(5*pi/6)*(cy+r*sin(t));};
border G12 (t=11*pi/12,15*pi/12) { x=cos(6*pi/6)*(cx+r*cos(t))-sin(6*pi/6)*(-cy-r*sin(t));
    y=sin(6*pi/6)*(cx+r*cos(t))+cos(6*pi/6)*(-cy-r*sin(t));};
border G13 (t=15*pi/12,11*pi/12) { x=cos(6*pi/6)*(cx+r*cos(t))-sin(6*pi/6)*(cy+r*sin(t));
```

```
        y=sin(6*pi/6)*(cx+r*cos(t))+cos(6*pi/6)*(cy+r*sin(t));};
border G14 (t=11*pi/12,15*pi/12) { x=cos(7*pi/6)*(cx+r*cos(t))-sin(7*pi/6)*(-cy-r*sin(t));
        y=sin(7*pi/6)*(cx+r*cos(t))+cos(7*pi/6)*(-cy-r*sin(t));};
border G15 (t=15*pi/12,11*pi/12) { x=cos(7*pi/6)*(cx+r*cos(t))-sin(7*pi/6)*(cy+r*sin(t));
        y=sin(7*pi/6)*(cx+r*cos(t))+cos(7*pi/6)*(cy+r*sin(t));};
border G16 (t=11*pi/12,15*pi/12) { x=cos(8*pi/6)*(cx+r*cos(t))-sin(8*pi/6)*(-cy-r*sin(t));
        y=sin(8*pi/6)*(cx+r*cos(t))+cos(8*pi/6)*(-cy-r*sin(t));};
border G17 (t=15*pi/12,11*pi/12) { x=cos(8*pi/6)*(cx+r*cos(t))-sin(8*pi/6)*(cy+r*sin(t));
        y=sin(8*pi/6)*(cx+r*cos(t))+cos(8*pi/6)*(cy+r*sin(t));};
border G18 (t=11*pi/12,15*pi/12) { x=cos(9*pi/6)*(cx+r*cos(t))-sin(9*pi/6)*(-cy-r*sin(t));
        y=sin(9*pi/6)*(cx+r*cos(t))+cos(9*pi/6)*(-cy-r*sin(t));};
border G19 (t=15*pi/12,11*pi/12) { x=cos(9*pi/6)*(cx+r*cos(t))-sin(9*pi/6)*(cy+r*sin(t));
        y=sin(9*pi/6)*(cx+r*cos(t))+cos(9*pi/6)*(cy+r*sin(t));};
border G20 (t=11*pi/12,15*pi/12) { x=cos(10*pi/6)*(cx+r*cos(t))-sin(10*pi/6)*(-cy-r*sin(t));
        y=sin(10*pi/6)*(cx+r*cos(t))+cos(10*pi/6)*(-cy-r*sin(t));};
border G21 (t=15*pi/12,11*pi/12) { x=cos(10*pi/6)*(cx+r*cos(t))-sin(10*pi/6)*(cy+r*sin(t));
        y=sin(10*pi/6)*(cx+r*cos(t))+cos(10*pi/6)*(cy+r*sin(t));};
border G22 (t=11*pi/12,15*pi/12) { x=cos(11*pi/6)*(cx+r*cos(t))-sin(11*pi/6)*(-cy-r*sin(t));
        y=sin(11*pi/6)*(cx+r*cos(t))+cos(11*pi/6)*(-cy-r*sin(t));};
border G23 (t=15*pi/12,11*pi/12) { x=cos(11*pi/6)*(cx+r*cos(t))-sin(11*pi/6)*(cy+r*sin(t));
        y=sin(11*pi/6)*(cx+r*cos(t))+cos(11*pi/6)*(cy+r*sin(t));};
border G24 (t=11*pi/12,15*pi/12) { x=cx+r*cos(t); y=-cy-r*sin(t);};

plot(G1(n)+G2(n)+G3(n)+G4(n)+G5(n)+G6(n)+G7(n)+G8(n)+G9(n)+G10(n)+G11(n)+G12(n)+G13(n)
    +G14(n)+G15(n)+G16(n)+G17(n)+G18(n)+G19(n)+G20(n)+G21(n)+G22(n)+G23(n)+G24(n));
mesh Th=buildmesh(G1(n)+G2(n)+G3(n)+G4(n)+G5(n)+G6(n)+G7(n)+G8(n)+G9(n)+G10(n)+G11(n)
    +G12(n)+G13(n)+G14(n)+G15(n)+G16(n)+G17(n)+G18(n)+G19(n)+G20(n)+G21(n)+G22(n)+G23(n)
    +G24(n),fixeborder=true);
plot(Th,wait=true,fill=true);

fespace Vh(Th,P2,periodic=[[2,(cos(pi/12))*(y-y02)-(sin(pi/12))*(x-x02)],
    [9,(cos(3*pi/12))*(y-y09)-(sin(-3*pi/12))*(x-x09)],
    [4,(cos(3*pi/12))*(y-y04)-(sin(3*pi/12))*(x-x04)],
    [11,(cos(pi/12))*(y-y11)-(sin(-pi/12))*(x-x11)],
    [6,(cos(5*pi/12))*(y-y06)-(sin(5*pi/12))*(x-x06)],
    [13,(cos(pi/12))*(y-y13)-(sin(pi/12))*(x-x13)],
    [8,(cos(7*pi/12))*(y-y08)-(sin(7*pi/12))*(x-x08)],
    [15,(cos(3*pi/12))*(y-y15)-(sin(3*pi/12))*(x-x15)],
    [10,(cos(9*pi/12))*(y-y10)-(sin(9*pi/12))*(x-x10)],
    [17,(cos(5*pi/12))*(y-y17)-(sin(5*pi/12))*(x-x17)],
    [12,(cos(11*pi/12))*(y-y12)-(sin(11*pi/12))*(x-x12)],
    [19,(cos(7*pi/12))*(y-y19)-(sin(7*pi/12))*(x-x19)],
    [14,(cos(13*pi/12))*(y-y14)-(sin(13*pi/12))*(x-x14)],
    [21,(cos(9*pi/12))*(y-y21)-(sin(9*pi/12))*(x-x21)],
    [16,(cos(15*pi/12))*(y-y16)-(sin(15*pi/12))*(x-x16)],
    [23,(cos(11*pi/12))*(y-y23)-(sin(11*pi/12))*(x-x23)],
    [18,(cos(17*pi/12))*(y-y18)-(sin(17*pi/12))*(x-x18)],
    [1,(cos(13*pi/12))*(y-y01)-(sin(13*pi/12))*(x-x01)],
    [20,(cos(19*pi/12))*(y-y20)-(sin(19*pi/12))*(x-x20)],
    [3,(cos(15*pi/12))*(y-y03)-(sin(15*pi/12))*(x-x03)],
    [22,(cos(21*pi/12))*(y-y22)-(sin(21*pi/12))*(x-x22)],
    [5,(cos(17*pi/12))*(y-y05)-(sin(17*pi/12))*(x-x05)],
    [24,(cos(23*pi/12))*(y-y24)-(sin(23*pi/12))*(x-x24)],
    [7,(cos(19*pi/12))*(y-y07)-(sin(19*pi/12))*(x-x07)]]);

Vh u1, u2;
real sigma = 0.00001;
varf op(u1,u2)= int2d(Th)(dx(u1)*dx(u2)+ dy(u1)*dy(u2)- sigma*(4*u1*u2)/(1-x^2-y^2)^2);
varf b([u1],[u2]) = int2d(Th)((4*u1*u2)/(1-x^2-y^2)^2);

matrix OP= op(Vh,Vh,solver=Crout,factorize=1);
matrix B=b(Vh,Vh,solver=CG,eps=1e-20);

real[int] ev(nev);
```

```
Vh[int] eV(nev);

int k=EigenValue(OP,B,sym=true,sigma=sigma,value=ev,vector=eV, tol=1e-10,maxit=0,ncv=0);
for (int i=0;i<k;i++) {
u1=eV[i];
real gg = int2d(Th)(dx(u1)*dx(u1) + dy(u1)*dy(u1));  real mm= int2d(Th)(u1*u1) ;

ofstream Eva(fnm, append);
Eva <<   ev[i] << "\n";

}
```

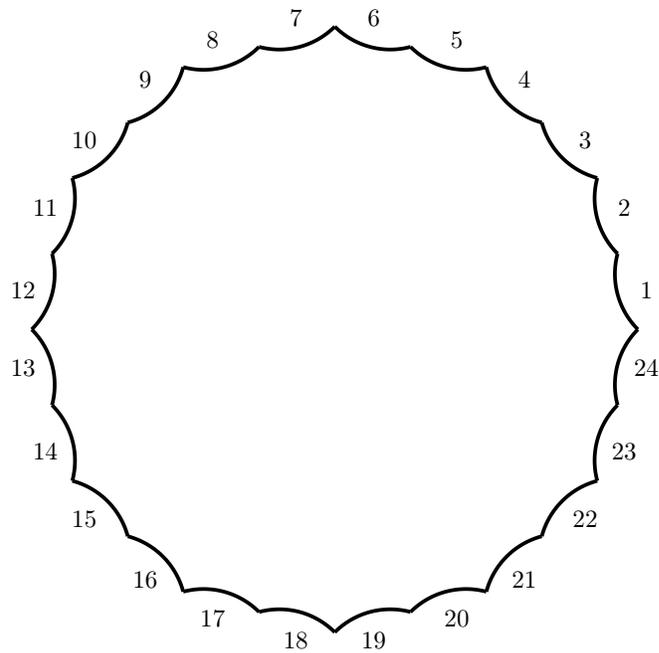

Figure D.1: Fundamental polygon for the surface $M3$

We also include the code for the Fermat quartic:

```
int n=30;
int nev=100;

string fnm="fermat_quartic-Eigenvalues-n30-nev100.txt"; //file name for saving

real cx = 1.0044;
real cy = 0.2275;
real r = 0.2463;

real s=0.7881;
real theta=0.1698;
real x01=s*cos(theta);
real y01=s*sin(theta);
real x06=s*cos(3*pi/4-theta);
real y06=s*sin(3*pi/4-theta);
```

```
real x03=s*cos(pi/4+theta);
real y03=s*sin(pi/4+theta);
real x08=s*cos(pi-theta);
real y08=s*sin(pi-theta);
real x05=s*cos(pi/2+theta);
real y05=s*sin(pi/2+theta);
real x10=s*cos(5*pi/4-theta);
real y10=s*sin(5*pi/4-theta);
real x07=s*cos(3*pi/4+theta);
real y07=s*sin(3*pi/4+theta);
real x12=s*cos(6*pi/4-theta);
real y12=s*sin(6*pi/4-theta);
real x09=s*cos(pi+theta);
real y09=s*sin(pi+theta);
real x14=s*cos(7*pi/4-theta);
real y14=s*sin(7*pi/4-theta);
real x11=s*cos(5*pi/4+theta);
real y11=s*sin(5*pi/4+theta);
real x16=s*cos(theta);
real y16=s*sin(-theta);
real x13=s*cos(6*pi/4+theta);
real y13=s*sin(6*pi/4+theta);
real x02=s*cos(pi/4-theta);
real y02=s*sin(pi/4-theta);
real x15=s*cos(7*pi/4+theta);
real y15=s*sin(7*pi/4+theta);
real x04=s*cos(pi/2-theta);
real y04=s*sin(pi/2-theta);

border G1 (t=11*pi/8,7*pi/8) { x=cx+r*cos(t); y=cy+r*sin(t);};
border G2 (t=7*pi/8,11*pi/8) { x=cos(pi/4)*(cx+r*cos(t))-sin(pi/4)*(-cy-r*sin(t));
    y=sin(pi/4)*(cx+r*cos(t))+cos(pi/4)*(-cy-r*sin(t));};
border G3 (t=11*pi/8,7*pi/8) { x=cos(pi/4)*(cx+r*cos(t))-sin(pi/4)*(cy+r*sin(t));
    y=sin(pi/4)*(cx+r*cos(t))+cos(pi/4)*(cy+r*sin(t));};
border G4 (t=7*pi/8,11*pi/8) { x=cos(pi/2)*(cx+r*cos(t))-sin(pi/2)*(-cy-r*sin(t));
    y=sin(pi/2)*(cx+r*cos(t))+cos(pi/2)*(-cy-r*sin(t));};
border G5 (t=11*pi/8,7*pi/8) { x=cos(pi/2)*(cx+r*cos(t))-sin(pi/2)*(cy+r*sin(t));
    y=sin(pi/2)*(cx+r*cos(t))+cos(pi/2)*(cy+r*sin(t));};
border G6 (t=7*pi/8,11*pi/8) { x=cos(3*pi/4)*(cx+r*cos(t))-sin(3*pi/4)*(-cy-r*sin(t));
    y=sin(3*pi/4)*(cx+r*cos(t))+cos(3*pi/4)*(-cy-r*sin(t));};
border G7 (t=11*pi/8,7*pi/8) { x=cos(3*pi/4)*(cx+r*cos(t))-sin(3*pi/4)*(cy+r*sin(t));
    y=sin(3*pi/4)*(cx+r*cos(t))+cos(3*pi/4)*(cy+r*sin(t));};
border G8 (t=7*pi/8,11*pi/8) { x=cos(pi)*(cx+r*cos(t))-sin(pi)*(-cy-r*sin(t));
    y=sin(pi)*(cx+r*cos(t))+cos(pi)*(-cy-r*sin(t));};
border G9 (t=11*pi/8,7*pi/8) { x=cos(pi)*(cx+r*cos(t))-sin(pi)*(cy+r*sin(t));
    y=sin(pi)*(cx+r*cos(t))+cos(pi)*(cy+r*sin(t));};
border G10 (t=7*pi/8,11*pi/8) { x=cos(5*pi/4)*(cx+r*cos(t))-sin(5*pi/4)*(-cy-r*sin(t));
    y=sin(5*pi/4)*(cx+r*cos(t))+cos(5*pi/4)*(-cy-r*sin(t));};
border G11 (t=11*pi/8,7*pi/8) { x=cos(5*pi/4)*(cx+r*cos(t))-sin(5*pi/4)*(cy+r*sin(t));
    y=sin(5*pi/4)*(cx+r*cos(t))+cos(5*pi/4)*(cy+r*sin(t));};
border G12 (t=7*pi/8,11*pi/8) { x=cos(6*pi/4)*(cx+r*cos(t))-sin(6*pi/4)*(-cy-r*sin(t));
    y=sin(6*pi/4)*(cx+r*cos(t))+cos(6*pi/4)*(-cy-r*sin(t));};
border G13 (t=11*pi/8,7*pi/8) { x=cos(6*pi/4)*(cx+r*cos(t))-sin(6*pi/4)*(cy+r*sin(t));
    y=sin(6*pi/4)*(cx+r*cos(t))+cos(6*pi/4)*(cy+r*sin(t));};
border G14 (t=7*pi/8,11*pi/8) { x=cos(7*pi/4)*(cx+r*cos(t))-sin(7*pi/4)*(-cy-r*sin(t));
    y=sin(7*pi/4)*(cx+r*cos(t))+cos(7*pi/4)*(-cy-r*sin(t));};
border G15 (t=11*pi/8,7*pi/8) { x=cos(7*pi/4)*(cx+r*cos(t))-sin(7*pi/4)*(cy+r*sin(t));
    y=sin(7*pi/4)*(cx+r*cos(t))+cos(7*pi/4)*(cy+r*sin(t));};
border G16 (t=7*pi/8,11*pi/8) { x=cx+r*cos(t); y=-cy-r*sin(t);};

plot(G1(n)+G2(n)+G3(n)+G4(n)+G5(n)+G6(n)+G7(n)+G8(n)+G9(n)+G10(n)+G11(n)+G12(n)+G13(n)
    +G14(n)+G15(n)+G16(n));
mesh Th=buildmesh(G1(n)+G2(n)+G3(n)+G4(n)+G5(n)+G6(n)+G7(n)+G8(n)+G9(n)+G10(n)+G11(n)
    +G12(n)+G13(n)+G14(n)+G15(n)+G16(n),fixeborder=true);
```

```
plot(Th,wait=true,fill=true);

fespace Vh(Th,P2,periodic=[[1,cos(pi/8)*(y-y01)-sin(pi/8)*(x-x01)],
    [6,cos(3*pi/8)*(y-y06)-sin(-3*pi/8)*(x-x06)],
    [3,cos(3*pi/8)*(y-y03)-sin(3*pi/8)*(x-x03)],
    [8,cos(pi/8)*(y-y08)-sin(-pi/8)*(x-x08)],
    [5,cos(5*pi/8)*(y-y05)-sin(5*pi/8)*(x-x05)],
    [10,cos(pi/8)*(y-y10)-sin(pi/8)*(x-x10)],
    [7,cos(7*pi/8)*(y-y07)-sin(7*pi/8)*(x-x07)],
    [12,cos(3*pi/8)*(y-y12)-sin(3*pi/8)*(x-x12)],
    [9,cos(9*pi/8)*(y-y09)-sin(9*pi/8)*(x-x09)],
    [14,cos(5*pi/8)*(y-y14)-sin(5*pi/8)*(x-x14)],
    [11,cos(11*pi/8)*(y-y11)-sin(11*pi/8)*(x-x11)],
    [16,cos(7*pi/8)*(y-y16)-sin(7*pi/8)*(x-x16)],
    [13,cos(13*pi/8)*(y-y13)-sin(13*pi/8)*(x-x13)],
    [2,cos(9*pi/8)*(y-y02)-sin(9*pi/8)*(x-x02)],
    [15,cos(15*pi/8)*(y-y15)-sin(15*pi/8)*(x-x15)],
    [4,cos(11*pi/8)*(y-y04)-sin(11*pi/8)*(x-x04)]]);

Vh u1, u2;
real sigma = 0.00001;
varf op(u1,u2)= int2d(Th)(dx(u1)*dx(u2)+ dy(u1)*dy(u2)- sigma*(4*u1*u2)/(1-x^2-y^2)^2);
varf b([u1],[u2]) = int2d(Th)((4*u1*u2)/(1-x^2-y^2)^2);

matrix OP= op(Vh,Vh,solver=Crout,factorize=1);
matrix B=b(Vh,Vh,solver=CG,eps=1e-20);

real[int] ev(nev);
Vh[int] eV(nev);

int k=EigenValue(OP,B,sym=true,sigma=sigma,value=ev,vector=eV, tol=1e-10,maxit=0,ncv=0);
for (int i=0;i<k;i++) {
u1=eV[i];
real gg = int2d(Th)(dx(u1)*dx(u1) + dy(u1)*dy(u1)); real mm= int2d(Th)(u1*u1) ;

ofstream Eva(fnm, append);
Eva <<  ev[i] << "\n";

}
```

**Remark D.1** *The observant reader will recognize $\lambda_2$ in Table D.3 - it is roughly equal to the first positive eigenvalue of the Bolza surface! Moreover, it appears to have the same multiplicity. In fact, if one takes the subset of eigenvalues in Table D.3 with multiplicity 1, 2 or 3, one finds that they match the subset of eigenvalues of Table C.1, and that the multiplicities are the same. This is no coincidence. Recall the presentation of the automorphism group given in Chapter 1:*

$$\left\langle a,b,c \,|\, a^8 = b^3 = c^2 = (ab)^2 = (a^2b^2)^3 = (a^4b^2)^3 = acac = bcbc = e \right\rangle,$$

*where $a$ is a rotation of order 8 around the centre of the fundamental 16-gon, $b$ is a rotation of order 3 around one of the 32 $(4,4,4)$ triangles that can be used to tessellate the surface (note this is twice as many as in the case of the Bolza surface), and $c$ is a reflection in the real line. As we have done with the Bolza surface and Klein quartic, we can input this presentation into GAP to get information about the structure of the group*

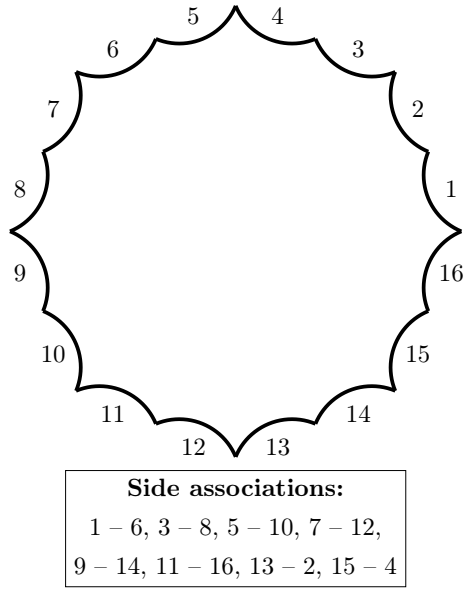

Figure D.2: Fundamental polygon for the Fermat quartic

and information about its representation theory. The code below indicates that the full automorphism group is isomorphic to $(((\mathbb{Z}_4 \times \mathbb{Z}_4) \rtimes \mathbb{Z}_3) \rtimes \mathbb{Z}_2) \rtimes \mathbb{Z}_2$.

```
gap> g:=FreeGroup(3);
<free group on the generators [ f1, f2, f3 ]>
gap> f:=g/[g.1^8,g.2^3,g.3^2,(g.1*g.2)^2,(g.1^2*g.2^2)^3,(g.1^4*g.2^2)^3,
    g.1*g.3*g.1*g.3,g.2*g.3*g.2*g.3];
<fp group on the generators [ f1, f2, f3 ]>
gap> Order(f);
192
gap> StructureDescription(f);
"(((C4xC4):C3):C2):C2"
gap> RequirePackage("repsn");
─────────────────────────────────────────────────
Repsn for Constructing Representations of Finite Groups
                Version 3.0.2

                   Written by
                Vahid Dabbaghian
─────────────────────────────────────────────────
true
gap> Display(CharacterTable(f));
CT1

     2  6  4  4  1  5  6  4  3  4  3  4  5  1  4
     3  1  .  .  1  .  .  1  .  .  .  .  1  .

        1a 2a 2b 3a 4a 2c 2d 8a 4b 8b 4c 4d 6a 2e
    2P  1a 1a 1a 3a 2c 1a 1a 4d 2c 4a 2c 2c 3a 1a
    3P  1a 2a 2b 1a 4a 2c 2d 8a 4b 8b 4c 4d 2d 2e
    5P  1a 2a 2b 3a 4a 2c 2d 8a 4b 8b 4c 4d 6a 2e
```

```
        7P  1a  2a  2b  3a  4a  2c  2d  8a  4b  8b  4c  4d  6a  2e

X.1          1   1   1   1   1   1   1   1   1   1   1   1   1   1
X.2          1  -1  -1   1   1   1   1  -1  -1  -1  -1   1   1   1
X.3          1  -1   1   1   1   1  -1  -1  -1   1   1   1  -1  -1
X.4          1   1  -1   1   1   1  -1   1   1  -1  -1   1  -1  -1
X.5          2   .   .  -1   2   2  -2   .   .   .   .   2   1  -2
X.6          2   .   .  -1   2   2   2   .   .   .   .   2  -1   2
X.7          3  -1  -1   .  -1   3   3   1  -1   1  -1  -1   .  -1
X.8          3  -1   1   .  -1   3  -3   1  -1  -1   1  -1   .   1
X.9          3   1  -1   .  -1   3  -3  -1   1   1  -1  -1   .   1
X.10         3   1   1   .  -1   3   3  -1   1  -1   1  -1   .  -1
X.11         6  -2   .   .   2  -2   .   .   2   .   .  -2   .   .
X.12         6   .  -2   .  -2  -2   .   .   .   .   2   2   .   .
X.13         6   .   2   .  -2  -2   .   .   .   .  -2   2   .   .
X.14         6   2   .   .   2  -2   .   .  -2   .   .  -2   .   .
```
```
gap> chi:=Irr(f);;
gap> IrreducibleAffordingRepresentation(chi[1]);
[ f1, f2, f3 ] -> [ [ [ 1 ] ], [ [ 1 ] ], [ [ 1 ] ] ]
gap> IrreducibleAffordingRepresentation(chi[2]);
[ f1, f2, f3 ] -> [ [ [ -1 ] ], [ [ 1 ] ], [ [ -1 ] ] ]
gap> IrreducibleAffordingRepresentation(chi[3]);
[ f1, f2, f3 ] -> [ [ [ -1 ] ], [ [ 1 ] ], [ [ 1 ] ] ]
gap> IrreducibleAffordingRepresentation(chi[4]);
[ f1, f2, f3 ] -> [ [ [ 1 ] ], [ [ 1 ] ], [ [ -1 ] ] ]
gap> IrreducibleAffordingRepresentation(chi[5]);
[ f1, f2, f3 ] -> [ [ [ 0, -1 ], [ -1, 0 ] ], [ [ E(3)^2, 0 ], [ 0, E(3) ] ],
  [ [ 0, 1 ], [ 1, 0 ] ] ]
gap> IrreducibleAffordingRepresentation(chi[6]);
[ f1, f2, f3 ] -> [ [ [ 0, 1 ], [ 1, 0 ] ], [ [ E(3)^2, 0 ], [ 0, E(3) ] ],
  [ [ 0, 1 ], [ 1, 0 ] ] ]
gap> IrreducibleAffordingRepresentation(chi[7]);
[ f1, f2, f3 ] -> [ [ [ 0, 1, 0 ], [ -1, 0, 0 ], [ 0, 0, 1 ] ],
  [ [ 0, 1, 0 ], [ 0, 0, 1 ], [ 1, 0, 0 ] ],
  [ [ 0, -1, 0 ], [ -1, 0, 0 ], [ 0, 0, -1 ] ] ]
gap> IrreducibleAffordingRepresentation(chi[8]);
[ f1, f2, f3 ] -> [ [ [ 0, 1, 0 ], [ -1, 0, 0 ], [ 0, 0, 1 ] ],
  [ [ 0, 1, 0 ], [ 0, 0, 1 ], [ 1, 0, 0 ] ],
  [ [ 0, 1, 0 ], [ 1, 0, 0 ], [ 0, 0, 1 ] ] ]
gap> IrreducibleAffordingRepresentation(chi[9]);
[ f1, f2, f3 ] -> [ [ [ 0, -1, 0 ], [ 1, 0, 0 ], [ 0, 0, -1 ] ],
  [ [ 0, 1, 0 ], [ 0, 0, 1 ], [ 1, 0, 0 ] ],
  [ [ 0, -1, 0 ], [ -1, 0, 0 ], [ 0, 0, -1 ] ] ]
gap> IrreducibleAffordingRepresentation(chi[10]);
[ f1, f2, f3 ] -> [ [ [ 0, -1, 0 ], [ 1, 0, 0 ], [ 0, 0, -1 ] ],
  [ [ 0, 1, 0 ], [ 0, 0, 1 ], [ 1, 0, 0 ] ],
  [ [ 0, 1, 0 ], [ 1, 0, 0 ], [ 0, 0, 1 ] ] ]
gap> IrreducibleAffordingRepresentation(chi[11]);
[ f1, f2, f3 ] ->
[ [ [ 0, 1, 0, 0, 0, 0 ], [ E(4), 0, 0, 0, 0, 0 ], [ 0, 0, E(4), 0, 0, 0 ],
      [ 0, 0, 0, 0, -E(4), 0 ], [ 0, 0, 0, 1, 0, 0 ],
      [ 0, 0, 0, 0, 0, -E(4) ] ],
  [ [ 0, 0, 1, 0, 0, 0 ], [ -1, 0, 0, 0, 0, 0 ], [ 0, -1, 0, 0, 0, 0 ],
      [ 0, 0, 0, 0, -1, 0 ], [ 0, 0, 0, 0, 0, -1 ], [ 0, 0, 0, 1, 0, 0 ] ],
```

```
      [ [ 0, 0, 0, 1, 0, 0 ], [ 0, 0, 0, 0, 1, 0 ], [ 0, 0, 0, 0, 0, 1 ],
          [ 1, 0, 0, 0, 0, 0 ], [ 0, 1, 0, 0, 0, 0 ], [ 0, 0, 1, 0, 0, 0 ] ] ]
gap> IrreducibleAffordingRepresentation(chi[12]);
[ f1, f2, f3 ] ->
[ [ [ 0, 0, 0, 1, 0, 0 ], [ 0, 0, 0, 0, 1, 0 ], [ 0, 0, 0, 0, 0, 1 ],
      [ -E(4), 0, 0, 0, 0, 0 ], [ 0, 1, 0, 0, 0, 0 ], [ 0, 0, E(4), 0, 0, 0 ]
    ], [ [ 0, 1, 0, 0, 0, 0 ], [ 0, 0, 1, 0, 0, 0 ], [ 1, 0, 0, 0, 0, 0 ],
      [ 0, 0, 0, 0, 0, -E(4) ], [ 0, 0, 0, E(4), 0, 0 ], [ 0, 0, 0, 0, 1, 0 ]
    ],
  [ [ 0, 0, -1, 0, 0, 0 ], [ 0, -1, 0, 0, 0, 0 ], [ -1, 0, 0, 0, 0, 0 ],
      [ 0, 0, 0, 0, 0, E(4) ], [ 0, 0, 0, 0, -1, 0 ],
      [ 0, 0, 0, -E(4), 0, 0 ] ] ]
gap> IrreducibleAffordingRepresentation(chi[13]);
[ f1, f2, f3 ] ->
[ [ [ 0, 0, 0, 1, 0, 0 ], [ 0, 0, 0, 0, 1, 0 ], [ 0, 0, 0, 0, 0, 1 ],
      [ -E(4), 0, 0, 0, 0, 0 ], [ 0, 1, 0, 0, 0, 0 ], [ 0, 0, E(4), 0, 0, 0 ]
    ], [ [ 0, 1, 0, 0, 0, 0 ], [ 0, 0, 1, 0, 0, 0 ], [ 1, 0, 0, 0, 0, 0 ],
      [ 0, 0, 0, 0, 0, -E(4) ], [ 0, 0, 0, E(4), 0, 0 ], [ 0, 0, 0, 0, 1, 0 ]
    ], [ [ 0, 0, 1, 0, 0, 0 ], [ 0, 1, 0, 0, 0, 0 ], [ 1, 0, 0, 0, 0, 0 ],
      [ 0, 0, 0, 0, 0, -E(4) ], [ 0, 0, 0, 0, 1, 0 ], [ 0, 0, 0, E(4), 0, 0 ]
    ] ]
gap> IrreducibleAffordingRepresentation(chi[14]);
[ f1, f2, f3 ] ->
[ [ [ 0, 1, 0, 0, 0, 0 ], [ E(4), 0, 0, 0, 0, 0 ], [ 0, 0, -E(4), 0, 0, 0 ],
      [ 0, 0, 0, 0, -E(4), 0 ], [ 0, 0, 0, 1, 0, 0 ], [ 0, 0, 0, 0, 0, E(4) ]
    ], [ [ 0, 0, 1, 0, 0, 0 ], [ 1, 0, 0, 0, 0, 0 ], [ 0, 1, 0, 0, 0, 0 ],
      [ 0, 0, 0, 0, 1, 0 ], [ 0, 0, 0, 0, 0, 1 ], [ 0, 0, 0, 1, 0, 0 ] ],
  [ [ 0, 0, 0, 1, 0, 0 ], [ 0, 0, 0, 0, 1, 0 ], [ 0, 0, 0, 0, 0, 1 ],
      [ 1, 0, 0, 0, 0, 0 ], [ 0, 1, 0, 0, 0, 0 ], [ 0, 0, 1, 0, 0, 0 ] ] ]
```

| Eigenvalue | Numerical multiplicity |
|:---:|:---:|
| 2.6767 | 8 |
| 6.61848 | 7 |
| 10.8637 | 6 |
| 12.1775 | 8 |
| 17.2397 | 7 |
| 21.9563 | 7 |
| 24.0649 | 8 |
| 25.9085 | 6 |
| 30.7817 | 6 |
| 36.4369 | 8 |
| 37.4073 | 8 |
| 41.4847 | 6 |
| 44.8658 | 8 |
| 49.0165 | 6 |
| 50.6025 | 6 |
| 57.7794 | 8 |
| 57.8082 | 7 |
| 61.9464 | 7 |
| 67.7629 | 8 |
| 69.9676 | 6 |
| 71.0138 | 8 |
| 75.6357 | 6 |
| 81.8839 | 6 |
| 82.6094 | 7 |
| 89.4299 | 7 |
| 90.5087 | 8 |
| 92.7509 | 8 |
| 96.8869 | 6 |

Table D.1: First positive eigenvalues of the Klein Quartic

| Eigenvalue | Numerical multiplicity |
|---|---|
| 2.45206 | 4 |
| 2.8381 | 3 |
| 4.35754 | 3 |
| 7.61965 | 4 |
| 8.45728 | 1 |
| 10.126 | 4 |
| 10.1635 | 3 |
| 11.76 | 3 |
| 14.1724 | 4 |
| 14.2628 | 2 |
| 17.6818 | 2 |
| 18.812 | 1 |
| 19.2547 | 3 |
| 21.2083 | 4 |
| 22.827 | 4 |
| 24.166 | 3 |
| 24.1744 | 4 |
| 26.9546 | 3 |
| 28.0159 | 1 |
| 28.3388 | 3 |
| 30.7729 | 4 |
| 35.8649 | 3 |
| 37.1106 | 3 |
| 37.3413 | 4 |
| 37.976 | 4 |
| 38.6914 | 3 |
| 38.9579 | 2 |
| 42.7433 | 3 |
| 44.0167 | 4 |
| 44.6724 | 4 |
| 46.6306 | 2 |
| 49.448 | 4 |

Table D.2: First positive eigenvalues of the surface $M3$

| Eigenvalue | Numerical multiplicity |
|---|---|
| 2.4877 | 6 |
| 3.83889 | 3 |
| 7.16032 | 6 |
| 8.24958 | 2 |
| 11.4272 | 6 |
| 13.3096 | 6 |
| 15.049 | 3 |
| 18.6587 | 3 |
| 19.3559 | 6 |
| 23.0736 | 1 |
| 24.7216 | 6 |
| 27.0788 | 6 |
| 27.5915 | 6 |
| 28.0796 | 3 |
| 32.6738 | 1 |
| 36.1647 | 6 |
| 36.2386 | 2 |
| 38.8039 | 6 |
| 40.1168 | 3 |
| 44.0176 | 3 |
| 44.412 | 6 |
| 45.5628 | 6 |

Table D.3: First positive eigenvalues of the Fermat quartic

# Appendix E

# FreeFEM++ Code for Genus 4

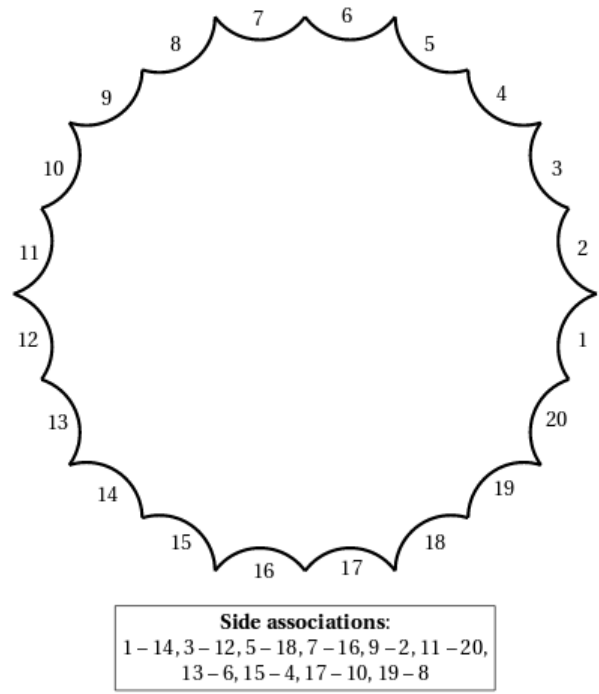

Figure E.1: Fundamental polygon for Bring's surface

The most symmetric surface of genus 4 is Bring's surface (also referred to as Bring's curve) [71]. It can be obtained by associating the sides of a hyperbolic icosagon of area $12\pi$, as shown in Figure E.1. Similarly to the Bolza surface and Klein quartic, we conjecture that because it has maximal symmetry among hyperbolic surfaces of genus 4, it also maximizes the first positive eigenvalue of the Laplace-Beltrami operator among surfaces of this class. In this appendix, we provide a finite element code for computing the first 100 eigenvalues of Bring's surface, along with the output, shown in Table E.1, which lists the eigenvalues and their predicted multiplicities based on the numerical results. The FreeFEM++ code follows the same style as that in earlier appendices; please refer to Appendix C for comments on how to interpret it.

The fundamental polygon for Bring's surface can be tessellated by 240 (2, 4, 5) triangles. The group actions that provide maps between the triangles in this tessellation give the 240 symmetries of Bring's surface. Note that this means Bring's surface is not a Hurwitz surface (recall Theorem 1.30), and that there does not exist a Hurwitz surface of genus 4. The generators of the group satisfy the following relations

$$\langle r,\, s,\, t \,|\, r^5 = s^2 = t^2 = rtrt = stst = (rs)^4 = (sr^3sr^2)^2 = e \rangle,$$

where $r$ is rotation of order 5 about the center of the fundamental icosagon, $s$ is rotation of order 2 about a vertex where 4 triangles meet in the tessellation by (2, 4, 5) triangles, and $t$ is a reflection in the real line (see [71] for more details on this group). The final piece of code in this appendix shows analysis of this finitely presented group in GAP. In particular, we see that it is isomorphic to $S_5 \times \mathbb{Z}_2$, and show its character table.

```
int n = 30;
int nev = 100;

string fnm="brings-surface-Eigenvalues-n30-nev100.txt"; //file name for saving

real cx = 1.00156;
real cy = 0.17178;
real r = 0.180621;

real s = 0.837843;
real theta = 0.138855;
real x01 = s*cos(theta);
real y01 = s*sin(theta);
real x02 = s*cos(pi/5-theta);
real y02 = s*sin(pi/5-theta);
real x03 = s*cos(pi/5+theta);
real y03 = s*sin(pi/5+theta);
real x04 = s*cos(2*pi/5-theta);
real y04 = s*sin(2*pi/5-theta);
real x05 = s*cos(2*pi/5+theta);
real y05 = s*sin(2*pi/5+theta);
real x06 = s*cos(3*pi/5-theta);
real y06 = s*sin(3*pi/5-theta);
real x07 = s*cos(3*pi/5+theta);
real y07 = s*sin(3*pi/5+theta);
real x08 = s*cos(4*pi/5-theta);
real y08 = s*sin(4*pi/5-theta);
real x09 = s*cos(4*pi/5+theta);
real y09 = s*sin(4*pi/5+theta);
real x10 = s*cos(5*pi/5-theta);
real y10 = s*sin(5*pi/5-theta);
real x11 = s*cos(5*pi/5+theta);
real y11 = s*sin(5*pi/5+theta);
real x12 = s*cos(6*pi/5-theta);
real y12 = s*sin(6*pi/5-theta);
real x13 = s*cos(6*pi/5+theta);
real y13 = s*sin(6*pi/5+theta);
real x14 = s*cos(7*pi/5-theta);
real y14 = s*sin(7*pi/5-theta);
```

```
real x15 = s*cos(7*pi/5+theta);
real y15 = s*sin(7*pi/5+theta);
real x16 = s*cos(8*pi/5-theta);
real y16 = s*sin(8*pi/5-theta);
real x17 = s*cos(8*pi/5+theta);
real y17 = s*sin(8*pi/5+theta);
real x18 = s*cos(9*pi/5-theta);
real y18 = s*sin(9*pi/5-theta);
real x19 = s*cos(9*pi/5+theta);
real y19 = s*sin(9*pi/5+theta);
real x20 = s*cos(theta);
real y20 = s*sin(-theta);

border G1 (t=7*pi/5,4*pi/5) { x=cx+r*cos(t); y=cy+r*sin(t);};
border G2 (t=4*pi/5,7*pi/5) { x=cos(pi/5)*(cx+r*cos(t))-sin(pi/5)*(-cy-r*sin(t));
    y=sin(pi/5)*(cx+r*cos(t))+cos(pi/5)*(-cy-r*sin(t));};
border G3 (t=7*pi/5,4*pi/5) { x=cos(pi/5)*(cx+r*cos(t))-sin(pi/5)*(cy+r*sin(t));
    y=sin(pi/5)*(cx+r*cos(t))+cos(pi/5)*(cy+r*sin(t));};
border G4 (t=4*pi/5,7*pi/5) { x=cos(2*pi/5)*(cx+r*cos(t))-sin(2*pi/5)*(-cy-r*sin(t));
    y=sin(2*pi/5)*(cx+r*cos(t))+cos(2*pi/5)*(-cy-r*sin(t));};
border G5 (t=7*pi/5,4*pi/5) { x=cos(2*pi/5)*(cx+r*cos(t))-sin(2*pi/5)*(cy+r*sin(t));
    y=sin(2*pi/5)*(cx+r*cos(t))+cos(2*pi/5)*(cy+r*sin(t));};
border G6 (t=4*pi/5,7*pi/5) { x=cos(3*pi/5)*(cx+r*cos(t))-sin(3*pi/5)*(-cy-r*sin(t));
    y=sin(3*pi/5)*(cx+r*cos(t))+cos(3*pi/5)*(-cy-r*sin(t));};
border G7 (t=7*pi/5,4*pi/5) { x=cos(3*pi/5)*(cx+r*cos(t))-sin(3*pi/5)*(cy+r*sin(t));
    y=sin(3*pi/5)*(cx+r*cos(t))+cos(3*pi/5)*(cy+r*sin(t));};
border G8 (t=4*pi/5,7*pi/5) { x=cos(4*pi/5)*(cx+r*cos(t))-sin(4*pi/5)*(-cy-r*sin(t));
    y=sin(4*pi/5)*(cx+r*cos(t))+cos(4*pi/5)*(-cy-r*sin(t));};
border G9 (t=7*pi/5,4*pi/5) { x=cos(4*pi/5)*(cx+r*cos(t))-sin(4*pi/5)*(cy+r*sin(t));
    y=sin(4*pi/5)*(cx+r*cos(t))+cos(4*pi/5)*(cy+r*sin(t));};
border G10 (t=4*pi/5,7*pi/5) { x=cos(pi)*(cx+r*cos(t))-sin(pi)*(-cy-r*sin(t));
    y=sin(pi)*(cx+r*cos(t))+cos(pi)*(-cy-r*sin(t));};
border G11 (t=7*pi/5,4*pi/5) { x=cos(pi)*(cx+r*cos(t))-sin(pi)*(cy+r*sin(t));
    y=sin(pi)*(cx+r*cos(t))+cos(pi)*(cy+r*sin(t));};
border G12 (t=4*pi/5,7*pi/5) { x=cos(6*pi/5)*(cx+r*cos(t))-sin(6*pi/5)*(-cy-r*sin(t));
    y=sin(6*pi/5)*(cx+r*cos(t))+cos(6*pi/5)*(-cy-r*sin(t));};
border G13 (t=7*pi/5,4*pi/5) { x=cos(6*pi/5)*(cx+r*cos(t))-sin(6*pi/5)*(cy+r*sin(t));
    y=sin(6*pi/5)*(cx+r*cos(t))+cos(6*pi/5)*(cy+r*sin(t));};
border G14 (t=4*pi/5,7*pi/5) { x=cos(7*pi/5)*(cx+r*cos(t))-sin(7*pi/5)*(-cy-r*sin(t));
    y=sin(7*pi/5)*(cx+r*cos(t))+cos(7*pi/5)*(-cy-r*sin(t));};
border G15 (t=7*pi/5,4*pi/5) { x=cos(7*pi/5)*(cx+r*cos(t))-sin(7*pi/5)*(cy+r*sin(t));
    y=sin(7*pi/5)*(cx+r*cos(t))+cos(7*pi/5)*(cy+r*sin(t));};
border G16 (t=4*pi/5,7*pi/5) { x=cos(8*pi/5)*(cx+r*cos(t))-sin(8*pi/5)*(-cy-r*sin(t));
    y=sin(8*pi/5)*(cx+r*cos(t))+cos(8*pi/5)*(-cy-r*sin(t));};
border G17 (t=7*pi/5,4*pi/5) { x=cos(8*pi/5)*(cx+r*cos(t))-sin(8*pi/5)*(cy+r*sin(t));
    y=sin(8*pi/5)*(cx+r*cos(t))+cos(8*pi/5)*(cy+r*sin(t));};
border G18 (t=4*pi/5,7*pi/5) { x=cos(9*pi/5)*(cx+r*cos(t))-sin(9*pi/5)*(-cy-r*sin(t));
    y=sin(9*pi/5)*(cx+r*cos(t))+cos(9*pi/5)*(-cy-r*sin(t));};
border G19 (t=7*pi/5,4*pi/5) { x=cos(9*pi/5)*(cx+r*cos(t))-sin(9*pi/5)*(cy+r*sin(t));
    y=sin(9*pi/5)*(cx+r*cos(t))+cos(9*pi/5)*(cy+r*sin(t));};
border G20 (t=4*pi/5,7*pi/5) { x=cx+r*cos(t); y=-cy-r*sin(t);};

plot(G1(n)+G2(n)+G3(n)+G4(n)+G5(n)+G6(n)+G7(n)+G8(n)+G9(n)+G10(n)+G11(n)+G12(n)+G13(n)
    +G14(n)+G15(n)+G16(n)+G17(n)+G18(n)+G19(n)+G20(n));mesh
Th=buildmesh(G1(n)+G2(n)+G3(n)+G4(n)+G5(n)+G6(n)+G7(n)+G8(n)+G9(n)+G10(n)+G11(n)+G12(n)
    +G13(n)+G14(n)+G15(n)+G16(n)+G17(n)+G18(n)+G19(n)+G20(n),fixeborder=true);
plot(Th,wait=true,fill=true);

fespace Vh(Th,P2,periodic=[[1,(cos(pi/10))*(y-y01)-(sin(pi/10))*(x-x01)],
    [14,(cos(3*pi/10))*(y-y14)-(sin(3*pi/10))*(x-x14)],
    [3,(cos(3*pi/10))*(y-y03)-(sin(3*pi/10))*(x-x03)],
    [12,(cos(pi/10))*(y-y12)-(sin(pi/10))*(x-x12)],
    [5,(cos(5*pi/10))*(y-y05)-(sin(5*pi/10))*(x-x05)],
    [18,(cos(7*pi/10))*(y-y18)-(sin(7*pi/10))*(x-x18)],
```

```
    [7,(cos(7*pi/10))*(y-y07)-(sin(7*pi/10))*(x-x07)],
    [16,(cos(5*pi/10))*(y-y16)-(sin(5*pi/10))*(x-x16)],
    [9,(cos(9*pi/10))*(y-y09)-(sin(9*pi/10))*(x-x09)],
    [2,(cos(11*pi/10))*(y-y02)-(sin(11*pi/10))*(x-x02)],
    [11,(cos(11*pi/10))*(y-y11)-(sin(11*pi/10))*(x-x11)],
    [20,(cos(9*pi/10))*(y-y20)-(sin(9*pi/10))*(x-x20)],
    [13,(cos(13*pi/10))*(y-y13)-(sin(13*pi/10))*(x-x13)],
    [6,(cos(15*pi/10))*(y-y06)-(sin(15*pi/10))*(x-x06)],
    [15,(cos(15*pi/10))*(y-y15)-(sin(15*pi/10))*(x-x15)],
    [4,(cos(13*pi/10))*(y-y04)-(sin(13*pi/10))*(x-x04)],
    [17,(cos(17*pi/10))*(y-y17)-(sin(17*pi/10))*(x-x17)],
    [10,(cos(19*pi/10))*(y-y10)-(sin(19*pi/10))*(x-x10)],
    [19,(cos(19*pi/10))*(y-y19)-(sin(19*pi/10))*(x-x19)],
    [8,(cos(17*pi/10))*(y-y08)-(sin(17*pi/10))*(x-x08)]]);

Vh u1, u2;
real sigma = 0.00001;
varf op(u1,u2)= int2d(Th)(dx(u1)*dx(u2)+ dy(u1)*dy(u2)- sigma* u1*u2);
varf b([u1],[u2]) = int2d(Th)((4*u1*u2)/(1-x^2-y^2)^2);

matrix OP= op(Vh,Vh,solver=Crout,factorize=1);
matrix B=b(Vh,Vh,solver=CG,eps=1e-20);

real[int] ev(nev);
Vh[int] eV(nev);

int k=EigenValue(OP,B,sym=true,sigma=sigma,value=ev,vector=eV, tol=1e-10,maxit=0,ncv=0);
for (int i=0;i<k;i++) {
u1=eV[i];
real gg = int2d(Th)(dx(u1)*dx(u1) + dy(u1)*dy(u1));  real mm= int2d(Th)(u1*u1) ;

//cout<<"----"<< i<<""<<ev[i]<<"err="
//<<int2d(Th)(4/(1-x^2-y^2)^2) << " --- "<<endl; plot(eV[i],cmm="Eigen Vector "+i+" valeur ="
    + ev[i] ,wait=1,value=1);

ofstream Eva(fnm, append);
Eva <<  ev[i] << "\n";

}
```

The following shows the GAP analysis of the symmetry group of Bring's surface. Please refer to Appendix A for details on how to interpret this code.

```
gap> f:=FreeGroup(3);
<free group on the generators [ f1, f2, f3 ]>
gap> g:=f/[f.1^5,f.2^2,f.3^2,f.1*f.3*f.1*f.3,f.2*f.3*f.2*f.3,
    (f.1*f.2)^4,(f.2*f.1^3*f.2*f.1^2)^2];
<fp group on the generators [ f1, f2, f3 ]>
gap> Order(g);
240
gap> StructureDescription(g);
"C2_x_S5"
gap> RequirePackage("repsn");;
gap> Display(CharacterTable(g));
CT1

    2  4  3  4  2  3  2  1  3  2  3  2   1  4  4
    3  1  1  .  1  1  1  .  .  1  .  1   .  .  1
    5  1  .  .  .  .  1  .  .  .  .  1   .  1
```

```
          1a  2a  2b  6a  2c  3a  5a  4a  6b  4b  6c  10a  2d  2e
      2P  1a  1a  1a  3a  1a  3a  5a  2d  3a  2d  3a   5a  1a  1a
      3P  1a  2a  2b  2c  2c  1a  5a  4a  2a  4b  2e  10a  2d  2e
      5P  1a  2a  2b  6a  2c  3a  1a  4a  6b  4b  6c   2e  2d  2e
      7P  1a  2a  2b  6a  2c  3a  5a  4a  6b  4b  6c  10a  2d  2e

X.1        1   1   1   1   1   1   1   1   1   1   1    1   1   1
X.2        1  −1  −1   1   1   1   1  −1  −1   1  −1   −1   1  −1
X.3        1  −1   1  −1  −1   1   1  −1  −1  −1   1    1   1   1
X.4        1   1  −1  −1  −1   1   1   1   1  −1  −1   −1   1  −1
X.5        4  −2   .   1  −2   1  −1   .   1   .   1   −1   .   4
X.6        4   2   .  −1   2   1  −1   .  −1   .   1   −1   .   4
X.7        4   2   .   1  −2   1  −1   .  −1   .  −1    1   .  −4
X.8        4  −2   .  −1   2   1  −1   .   1   .  −1    1   .  −4
X.9        5   1   1   1   1  −1   .  −1   1  −1  −1    .   1   5
X.10       5  −1   1  −1  −1  −1   .   1  −1   1  −1    .   1   5
X.11       5  −1  −1   1   1  −1   .   1  −1  −1   1    .   1  −5
X.12       5   1  −1  −1  −1  −1   .  −1   1   1   1    .   1  −5
X.13       6   .  −2   .   .   .   1   .   .   .   .    1  −2   6
X.14       6   .   2   .   .   .   1   .   .   .   .   −1  −2  −6
gap> ConjugacyClasses(g);
[ <identity ...>^G, f1^2*(f1^2*f2)^2*f1*f2*f1*f3^G, f3^G,
  f1^2*(f1^2*f2)^2*f1*f2*f1^G, f2^G, f1^3*(f1*f2)^2*f1^4^G, f1^G,
  f1^2*(f1^2*f2)^2*f1*f2*f3^G, f2*f1^4*f3^G, f2*f1^G,
  f1^3*(f1*f2)^2*f1^3*f3^G, f1^3*(f1*f2)^2*f1^2*f3^G, f1*(f1*f2)^2*f1^4^G,
  (f1*f2*f1)^3*f2*f3^G ]
```

| Eigenvalue | Numerical multiplicity |
|---|---|
| 1.91556 | 6 |
| 2.78954 | 5 |
| 5.89883 | 5 |
| 7.3379 | 5 |
| 8.26099 | 4 |
| 8.55338 | 6 |
| 13.2289 | 6 |
| 13.5082 | 4 |
| 15.1462 | 4 |
| 17.0459 | 6 |
| 17.7721 | 5 |
| 19.0882 | 6 |
| 22.5953 | 4 |
| 24.3521 | 5 |
| 24.6075 | 4 |
| 28.0573 | 6 |
| 28.5767 | 6 |
| 29.6861 | 4 |
| 31.4013 | 5 |

Table E.1: First positive eigenvalues of Bring's surface

# Appendix F

# Mathematica Code for Bolza Surface

The following block of code helps us to choose the optimal value of $L$ to define the spectral function $h(t)$ in Section 3.4. We first define "lemult" as the ordered list of lengths of geodesics in the length spectrum, using the formula in Section 3.2, along with their multiplicities. Note that only the multiplicity of the first length in the spectrum is precise; all others have been calculated numerically. Also note that whilst the multiplicity of the systole of the Bolza surface is 12, here we are concerned about moving along the geodesic in both a positive and negative direction, so we double the multiplicity to take this into account. We then use "les" to flatten this into one list, rather than a list of pairs, where each length appears multiple times.

We define the spectral function
$$h(x, L) = \left(\frac{\sin Lx}{Lx}\right)^4,$$
along with its Fourier transform $g(\xi, L)$, and the derivative of the Fourier transform, $gd(x, L)$. Note that all of these are functions of $L$ as well as the usual variable.

For a given value of $L$, the function "multbound" calculates the bound on the multiplicity of the first eigenvalue, using the method described in Section 3.4. Finally, we can use Mathematica's "DiscretePlot" function to plot the bounds on $m$ given by "multbound(L)" for values of $L$ between 0.4 and 1, with a step size of 0.01. The output is shown in Figure F.1, where we see that $L = \frac{93}{100}$ gives the lowest bound.

```
lemult = {{2 ArcCosh[1 + Sqrt[2]], 24}, {2 ArcCosh[3 + 2 Sqrt[2]],
    24}, {2 ArcCosh[5 + 3 Sqrt[2]], 48}, {2 ArcCosh[7 + 5 Sqrt[2]],
    96}, {2 ArcCosh[9 + 6 Sqrt[2]], 48}, {2 ArcCosh[9 + 7 Sqrt[2]],
    48}, {2 ArcCosh[11 + 8 Sqrt[2]], 8}, {2 ArcCosh[13 + 9 Sqrt[2]],
    96}};

les = Flatten[
   Table[Table[lemult[[i]][[1]], {k, 1, lemult[[i]][[2]]}], {i, 1,
     8}]];
```

```
h[x_, L_] = (Sin[L x]/(L x))^4
g[\[Xi]_, L_] = FourierTransform[h[t, L], t, \[Xi]]/Sqrt[2 Pi]
Zer[x_] = 0;
gd[x_, L_] = D[g[x, L], x] /. {Sign'->Zer}

multbound[L_]:=(h[x_]=h[x,L];
    g[\[Xi]_]=g[\[Xi],L];
    gd[x_]=gd[x,L];
    RightSide=-2NIntegrate[gd[x]/Sinh[x/2],{x,0,4L}]+
    Sum[Sum[les[[i]]g[k*les[[i]]]/(2 Sinh[k*les[[i]]/2]),{i,1,Length[les]}],{k,1,10}];
    (RightSide-h[I/2])/h[Sqrt[116469/28089-1/4]])

DiscretePlot[multbound[L],{L,0.4,1.0,0.01}]
```

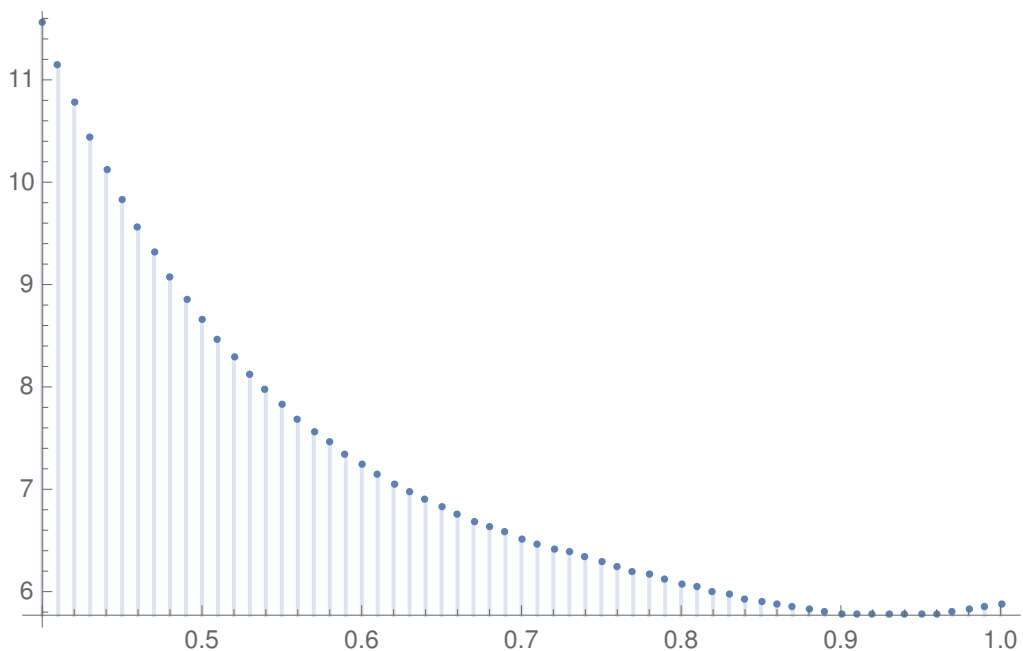

Figure F.1: Output of the function "multbound" for different values of $L$

We also include the code used to produce Figure 3.19. Here we construct a test function to satisfy the hypothesis of Corollary 3.17, and plot its nodal line as a function of $t$. We superimpose this plot of the nodal line against a plot of the boundary of the annulus considered in Section 3.4.1 to show that the test function is positive inside the desired region.

```
(* Eigenvalue *)
\[Lambda] = 5.;

(* Length of core geodesic *)
ell = 2. ArcCosh[1/2 Csc[\[Pi]/8]];

(* Hypergeometric functions with eigenvalue \[Lambda]. For efunc1, \
```

```
(* the m variable gives the restriction of the eigenfunction to the mth \
Fourier node *)
efunc0[\[Rho]_] =
  Re[Evaluate[
    Sinh[\[Rho]] Hypergeometric2F1[(3/2 + I Sqrt[\[Lambda] - 1/4])/
       2, (2 - (1/2 + I Sqrt[\[Lambda] - 1/4]))/2,
      3/2, -Sinh[\[Rho]]^2] ]];
efunc1[\[Rho]_, m_] =
  Re[Evaluate[
    Sinh[\[Rho]] Cosh[\[Rho]]^(2 Pi I m/
         ell) Hypergeometric2F1[(3/2 + I Sqrt[\[Lambda] - 1/4])/2 +
       I m Pi/ell, (2 - (1/2 + I Sqrt[\[Lambda] - 1/4]))/2 +
       Pi I m /ell, 3/2, -Sinh[\[Rho]]^2] ]];

(* Nodal Line as a function of t. *)
\[Rho]zero[
   t_] := (\[Rho] /.
    FindRoot[
     efunc0[\[Rho]] - (0.038 Sin[2 Pi (t)/ell - Pi/2]) efunc1[\[Rho],
         1] - (0.0002 Cos[2 2 Pi (t)/ell]) efunc1[\[Rho], 2], {\[Rho],
       0.5, 1.9}]);
nodal = Plot[{\[Rho]zero[t]}, {t, 0, ell}];

(* Parametric plot of boundary of annulus *)
x1 = Cos[2 Pi t/ell];
y1 = Sin[2 Pi t/ell];
x2 = 4.61158 Cos[2 Pi t/ell];
y2 = 4.61158 Sin[2 Pi t/ell];
x3 = 1.0986840182524393 + 0.4550896306919397 Cos[2 Pi t/ell];
y3 = 0.4550896306919397 Sin[2 Pi t/ell];
x3 = 1.0986840182524393 + 0.4550896306919397 Cos[2 Pi t/ell];
y3 = 0.4550896306919397 Sin[2 Pi t/ell];
x4 = 2.5537741953480246 + 1.3820929243421862 Cos[2 Pi t/ell];
y4 = 1.3820929243421862 Sin[2 Pi t/ell];
x5 = 5.066671210777469 + 2.0986830534021754 Cos[2 Pi t/ell];
y5 = 2.0986830534021754 Sin[2 Pi t/ell];
boundary =
  ParametricPlot[{{Log[Sqrt[y1^2 + x1^2]],
     ArcSinh[x1/y1]}, {Log[Sqrt[y2^2 + x2^2]],
     ArcSinh[x2/y2]}, {Log[Sqrt[y3^2 + x3^2]],
     ArcSinh[x3/y3]}, {Log[Sqrt[y4^2 + x4^2]],
     ArcSinh[x4/y4]}, {Log[Sqrt[y5^2 + x5^2]], ArcSinh[x5/y5]}}, {t,
    0.0001, ell}, AspectRatio -> Automatic, PlotRange -> All];

(* Superimpose nodal line plot onto boundary of annulus *)
\
Show[nodal, boundary, PlotRange -> {{0, 1.6}, {1, 2.3}}]
```

# Appendix G

# Alternative approach to proving Conjecture 3.12

In this appendix, we outline another method that could be used to prove Conjecture 3.12. As in Section 3.4.1, consider the annulus $\Omega$ in Figure 3.18 with Dirichlet boundary conditions. In particular, we can take $\Omega$ to be a subspace of a hyperbolic cylinder with $\alpha$ as its core geodesic. We have a Dirichlet boundary condition on $\alpha$, so one method that we could use to bound $\lambda_1(\Omega)$ would be to compare it to the first eigenvalue of the symmetric half cylinder with core geodesic $\alpha$, having the same area as $\Omega$, that is, $\pi$. Here we make the following conjecture, which is essentially a Faber-Krahn inequality for half-cylinders.

**Conjecture G.1** *Let $\Omega_0$ be the symmetric half-cylinder with core geodesic $\alpha$ of length*

$$l(\alpha) = 2\operatorname{arccosh}\left(\frac{\csc\left(\frac{\pi}{8}\right)}{2}\right),$$

*and area $\pi$. Then*

$$\lambda_1(\Omega) \geq \lambda_1(\Omega_0).$$

**Remark G.2** *Recall from Lemma 2.35 that annuli in hyperbolic cylinders satisfy an isoperimetric inequality. This is a similar hypothesis to that in the classical Faber-Krahn inequality (compare Theorem 1.14), so it may be possible to prove Conjecture G.1 in a similar way.*

If we assume Conjecture G.1 holds, we consider the first odd eigenfunction on $\Omega_0$, that is,

$$\phi(s) = \sinh(l(\alpha)) \, {}_2\mathbf{F}_1\left(\frac{1+s}{2}, \frac{2-s}{2}; \frac{3}{2}; -\sinh^2(l(\alpha))\right)$$

as a function of $t$, with $s = \frac{1}{2} + it$. The corresponding eigenvalue $\lambda = s(1-s)$ is obtained by taking the first zero of $\phi(s)$. This eigenfunction is odd on a symmetric cylinder of area $2\pi$ with respect to the core geodesic $\alpha$, so will satisfy a Dirichlet boundary condition on $\alpha$, as well as its other boundary geodesic, that is, it satisfies Dirichlet boundary conditions on a symmetric half-cylinder of area $\pi$.

To get a lower bound on $\lambda$, we analyze the differential equation, loosely following the method of Appendix A in [77]. As a function of $\rho$, we know that

$$\phi(\rho) = \sinh(\rho) \, _2F_1\left(\frac{1+s}{2}, \frac{2-s}{2}; \frac{3}{2}; -\sinh^2(\rho)\right)$$

satisfies

$$\left(-\frac{1}{\cosh(\rho)}\frac{d}{d\rho}\cosh(\rho)\frac{d}{d\rho} - \lambda\right)\phi(\rho) = 0$$

with initial condition $\phi(0) = 0$ and $\frac{d}{d\rho}\phi(0) = 1$, where $\lambda = s(1-s)$ (recall our definition of the Laplacian on a hyperbolic cylinder in Section 2.3.2). We substitute

$$\varphi = 2\arctan\left(\tanh\left(\frac{\rho}{2}\right)\right)$$

to get

$$\left(-\frac{d^2}{d\varphi^2} - \frac{\lambda}{\cos^2(\varphi)}\right)\Psi(\varphi) = 0, \tag{G.1}$$

where $\Psi(\varphi) = \phi(\rho(\varphi))$ satisfies $\Psi(0) = 0$ and $\frac{d}{d\rho}\Psi(0) = 1$. We prove the following

**Theorem G.3** *Suppose we have $g \in C^2[0, x_0]$ such that $g \geq 0$ on $[0, x_0]$, $g > 0$ on $(0, x_0)$, and*

$$\left(-\frac{g''(x)}{g(x)} - \frac{\lambda}{\cos^2(x)}\right) \geq 0 \tag{G.2}$$

*on $[0, x_0]$. Then $f \in C^2[0, x_0]$ that is a solution to Equation (G.1), that is, $f$ solves*

$$-f''(x) - \frac{\lambda}{\cos^2(x)}f(x) = 0, \quad f(0) = 0, \quad f'(0) = 1,$$

*has no zero in $(0, x_0)$.*

**Proof:** The proof is done by contradiction. Let $r_0 \in (0, x_0)$ be the first positive zero of $f$. Then, on $[0, r_0)$ we have $f \geq 0$, $g \geq 0$, and on $(0, r_0)$, $g > 0$. We also have $f > 0$ on $(0, r_0)$. This is because $f'(0) = 1$. The definition of the derivative implies that $f$ must be positive for some value close to zero, in particular

$$1 = f'(0) = \lim_{h \to 0} \frac{f(0+h)}{h},$$

so $f$ will be positive until it is zero again at $r_0$. Consider the Wronskian

$$F(x) = f'(x)g(x) - g'(x)f(x),$$

noting that $F(0) = 0$, and its derivative

$$F'(x) = f''(x)g(x) - g''(x)f(x)$$

$$= \left(-\frac{\lambda}{\cos^2(x)} - \frac{g''(x)}{g(x)}\right) g(x) f(x).$$

Then $F'(x) \geq 0$, that is, it is increasing on $(0, r_0)$, which implies $F(r_0) \geq 0$. On the other hand, consider the limit from below of

$$f'(r_0) = \lim_{h \to 0^-} \frac{f(r_0 + h)}{h},$$

that is, for negative $h$. We know that $f(r_0 - h)$ is positive for small $h$, since $f(x)$ is positive on $(0, r_0)$, so $f'(r_0) < 0$. With $f'(r_0) < 0$ and $g(r_0) > 0$, we have $F(r_0) < 0$, a contradiction. $\square$

Choosing a suitable function $g(x)$ that satisfies the hypothesis of Theorem G.3, we can show that $\phi(\rho)$ does not have a zero on the interval $(0, \rho)$ for certain $\lambda$. The "certain $\lambda$" condition will give us our bound on $\lambda_1(\Omega)$. First, we need to change variables from $\rho$ to $\varphi$, that is

$$\varphi = 2 \arctan\left(\tanh\left(\frac{\operatorname{arcsinh}\left(\frac{\pi}{2\operatorname{arccosh}\left(\frac{1}{2}\csc\left(\frac{\pi}{8}\right)\right)}\right)}{2}\right)\right) \approx 1.11795903.$$

Take $g(x) = x - \frac{4x^3}{5}$. We have to show that it is positive on $(0, \varphi)$, and find the values of $\lambda$ for which Equation (G.2) is satisfied on $[0, \varphi]$. Let us start by showing $g(x) > 0$. We have

$$g'(x) = 1 - \frac{12x^2}{5},$$

which is zero if $x = \pm\sqrt{\frac{5}{12}}$. The second derivative

$$g''(x) = -\frac{24x}{5}$$

tells us that $g(x)$ has a maximum at $\sqrt{\frac{5}{12}}$ and a minimum at $-\sqrt{\frac{5}{12}}$. Now,

$$g\left(\sqrt{\frac{5}{12}}\right) = \sqrt{\frac{5}{12}}\left(\frac{2}{3}\right),$$

so the maximum value is positive. The roots of $g(x)$ are at $x = 0$ and $x = \pm\sqrt{\frac{5}{4}}$. Note that $\sqrt{\frac{5}{4}}$ is outside the interval $[0, \varphi]$, so $g(0) = 0$ and for all other values on the interval it is positive.

Now we have to find the values of $\lambda > 0$ for which

$$\frac{24}{5} \geq \frac{\lambda\left(1 - \frac{4x^2}{5}\right)}{\cos^2(x)}.$$

Cross multiplying here is fine, since $1 - \frac{4x^2}{5} > 0$ on $[0, \varphi]$. Let us try to maximize

$$h(x) = \sec^2(x)\left(1 - \frac{4x^2}{5}\right).$$

The (diagonal) Padé approximant of $\sec^2(x)$ of order 4 is

$$\frac{\frac{13x^4}{945} + \frac{11x^2}{63} + 1}{\frac{163x^4}{945} - \frac{52x^2}{63} + 1}.$$

We use this expression to approximate $\sec^2(x)$ on $[0, \varphi]$; we can get an error bound on our approximant using Taylor's remainder theorem. First, we maximize

$$\tilde{h}(x) = \left(\frac{\frac{13x^4}{945} + \frac{11x^2}{63} + 1}{\frac{163x^4}{945} - \frac{52x^2}{63} + 1}\right)\left(1 - \frac{4x^2}{5}\right).$$

We have

$$\tilde{h}'(x) = -\frac{2x\left(8476x^8 - 81120x^6 - 798345x^4 + 2664900x^2 - 893025\right)}{5\left(163x^4 - 780x^2 + 945\right)^2}.$$

This clearly has a root at zero, and to find the other roots, we need to find the roots of the order 8 polynomial in parentheses. Notice that this is quartic in $x^2$, so we may use Ferrari's method to find its roots and then take their square roots. For a general quartic polynomial

$$ax^4 + bx^3 + cx^2 + dx + e,$$

the four roots are given by

$$x_{1,2} = -\frac{b}{4a} \pm \frac{1}{2}\sqrt{-2p + \frac{q}{S} - 4S^2} - S \quad \text{and}$$

$$x_{3,4} = -\frac{b}{4a} \pm \frac{1}{2}\sqrt{-2p - \frac{q}{S} - 4S^2} + S,$$

where

$$p = \frac{8ac - 3b^2}{8a^2},$$

$$q = \frac{8a^2d - 4abc + b^3}{8a^3},$$

$$\Delta_0 = 12ae - 3bd + c^2,$$

$$\Delta_1 = -72ace + 27ad^2 + 27b^2e - 9bcd + 2c^3,$$

$$Q = \sqrt[3]{\frac{\sqrt{\Delta_1^2 - 4\Delta_0^3} + \Delta_1}{2}}, \quad \text{and}$$

$$S = \frac{1}{2}\sqrt{\frac{1}{3a}\left(\frac{\Delta_0}{Q} + Q\right) - \frac{2p}{3}}.$$

We find that for $a = 8476$, $b = -81120$, and so on, the only positive real root of the polynomial whose square root lies in $[0, \varphi]$ is $x_1$. Therefore, the turning points of $\tilde{h}'(x)$ that we need to consider are at $0$ and $\sqrt{x_1}$. The second derivative

$$-\frac{2\left(1381588x^{12} - 19833840x^{10} + 652299885x^8 - 2085793200x^6 - 674616600x^4 + 5465313000x^2 - 843908625\right)}{5\left(163x^4 - 780x^2 + 945\right)^3}$$

is positive at $0$ and negative at $x_1$, so $\tilde{h}(x)$ takes its maximum value at $x_1$. We do not show the lengthy closed form here, but it is roughly equal to $1.04542$. To find out how accurate this is, we can compare the Taylor series of $\sec^2(x)$,

$$1 + x^2 + \frac{2x^4}{3} + \frac{17x^6}{45} + \frac{62x^8}{315} + \frac{1382x^{10}}{14175} + O\left(x^{11}\right),$$

with the Taylor series of our approximant, that is,

$$1 + x^2 + \frac{2x^4}{3} + \frac{17x^6}{45} + \frac{62x^8}{315} + \frac{28963x^{10}}{297675} + O\left(x^{11}\right).$$

We see that these agree for the first five terms, and so the error is the maximum value of

$$\frac{1382x^{10}}{14175} - \frac{28963x^{10}}{297675} = \frac{59x^{10}}{297675}$$

on $[0, \varphi]$. This occurs at $\varphi$, since $\sec^2(x)$ is increasing on this interval, so our error is

$$\delta x_1 = \frac{59\varphi^{10}}{297675} \approx 0.000604461.$$

Now

$$\frac{24}{5} \geq \lambda\left(\tilde{h}(x_1 - \delta x_1)\right)$$

is satisfied for positive $\lambda$ when

$$\lambda \leq \frac{24}{5\left(\tilde{h}(x_1 - \delta x_1)\right)} \approx 4.591456764,$$

that is, for any such $\lambda$, $\phi(\rho)$ does not have a zero on $(0, \rho)$. Hence

$$\lambda_1(\Omega_i) > \frac{24}{5\left(\tilde{h}(x_1 - \delta x_1)\right)} \approx 4.591456764, \qquad i \in \{11, 12, 13\}.$$